# EQUIVARIANT LOCALIZATION IN FACTORIZATION HOMOLOGY AND APPLICATIONS IN MATHEMATICAL PHYSICS I: FOUNDATIONS

## DYLAN BUTSON

ABSTRACT. We develop a theory of equivariant factorization algebras on varieties with an action of a connected algebraic group $G$, extending the definitions of Francis-Gaitsgory [FG11] and Beilinson-Drinfeld [BD04] to the equivariant setting. We define an equivariant analogue of factorization homology, valued in modules over $H_G^\bullet(\mathrm{pt})$, and in the case $G = (\mathbb{C}^\times)^n$ we prove an equivariant localization theorem for factorization homology, analogous to the classical localization theorem [AB95]. We establish a relationship between $\mathbb{C}^\times$ equivariant factorization algebras and filtered quantizations of their restrictions to the fixed point subvariety. These results provide a model for predictions from the physics literature about the $\Omega$-background construction introduced in [Nek03], interpreting factorization $\mathbb{E}_n$ algebras as observables in mixed holomorphic-topological quantum field theories.

In the companion paper [But20b], we develop tools to give geometric constructions of factorization $\mathbb{E}_n$ algebras, and apply them to define those corresponding to holomorphic-topological twists of supersymmetric gauge theories in low dimensions. Further, we apply our above results in these examples to give an account of the predictions of [CG18] as well as [BLL+15], and explain the relation between these constructions from this perspective.

## CONTENTS

## 1. INTRODUCTION

We begin with some overarching remarks about the background and motivation for the present work, as well as its companion paper and formal sequel [But20b], which we call Parts I and II, respectively. In Section 1.3 we give a detailed overview and summary of Part I. In Section 1.4 we outline the complementary results of Part II, and in Section 1.5 we give a brief description of the future directions of this project. Together, these results will comprise the Author's thesis.

1.1. **Background: Factorization algebras, representation theory, and quantum field theory.** Factorization algebras were introduced by Beilinson and Drinfeld in [BD04] as a model for algebras of observables in two dimensional chiral conformal quantum field theories, defined in the language of algebraic geometry. Factorization algebras in this setting generalize vertex algebras to global objects defined over algebraic curves, vaguely analogous to sheaves on them. A generalization of the theory of factorization algebras to higher dimensional varieties was also given in [FG11], analogously modeling holomorphic quantum field theories in higher dimensions, which by definition generalize the holomorphic behaviour of observables in chiral conformal field theories in two real dimensions. From the beginning, the development of this theory was motivated by the essential connection between chiral conformal field theory and representation theory of affine Lie algebras.

An analogue of factorization algebras defined over smooth manifolds in the language of algebraic topology was proposed by Lurie in [Lur08], as an example of a class of extended topological field theories in the mathematical sense defined therein, and pursued by Ayala, Francis, Lurie, and collaborators in [AF15, AFR15, AFT16, Lur09a, Lur12]. Factorization algebras in the topological setting analogously generalize algebras over the little $n$-discs operad, and again describe the algebras of observables in topological quantum field theories of dimension $n$. In the case $n = 1$ these are equivalent to usual (homotopy) associative algebras, a central topic of study in classical representation theory.

Thus, there is a natural dictionary between predictions of quantum field theory or string theory, which have led to groundbreaking ideas in a variety of areas of mathematics, and statements in representation theory phrased in terms of factorization algebras. This dictionary is both a primary motivation and the main source of new ideas for this series of papers.

The work of Costello [Cos11] and Costello-Gwilliam [CG16] established such a dictionary in a more analytic context, constructing a variant of factorization algebras defined over smooth manifolds in terms of the differential geometric input data of a Lagrangian classical field theory satisfying certain ellipticity requirements together with a choice of renormalization scheme. These ideas were very influential for the present series of papers, and have led to many other developments following this paradigm [Cos13, CS15, BY16, GW18, ES19, SW19, ESW20].

The present series of papers also closely follows the program of Ben-Zvi, Nadler, and collaborators, which gives approaches to many facets of geometric representation theory in terms of extended topological field theory and derived algebraic geometry [BZN09, BZFN10, BZN13, BZG17, BZN18]. In particular, the use of sheaf theory in constructing extended topological field theories from geometry is a central theme of the present series of papers, which is borrowed from *loc. cit.*. Further, the derived stacks and sheaf theories defined on them which are relevant for our constructions can often be predicted from statements about the shifted symplectic geometry of the spaces of solutions to the Euler-Lagrange equations in the relevant classical field theories. This relies on a family of ideas about functoriality of shifted geometric quantization, closely related to those in *loc. cit.*, which I learned from Pavel Safronov [Saf20].

Finally, the circle of ideas and mathematical technology around the local geometric Langlands correspondence [ABC+18], derived geometric Satake correspondence [BF08], and Coulomb branch construction [BFM05, BFN18, BFN19b], provided a collection of mathematically well-understood examples and established techniques which were crucial for the technical underpinning for the present series of papers. In particular, we follow the sheaf theory foundations given in [GR14a, GR14b, Gai15, GR17a, GR17b, Ras15a, Ras15b, Ras20b] and references therein. These ideas can naturally be interpreted in certain holomorphic-topological twists of supersymmetric quantum field theories, as we explain below. These interpretations have also been studied in a more mathematical context, for example in [EY18, BZN18, EY19, EY20, RY19], and I have benefited greatly from ongoing discussions with Justin Hilburn and Philsang Yoo about these ideas. In particular, the forthcoming papers [HY], [GY], and [HR] will also contain some of their ideas that we follow in the present work.

In terms of the various perspectives we have just discussed, we can summarize an underlying goal of this series of papers as follows:

We develop a dictionary between factorization algebras and quantum field theory in the mixed holomorphic-topological setting, using a synthesis of the chiral and topological variants of factorization algebras; examples of interest are given by factorization compatible sheaf theory constructions, motivated by shifted geometric quantization of spaces of solutions to equations of motion in supersymmetric gauge theories, and using tools from geometric representation theory and derived algebraic geometry.

## 1.2. Motivation: Holomorphic-topological twists of supersymmetric quantum field theories and $\Omega$-backgrounds.
The more broad goal of this series of papers is to use this dictionary to formulate and prove results from a particular family of interconnected predictions of string theory, at the intersections of affine representation theory [KW07, GW09, Gai18, BPRR15], enumerative geometry [AGT10, Nek16, NP17, GR19], low-dimensional topology [GGP16, DGP18, Wit12], and integrable systems [Nek03, NS10, NW10]. These ideas are centred around the six dimensional $\mathcal{N} = (2,0)$ superconformal field theory, sometimes called "theory X", which is an elusive, non-Lagrangian quantum field theory that morally describes fluctuations of M5 branes in M theory (which we remind the reader are geometric objects supported on six dimensional spaces). This theory is considered on a spacetime of the form $C \times M$, for $C$ a smooth algebraic curve and $M$ a smooth four manifold, and this gives rise to natural predictions relating chiral factorization algebras over the curve $C$ with the differential topology of the four manifold $M$, or the enumerative geometry of sheaves in the case $M = S$ is a smooth algebraic surface over $\mathbb{C}$.

As an intermediate step, we establish analogous predictions from three and four dimensional gauge theories following [BDG17, CG18, BLL+15], which correspondingly relate to the representation theory of classical Lie algebras, and of quantizations of symplectic singularities more generally [BDGH16], as well as to more classical aspects of enumerative geometry [BDG+18] and integrable systems [NS09, CWY18].

Similar ideas have been studied extensively in mathematics already in both of the above contexts, often explicitly motivated by the same physics considerations; for example [FG06, Ara18, Bra04, SV13, MO19, BFN14, Neg17, RSYZ19, FG20, BD99, BFN18, BZG17, Cos13] are a few which have been influential in our understanding of this family of ideas, ordered roughly corresponding to the physics references above.

The preceding predictions are nominally phrased in terms of string theory and supersymmetric quantum field theory, which are notoriously difficult to understand and often not yet defined mathematically, but an important common feature of these results from the physical perspective is that they often factor through mixed holomorphic-topological twists of the relevant quantum field theories. As a result, these theories are expected to be amenable to descriptions in terms of algebraic geometry and topology, and in particular the algebras of observables of these theories are expected to correspond to objects in the synthesis of chiral and topological factorization algebras mentioned above that we study in the present work. This is the fundamental reason for the effectiveness of the mathematical tools considered in the present work in the relevant physics context.

However, there is another salient feature of many of the physical constructions and corresponding mathematical interpretations mentioned above, which has not been codified mathematically in our explanation so far: in the seminal paper [Nek03], Nekrasov introduced a construction in quantum field theory called an $\Omega$-*background*, an additional structure on a partially topological quantum field theory which (when it exists) deforms the given theory in a way that enforces rotational equivariance with respect to a fixed $S^1$ action on the underlying spacetime. A primary consequence is that cohomological calculations in $\Omega$-deformed topological field theories are given by the analogous calculations in equivariant cohomology.

Moreover, motivated by the localization theorem in equivariant cohomology, we expect these calculations should in some sense localize to the fixed points of the underlying $S^1$ action, after passing to an appropriate localization $\mathbb{K}[\varepsilon][f^{-1}]$ of the base ring $\mathbb{K}[\varepsilon] := H^{\bullet}_{S^1}(\mathrm{pt}; \mathbb{K})$. In fact, calculations in the algebras of observables of the $\Omega$-deformed theories localize to calculations in (families over $\mathbb{K}[\varepsilon][f^{-1}]$ of) algebras of observables over the fixed point locus. Furthermore, such families of algebras of observables have been observed in [NS09, NW10] to define filtered quantizations of the algebra specialized at the central fibre over $\mathbb{K}[\varepsilon]$.

In the present work, formally Part I of the series, we establish the foundations of the theory of equivariant factorization algebras in the mixed chiral-topological setting. Moreover, in this language we give an account of the equivariant localization and quantization phenomena associated with the $\Omega$-background construction in holomorphic-topological quantum field theory described above. We give an overview of these results presently in Section 1.3.

In the companion paper [But20b], formally Part II of the series, we develop methods for constructing examples of equivariant factorization algebras corresponding to holomorphic-topological twists of supersymmetric gauge theories, and apply the results of Part I in these examples. We give a preview of the results of Part II in Section 1.4 below.

## 1.3. **Overview of Part I.** In this subsection, we give an overview of the results of the present work.

1.3.1. *Overview of Chapter 1.* The first chapter recalls the basics of the theory of algebraic factorization algebras and its relation to vertex algebras, following [BD04] and [FG11]. None of this material is original, but we hope that the relatively concrete summary given here will help make the subject more accessible for the reader. We postpone a detailed overview until Section 3.

1.3.2. *Overview of Chapter 2.* In Chapter 2, we begin by establishing the elementary foundations of the theory of equivariant factorization algebras $A \in \mathrm{Alg}^{\mathrm{fact}}(X)^G$ on algebraic varieties $X$ with the action of a connected algebraic group $G$. There is a key vector space valued invariant of factorization algebras called factorization homology, which defines a functor

$$\int_X : \mathrm{Alg}^{\mathrm{fact}}(X) \to \mathrm{Vect} \ ,$$

analogous to sheaf cohomology. The factorization homology of factorization algebras generalizes the spaces of conformal blocks of vertex algebras and Hochschild homology of associative algebras. In Section 18, we define an equivariant analogue of factorization homology

$$\int_X^G : \mathrm{Alg}^{\mathrm{fact}}(X)^G \to \mathrm{H}_G^\bullet(\mathrm{pt})\text{-Mod} \ ,$$

and in the case $G = (\mathbb{C}^\times)^n$, we prove an equivariant localization theorem in this context:

*Theorem* 1.3.1. Let $A \in \mathrm{Alg}^{\mathrm{fact}}(X)^G$ be an equivariant factorization algebra. The natural map

$$\int_{X^G}^G \iota^! A \xrightarrow{\ \cong\ } \int_X^G A$$

induces an equivalence over the localization $\mathrm{H}_G^\bullet(\mathrm{pt})[f_k^{-1}]$.

Once correctly formulated, the proof of this statement follows straightforwardly from the results of [GKM97]. Nonetheless, it provides an important link between higher dimensional factorization algebras on $X$, which are often subtle to understand algebraically due to their homotopical nature, and lower dimensional factorization algebras on $X^G$, which can be identified with more familiar objects such as associative algebras or vertex algebras.

Next, we carry out a basic study of the algebraic structure of equivariant factorization algebras in the simplest examples, explain relations to algebras over variants of the framed little $n$-disks operad, and give an account in this language of the relationship to deformation quantization predicted in the physics literature, as described above. The latter proceeds as follows:

In general, the restriction of an equivariant factorization algebra

$$\iota^! A \in \mathrm{Alg}^{\mathrm{fact}}(X^G)^G \cong \mathrm{Alg}^{\mathrm{fact}}(X^G)_{/\mathrm{H}_G^\bullet(\mathrm{pt})}$$

defines a family of factorization algebras on $X^G$ parameterized by $\mathrm{H}_G^\bullet(\mathrm{pt})$, since $G$ acts trivially on $X^G$. The case when $G = \mathbb{C}^\times$ corresponds to the usual $\Omega$-background construction, and we show that a $\mathbb{C}^\times$ equivariant factorization algebra $A \in \mathrm{Alg}^{\mathrm{fact}}(X)^{\mathbb{C}^\times}$ induces a family of factorization algebras $\iota^! A \in \mathrm{Alg}^{\mathrm{fact}}(X^{\mathbb{C}^\times})_{/\mathbb{K}[\varepsilon]}$ over $\mathbb{K}[\varepsilon] = \mathrm{H}_{\mathbb{C}^\times}^\bullet(\mathrm{pt})$, which defines a filtered quantization of the central fibre. For simplicity, we consider factorization algebras which are $\mathbb{G}_a$ equivariant, or equivalently topological, along $\mathbb{A}^1$

$$(1.3.1) \qquad \mathrm{Alg}^{\mathrm{fact}}(X \times \mathbb{A}^1)^{\mathbb{G}_a} \cong \mathrm{Alg}_{\mathbb{E}_2}^{\mathrm{fact}}(X) \qquad \text{so that} \qquad \mathrm{Alg}^{\mathrm{fact}}(X \times \mathbb{A}^1)^{\mathbb{G}_a \rtimes \mathbb{G}_m} \cong \mathrm{Alg}_{\mathbb{E}_2^{S^1}}^{\mathrm{fact}}(X)$$

the additional $\mathbb{G}_m$ equivariance is equivalent to a framed, or $S^1$ equivariant, enhancement of the $\mathbb{E}_2$ structure.

In Section 23, we explain an application of the Goresky-Kottwitz-MacPherson Koszul duality result [GKM97] to equivariant operads, in the sense of [SW03]. In this example, it gives an equivalence between $S^1$ equivariant $\mathbb{E}_2$ algebras and algebras over (a two-periodic variant of) the zeroth Beilinson-Drinfeld operad $\mathbb{BD}_0^u$

$$(1.3.2)$$

$$\mathrm{Alg}_{\mathbb{E}_2^{S^1}}(\mathrm{Perf}_{\mathbb{K}}) \cong \mathrm{Alg}_{\mathbb{BD}_0^u}(\mathrm{D}_{\mathrm{fg}}^b(\mathbb{K}[u])) \qquad \text{and similarly} \qquad \mathrm{Alg}_{\mathbb{E}_{n+2}^{S^1}}(\mathrm{Perf}_{\mathbb{K}}) \cong \mathrm{Alg}_{\mathbb{BD}_n^u}(\mathrm{D}_{\mathrm{fg}}^b(\mathbb{K}[u])) \ ,$$

where $\mathbb{K}[u] = H_{S^1}^\bullet(\mathrm{pt})$. In general, algebras over the operad $\mathbb{BD}_n^u$ define graded quantizations of $\mathbb{P}_{n+2}$ algebras to $\mathbb{E}_n$ algebras, over the base ring $\mathbb{K}[u]$. A similar result was announced in [BBZB+20] as to appear in [BZN]. The preceding equivalences also extend to factorization objects, so that in summary we have:

*Theorem* 1.3.2. There are equivalences of categories

$$\mathrm{Alg}^{\mathrm{fact}}(X \times \mathbb{A}^1)^{\mathbb{G}_a \rtimes \mathbb{G}_m} \cong \mathrm{Alg}_{\mathbb{E}_2^{S^1}}^{\mathrm{fact}}(X) \cong \mathrm{Alg}_{\mathbb{BD}_0^u}^{\mathrm{fact}}(X) \ ,$$

such that the latter intertwines the functors of forgetting the $\mathbb{E}_2^{S^1}$ structure and taking the homology $\mathbb{P}_2$ algebra, with restriction to the generic and central fibres $\{u = 1\}$ and $\{u = 0\}$, respectively.

The latter category is equivalent to that of (two-periodic) filtered quantizations of (shifted) Coisson algebras to chiral factorization algebras on $X$, by the chiral Poisson additivity theorem of Rozenblyum. Thus, equivariant factorization algebras $A \in \mathrm{Alg}^{\mathrm{fact}}(X \times \mathbb{A}^1)^{\mathbb{G}_a \rtimes \mathbb{G}_m}$ induce quantizations $\iota^! A \in \mathrm{Alg}^{\mathrm{fact}}(X)_{/\mathbb{K}[\varepsilon]}$.

Finally, we explain the manifestation in this language of the physical principle of equivariant cigar reduction, which plays a central role in our applications of interest in Part II [But20b], as we explain below.

Consider an $S^1$ equivariant factorization $\mathbb{E}_2$ algebra

$$A \in \mathrm{Alg}_{\mathbb{E}_2^{S^1}}^{\mathrm{fact}}(X) \qquad \text{and define} \qquad A_0 := \mathrm{oblv}_{\mathbb{E}_2^{S^1}}^{\mathbb{E}_0} A \ \in \mathrm{Alg}^{\mathrm{fact}}(X) \ .$$

Note that $A_0$ is canonically a module over $A$ in the $\mathbb{E}_2$ sense, so that there is a module structure

$$(1.3.3) \qquad A_0 \in \mathrm{CH}_\bullet(A)\text{-Mod}(\mathrm{Alg}^{\mathrm{fact}}(X)) \qquad \text{or equivalently a map} \qquad \mathrm{CH}_\bullet(A) \to \mathrm{CH}^\bullet(A_0)$$

in the category $\mathrm{Alg}_{\mathbb{E}_1}^{\mathrm{fact}}(X)$ of factorization $\mathbb{E}_1$ algebras. In these terms, we have the following additional structure relating the factorization $\mathbb{E}_2^{S^1}$ algebra and the corresponding factorization $\mathbb{BD}_0^u$ algebra, which has a geometric interpretation in physics as the *equivariant cigar reduction principle*, pictured in Figure 1:

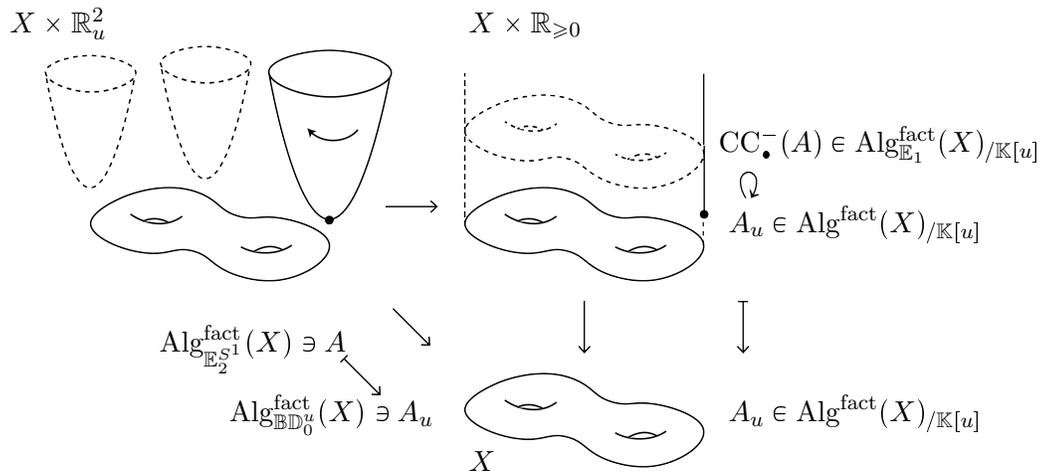

FIGURE 1. The equivariant cigar reduction principle

*Proposition* 1.3.3. The family of factorization algebras $A_u \in \text{Alg}^{\text{fact}}(X)_{/\mathbb{K}[u]}$ underlying the factorization $\mathbb{BD}_0^u$ algebra corresponding to $A \in \text{Alg}^{\text{fact}}_{\mathbb{E}_2^{S^1}}(X)$ under Theorem 1.3.2, admits a canonical module structure

$$A_u \in \text{CC}_\bullet^-(A)\text{-Mod}(\text{Alg}^{\text{fact}}(X)_{/\mathbb{K}[u]}) \qquad \text{such that} \qquad A_u|_{\{u=0\}} = A_0 \in \text{CH}_\bullet(A)\text{-Mod}(\text{Alg}^{\text{fact}}(X)) \ ,$$

its restriction to the central fibre agrees with the module structure of Equation 1.3.3, where $\text{CC}_\bullet^-(A) \in \text{Alg}^{\text{fact}}_{\mathbb{E}_1}(X)_{/\mathbb{K}[u]}$ denotes the negative cyclic chains on $A$, considered as a family of factorization $\mathbb{E}_1$ algebras over $\mathbb{K}[u] = \text{H}^\bullet_{S^1}(\text{pt})$ with central fibre $\text{CH}_\bullet(A) \in \text{Alg}^{\text{fact}}_{\mathbb{E}_1}(X)$ .

In Part II of this series of papers [But20b], the preceding proposition provides an explanation of the relationship between the construction of chiral algebras corresponding to four dimensional $\mathcal{N} = 2$ superconformal field theories in [BLL$^+$15], which we give a mathematical account of in terms of equivariant factorization algebras, and the construction of boundary chiral algebras for (holomorphic-)topological twists of three dimensional $\mathcal{N} = 4$ theories following [CG18], which is the other central topic of Part II.

1.4. **Preview of Part II.** In the companion paper and formal sequel [But20b] to the present work, we develop methods to give geometric constructions of equivariant factorization $\mathbb{E}_n$ algebras corresponding to holomorphic-topological twists of supersymmetric gauge theories equipped with an $\Omega$ background.

The first main example of interest is the three dimensional A model gauge theory which occurs as a topological twist of three dimensional $\mathcal{N} = 4$ supersymmetric gauge theory. A factorization $\mathbb{E}_1$ algebra $\mathcal{A}(G, N) \in \text{Alg}^{\text{fact}}_{\mathbb{E}_1}(C)$ describing the local observables of the three dimensional $A$ model gauge theory on $C \times \mathbb{R}$ with gauge group $G$ and matter representation $T^\vee N$ was introduced in [BFN18]. Moreover, it is explained in *loc. cit.* that this construction also gives a filtered quantization of a graded Poisson algebra, which in good cases describes a quantization of the symplectic singularity which is dual to $T^\vee N$ in the sense of symplectic duality [BLPW14], or three dimensional mirror symmetry [IS96]; the latter was the original motivation for the construction.

The relationship between these results is an example of the equivalence of Theorem 1.3.1:

*Theorem* 1.4.1. [BFN18] For $C = \mathbb{A}^1$, the factorization $\mathbb{E}_1$ algebra

$$\mathcal{A}(G, N) \in \text{Alg}^{\text{fact}}_{\mathbb{E}_1}(\mathbb{A}^1)^{\mathbb{G}_a \rtimes \mathbb{G}_m} \cong \text{Alg}_{\mathbb{E}_2^{S^1}}(\text{Vect}_\mathbb{K}) \cong \text{Alg}_{\mathbb{BD}_1^u}(\text{D}(\mathbb{K}[u])) \ ,$$

admits a canonical $\mathbb{G}_a \rtimes \mathbb{G}_m$ equivariant structure and thus, under the equivalence Theorem 1.3.1, defines a filtered quantization of a (2-shifted) Poisson algebra to an associative (or $\mathbb{E}_1$) algebra.

Concretely, passing to $\mathbb{G}_m$ equivariant (with respect to loop rotation) Borel-Moore homology in the definition of $\mathcal{A}(G, N)$ in [BFN18] gives a quantization of the homology $\mathbb{P}_3$ algebra, which they view as a graded commutative algebra with Poisson bracket of degree $-2$. We also explain an analogous construction of the three dimensional $B$ model in [But20b], which gives a filtered quantization of $T^\vee Y$ itself by this mechanism.

The next main topic is the factorization algebra $\mathcal{D}^{\text{ch}}(Y) \in \text{Alg}^{\text{fact}}(C)$ on $C$ of chiral differential operators on $Y = N/G$, and its relationship to the three dimensional $A$ model above, culminating in a proof in this language of the prediction of Costello-Gaiotto from [CG18] that the three dimensional A model admits a boundary condition with local observables described by $\mathcal{D}^{\text{ch}}(Y)$.

*Remark* 1.4.2. For $Y$ a scheme the construction of $\mathcal{D}^{\text{ch}}(Y)$ requires a trivialization of the determinant gerbe [KV06], which for $Y = N/G$ we identify with a lift of the $G_\mathcal{O}$ action on $\mathcal{D}^{\text{ch}}(N)$ to an action of $\hat{\mathfrak{g}}$ at level $-$Tate. Physically, this corresponds to the requirement that the corresponding four dimensional $\mathcal{N} = 2$ theory is superconformal.

Our formulation of the prediction of interest from [CG18] is, under the hypotheses of the preceding remark, the following:

*Theorem* 1.4.3. [But20b] The chiral differential operators on $Y = N/G$ admits a canonical module structure

$$\mathcal{D}^{\mathrm{ch}}(Y) \ \in \ \mathcal{A}(G, N)\text{-Mod}(\mathrm{Alg}^{\mathrm{fact}}(C))$$

over the factorization $\mathbb{E}_1$ algebra $\mathcal{A}(G, N) \in \mathrm{Alg}^{\mathrm{fact}}_{\mathbb{E}_1}(C)$ on $C$ constructed in [BFN18].

This result corresponds to the statement that the three dimensional $A$ model to $Y$ admits a chiral boundary condition, so that the algebra of local observables $\mathcal{A}(G, N)$ of the three dimensional theory on $C \times \mathbb{R}_{\geqslant 0}$ acts on the chiral algebra $\mathcal{D}^{\mathrm{ch}}(Y)$ of boundary observables on $C$, as pictured on the right.

We also construct a family of factorization $\mathbb{E}_1$ algebras $\mathcal{C}(Y)^{\hbar} \in \mathrm{Alg}^{\mathrm{fact}}_{\mathbb{E}_1}(X)_{/\mathbb{K}[\hbar]}$ over $\mathbb{K}[\hbar]$, with generic fibre $\mathcal{C}(Y)^{\hbar}|_{\{\hbar=1\}} = \mathcal{A}(G, N)$ the factorization $\mathbb{E}_1$ algebra of the three dimensional A model, and central fibre that of the

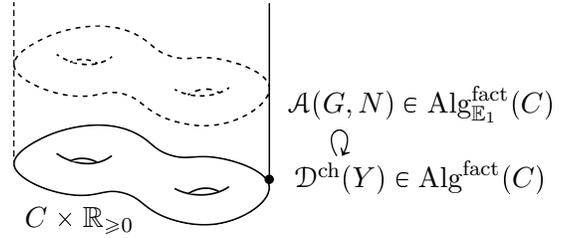

$\mathcal{A}(G, N) \in \mathrm{Alg}^{\mathrm{fact}}_{\mathbb{E}_1}(C)$

$\mathcal{D}^{\mathrm{ch}}(Y) \in \mathrm{Alg}^{\mathrm{fact}}(C)$

$C \times \mathbb{R}_{\geqslant 0}$

holomorphic-B twist. Moreover, we show that the module structure on $\mathcal{D}^{\mathrm{ch}}(Y)$ over $\mathcal{A}(G, N)$ given above extends to one over $\mathcal{C}(Y)^{\hbar}$,

$$(1.4.1) \qquad \qquad \mathcal{D}^{\mathrm{ch}}(Y)_{\hbar} \in \mathcal{C}(Y)^{\hbar}\text{-Mod}(\mathrm{Alg}^{\mathrm{fact}}(C)_{/\mathbb{K}[\hbar]}) \ ,$$

where $\mathcal{D}^{\mathrm{ch}}(Y)_{\hbar} \in \mathrm{Alg}^{\mathrm{fact}}(C)_{/\mathbb{K}[\hbar]}$ is the filtered quantization of chiral differential operators to $Y$.

Further, we use the theory of equivariant factorization algebras developed in the present work to give a mathematical account of the construction of chiral algebras corresponding to 4d $\mathcal{N} = 2$ superconformal gauge theories introduced in [BLL$^+$15]. For $Y$ satisfying the hypotheses of Remark 1.4.2, we have:

*Theorem* 1.4.4. [But20b] There is a natural factorization $\mathbb{E}_2^{S^1}$ algebra $\mathcal{F}(Y) \in \mathrm{Alg}^{\mathrm{fact}}_{\mathbb{E}_2^{S^1}}(C)$ such that

$$\mathcal{F}(Y) \mapsto \mathcal{D}^{\mathrm{ch}}(Y)_u \qquad \text{under the equivalence} \qquad \mathrm{Alg}^{\mathrm{fact}}_{\mathbb{E}_2^{S^1}}(C) \cong \mathrm{Alg}^{\mathrm{fact}}_{\mathbb{BD}_0^u}(C)$$

of Theorem 1.3.2, where $\mathcal{D}^{\mathrm{ch}}(Y)_u \in \mathrm{Alg}^{\mathrm{fact}}_{\mathbb{BD}_0^u}(X)$ is the (two-periodic) filtered quantization of the factorization algebra of chiral differential operators to $Y$.

Further, we explain the relation of this construction with our formulation of the predictions of [CG18], via the equivariant cigar reduction principle described in Proposition 1.3.3. Let $\mathcal{C}(Y)^u \in \mathrm{Alg}^{\mathrm{fact}}_{\mathbb{E}_1}(C)_{/\mathbb{K}[u]}$ denote the two-periodic variant of the family of factorization $\mathbb{E}_1$ described above:

*Theorem* 1.4.5. [But20b] There is an equivalence of families of factorization $\mathbb{E}_1$ algebras on $X$

$$\mathrm{CC}^-_\bullet(\mathcal{F}(Y)) \xrightarrow{\cong} \mathcal{C}(Y)^u \ \in \ \mathrm{Alg}^{\mathrm{fact}}_{\mathbb{E}_1}(C)_{/\mathbb{K}[u]} \ ,$$

such that under the equivalence of Theorem 1.3.2, the module structure of the preceding proposition

$$\mathcal{F}(Y)_u \in \mathrm{CC}^-_\bullet(\mathcal{F}(Y))\text{-Mod}(\mathrm{Alg}^{\mathrm{fact}}(C)_{/\mathbb{K}[u]}) \qquad \text{identifies with} \qquad \mathcal{D}^{\mathrm{ch}}(Y)_u \in \mathcal{C}(Y)^u\text{-Mod}(\mathrm{Alg}^{\mathrm{fact}}(C)_{/\mathbb{K}[u]}) \ ,$$

equipped with the module structure recalled in Equation 1.4.1 above.

The preceding theorem identifies the equivariant $S^1$ reduction $\mathrm{CC}^-_\bullet(\mathcal{F}(Y))$ of the holomorphic-B twist of four dimensional $\mathcal{N}=2$ gauge theory with the deformation $\mathcal{C}(Y)^u$ from the holomorphic-B twist to the A twist of three dimensional $\mathcal{N}=4$ gauge theory. It also identifies the family of boundary conditions for the former induced by the 'cigar tip' as in Proposition 1.3.3, with the family of boundary conditions for the latter whose local observables are the filtered quantization of the factorization algebra of chiral differential operators, as in Equation 1.4.1.

### 1.5. Preview of future directions.

The other main intended application of the results developed in the present work is to establish a variant of the AGT conjecture in the factorization setting, and construct an approximation to the conjectural vertex algebras $\mathrm{VOA}[M_4]$ introduced in [GGP16, DGP18, FG20] via factorization homology. Analogous to the $\mathbb{E}_2^{S^1}$ enhancement of the chiral differential operators claimed in Theorem 1.4.4 above, we expect the following:

*Proposal* 1.5.1. There is a canonical $\mathbb{E}_4^{S^1 \times S^1}$ enhancement of the principal affine W-algebra

$$\tilde{\mathcal{W}}(\mathfrak{gl}_r) \in \mathrm{Alg}^{\mathrm{fact}}(C \times \mathbb{A}^2)^{\mathbb{G}_a^2 \rtimes \mathbb{G}_m^2} \cong \mathrm{Alg}^{\mathrm{fact}}_{\mathbb{E}_4^{S^1 \times S^1}}(C)$$

as a factorization algebra on any smooth algebraic curve $C$.

The structures on $\mathbb{G}_m$ equivariant factorization algebras outlined above and established in the present work are thus expected for $\mathcal{W}(\mathfrak{gl}_r)$ in two distinct, compatible ways; Figure 2 is the analogue of Figure 1 in this setting.

In fact, we expect that $\mathcal{W}(\mathfrak{gl}_r) \in \mathrm{Alg}^{\mathrm{fact}}_{\mathbb{E}_4^{\mathrm{fr}}}(X)$ is framed, and assuming this we can make the following definitions:

*Definition* 1.5.2. Let $M_4$ and $N_3$ be oriented manifolds of dimensions 4 and 3. We define

$$\mathcal{W}(M_4, \mathfrak{gl}_r) = \int_{M_4} \tilde{\mathcal{W}}(\mathfrak{gl}_r) \in \mathrm{Alg}^{\mathrm{fact}}(C) \qquad \text{and} \qquad \mathcal{A}(N_3, \mathfrak{gl}_r) = \int_{N_3} \tilde{\mathcal{W}}(\mathfrak{gl}_r) \in \mathrm{Alg}^{\mathrm{fact}}_{\mathbb{E}_1}(C) \ .$$

By the tensor excision theorem for factorization homology proved in [Lur09a, AF15], these factorization algebras would necessarily satisfy the following 'gluing construction':

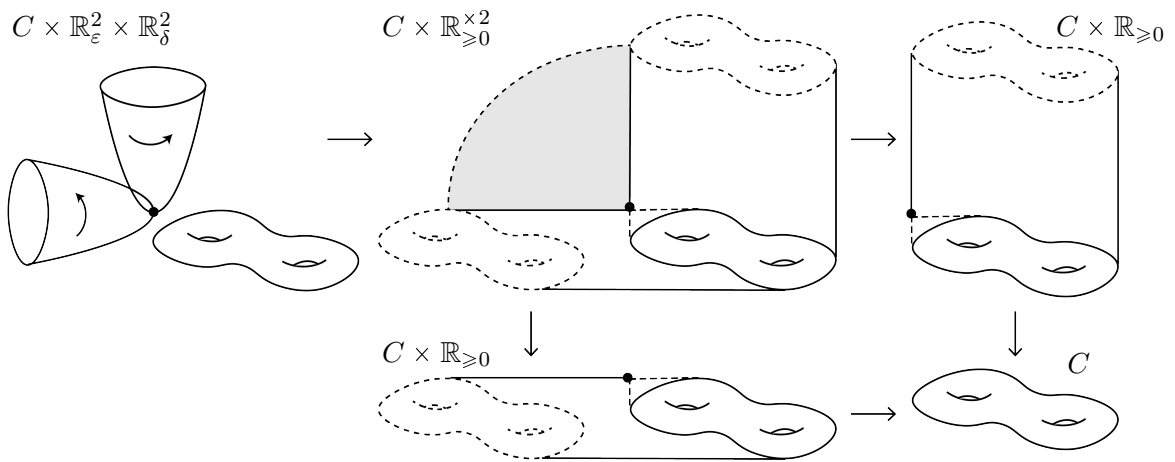

FIGURE 2. The $S^1 \times S^1$ equivariant $\mathbb{E}_4$ enhancement of $\mathcal{W}(\mathfrak{gl}_r)$

*Corollary* 1.5.3. For $M_4 = M_4^- \cup M_4^+$ a collar-gluing presentation with $M_4^- \cap M_4^+ \cong N_3 \times \mathbb{R}$ we have

$$\mathcal{W}(M_4^-, \mathfrak{gl}_r), \mathcal{W}(M_4^+, \mathfrak{gl}_r) \in \mathcal{A}(N_3, \mathfrak{gl}_r)\text{-Mod}(\text{Alg}^{\text{fact}}(C))$$

and moreover there is a canonical equivalence

$$(1.5.1) \qquad \mathcal{W}(M_4^-, \mathfrak{gl}_r) \underset{\mathcal{A}(N_3, \mathfrak{gl}_r)}{\otimes} \mathcal{W}(M_4^+, \mathfrak{gl}_r) \xrightarrow{\cong} \mathcal{W}(M_4, \mathfrak{gl}_r) .$$

*Remark* 1.5.4. The vertex algebras $\mathcal{W}(M_4, \mathfrak{gl}_r)$ are an approximation to the conjectural algebras VOA$[M_4]$ proposed in [GGP16, DGP18, FG20], which are meant to encode rich information about the differential topology of $M_4$ (or the enumerative geometry of sheaves on $S$, in the case that $M_4$ is given by a smooth algebraic surface $S$). In general, we expect there is a map

$$\mathcal{W}(M_4, \mathfrak{gl}_r) \to \text{VOA}[M_4] ,$$

but the latter is typically much larger: in the setting of the preceding corollary, an analogous gluing construction has been conjectured for VOA$[M_4]$, but with the tensor product replaced by a 'vertex algebra extension', of which the tensor product formula of Equation 1.5.1 gives only the first summand. On $M_4 = \mathbb{R}^4$ both constructions give the principal affine $\mathcal{W}$ algebra, but its more exotic gluing construction makes VOA$[M_4]$ a much more interesting invariant in general.

From the perspective of the more common constructions of $\mathcal{W}$ algebras associated to algebraic surfaces in terms of enumerative geometry, this discrepancy corresponds to the fact that the algebra of modes $\mathcal{W}(M, \mathfrak{gl}_r)^{\text{as}}$ is built from Hecke modifications on the moduli of instantons (or torsion free sheaves in the case of an algebraic surface $S$) supported at points, while the conjectural VOA$[M_4]$ should also encode modifications along two dimensional submanifolds (or algebraic curve classes).

For $S$ a smooth, toric algebraic surface, the equivariant localization formula for factorization homology given in Theorem 1.3.1 applied to calculate $\mathcal{W}(S, \mathfrak{gl}_r)$ gives the following:

*Corollary* 1.5.5. The natural map defines an equivalence

$$(1.5.2) \qquad \bigotimes_{s \in S^G} \mathcal{W}(\mathfrak{gl}_r)_{\varepsilon_s, \delta_s} \xrightarrow{\cong} \mathcal{W}(S, \mathfrak{gl}_r)_{\varepsilon, \delta} \qquad \text{as modules over } \mathbb{K}[\varepsilon, \delta][f_k^{-1}].$$

This is illustrated in Figure 3 below in the case that $S = \mathbb{P}^2$ with the usual $(\mathbb{C}^\times)^2$ action.

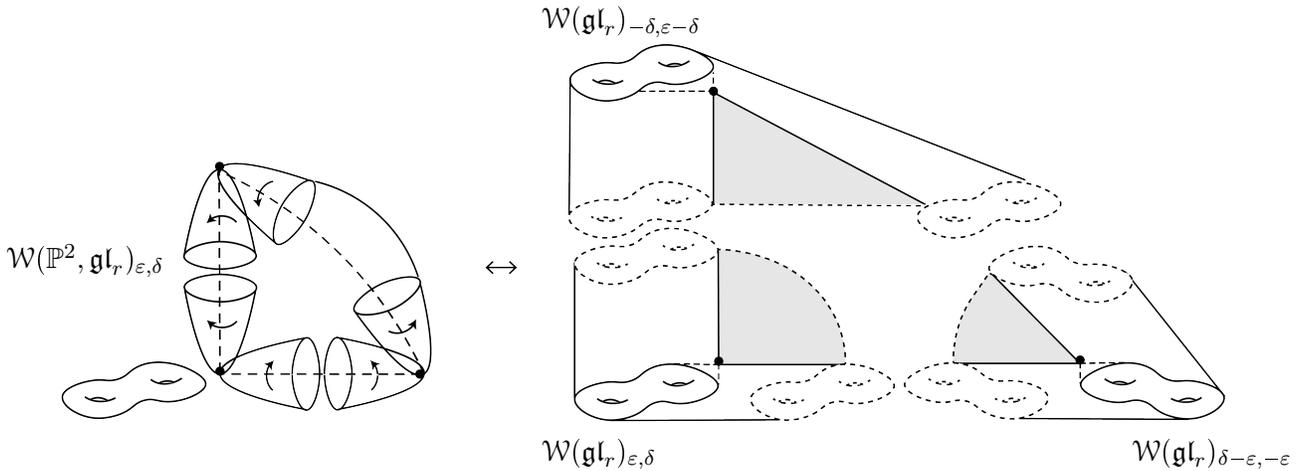

FIGURE 3. Equivariant localization formula for $\mathcal{W}(\mathbb{P}^2, \mathfrak{gl}_r)$

Moreover, toric surfaces $S$ are equivariantly formal [GKM97], and for such spaces there is a refinement of the classical localization theorem [CS74, GKM97]. The goal of our in progress project is to use an analogous enhancement in factorization homology to calculate $\mathcal{W}(S, \mathfrak{gl}_r)$, refining the description of Equation 1.5.2 above. In fact, we also have a proposal for a similar construction (and partial generalization) of the full VOA$[M_4]$ in the case that $M_4 = S$ is a smooth toric algebraic surface occuring as a reduced divisor in a toric Calabi-Yau threefold, though this is the topic of a separate project.

### 1.6. Acknowledgments.

I am especially thankful to Kevin Costello, for first introducing me to so many of the ideas underlying this series of papers, as well as to Sam Raskin, for patient explanations about various more technical aspects. I would also like to thank David Ben-Zvi, Alexander Braverman, Davide Gaiotto, Justin Hilburn, Surya Raghavendran, Pavel Safronov, Brian Williams, and Philsang Yoo for useful discussions.

## 2. Conventions

### 2.1. General conventions and notation.

The required notations are introduced, together with the relevant definitions, in the many appendices to the present work. In this subsection, we briefly recall some of the most commonly used notation.

We fix a base field $\mathbb{K}$ of characteristic zero, which we occasionally assume to be given by the complex numbers $\mathbb{K} = \mathbb{C}$. We typically use $X$ to denote a smooth algebraic variety over $\mathbb{K}$ of dimension $d_X$. In such cases, we let $\mathcal{O}_X$ denote the sheaf of regular functions, $\mathcal{D}_X$ the sheaf of differential operators, and $\omega_X = \Omega_X^{d_X}[d_X]$ the dualizing sheaf of $X$

Moreover, we let $D(X)$ denote the category of $D$ modules on $X$; the theory of such is reviewed in Appendix A.2. In particular, given a map $f : X \to Y$ of algebraic varieties, there are induced functors, denoted by

$$f_* : D(X) \to D(Y) \qquad \text{and} \qquad f^! : D(Y) \to D(X) \ ,$$

given by the pullback and pushforward of $D$ modules.

We let $\mathrm{D}(A)$ denote the derived category of modules over an algebra $A$, and similarly $\mathrm{D}^b(A)$, $\mathrm{D}^+(A)$, and $\mathrm{D}_{\mathrm{fg}}(A)$ its bounded, bounded above, and finitely generated variants.

We let $\mathrm{Op}(\mathcal{C})$ denote the category of (by default symmetric) operads internal to a symmetric monoidal category $\mathcal{C}$; the theory of such is reviewed in Appendix C.1. For $\mathcal{O}, \mathcal{O}' \in \mathrm{Op}(\mathcal{C})$, we let $\mathrm{Alg}_{\mathcal{O}}(\mathcal{O}') = \mathrm{Hom}_{\mathrm{Op}(\mathcal{C})}(\mathcal{O}, \mathcal{O}')$ denote the category of $\mathcal{O}$ algebras internal to the operad $\mathcal{O}'$. We let Ass, Comm, Lie, and $\mathbb{P}_1$ denote the associative, commutative, Lie, and Poisson operads. Further, we let $\mathbb{E}_n$ denote the little $n$-disks operad, $\mathbb{P}_n$ the corresponding shifted Poisson operad, and $\mathbb{BD}_n$ the Beilinson-Drinfeld operad, as recalled in appendices C.4, C.5, and C.6, respectively.

Throughout, whenever possible we work in the framework of cocomplete DG categories and stable infinity categories, following [Lur12] for example.

*Warning* 2.1.1. We often fall short of a complete explanation at the level of homotopical precision typical in the study of such categories. This series of papers is ultimately about concrete *objects*, like vector spaces or sheaves, and the algebraic structures they carry, and we carefully establish our results along these lines. However, for the purposes of exposition we occasionally use more categorical language in situations where we have not established a precise enough setting for their careful interpretation. We hope that this will not be a cause of confusion for the reader, and we have included many similar warnings throughout the text.

## 2.2. Notation around partitions and diagonal embeddings.

We record here various notations related to certain index categories of finite sets, which are used throughout the text.

Let fSet denote the category of (possibly empty) finite sets with arbitrary maps of sets $\pi : I \to J$, and fSet$^{\mathrm{surj}}$ denote the category of non-empty finite sets with surjective maps $\pi : I \twoheadrightarrow J$.

*Remark* 2.2.1. The category fSet$^{\mathrm{surj}}$ parameterizes a diagram whose colimit defines the moduli space of non-empty finite subsets Ran$_X$ of a space $X$, while the category fSet is used analogously to describe the moduli space of possibly empty finite subsets Ran$_{X,\mathrm{un}}$ of a space $X$; see definitions 4.1.2 and 4.3.1.

These conventions are also be used in a closely related way in the definition of operads, recalled in Appendix C.1; the relation is for example apparent in the results of Section 20.

We identify a surjective map $\pi : I \twoheadrightarrow J$, a morphism in fSet$^{\mathrm{surj}}$, with a $J$-coloured partition of $I$, given by

$$ I = \bigsqcup_{j \in J} I_j \qquad \text{where} \qquad I_j := \pi^{-1}(j) \ . $$

A general map $\pi : I \to J$, a morphism in fSet, is similarly equivalent to a $J$-coloured partition of $I$, given by

$$ I = \bigsqcup_{j \in J} I_j = \bigsqcup_{j \in \mathrm{im}(\pi)} I_j \sqcup \bigsqcup_{j \in I_\pi} \emptyset_j \qquad \text{where} \qquad I_j := \pi^{-1}(j) \quad \text{and} \quad I_\pi = J \backslash \mathrm{im}(\pi) \ . $$

Here the subsets of $I$ corresponding to certain colours $j \in J$ are allowed to be empty, while still recorded in this combinatorial model. We also encode this data equivalently by the induced surjection

$$ \alpha_\pi := \pi \times \mathbb{1}_{I_\pi} : I \sqcup I_\pi \twoheadrightarrow \mathrm{im}(\pi) \sqcup I_\pi = J \ . $$

For each $\pi : I \to J$, we define a corresponding diagonal map

$$ \Delta(\pi) : X^J \to X^I \qquad (x_j)_{j \in J} \mapsto (x_{\pi(i)})_{i \in I} \ . $$

If $\pi : I \twoheadrightarrow J$ is a surjection, then $\Delta(\pi) : X^J \hookrightarrow X^I$ is a closed embedding, while if $\pi : I \hookrightarrow J$ is injective, $\Delta(\pi) : X^J \twoheadrightarrow X^I$ is smooth over $X^I$. In general, $\Delta(\pi)$ is smooth over its image with fibre $X^{I_\pi}$. In the case $\pi : I \to \{\mathrm{pt}\}$, we use the notation $\Delta^{(I)} = \Delta(\pi) : X \to X^I$ for the small diagonal embedding.

For each surjection $\pi : I \twoheadrightarrow J$ we define a corresponding diagonal complement

$$ j(\pi) : U(\pi) \hookrightarrow X^I \qquad U(\pi) = \{(x_i)_{i \in I} | x_i \neq x_j \text{ if } \pi(x_i) \neq \pi(x_j)\} \ . $$

For example, for $\pi : I \to \{\mathrm{pt}\}$ this gives $U(\pi) = X^I$, while for $\pi = \mathbb{1}_I : I \to I$ this gives $U(\pi)$ equal to the complement of all partial diagonals in $X^I$. Note that $U(\pi)$ is not in general complementary to the image of $\Delta(\pi)$.

In the case $\pi = \mathbb{1}_I : I \to I$, we use the notation $U^{(I)} = U(\mathbb{1}_I)$ and $j^{(I)} = j(\mathbb{1}_I) : U^{(I)} \hookrightarrow X^I$ for the complement to the union of all partial diagonals in $X^I$.

## 2.3. References to the companion paper.

As we have mentioned, this paper is the first among at least two papers in this series, and we will systematically refer to propositions, definitions, etc. from each of the companion papers in the other. In this paper, references to the present text are given by red hyperlinks, such as 4.2.3 which refers to the remark of that enumeration. References to Part II are given by turquoise hyperlinks (which link to the companion pdf if both files are in the same folder) and their enumeration is prefaced by a II, such as II-4.2.3.

# Chapter 1

# Factorization Algebras, Chiral Algebras, and Vertex Algebras

### 3. Overview of Chapter 1

In this chapter, we review the theory of factorization algebras in the setting of algebraic geometry, following [BD04, FG11, Ras15a], and explain the relationship of these objects to vertex algebras. We also give an exposition of a few more specialized topics in the theory of chiral factorization algebras which will be required in this series of papers. None of this material is original, but we hope that the relatively concrete summary given here will help make the subject more accessible for the reader.

3.1. **General Overview: Local observables, the Ran space, and factorization algebras.** Recall from Section 1.1 that factorization algebras are algebraic objects defined over algebraic varieties or manifolds, vaguely analogously to sheaves on them, which describe the algebras of observables of quantum field theories defined on these spaces. Factorization algebras can be defined in the language of algebraic geometry or topology, and the resulting objects generalize vertex algebras and usual (associative) algebras, respectively. Indeed, these can be interpreted as the algebras of observables of chiral conformal field theories in real dimension two, or topological field theories in real dimension one, respectively.

Higher dimensional factorization algebras on algebraic varieties describe the local observables of holomorphic field theories, which by definition generalize the holomorphic behaviour of observables in chiral conformal field theories in two real dimensions. Similarly, higher dimensional topological factorization algebras generalize algebras over the little n-disks operad $\mathbb{E}_n$, which are reviewed in Appendix C.4. The results of this paper emphasize applications using factorization algebras in the setting of algebraic geometry, which we will sometimes call the *chiral* setting. We now give a heuristic overview of the definition of factorization algebras in the context of algebraic geometry:

The starting point for the definition of factorization algebras on an algebraic variety $X$ is a space $\mathrm{Ran}_X$ called the Ran space of $X$, which is by definition the moduli space of non-empty, finite subsets $\{x_i\} \subset X$, which we record heuristically as

$$\mathrm{Ran}_X = \{ \ \{x_i\} \subset X \ \} \ .$$

A factorization algebra is meant to describe the algebra of local operators of a quantum field theory, and this data can naturally be interpreted as defining a sheaf over the space $\mathrm{Ran}_X$: to each point $\mathbf{x} = \{x_i\} \in \mathrm{Ran}_X$ corresponding to a finite subset $\{x_i\} \subset X$, assign the vector space

$$\mathcal{A}_{\{x_i\}} := \{ \ \text{observables local to the collection of points } \{x_i\} \ \} \ .$$

In particular, on the copy of $X \hookrightarrow \mathrm{Ran}_X$ corresponding the singleton subsets $\{x\} \subset X$, the vector spaces $\mathcal{A}_x := \mathcal{A}_{\{x\}}$ are just the usual family of vector spaces describing the spaces of local operators over each point $x \in X$, which defines the vector space underlying the corresponding vertex algebra in the case that $X$ is of complex dimension one.

The fact that these vector spaces glue together to define a sheaf over $\mathrm{Ran}_X$ encodes the condition that given two distict points $x_0, x_1 \in X$, as $x_1$ approaches $x_0$ there should be a gluing map

$$\mathcal{A}_{\{x_0, x_1\}} \rightsquigarrow \mathcal{A}_{x_0} \ .$$

As we will see, the correct data to prescribe this sheaf over $\mathrm{Ran}_X$ is not a map of the underlying vector spaces, but rather a meromorphic family of products with poles concentrated along the diagonal $x_0 = x_1$, as expected from operator product expansions in two dimensional chiral conformal field theories.

There is one more natural condition to impose on this sheaf $\mathcal{A}$ on $\mathrm{Ran}_X$, called factorizability, which corresponds to the notion of locality in quantum field theory: for distinct points $x_0 \neq x_1$ away from the diagonal, the space of observables local to the set $\{x_0, x_1\}$ should be equivalent to the tensor product

$$\mathcal{A}_{\{x_0,x_1\}} \cong \mathcal{A}_{x_0} \otimes \mathcal{A}_{x_1} \ ,$$

of the space of observables local to $x_0$ with the space of observables local to $x_1$. In the limit as $x_1$ approaches $x_0$, we can combine this identification with the gluing map above to obtain a map

$$\mathcal{Y}_{x_0} : \mathcal{A}_{x_0} \otimes \mathcal{A}_{x_0} \rightsquigarrow \mathcal{A}_{x_0}$$

which describes the operator product expansion of observables in an infinitesimal neighbourhood of each point $x_0 \in X$. In the case that $X$ is of complex dimension one, this will recover the usual vertex operator structure map of a vertex algebra.

The first several sections of this Chapter are devoted to filling in the mathematical content of this heuristic description, carefully identifying the resulting algebraic structures, and comparing them with the more concrete structures in the theory of vertex algebras. The latter sections treat more specialized topics which will be required in the present work and its sequel Part II [But20b].

3.2. **Summary.** We now give a summary of each of the sections in this chapter. None of the results of this chapter are new: we follow [BD04] and [FG11] throughout.

3.2.1. *The* Ran *space and the category* $D(\mathrm{Ran}_X)$. In Section 4, we give a more mathematically detailed description of the space $\mathrm{Ran}_X$, and define the category of $D$ modules $D(\mathrm{Ran}_X)$ on $\mathrm{Ran}_X$, which is the appropriate variant of sheaf on $\mathrm{Ran}_X$ to define factorization algebras corresponding to usual holomorphic field theory.

3.2.2. *Monoidal structures on* $D(\mathrm{Ran}_X)$. In Section 5, we recall several monoidal structures on the category of $D$ modules on $\mathrm{Ran}_X$ and their basic properties.

3.2.3. *Factorization algebras.* In Section 6, we recall the definition of factorization algebras in terms of the monoidal structures of the preceding section.

3.2.4. *Chiral algebras.* In Section 7, we recall the definition of chiral algebras in terms of the monoidal structures of Section 5. These are equivalent to factorization algebras, as we recall in Section 14 below.

3.2.5. *Chiral algebras and vertex algebras.* In Section 8 we recall the equivalence between weakly $\mathbb{G}_a$ equivariant chiral algebras on $\mathbb{A}^1$ and vertex algebras.

3.2.6. *OPE algebras.* In Section 9 we recall the notion of OPE algebras, which gives the most direct generalization of the operator product expansion map of vertex algebras to a global operation defined over algebraic curves. Further, we recall that these are equivalent to chiral algebras.

3.2.7. *Chiral algebras and topological associative algebras.* In Section 11, we recall the 'algebra of modes' construction, which from a chiral algebra $A$ defines a topological associative algebra $A_x^{\mathrm{as}}$ at each point $x \in X$.

3.2.8. *Lie$^*$, Comm$^!$, and Coisson algebras.* In Section 10 we recall the analogues of Lie algebras, commutative algebras, and Poisson algebras in the chiral setting, under the analogy that chiral algebras correspond to usual associative algebras. Further, in the weakly $\mathbb{G}_a$ equivariant case on $\mathbb{A}^1$, we identify these with vertex Lie algebras, commuative vertex algebras, and vertex Poisson algebras, respectively.

3.2.9. *Chiral enveloping algebras.* In Section 12, we recall the analogue of the universal enveloping algebra construction in the chiral setting, which for a Lie$^*$ algebra $L$ constructs a chiral algebra $\mathcal{U}^{\mathrm{ch}}(L)$ which satisfies the analogous universal property. We also explain the translation of this construction to the vertex algebra setting. We also recall a variant of the construction that is twisted by an appropriate notion of central extension of the Lie$^*$ algebra $L$, and its analogue in the vertex algebra setting. This allows us to recover many familiar examples such as the affine Kac-Moody and Virasoro vertex algebras.

3.2.10. *BRST reduction of chiral algebras and vertex algebras.* In Section 13, we outlined the theory of BRST reduction of chiral algebras following [BD04], and relate it to the notion of BRST reduction of vertex algebras described in [FBZ04].

3.2.11. *Francis-Gaitsgory chiral Koszul duality.* In Section 14, we give a concrete explanation of the equivalence of chiral algebras and factorization algebras in complex dimension one, and outline the general proof of this fact given in [FG11].

## 4. The Ran space and the category $D(\mathrm{Ran}_X)$

Let $X \in \mathrm{Sch}_{\mathbb{K}}$ be a scheme over $\mathbb{K}$; we will primarily be interested in the case when $X$ is a smooth, finite type variety over $\mathbb{K} = \mathbb{C}$.

4.1. **The Ran space.** Following the heuristic descriptions in Subsection 3.1, we wish to define the moduli space of non-empty finite subsets of $X$, called the Ran space $\mathrm{Ran}_X$ of $X$. The basic idea for constructing $\mathrm{Ran}_X$ is that it should be glued together from various powers $X^I$ of $X$ as follows: it should contain a copy of $X$ corresponding to the space of one point subsets of $X$, glued diagonally into a copy of $\mathrm{Sym}^2 X$ corresponding to the space of two point subsets, which is further glued diagonally into a quotient of $\mathrm{Sym}^3 X$ (which identifies e.g. $(x, x, y)$ with $(x, y, y)$, as both are representatives of $\{x, y\} \subset X$) corresponding to the space of three point subsets, and so on. We summarize this gluing procedure in the putative definition

$$(4.1.1) \qquad \mathrm{Ran}_X = \mathrm{colim} \left[ X \xrightarrow{\Delta} X^2 \underset{S_2}{\overset{}{\rightrightarrows}} X^3 \underset{S_3}{\overset{}{\Rrightarrow}} \ldots \right] .$$

We now formalize this construction:

Recall from 2.2 that $\mathrm{fSet}^{\mathrm{surj}}$ denotes the category with objects given by non-empty finite sets $I$ and morphisms given by surjections $\pi : I \twoheadrightarrow J$. There is a natural $\mathrm{fSet}^{\mathrm{surj,op}}$ diagram in $\mathrm{Sch}_{\mathrm{ft}}$, given by

$$(4.1.2) \qquad X^\bullet : \mathrm{fSet}^{\mathrm{surj,op}} \to \mathrm{Sch} \qquad \text{by} \qquad \begin{cases} I \mapsto X^I \\ [\pi : I \twoheadrightarrow J] \mapsto [\Delta(\pi) : X^J \hookrightarrow X^I] \end{cases} ,$$

where $\Delta(\pi) : X^J \to X^I$ is the corresponding diagonal embedding, defined in 2.2.

*Remark* 4.1.1. The bijections in fSet$^{\text{surj,op}}$ induce an $S_n$ automorphism group at each object $I$ with $|I| = n$.

We now state the formal definition of the space Ran$_X$. We encourage the reader not familiar with the technical notions mentioned to ignore them; as we explain in Remark 4.1.4 below, our exposition in the remainder of the text is formally independent of the actual definition of Ran$_X$.

*Definition* 4.1.2. The Ran space of $X$ is the pseudo indscheme presented by

$$\text{Ran}_X = \operatorname*{colim}_{I \in \text{fSet}^{\text{surj,op}}} X^I = \text{colim} \left[ \dots \leftarrow X^I \xleftarrow{\Delta(\pi)} X^J \leftarrow \dots \right],$$

the colimit of the diagram in Equation 4.1.2, evaluated after composing with the Yoneda embedding Sch $\hookrightarrow$ PreStk.

*Remark* 4.1.3. A pseudo indscheme is a prestack presented as a (not necessarily filtered) colimit of schemes under closed embeddings. The maps in the diagram of Equation 4.1.2 are evidently closed embeddings, but the index category fSet$^{\text{surj,op}}$ is not filtered, so Ran$_X$ does not define an ind scheme in the usual sense.

*Remark* 4.1.4. The category $D(\text{Ran}_X)$ will be defined in 4.2.1 below, without explicit reference to the preceeding definition of Ran$_X$; see remarks 4.2.2 and 4.2.4. In general, the foundational definitions and results in this paper will be stated in terms of the category $D(\text{Ran}_X)$, in a way that is similarly formally independent of the definition of the underlying space Ran$_X$.

*Remark* 4.1.5. There exist canonical maps $\Delta^I : X^I \to \text{Ran}_X$ with images $\text{Ran}_{X, \leqslant n}$ for $|I| = n$ that define a filtration of Ran$_X$ by finite dimensional subschemes, and corresponding stratification of Ran$_X$ by subschemes $\text{Ran}_{X,n}$ parameterizing subsets of $X$ of cardinality $n$.

## 4.2. The category $D(\mathbf{Ran}_X)$ of $D$ modules on $\mathbf{Ran}_X$.

Following further the discussion in 3.1, our initial object of interest is the category of $D$ modules on the Ran space Ran$_X$ of $X$. In terms of the heuristic summarized in Equation 4.1.1, a $D$ module $\mathcal{A} \in D(\text{Ran}_X)$ on Ran$_X$ should be specified by:

- a $D$ module $A \in D(X)$ on $X$;
- a $D$ module $A_2 \in D(X^2)$ on $X^2$, an $S_2$ equivariant structure on $A_2$, and an isomorphism $\Delta^! A_2 \cong A$;
- a $D$ module $A_3 \in D(X^3)$ on $X^3$, an $S_3$ equivariant structure on $A_3$, and isomorphisms $\Delta(\pi)^! A_3 \cong A_2$ for each $\Delta(\pi) : X^2 \hookrightarrow X^3$ corresponding to $\pi : \{1, 2, 3\} \twoheadrightarrow \{1, 2\}$;

and so on. Thus, following [BD04] and [FG11], we make the following definition:

*Definition* 4.2.1. An object $\mathcal{A} \in D(\text{Ran}_X)$ is an assignment

$$I \mapsto A_I \in D(X^I) \qquad\qquad [\pi : I \twoheadrightarrow J] \mapsto [\Delta(\pi)^! A_I \xrightarrow{\cong} A_J]$$

defined for each finite set $I \in$ fSet and surjection $\pi : I \twoheadrightarrow J$.

A morphism $f : \mathcal{A} \to \mathcal{B}$ between $\mathcal{A}, \mathcal{B} \in D(\text{Ran}_X)$ is given by an assignment

$$I \mapsto [f_I : A_I \to B_I] \qquad [\pi : I \twoheadrightarrow J] \mapsto \begin{array}{ccc} \Delta(\pi)^! A_I & \longrightarrow & A_J \\ {\scriptstyle \Delta(\pi)^!(f_I)} \Big\downarrow & {\scriptstyle \sim} & \Big\downarrow {\scriptstyle f_J} \\ \Delta(\pi)^! B_I & \longrightarrow & B_J \end{array}$$

defined for each finite set $I \in$ fSet and surjection $\pi : I \twoheadrightarrow J$.

An object $\mathcal{A} \in D(\mathrm{Ran}_X)$ is called coherent, regular holonomic, ... if $A_I \in D(X^I)$ is such for each $I \in \mathrm{fSet}$.

*Remark* 4.2.2. The preceding definition is evidently independent of the definition of $\mathrm{Ran}_X$ in 4.1.2, in keeping with Remark 4.1.4.

*Remark* 4.2.3. The preceding definition of the category of $D$ modules on $\mathrm{Ran}_X$ can be formalized

$$D(\mathrm{Ran}_X) = \lim_{I \in \mathrm{fSet}^{\mathrm{surj}}} D^!(X^I) = \lim \left[ \ldots \to D(X^I) \xrightarrow{\Delta(\pi)^!} D(X^J) \to \ldots \right] ,$$

where $D^! : \mathrm{Sch}^{\mathrm{op}} \to \mathrm{DGCat}_{\mathrm{cont}}$ is as defined in Appendix II-B.6.

*Remark* 4.2.4. In [Gai12] and [Ras15a], the pseudo indscheme $\mathrm{Ran}_X$ is defined as above, and the definition of its $D$ module category is given in terms of a general definition of the category of $D$ modules on such spaces. A proof of the equivalence of this approach with the above is given in Section 8 of [Ras15a], for example.

*Remark* 4.2.5. There are canonical functors

$$\Delta^{I,!} : D(\mathrm{Ran}_X) \to D(X^I) \qquad \text{and} \qquad \Delta^I_* : D(X^I) \to D(\mathrm{Ran}_X) ,$$

which correspond heuristically to pullback and pushforward along the diagonal embedding $\Delta^I : X^I \hookrightarrow \mathrm{Ran}_X$, as the notation suggests.

*Remark* 4.2.6. The functors of 4.2.5 for $|I| = 1$ are denoted $\Delta^{\mathrm{main}} := \Delta^I$. In this case, the functors $\Delta^{\mathrm{main}}_*, \Delta^{\mathrm{main},!}$ define inverse equivalences $D(X) \cong D(\mathrm{Ran}_X)_X$, by Kashiwara's lemma, where $D(\mathrm{Ran}_X)_X$ denotes the full subcategory on $D$ modules supported on the main diagonal.

*Example* 4.2.7. There is a canonical object $\omega_{\mathrm{Ran}_X} \in D(\mathrm{Ran}_X)$ defined by the assignment

$$I \mapsto \left( \omega_{X^I} \in D(X^I) \right) \qquad [\pi : I \twoheadrightarrow J] \mapsto \left[ \Delta(\pi)^! \omega_{X^I} \xrightarrow{\cong} \omega_{X^J} \right] .$$

4.3. **Unital $D$ modules on $\mathrm{Ran}_X$.** We introduce a 'unital' variant of the above notion of $D$ module on $\mathrm{Ran}_X$. The additional data of the unital structure will be used to define the notion of unital factorization algebra, and this data will correspond to the the unit of a vertex algebra under the equivalence of Section 8, and similarly the unital structure on the $\mathbb{E}_n$ operad under the equivalence of Section 20.

*Definition* 4.3.1. An object $\mathcal{A} \in D(\mathrm{Ran}_{X,\mathrm{un}})$ is an assignment

$$I \mapsto A_I \in D(X^I) \qquad [\pi : I \to J] \mapsto [\Delta(\pi)^! A_I \to A_J]$$

defined for each (possibly empty) finite set $I \in \mathrm{fSet}$ and map $\pi : I \to J$, such that the maps corresponding to surjections $\pi : I \twoheadrightarrow J$ are isomorphisms.

A morphism $f : \mathcal{A} \to \mathcal{B}$ between $\mathcal{A}, \mathcal{B} \in D(\mathrm{Ran}_X)$ is given by an assignment

$$I \mapsto [f_I : A_I \to B_I] \qquad [\pi : I \to J] \mapsto \begin{array}{ccc} \Delta(\pi)^! A_I & \longrightarrow & A_J \\ {\scriptstyle \Delta(\pi)^!(f_I)} \downarrow \ \ \Big/ {\scriptstyle \sim} & & \downarrow {\scriptstyle f_J} \\ \Delta(\pi)^! B_I & \longrightarrow & B_J \end{array}$$

defined for each $I \in \mathrm{fSet}$ and $\pi : I \to J$.

*Remark* 4.3.2. For a non surjective map $\pi : I \to J$, the factorization

$$I \xrightarrow{\bar{\pi}} \operatorname{im}(\pi) \hookrightarrow J \qquad \text{induces an identification} \qquad \Delta(\pi)^! A_I = \Delta(\bar{\pi})^! A_I \boxtimes \omega_{X^{I_\pi}} \ .$$

In particular, for $\pi : I \hookrightarrow J$ injective, or further in particular for $\pi : \varnothing \to J$, the assignment of the preceding definition is required to give maps

$$(4.3.1) \qquad\qquad A_I \boxtimes \omega_{X^{I_\pi}} \to A_J \qquad \text{and} \qquad \omega_{X^J} \to A_J \ ,$$

respectively.

*Remark* 4.3.3. The structure maps which are not necessarily invertible are evidently not interpretable as gluing data for a sheaf on a usual space. However, we will describe an analogous geometric interpretation of $D(\operatorname{Ran}_{X,\operatorname{un}})$ as the category of $D$ modules on a *lax* prestack $\operatorname{Ran}_{X,\operatorname{un}}$ in Section II-2.2, following [Ras15a].

*Remark* 4.3.4. There is evidently a natural forgetful functor $D(\operatorname{Ran}_{X,\operatorname{un}}) \to D(\operatorname{Ran}_X)$ defined by restricting the above assignments along the inclusion $\mathrm{fSet}^{\mathrm{surj}} \hookrightarrow \mathrm{fSet}$.

*Example* 4.3.5. There is a canonical object $\omega_{\operatorname{Ran}_{X,\operatorname{un}}} \in D(\operatorname{Ran}_{X,\operatorname{un}})$ defined by the assignment

$$I \mapsto \big( \omega_{X^I} \in D(X^I) \big) \qquad [\pi : I \to J] \mapsto \Big[ \Delta(\pi)^! \omega_{X^I} \xrightarrow{\cong} \omega_{X^J} \Big] \ .$$

## 5. Monoidal structures on $D(\operatorname{Ran}_X)$ and operad structures on $D(X)$

### 5.1. The $\otimes^!$ monoidal structure.
The $\otimes^!$ monoidal structure on $D(\operatorname{Ran}_X)$ is in principle just the symmetric monoidal structure defined on $D(Y)$ for any scheme $Y$ as Definition A.2.6; in keeping with Remark 4.1.4, we make the following formal definition:

*Definition* 5.1.1. The $\otimes^!$ monoidal structure on $D(\operatorname{Ran}_X)$ is defined by

$$\otimes^! : D(\operatorname{Ran}_X)^{\times 2} \to D(\operatorname{Ran}_X) \qquad (\mathcal{A}, \mathcal{B}) \mapsto (\mathcal{A} \otimes^! \mathcal{B})_I = A_I \otimes^! B_I \quad \in \quad D(X^I) \ ,$$

$$\otimes^!_{j \in J} : D(\operatorname{Ran}_X)^J \to D(\operatorname{Ran}_X) \qquad (\mathcal{A}^j) \mapsto (\otimes^!_{j \in J} \mathcal{A}^j)_I = \otimes^!_{j \in J} A_I^j \quad \in \quad D(X^I) \ .$$

together with the coherence isomorphisms given by the analogous product of those for $\mathcal{A}$ and $\mathcal{B}$.

*Proposition* 5.1.2. The functors $\Delta_*^I : D(X^I) \to D(\operatorname{Ran}_X)$ and $\Delta^{I,!} : D(\operatorname{Ran}_X) \to D(X^I)$ of 4.2.5 are symmetric monoidal with respect to the $\otimes^!$ monoidal structure on each category.

*Corollary* 5.1.3. The functors $\Delta_*^{\mathrm{main}} : D(X) \to D(\operatorname{Ran}_X)$ and $\Delta^{\mathrm{main},!} : D(\operatorname{Ran}_X) \to D(X)$ define a natural symmetric monoidal equivalence $D(X)^! \cong D(\operatorname{Ran}_X)_X^!$.

### 5.2. The $\otimes^*$ monoidal structure.
The space $\operatorname{Ran}_X$ is naturally a commutative monoid object under the operation $\cup : \operatorname{Ran}_X^{\times 2} \to \operatorname{Ran}_X$ of union of finite subsets; this structure is discussed more formally in Section II-2.3, but only motivational in this section. The union map equips $D(\operatorname{Ran}_X)$ with an additional monoidal structure, as follows:

The $\otimes^*$ monoidal structure on $D(\operatorname{Ran}_X)$ is defined geometrically by push forward along the union map of the exterior product

$$\cup_* \circ \boxtimes : D(\operatorname{Ran}_X)^{\times 2} \to D(\operatorname{Ran}_X) \ ,$$

and similarly for higher arity monoidal products. The union map is presented as the colimit of maps $\operatorname{Ran}_{X, \leqslant n_1} \times \operatorname{Ran}_{X, \leqslant n_2} \to \operatorname{Ran}_{X, \leqslant n_1 + n_2}$ on bounded cardinality subsets, induced by the maps $X^I \times X^J \to X^{I \sqcup J}$. In keeping with Remark 4.1.4, we make the following formal definition:

*Definition* 5.2.1. The monoidal structure $\otimes^* : D(\mathrm{Ran}_X)^{\times 2} \to D(\mathrm{Ran}_X)$ and its higher arity components $\otimes^* : D(\mathrm{Ran}_X)^J \to D(\mathrm{Ran}_X)$ are presented by the functors

$$D(X^{I_1}) \times D(X^{I_2}) \to D(X^{I_1 \sqcup I_2}) \qquad (M_I, M_J) \mapsto M_I \boxtimes M_J \ ,$$
$$\times_{j \in J} D(X^{I_j}) \to D(X^I) \qquad (M_{I_j}) \mapsto \boxtimes_{j \in J} M_{I_j} \ .$$

It is difficult to write an explicit expression for the above monoidal product in general. However, we have the following results from [FG11] and [GL19]:

*Proposition* 5.2.2. For each $\mathcal{A}, \mathcal{B}, \mathcal{A}^j \in D(\mathrm{Ran}_X)$, there are canonical maps

$$\bigoplus_{\pi : I \,\twoheadrightarrow\, \{1,2\}} A_{I_1} \boxtimes B_{I_2} \to (\mathcal{A} \otimes^* \mathcal{B})_I \qquad \text{and} \qquad \bigoplus_{\pi : I \,\twoheadrightarrow\, J} \boxtimes_{j \in J} A^j_{I_j} \to (\otimes^*_{j \in J} \mathcal{A}^j)_I \ .$$

*Proposition* 5.2.3. Let $\mathcal{A}, \mathcal{B} \in D(\mathrm{Ran}_X)$. Then there is a canonical equivalence

$$j^{(I)*}(\mathcal{A} \otimes^* \mathcal{B})_I \cong j^{(I)*}\left( \bigoplus_{I = I_0 \cup I_1} \Delta^!_{I_0, I_1} (A_{I_0} \boxtimes B_{I_1}) \right) \ ,$$

where $j^{(I)} : U^{(I)} \to X^I$ is the complement of the union of all partial diagonals, and $\Delta_{I_0, I_1} : X^I \hookrightarrow X^{I_0} \times X^{I_1}$ is the diagonal embedding corresponding to the union map $I_0 \sqcup I_1 \twoheadrightarrow I_0 \cup I_1 = I$. The analogous statement also holds for higher arity tensor products.

*Remark* 5.2.4. The essential image $D(\mathrm{Ran}_X)_X$ of the inclusion $\Delta_* : D(X) \to D(\mathrm{Ran}_X)$ is evidently not closed under $\otimes^*$, so that $\otimes^*$ does not restrict to a monoidal structure on $D(X)$. However, it restricts to define an operad (or 'pseudo tensor', in the language of [BD04]) structure on $D(X)$, as in Example C.1.11.

*Definition* 5.2.5. The $\otimes^*$ operad structure on $D(X)$ is defined by the inclusion $D(X) \hookrightarrow D(\mathrm{Ran}_X)^*$, as in Example C.1.11.

*Corollary* 5.2.6. The functors $\Delta^{\mathrm{main}}_* : D(X) \to D(\mathrm{Ran}_X)$ and $\Delta^{\mathrm{main},!} : D(\mathrm{Ran}_X) \to D(X)$ of Remark 4.2.5 define an equivalence of operads between $D(\mathrm{Ran}_X)^*_X$ and $D(X)^*$.

*Example* 5.2.7. The multilinear operations in $D(X)^*$ are given by

$$\mathrm{Hom}_{D(X)^*}(\{M_i\}, L) = \mathrm{Hom}_{D(X^I)}(\boxtimes_{i \in I} M_i, \Delta^I_* L) \ .$$

For $b \in \mathrm{Hom}_{D(X)^*}(\{L, L\}, L)$, the composition $b \circ (b \otimes \mathbb{1}) \in \mathrm{Hom}_{D(X)^*}(\{L, L, L\}, L)$ is defined by

$$L \boxtimes L \boxtimes L \xrightarrow{b \boxtimes \mathbb{1}} \Delta_* L \boxtimes L \cong \Delta^{12,3}_*(L \boxtimes L) \xrightarrow{\Delta^{12,3}_*(b)} \Delta^{12,3}_X \Delta^{(2)}_* L = \Delta^{(3)}_* L \ ,$$

where $\Delta^{12,3} : X^2 \to X^3$ is defined by $(x, y) \mapsto (x, x, y)$.

## 5.3. The $\otimes^{\mathbf{ch}}$ tensor structure.

The space $\mathrm{Ran}_X$ has an additional monoid structure in the correspondence category, corresponding to the operation of disjoint union of finite subsets

$$\mathrm{Ran}_X^{\times 2} \xleftarrow{j_{\mathrm{disj}}} (\mathrm{Ran}_X^{\times 2})_{\mathrm{disj}} \xrightarrow{\sqcup} \mathrm{Ran}_X \ ;$$

again, this structure is discussed more formally in Section II-2.3, but only motivational in this section. The disjoint union correspondence equips $D(\mathrm{Ran}_X)$ with an another additional tensor structure.

The $\otimes^{\mathrm{ch}}$ tensor structure is defined geometrically as the exterior product, pulled back along $j_{\mathrm{disj}}$, and pushed forward along $\sqcup$, so that we define $\sqcup_* \circ j^!_{\mathrm{disj}} \circ \boxtimes : D(\mathrm{Ran}_X)^{\times 2} \to D(\mathrm{Ran}_X)$, and similarly for higher arity tensor products. This composition can again be presented as the colimit of maps of the finite cardinality subspaces, induced by maps on powers of $X$. Thus, as in the discussion preceeding Definition 5.2.1 and in keeping with Remark 4.1.4, we make the following definition:

*Definition* 5.3.1. The tensor structure $\otimes^{\mathrm{ch}} : D(\mathrm{Ran}_X)^{\times 2} \to D(\mathrm{Ran}_X)$ and its higher arity components $\otimes^{\mathrm{ch}} : D(\mathrm{Ran}_X)^J \to D(\mathrm{Ran}_X)$ are presented by the functors

$$D(X^{I_1}) \times D(X^{I_2}) \to D(X^{I_1 \sqcup I_2}) \qquad (M_I, M_J) \mapsto j_* j^! (M_I \boxtimes M_J) \ ,$$
$$\times_{j \in J} D(X^{I_j}) \to D(X^I) \qquad (M_{I_j}) \mapsto j(\pi)_* j(\pi)^! (\boxtimes_{j \in J} M_{I_j}) \ ,$$

where $j : U \hookrightarrow X^2$ is the compliment of the diagonal, and $j(\pi) : U(\pi) \hookrightarrow X^I$ is the partial diagonal complement determined by $\pi : I \twoheadrightarrow J$ as defined in Section 2.2.

*Remark* 5.3.2. The functors $j^! = j^*$ are canonically equivalent, and always defined for open embeddings. We use the notation $j^*$ throughout for the various open embeddings in this setting, in keeping with the notation of [BD04] and [FG11].

*Proposition* 5.3.3. For each $\mathcal{A}, \mathcal{B}, \mathcal{A}^j \in D(\mathrm{Ran}_X)$, there are canonical equivalences

$$(\mathcal{A} \otimes^{\mathrm{ch}} \mathcal{B})_I \cong \bigoplus_{\pi : I \twoheadrightarrow \{1,2\}} j(\pi)_* j(\pi)^* (A_{I_1} \boxtimes B_{I_2}) \qquad \text{and} \qquad (\otimes^{\mathrm{ch}}_{j \in J} \mathcal{A}^j)_I \cong \bigoplus_{\pi : I \twoheadrightarrow J} j(\pi)_* j(\pi)^* (\boxtimes_{j \in J} \mathcal{A}^j_{I_j}) \ .$$

*Remark* 5.3.4. Note that any $n$-fold tensor product vanishes when restricted to $X^I$ for $|I| < n$.

*Proposition* 5.3.5. For each $\mathcal{A}, \mathcal{B}, \mathcal{A}^j \in D(\mathrm{Ran}_X)$ there are natural maps

$$\mathcal{A} \otimes^* \mathcal{B} \to \mathcal{A} \otimes^{\mathrm{ch}} \mathcal{B} \qquad \text{and} \qquad \otimes^*_{j \in J} \mathcal{A}^j \to \otimes^{\mathrm{ch}}_{j \in J} \mathcal{A}^j \ ,$$

so that the identity defines an oplax symmetric monoidal functor $D(\mathrm{Ran}_X)^* \to D(\mathrm{Ran}_X)^{\mathrm{ch}}$.

*Remark* 5.3.6. As for $\otimes^*$, the essential image $D(\mathrm{Ran}_X)_X$ of the inclusion $\Delta_* : D(X) \to D(\mathrm{Ran}_X)$ is evidently not closed under $\otimes^{\mathrm{ch}}$, so that $\otimes^{\mathrm{ch}}$ does not restrict to a monoidal structure on $D(X)$, but still defines an operad structure on $D(X)$, as in Example C.1.11.

*Definition* 5.3.7. The $\otimes^{\mathrm{ch}}$ operad structure on $D(X)$ is defined by the inclusion $D(X) \hookrightarrow D(\mathrm{Ran}_X)^{\mathrm{ch}}$, as in Example C.1.11.

*Corollary* 5.3.8. The functors $\Delta^{\mathrm{main}}_* : D(X) \to D(\mathrm{Ran}_X)$ and $\Delta^{\mathrm{main},!} : D(\mathrm{Ran}_X) \to D(X)$ of 4.2.5 define an equivalence of operads between $D(\mathrm{Ran}_X)^{\mathrm{ch}}_X$ and $D(X)^{\mathrm{ch}}$.

*Example* 5.3.9. The multilinear operations in $D(X)^{\mathrm{ch}}$ are given by

$$\mathrm{Hom}_{D(X)^{\mathrm{ch}}}(\{M_i\}, L) = \mathrm{Hom}_{D(X^I)}(j^{(I)}_* j^{(I),*}(\boxtimes_{i \in I} M_i), \Delta^{(I)}_* L) \ ,$$

where $j^{(I)} : U^{(I)} \hookrightarrow X^I$ is the compliment of the union of all partial diagonals, as defined in 2.2.

For $\mu \in \mathrm{Hom}_{D(X)^{\mathrm{ch}}}(\{L, L\}, L)$, the composition $\mu \circ (\mu \otimes \mathbb{1}) \in \mathrm{Hom}_{D(X)^{\mathrm{ch}}}(\{A, A, A\}, A)$ is defined by

$$j^{(3)}_* (j^{(3)})^* (A \boxtimes A \boxtimes A) = j^{12,3}_* j^{12,3*} (j_* j^* (A \boxtimes A) \boxtimes A) \xrightarrow{j^{12,3}_* j^{12,3*}(\mu \boxtimes \mathbb{1})} j^{12,3}_* j^{12,3*} (\Delta_*(A) \boxtimes A)$$

$$= \Delta^{12,3}_* (j_* j^* (A \boxtimes A)) \xrightarrow{\Delta^{12,3}_*(\mu)} \Delta^{12,3}_* \Delta^{(2)}_* A = \Delta^{(3)}_* A$$

where $\Delta^{12,3} : X^2 \to X^3$ is defined by $(x, y) \mapsto (x, x, y)$, $U^{12,3} = \{(x, y, z) | x, y \neq z\}$ and $j^{12,3} : U^{12,3} \hookrightarrow X^3$.

*Remark* 5.3.10. There is a natural map of operads $D(X)^* \to D(X)^{\mathrm{ch}}$ which is the identity on objects and arity 1 morphisms, and defined on higher arity morphisms by the maps

$$\mathrm{Hom}_{D(X^I)}(\boxtimes_{i \in I} M_i, \Delta_*^{(I)} L) \to \mathrm{Hom}_{D(X^I)}(j_*^{(I)} j^{(I),*}(\boxtimes_{i \in I} M_i), \Delta_*^{(I)} L)$$

induced by

$$\boxtimes_{i \in I} M_i \to j_*^{(I)} j^{(I),*}(\boxtimes_{i \in I} M_i) \ ,$$

where the latter are given by the unit of the $(j^{(I),*}, j_*^{(I)})$ adjunction. These are intertwined with the maps of Remark 5.3.5 via the equivalences of Corollaries 5.2.6 and 5.3.8.

## 6. Factorization algebras

Following the discussion in Section 1, a factorization algebra is given by the data of

- A sheaf $\mathcal{A} \in D(\mathrm{Ran}_X)$ on $\mathrm{Ran}_X$, and
- an isomorphism $j^{(I)*} A_I \cong j^{(I)*} A_1^{\boxtimes I}$ for each $I \in \mathrm{fSet}$.

Evidently, we should require compatibility of the latter isomorphisms with the gluing maps $\Delta(\pi)^! A_I \cong A_J$, but some care is required to carefully state the homotopy coherence data in this heuristic definition. In this section, we give a formal definition of (non-unital) factorization algebra, following [FG11].

*Example* 6.0.1. A non-unital cocommutative coalgebra object $\mathcal{A} \in \mathrm{CoComm}^{\mathrm{nu}}(D(\mathrm{Ran}_X)^{\mathrm{ch}})$ is an object $\mathcal{A} \in D(\mathrm{Ran}_X)$, together with a map $\mathcal{A} \to \mathcal{A} \otimes^{\mathrm{ch}} \mathcal{A}$ and a compatible collection of higher arity analogues. Concretely, the data of the map $\mu$ must be specified compatibly over $X^I$ for each $I$: For $|I| = 1$ there is no data as $(A \otimes^{\mathrm{ch}} A)_{\{1\}} = 0$, but for $|I| = 2$, the required map is

$$A_2 \to j_* j^*(A_1 \boxtimes A_1) \qquad \text{or equivalently} \qquad j^* A_2 \to j^*(A_1 \boxtimes A_1) \ .$$

More generally, the required maps $\mathcal{A} \to \otimes_{j \in J}^{\mathrm{ch}} \mathcal{A}$ are specified over each stratum by maps

$$A_I \to j(\pi)_* j(\pi)^*(\boxtimes_{j \in J} A_{I_j}) \qquad \text{or equivalently} \qquad j(\pi)^* A_I \to j(\pi)^*(\boxtimes_{j \in J} A_{I_j})$$

for each $I, J$ and $\pi : I \twoheadrightarrow J$.

Note that for $\pi = \mathbb{1}_I$, these are maps of the type required in the heuristic definition of factorization algebra above, except that they are not necessarily equivalences. Thus, we make the following definition:

*Definition* 6.0.2. A non-unital factorization algebra on $X$ is a non-unital cocommutative coalgebra object $\mathcal{A} \in D(\mathrm{Ran}_X)^{\mathrm{ch}}$ such that the induced maps

$$(6.0.1) \qquad\qquad j(\pi)^* A_I \xrightarrow{\cong} j(\pi)^*(\boxtimes_{j \in J} A_{I_j})$$

are equivalences for each $I, J$ and $\pi : I \twoheadrightarrow J$.

A unital factorization algebra on $X$ is an object $\mathcal{A} \in D(\mathrm{Ran}_{X,\mathrm{un}})$ with a non-unital factorization algebra structure on its image in $D(\mathrm{Ran}_X)$, and compatibility data with the unital structure on $\mathcal{A}$; see Remark 6.0.3.

Let $\mathrm{Alg}^{\mathrm{fact}}(X)$ denote the category of non unital factorization algebras, defined as the full subcategory of $\mathrm{CoComm}^{\mathrm{nu}}(D(\mathrm{Ran}_X)^{\mathrm{ch}})$. Similarly, let $\mathrm{Alg}_{\mathrm{un}}^{\mathrm{fact}}(X)$ denote the category of unital factorization algebras.

*Remark* 6.0.3. The precise statement of the compatibility data is slightly involved, so we defer the formal definition of the category of unital factorization algebras to Example II-2.4.4. For now, we give the following example of the compatibility data: we require commutativity of the diagram

$$(6.0.2) \qquad \begin{array}{ccc} A_1 \boxtimes \omega_X & \longrightarrow & j_*j^*(A_1 \boxtimes \omega_X) \\ \downarrow & & \downarrow \\ A_2 & \longrightarrow & j_*j^*(A_1 \boxtimes A_1) \end{array}$$

and its higher arity analogues, where the vertical arrows are those from Equation 4.3.1.

*Example* 6.0.4. The dualizing sheaf $\omega_{\mathrm{Ran}_X} \in D(\mathrm{Ran}_X)$ of Example 4.2.7 defines a non-unital factorization algebra $\omega_{\mathrm{Ran}_X} \in \mathrm{Alg}^{\mathrm{fact}}(X)$, with structure maps

$$\omega_{X^I} \to j(\pi)_* j(\pi)^*(\boxtimes_{j \in J} \omega_{X^{I_j}})$$

given by the unit of the $(j^*, j_*)$ adjunction under the identification $\omega_{X^I} \boxtimes \omega_{X^J} \cong \omega_{X^{I \sqcup J}}$.

Similarly, the dualizing sheaf $\omega_{\mathrm{Ran}_{X,\mathrm{un}}} \in D(\mathrm{Ran}_{X,\mathrm{un}})$ of Example 4.3.5 defines a non-unital factorization algebra $\omega_{\mathrm{Ran}_{X,\mathrm{un}}} \in \mathrm{Alg}^{\mathrm{fact}}_{\mathrm{un}}(X)$.

*Proposition* 6.0.5. The $\otimes^!$ monoidal structure on $D(\mathrm{Ran}_{X,\mathrm{un}})$ induces a symmetric monoidal structure on $\mathrm{Alg}^{\mathrm{fact}}_{\mathrm{un}}(X)$, such that $\omega_{\mathrm{Ran}_{X,\mathrm{un}}}$ is the tensor unit.

*Proof.* We postpone the proof until that of Proposition II-6.0.7, which is the first place it is essentially used. $\qquad\square$

Let $\mathrm{Alg}^{\mathrm{fact}}_{\mathrm{un}}(X)^{\otimes !}$ denote the symmetric monoidal category of unital factorization algebras in the $\otimes^!$ monoidal structure.

*Definition* 6.0.6. A non unital factorization algebra $\mathcal{A} \in \mathrm{Alg}^{\mathrm{fact}}(X)$ is called commutative if the inverse of the equivalence 6.0.1 extends to a map

$$\boxtimes_{j \in J} A_{I_j} \to A_I$$

for each $I, J$ and $\pi : I \twoheadrightarrow J$. A unital factorization algebra $\mathcal{A} \in \mathrm{Alg}^{\mathrm{fact}}_{\mathrm{un}}(X)$ is called commutative if it is commutative as a non unital factorization algebra.

## 7. Chiral algebras

In this section we define the category of chiral algebras, which is equivalent to the category of factorization algebras, as we recall in Section 14. A chiral algebra is more closely analogous to a global analogue of the notion of vertex algebra, and in Section 8 we exhibit an equivalence of categories between weakly translation invariant, unital chiral algebras on $\mathbb{A}^1$ and vertex algebras. Again, we follow the approach of [FG11], and in turn [BD04], throughout this section.

*Example* 7.0.1. A Lie algebra object $\mathcal{L} \in D(\mathrm{Ran}_X)^{\mathrm{ch}}$ is an object $\mathcal{L}$, together with a map $\mu : \mathcal{L} \otimes^{\mathrm{ch}} \mathcal{L} \to \mathcal{L}$ and its higher arity analogues. Concretely, the data of $\mu$ must be specified compatibly over $X^I$ for each $I$. For $|I| = 1$ there is no data since $(\mathcal{L} \otimes^{\mathrm{ch}} \mathcal{L})_{\{1\}} = 0$, but for $|I| = 2$ the required map is

$$\mu : j(\pi)_* j(\pi)^*(L_1 \boxtimes L_1) \to L_2 \ .$$

More generally, for each $J$ we require maps $\mu_J : \otimes^{\mathrm{ch}}_{j \in J} L \to L$ which are specified on $X^I$ by maps

$$\mu(\pi) : j(\pi)_* j(\pi)^*(\boxtimes_{j \in J} L_{I_j}) \to L_I$$

for each $I, J$ and $\pi : I \twoheadrightarrow J$.

*Definition* 7.0.2. A (non-unital) chiral algebra on $X$ is a Lie algebra object $\mathcal{L} \in \mathrm{Lie}(D(\mathrm{Ran}_X)^{\mathrm{ch}})$ such that underlying object $\mathcal{L} \in D(\mathrm{Ran}_X)_X$ is supported on the image of the main diagonal $\Delta^{\mathrm{main}} : X \to \mathrm{Ran}_X$.

A unital chiral algebra on $X$ is an object $\mathcal{L} \in D(\mathrm{Ran}_{X,\mathrm{un}})$ with a non unital chiral algebra structure on its image in $D(\mathrm{Ran}_X)$, and compatibility data with the unital structure on $\mathcal{L}$; see Remark 7.0.5.

Let $\mathrm{Alg}^{\mathrm{ch}}(X)$ denote the category of chiral algebras, defined as the full subcategory of $\mathrm{Lie}(D(\mathrm{Ran}_X)^{\mathrm{ch}})$. Similarly, let $\mathrm{Alg}^{\mathrm{ch}}_{\mathrm{un}}(X)$ denote the category of unital chiral algebras.

*Remark* 7.0.3. Concretely, the condition $\mathcal{L} \in D(\mathrm{Ran}_X)_X$ is the requirement that

$$L_I = \Delta^{(I)}_* A \qquad \text{for some} \qquad A \in D(X)$$

for each $I \in \mathrm{fSet}$, where $\Delta^{(I)} : X \hookrightarrow X^I$ is the small diagonal embedding, as defined in Section 2.2.

Evidently, for such objects the data of the chiral Lie algebra structure maps are given by maps

$$(7.0.1) \qquad \mu(\pi) : j(\pi)_* j(\pi)^* (\boxtimes_{j \in J} \Delta^{(I_j)}_* A) \to \Delta^{(I)}_* A$$

for each $I, J$ and $\pi : I \twoheadrightarrow J$. This is essentially the same statement as Corollary 7.0.7 below.

*Warning* 7.0.4. We change notation and write $A \in \mathrm{Alg}^{\mathrm{ch}}(X)$ instead of denoting it by the underlying object $\mathcal{L} = \Delta^{\mathrm{main}}_* A \in D(\mathrm{Ran}_X)_X$.

*Remark* 7.0.5. Again, the precise definition of the unit compatibility data is not given here. In Example II-2.4.4, we give the formal definition of a unital factorization algebra, and we formally define a unital structure on a chiral algebra to correspond to that on a factorization algebra under the equivalence $\mathrm{Alg}^{\mathrm{fact}}(X) \cong \mathrm{Alg}^{\mathrm{ch}}(X)$ of Section 14. For now, we give the following example of the compatibility data: we require commutativity data for the diagram

$$(7.0.2) \qquad \begin{array}{ccc} j_* j^* (A \boxtimes \omega_X) & \longrightarrow & \Delta_* A \\ \downarrow & & \| \\ j_* j^* (A \boxtimes A) & \longrightarrow & \Delta_* A \end{array}$$

and its higher arity analogues, where the vertical map is that induced by the unit map of Equation 4.3.1.

*Proposition* 7.0.6. The forgetful functor $\mathrm{Alg}^{\mathrm{ch}}_{\mathrm{un}}(X) \to \mathrm{Alg}^{\mathrm{ch}}(X)$ admits a left adjoint, defined on the underlying $D$ module by the functor $A \mapsto \omega_X \oplus A$.

The resulting unital chiral algebra $A \oplus \omega_X \in \mathrm{Alg}^{\mathrm{ch}}_{\mathrm{un}}(X)$ is called the free unital chiral algebra on $A$.

The following is a formal consequence of the definition of chiral algebra:

*Corollary* 7.0.7. A non unital chiral algebra is equivalent to a Lie algebra $A \in \mathrm{Alg}_{\mathrm{Lie}}(D(X)^{\mathrm{ch}})$ internal to the operad $D(X)^{\mathrm{ch}}$, by Corollary 5.3.8.

*Remark* 7.0.8. In particular, we can summarize the above maps 7.0.1 in terms of the Lie algebra structure maps

$$\mu_I \in \mathrm{Hom}_{D(X)^{\mathrm{ch}}}(\{A\}_{i \in I}, A) = \mathrm{Hom}_{D(X^I)}(j^{(I)}_* j^{(I),*} (\boxtimes_{i \in I} A), \Delta^{(I)}_* A) \ .$$

*Definition* 7.0.9. A non unital chiral algebra $A \in \mathrm{Alg}^{\mathrm{ch}}(X)$ is called commutative if the composition

$$\boxtimes_{j \in J} \Delta_*^{(I_j)} A \to j(\pi)_* j(\pi)^* (\boxtimes_{j \in J} \Delta_*^{(I_j)} A) \xrightarrow{\mu(\pi)} \Delta_*^{(I)} A$$

vanishes for each $I, J \in \mathrm{fSet}$ and $\pi : I \twoheadrightarrow J$, where the first map is the unit of the $(j^*, j_*)$-adjunction and the second is the chiral algebra structure map $\mu(\pi)$ from equation 7.0.1. A unital chiral algebra $A \in \mathrm{Alg}^{\mathrm{ch}}_{\mathrm{un}}(X)$ is called commutative if it is commutative as a non unital chiral algebra.

*Remark* 7.0.10. In terms of the corresponding Lie algebra operad $D(X)^{\mathrm{ch}}$, this is equivalent to the vanishing of the composition

$$\boxtimes_{i \in I} A \to j_*^{(I)} j^{(I),*} (\boxtimes_{i \in I} A) \to \Delta_*^{(I)} A$$

for each $I \in \mathrm{fSet}$.

## 8. FROM CHIRAL ALGEBRAS TO VERTEX ALGEBRAS

*Warning* 8.0.1. Throughout this section, all of the objects will be of cohomological degree zero (in the heart of the relevant t-structure) and all the functors non-derived, in contrast with our general conventions.

Let $X = \mathbb{A}^n$ be $n$ dimensional affine space and let $G = \mathbb{G}_a^n$ act on $\mathbb{A}^n$ by translation. Fix global coordinates $x_i$ on $\mathbb{A}^n$ and note the resulting identification of the algebra of differential operators and its translation invariant subalgebra $\Gamma(\mathbb{A}^n, D_{\mathbb{A}^n}) \cong \mathbb{K}[x_i, \partial_{x_i}]$ and $\Gamma(\mathbb{A}^n, D_{\mathbb{A}^n})^{\mathbb{G}_a^n} \cong \mathbb{K}[\partial_{x_i}]$. We work with the left $D$ module model for $D(X)$ throughout this section and the next.

Consider the category $D(\mathbb{A}^n)^{\mathbb{G}_a^n, w}$ of weakly translation equivariant $D$ modules on $\mathbb{A}^n$; see Section 15 for a breif review of equivariant D modules.

*Remark* 8.0.2. For each $M \in D(\mathbb{A}^n)^{\mathbb{G}_a^n, w}$ the space of translation invariant sections $\Gamma(\mathbb{A}^n, M)^{\mathbb{G}_a^n}$ is naturally a module for the translation invariant differential operators $\Gamma(\mathbb{A}^n, D_{\mathbb{A}^n})^{\mathbb{G}_a^n}$. The original object $M$ can be recovered from this data as

$$M = \Gamma(\mathbb{A}^n, M)^{\mathbb{G}_a^n} \otimes_{\Gamma(\mathbb{A}^n, D_{\mathbb{A}^n})^{\mathbb{G}_a^n}} \Gamma(\mathbb{A}^n, D_{\mathbb{A}^n}) \ .$$

*Remark* 8.0.3. Given a choice of closed point $\iota_0 : 0 \hookrightarrow \mathbb{A}^n$, there is a natural identification

$$M_0 := \iota_0^! M \cong \Gamma(\mathbb{A}^n, M)^{\mathbb{G}_a^n}$$

given by extending the element of the fibre to a translation invariant section over $\mathbb{A}^n$; such a section exists and is unique since $\mathbb{A}^n$ is a $\mathbb{G}_a^n$ torsor.

In particular, there is a canonical trivialization of $M$ as an $\mathcal{O}$ module, such that the action of differential operators is given by the usual action on functions $\mathcal{O}_{\mathbb{A}^n}$ together with the $\mathbb{K}[\partial_{x_i}]$ module structure on $M_0$:

$$M \cong M_0 \otimes_{\mathbb{K}[\partial_{x_i}]} \mathbb{K}[x_i, \partial_{x_i}] \cong M_0 \otimes_{\mathbb{K}} \mathbb{K}[x_i] \qquad \text{where} \qquad \partial_{x_i}(m \otimes f) = \partial_{x_i} m \otimes f + m \otimes \partial_{x_i} f$$

for each $m \in M_0$ and $f \in \mathbb{K}[x_i]$.

In summary, we have the following:

*Proposition* 8.0.4. There is an equivalence of categories

$$(8.0.1) \qquad D(\mathbb{A}^n)^{\mathbb{G}_a^n, w} \underset{\cong}{\overset{\cong}{\rightleftarrows}} D(\mathbb{K}[\partial_{x_i}]) \qquad \text{defined by} \qquad M \mapsto M_0 \qquad \mathbb{V} \otimes_{\mathbb{K}} \mathbb{K}[x_i] \leftarrow\!\shortmid \mathbb{V} \ .$$

Now, we consider the case when $X = \mathbb{A}^1$ and state the main result of this section, following Chapter 3.6 of [BD04]; see Definition 17.0.6 for the notion of weakly equivariant chiral algebra.

*Theorem* 8.0.5. There is an equivalence of categories between the category of weakly translation equivariant chiral algebras on $\mathbb{A}^1$ and that of vertex algebras

$$\mathrm{Alg}_{\mathrm{un}}^{\mathrm{ch}}(\mathbb{A}^1)^{\mathbb{G}_a,w} \underset{\cong}{\overset{\cong}{\rightleftarrows}} \mathrm{VOA}_{\mathbb{K}} \qquad \text{defined by} \qquad A \mapsto A_0 \qquad \mathbb{V} \otimes_{\mathbb{K}} \mathbb{K}[x] \hookleftarrow \mathbb{V},$$

where $A \in D(\mathbb{A}^1)^{\mathbb{G}_a,w}$ is such that $\Delta_*^{\mathrm{main}} A = \mathcal{L} \in D(\mathrm{Ran}_{\mathbb{A}^1})_{\mathbb{A}^1}^{\mathbb{G}_a,w}$ and $A_0 = \iota_0^! A$, as above.

*Proof.* Let $\mathcal{L} \in \mathrm{Alg}^{\mathrm{ch}}(\mathbb{A}^1)^{\mathbb{G}_a,w}$ be a weakly translation equivariant chiral algebra, $A \in D(\mathbb{A}^1)^{\mathbb{G}_a,w}$ be such that $\Delta_*^{\mathrm{main}} A = \mathcal{L} \in D(\mathrm{Ran}_{\mathbb{A}^1})^{\mathbb{G}_a,w}$, and let $A = \mathbb{V} \otimes_{\mathbb{K}} \mathbb{K}[x]$ for $\mathbb{V} \in D(\mathbb{K}[\partial])$ the corresponding complex of $\mathbb{K}[\partial]$ modules under the equivalence of Proposition 8.0.1 above.

The underlying vector space $\mathbb{V} \in \mathrm{Vect}_{\mathbb{K}}$ defines the state space of our putative vertex algebra, and the module structure $\rho : \mathbb{K}[\partial] \to \mathrm{End}_{\mathbb{K}}(\mathbb{V})$ is equivalent to the action of the generator $T := \rho(\partial) \in \mathrm{End}(\mathbb{V})$ defines the translation operator. The unit map $\omega_{\mathbb{A}^1} \to L$ defines a distinguished vector $\varnothing \in \mathbb{V}$ such that $T(\varnothing) = 0$, which defines the unit for the vertex algebra $\mathbb{V}$.

It remains to identify the weakly equivariant chiral algebra structure map

$$\mu \in \mathrm{Hom}_{D(\mathbb{A}^1)^{(\mathbb{G}_a,w),\mathrm{ch}}}(\{A,A\},A) = \mathrm{Hom}_{D(\mathbb{A}^2)^{\mathbb{G}_a,w}}(j_*j^*(A\boxtimes A),\Delta_*A)$$

with a vertex operator map $Y(\cdot,z) : \mathbb{V}^{\otimes 2} \to \mathbb{V}((z))$ satisfying the conditions in the definition of vertex algebra, recalled in Definition E.1.3. The chiral algebra structure map can be written explicitly in coordinates as a map

$$\mathbb{V}^{\otimes 2} \otimes_{\mathbb{K}} \mathbb{K}[x,y,(x-y)^{-1}] \xrightarrow{b} \mathbb{V} \otimes_{\mathbb{K}} \mathbb{K}[x] \otimes_{\mathbb{K}} \delta_{x-y}$$

which intertwines the action of $\mathbb{K}[x,y,\partial_x,\partial_y]$ on each side. Here $\delta_{x-y} \cong (x-y)^{-1}\mathbb{K}[(x-y)^{-1}]$ is the delta function $D$ module on the diagonal.

Now, we apply Proposition 9.1.4 of the next subsection, just as in Example 9.2.9 but in the special case $X = \mathbb{A}^1$, which implies there exist natural maps $\bar{\mu}, \tilde{\mu}$ such that the following commutes

$$(8.0.2) \qquad \begin{array}{ccc} \mathbb{V}^{\otimes 2} \otimes_{\mathbb{K}} \mathbb{K}[x,y] & \xrightarrow{\tilde{\mu}} & \mathbb{V} \otimes \mathbb{K}[x] \otimes_{\mathbb{K}} \mathbb{K}((x-y)) \\ {\scriptstyle \iota} \downarrow & \overset{\bar{\mu}}{\nearrow} & \downarrow {\scriptstyle q} \\ \mathbb{V}^{\otimes 2} \otimes_{\mathbb{K}} \mathbb{K}[x,y,(x-y)^{-1}] & \xrightarrow{\mu} & \mathbb{V} \otimes_{\mathbb{K}} \mathbb{K}[x] \otimes_{\mathbb{K}} \delta_{x-y} \end{array},$$

and moreover each of the three maps $\mu,\bar{\mu},\tilde{\mu}$ is uniquely determined by the others; this is just the diagram of Equation 9.2.2 in the case $X = \mathbb{A}^1$.

The maps $\mu,\bar{\mu},\tilde{\mu}$ are all maps of weakly $\mathbb{G}_a$ equivariant $D$ modules, and in particular $\tilde{\mu}$ is determined by its restriction to $(0,0) \in \mathbb{A}^2$, which defines

$$Y : \mathbb{V}^{\otimes 2} \to \mathbb{V} \otimes_{\mathbb{K}} \mathbb{K}((z))$$

for $z = (x-y)$, as desired. Moreover:

- The commutativity of the diagram 7.0.2 and its analogue with the arguments $\omega_X$ and $A$ interchanged are equivalent to the conditions $Y(\varnothing,z) = \mathbb{1}_V$ and $Y(a,z)(\varnothing) \in V[[z]]$ with $Y(a,z)(\varnothing)|_{z=0} = a$.
- The fact that $\tilde{\mu}$ is a map of $D$ modules is equivalent to the condition that $[T,Y(a,z)] = \partial_z Y(a,z)$.
- The Jacobi identity for the chiral Lie bracket $\mu$ is equivalent to the mutual locality or 'associativity' condition of the vertex algebra; writing the Jacobi identity in coordinates in terms of the identifications above gives the well-known Jacobi type formulation of the associativity axiom of vertex algebras.

Alternatively, the Theorem 9.2.10 of the following subsection in the case of $X = \mathbb{A}^1$, restricted to the subcategories of $\mathbb{G}_a$ equivariant objects, implies these results. $\qquad\square$

## 9. Operator product expansions

In this section we summarize the results of Section 3.8 of [BD04], a special case of which is recalled in the preceding section.

*Warning* 9.0.1. Throughout this section, all of the objects will be of cohomological degree zero (in the heart of the relevant t-structure) and all the functors non-derived, in contrast with our general conventions.

### 9.1. **Abstract preliminaries.** 

We begin by recalling some abstract preliminary material, which we recommend the reader skip, and return to only as necessary.

Let $P$ a smooth algebraic variety. A $D_P$-sheaf is a (not neccesarily quasi coherent) sheaf of modules for $D_P$ over $P$ in the etale topology. Let $\tilde{\mathcal{M}}(P) = D_{\mathrm{et}}(P)^\heartsuit$ denote the abelian category of $D_P$-sheaves, in which $D(P)^\heartsuit$ is a full subcategory.

*Construction* 9.1.1. Let $i : Z \hookrightarrow P$ be a closed embedding of a smooth subvariety with $J \subset \mathcal{O}_P$ the corresponding ideal sheaf. There is a functor $\iota^* : \tilde{\mathcal{M}}(P) \to \tilde{\mathcal{M}}(Z)$ defined by $\iota^* F = \iota^\cdot(F/J \cdot F)$, with $D_Z$ module structure defined as usual.

The functor $\iota^*$ admits a right adjoint $\hat{\iota}_* : \tilde{\mathcal{M}}(Z) \to \tilde{\mathcal{M}}(P)$. This functor is exact and fully faithful, and its image is the full subcategory on objects which are complete along $Z \hookrightarrow P$.

The analogue of the usual functor

$$\iota_* : \tilde{\mathcal{M}}(Z) \to \tilde{\mathcal{M}}(P) \qquad \text{is defined by} \qquad \iota_* G = \hat{\iota}_* G \otimes (\iota_* \omega_Z)^l \ .$$

This functor is exact, fully faithful, and agrees with the usual direct image functor when restricted to $D^l(X)$. Its image is the full subcategory on objects such that each local section is annihilated by some power of $J$.

*Construction* 9.1.2. Let $j : U \hookrightarrow P$ denote the complementary open embedding to $\iota : Z \hookrightarrow P$, and define

$$\tilde{\iota}_* : \tilde{\mathcal{M}}(Z) \to \tilde{\mathcal{M}}(P) \qquad \text{by} \qquad \tilde{\iota}_* G = \hat{\iota}_* G(\cdot) \otimes j_* \mathcal{O}_U \ .$$

For $Z$ of codimension 1 in $P$, the short exact sequence

$$\mathcal{O}_P \hookrightarrow j_* j^* \mathcal{O}_P \twoheadrightarrow \iota_* \iota^*(\mathcal{O}_P) \qquad \text{reduces to} \qquad \mathcal{O}_P \hookrightarrow j_*(\mathcal{O}_U) \twoheadrightarrow \iota_*(\mathcal{O}_Z)$$

so that for each $G \in \tilde{\mathcal{M}}(Z)$ we obtain the exact sequence

$$\hat{\iota}_* G \hookrightarrow \tilde{\iota}_* G \twoheadrightarrow \iota_* G \ . \tag{9.1.1}$$

*Example* 9.1.3. Let $Z$ be a smooth algebraic curve $X$, $P = X \times X$, and $\iota = \Delta : X \to X \times X$. Then

$$\hat{\Delta}_* G \cong \mathrm{p}_0^{-1} G \otimes_{\mathrm{p}_0^{-1} \mathcal{O}_X} (\mathrm{p}_0^{-1} \mathcal{O}_X \otimes_{\mathbb{K}} \mathbb{K}[[x - y]]) \cong G_0 \otimes_{\mathbb{K}} \mathbb{K}[[x - y]] \cong G_1 \otimes_{\mathbb{K}} \mathbb{K}[[x - y]]$$

where $x, y$ are local coordinates on $X$, $G_i = \mathrm{p}_i^{-1} G$, and similarly

$$\tilde{\Delta}_* G \cong G_0 \otimes_{\mathbb{K}} \mathbb{K}((x - y)) \cong G_1 \otimes_{\mathbb{K}} \mathbb{K}((x - y)) \tag{9.1.2}$$

In particular, for $Z = \mathbb{A}^1$ we have $\hat{\Delta}_* G \cong G \otimes_{\mathbb{K}[x]} \mathbb{K}[x][[x - y]] \cong G_1 \otimes_{\mathbb{K}} \mathbb{K}[[x - y]] \cong G_2 \otimes_{\mathbb{K}} \mathbb{K}[[x - y]]$ and $\tilde{\Delta}_* G$.

We now state the key lemma which is used in the comparision of chiral structure of equation 7.0.1 and the operator product expansion maps of Definition 9.2.1 below. The latter is used in the definition 9.2.8 of OPE algebras, which generalizes the notion of vertex algebras to global curves, and is essentially equivalent to the notion of a bundle of vertex algebras in [FBZ04]. This comparison is used in the proof of Theorem 9.2.10 below, which generalizes Theorem 8.0.5 of the previous section.

*Proposition* 9.1.4. Let $G \in \bar{\mathcal{M}}(Z)$, $F \in \bar{\mathcal{M}}(P)$ and suppose $Z$ is codimension 1. Then the natural maps

$$\mathrm{Hom}(F \otimes j_* \mathcal{O}_U, \tilde{\iota}_* G) \twoheadrightarrow \mathrm{Hom}(F, \tilde{\iota}_* G) \qquad \text{and} \qquad \mathrm{Hom}(F \otimes j_* \mathcal{O}_U, \tilde{\iota}_* G) \hookrightarrow \mathrm{Hom}(F \otimes j_* \mathcal{O}_U, \iota_* G) \ ,$$

are isomorphisms, where the former is given by precomposition with $i$ and the latter by postcomposition with $q$, as summarized in the diagram

$$(9.1.3)$$

$$
\begin{array}{ccc}
F & \xrightarrow{\tilde{\varphi}} & \tilde{\iota}_* G \\
{\scriptstyle i}\downarrow & \overset{\bar{\varphi}}{\nearrow} & \downarrow{\scriptstyle q} \\
F \otimes j_* \mathcal{O}_U & \xrightarrow{\varphi^l} & \iota_* G
\end{array}
$$

where $\tilde{\varphi}$ and $\varphi^l$ are the images of $\bar{\varphi}$ under the two equivalences.

## 9.2. Operator product expansions.
We now restict our attention again to the case that $X$ is a smooth algebraic curve.

*Definition* 9.2.1. Let $\{F_i\}, G \in \bar{\mathcal{M}}(X)$. The space of operator product expansion (OPE) operations is defined as

$$O_I(\{F_i\}, G) = \mathrm{Hom}(\boxtimes_I F_i, \tilde{\Delta}_*^{(I)}(G)) \ ,$$

where $\tilde{\Delta}_*^{(I)} := \hat{\Delta}_*^{(I)}(\cdot) \otimes j_*^{(I)} \mathcal{O}_{U^{(I)}} : \bar{\mathcal{M}}(X) \to \bar{\mathcal{M}}(X^I)$, and $\hat{\Delta}_*^{(I)} = (\Delta^{(I)})_*^{\check{}}$ as defined in 9.1.1.

*Remark* 9.2.2. Note that $\tilde{\Delta}_*^{(I)} \neq (\Delta_*^{(I)})_*^{\check{}}$ unless $|I| = 2$, where the latter is as defined in 9.1.2.

*Proposition* 9.2.3. For each $\pi : J \twoheadrightarrow I$, there are natural composition maps

$$\bigotimes_{i \in I} O_{J_i}(\{H_j\}, F_i) \otimes O_I(\{F_i\}, G) \to \mathrm{Hom}(\boxtimes_{j \in J} H_j, \tilde{\Delta}_*^{(\{I,J\})} G)$$

for any $\{H_i^j\}, \{F_i\}, G \in \mathcal{M}(X)$, where $\tilde{\Delta}^{(\{I,J\})} = (\hat{\Delta}(\pi)_* \circ \tilde{\Delta}_*^{(J)}(\cdot)) \otimes j(\pi)_* \mathcal{O}_{U(\pi)}$.

*Remark* 9.2.4. These maps do not define an operad stucture on $D(X)^{\heartsuit}$ with multilinear maps defined by OPE operations, as the composition of two OPE operations may fail to define another OPE operation.

*Remark* 9.2.5. The above composition maps are still associative in the appropriate sense: there is no ambiguity in preforming itterated compositions, although such compositions are valued in generalizations of the spaces of operations, as in the preceeding proposition.

*Remark* 9.2.6. There is a natural inclusion $\tilde{\Delta}_*^{(J)} G \hookrightarrow \tilde{\Delta}^{(\{I,J\})} G$ giving an inclusion

$$(9.2.1) \qquad O_J(\{H_j\}, G) = \mathrm{Hom}(\boxtimes_{j \in J} H_j, \tilde{\Delta}_*^{(J)} G) \hookrightarrow \mathrm{Hom}(\boxtimes_{j \in J} H_j, \tilde{\Delta}_*^{(\{I,J\})} G)$$

*Definition* 9.2.7. A collection of OPE operations $\gamma \otimes (\otimes_i \delta_i) \in \bigotimes_{i \in I} O_{J_i}(\{H_j\}, F_i) \otimes O_I(\{F_i\}, G)$ compose nicely if $\gamma(\delta_i)$ is in the image of the inclusion.

*Definition* 9.2.8. Let $\tilde{\mu} \in O_2(\{G, G\}, G)$ a binary ope operation. Then

- $\tilde{\mu}$ is called associative if both $\tilde{\mu} \otimes (\tilde{\mu} \otimes \mathbb{1}_G)$ and $\tilde{\mu} \otimes (\mathbb{1}_G \otimes \tilde{\mu})$ compose nicely, and their values in $O_3(\{G, G, G\}, G)$ coincide.
- $\tilde{\mu}$ is called commutative if it invariant under the natural transposition of factors isomorphism on $O_2(\{G, G\}, G)$.

An associative and commutative binary ope operation $\tilde{\mu} \in O_2(\{G, G\}, G)$ is called an ope algebra structure on $G$, and $(G, \tilde{\mu})$ is called an OPE algebra.

A unit 1 for an ope algebra $(G, \tilde{\mu})$ is a horizontal section $1 \in G$ such that for each $a \in G$ we have $\tilde{\mu}(1 \boxtimes a), \tilde{\mu}(a \boxtimes 1) \in \hat{\Delta}_* G \subset \tilde{\Delta}_* G$, and both project to $a \in G$ under $\hat{\Delta}_* G \to G$ the cokernel of the inclusion $J_\Delta \cdot \hat{\Delta}_* G \hookrightarrow \hat{\Delta}_* G$.

*Example* 9.2.9. Applying proposition 9.1.4 in the case $P = X \times X$, $Z = X$ and $\iota = \Delta$ with $F = L_1^l \boxtimes L_2^l$ and $G = M^l$, we obtain an isomorphism

$$O_2(L_1, L_2; M) = \operatorname{Hom}(L_1 \boxtimes L_2, \tilde{\Delta}_* M) \cong \operatorname{Hom}(j_* j^*(L_1^l \boxtimes L_2^l), \Delta_* M^l) \otimes_\mathbb{K} \lambda_2 = \operatorname{Hom}_{D(X)^{\operatorname{ch}}}(L_1, L_2; M) \otimes_\mathbb{K} \lambda_2 \,,$$

where $\lambda_2 = \mathbb{K}_{\operatorname{sign}} \in \mathbb{K}[S_2]$-Mod is the sign representation. More generally there is a canonical embedding, which is not an equivalence in keeping with Remark 9.2.4 above,

$$O_I(\{L_i^l\}, M^l) = \operatorname{Hom}(\boxtimes_i L_i^l, \tilde{\Delta}_*^{(I)} M^l) \hookrightarrow \operatorname{Hom}(j_*^{(I)} j^{(I)*}(\boxtimes_i L_i^l), \Delta_*^{(I)} M^l) \otimes \lambda_I = \operatorname{Hom}_{D(X)^{\operatorname{ch}}}(\{L_i\}, M) \otimes \lambda_I$$

where $\lambda_I = \omega_X^{\boxtimes I} \otimes \omega_{X^I}^{-1}$, recalling $j_* j^*(F) \cong F \otimes j_* \mathcal{O}_U$ for $F \in \bar{\mathcal{M}}(X)$.

In the case of interest when $L_1^l = L_2^l = M^l = A^l$, the analogue of the diagram 9.1.3 is given by

(9.2.2)
$$\begin{array}{ccc} A \boxtimes A & \xrightarrow{\tilde{\mu}} & \tilde{\Delta}_* A \\ {\scriptstyle i} \downarrow & {\scriptstyle \bar{\mu}} \nearrow & \downarrow {\scriptstyle q} \\ j_* j^*(A \boxtimes A) & \xrightarrow{\mu} & \Delta_* A \end{array} \quad .$$

Now, the main result of Chapter 3.8 of [BD04] is the following:

*Theorem* 9.2.10. The isomorphism

$$O_2(\{A^l, A^l\}, A^l) \cong P_2^{\operatorname{ch}}(\{A, A\}, A)$$

gives a bijection between the set of ope algebra structures on $A^l$ and the set of non-unital chiral algebra structures on $A$, such that a flat section $1 \in A^l$ defines a unit for an ope structure if and only if it defines a unit for the corresponding chiral algebra.

## 10. From Lie*, Comm$^!$ and Coisson to Lie, commutative, and Poisson vertex algebras

In Section 11, we recall that chiral algebras are closely related to topological associative algebras, following Section 3.6 of [BD04]. In this section, we recall the chiral analogues of Lie, commutative, and Poisson algebras, called Lie*, Comm$^!$, and Coisson algebras, respectfully. Further, we show that on $X = \mathbb{A}^1$ in the $\mathbb{G}_a$ equivariant setting, such objects are equivalent to Lie, commutative, and Poisson vertex algebras, respectively, in analogy with the results of Section 8. We follow sections 2.3, 2.5, and 2.6 of [BD04] throughout.

10.1. **Overview.** Recall that there are canonical maps of operads
(10.1.1)

Lie → Ass → Comm        and corresponding functors        $\mathrm{Alg}_{\mathrm{Comm}}(\mathcal{O}) \to \mathrm{Alg}_{\mathrm{Ass}}(\mathcal{O}) \to \mathrm{Alg}_{\mathrm{Lie}}(\mathcal{O})$

on algebras internal any operad $\mathcal{O} \in \mathrm{Op}(\mathrm{Vect}_{\mathbb{K}})$. These functors are just the inclusion of commutative algebras as a full subcategory of associative algebras, and the forgetful functor from associative to Lie algebras given by remembering only commutators, respectfully. This sequence is 'left exact', in the sense that the functor $\mathrm{Alg}_{\mathrm{Comm}}(\mathcal{O}) \to \mathrm{Alg}_{\mathrm{Ass}}(\mathcal{O})$ is the inclusion of the full subcategory on objects whose image under $\mathrm{Alg}_{\mathrm{Ass}}(\mathcal{O}) \to \mathrm{Alg}_{\mathrm{Lie}}(\mathcal{O})$ have trivial Lie structure maps.

Poisson algebras also arise naturally in this setting: given a one parameter family of associative algebras $A_\hbar \in \mathrm{Alg}_{\mathrm{Ass}}(\mathrm{D}(\mathbb{K}[\hbar]))$ with central fibre $A_0 = A_\hbar|_{\{\hbar=0\}} \in \mathrm{Alg}_{\mathrm{Comm}}(\mathrm{Vect}_{\mathbb{K}})$ a commutative algebra, there is a canonical lift $A_0 \in \mathrm{Alg}_{\mathbb{P}_1}(\mathrm{Vect}_{\mathbb{K}})$ of $A_0$ to a Poisson algebra. By definition, a Poisson algebra is a commutative algebra with a Lie bracket $\{\cdot, \cdot\} : A_0 \otimes_{\mathbb{K}} A_0 \to A_0$ that acts as a bi-derivation of the product, which in the setting at had is defined by

$$(10.1.2) \qquad\qquad \{\cdot, \cdot\} = \frac{1}{\hbar} [\cdot, \cdot]_\hbar |_{\{\hbar=0\}}$$

where $[\cdot, \cdot]_\hbar : A_\hbar \otimes_{\mathbb{K}} A_\hbar \to A_\hbar$ is the commutator in $A_\hbar$.

In the following, we define the categories $\mathrm{Lie}^*(X)$, $\mathrm{Comm}^!(X)$ and $\mathrm{Cois}(X)$ of $\mathrm{Lie}^*$, $\mathrm{Comm}^!$ and Coisson algebras on a variety $X$. Further, we define functors

$$(10.1.3) \qquad\qquad \mathrm{Comm}^!(X) \to \mathrm{Alg}^{\mathrm{ch}}(X) \to \mathrm{Lie}^*(X)$$

analogous to those of Equation 10.1.1: The former will be the inclusion of $\mathrm{Comm}^!(X)$ as the full subcategory of $\mathrm{Alg}^{\mathrm{ch}}(X)$ on objects whose image under the latter functor have trivial $\mathrm{Lie}^*$ structure maps, which are precisely the commutative chiral algebras of Definition 7.0.9.

Similarly, Coisson algebras are by definition $\mathrm{Comm}^!$ algebras with an analogously compatible $\mathrm{Lie}^*$ bracket. Further, given a one parameter family of chiral algebras $A_\hbar \in \mathrm{Alg}^{\mathrm{ch}}(X)_{/\mathbb{K}[\hbar]}$ with central fibre $A_0 = A_\hbar|_{\{\hbar=0\}} \in \mathrm{Comm}^!(X)$ commutative, there is a canonical lift $A_0 \in \mathrm{Cois}(X)$ of $A_0$ to a Coisson algebra, with $\mathrm{Lie}^*$ bracket defined analogously as the first order approximation to the associated family of $\mathrm{Lie}^*$ algebras.

10.2. **$\mathrm{Lie}^*$ algebras.** To begin, we recall the rudiments of the theory of $\mathrm{Lie}^*$ algebras.

*Definition* 10.2.1. A $\mathrm{Lie}^*$ algebra on $X$ is a Lie algebra object in the operad $D(X)^*$ of $D$ modules on $X$ with the $*$ operad structure of Definition 5.2.5.

Let $\mathrm{Lie}^*(X) = \mathrm{Alg}_{\mathrm{Lie}}(D(X)^*)$ denote the category of $\mathrm{Lie}^*$ algebras on $X$.

*Remark* 10.2.2. Concretely, a $\mathrm{Lie}^*$ algebra on $X$ is given by a $D$ module $L \in D(X)$ together with iterated Lie bracket maps

$$b_I : \boxtimes_{i \in I} L \to \Delta_*^{(I)} L \ ,$$

in $D(X^I)$ for each $I \in \mathrm{fSet}$.

*Example* 10.2.3. Suppose that $L = \tilde{L}_{\mathcal{D}} := \tilde{L} \otimes_{\mathcal{O}_X} \mathcal{D}_X \in D(X)$ is an induced $D$ module on $\tilde{L} \in \mathrm{QCoh}(X)$. Then a $\mathrm{Lie}^*$ structure on $L$ is equivalent to a skew-symmetric bidifferential operator $b \in \mathrm{PDiff}(\tilde{L}, \tilde{L}; \tilde{L})$ satisfying the Jacobi identity with respect to the usual composition of polydifferential operators.

*Example* 10.2.4. Let $\mathfrak{g}$ be a finite type Lie algebra. Then $L = \mathfrak{g} \otimes_{\mathbb{K}} \mathcal{D}_X = \mathfrak{g} \otimes_{\mathbb{K}} \mathcal{O}_X \otimes_{\mathcal{O}_X} \mathcal{D}_X \in \mathrm{Lie}^*(X)$ is naturally a Lie* algebra under the $\mathcal{D}_X$ linear extension of the Lie bracket map on $\mathfrak{g}$. The corresponding bidifferential operator in this example is the $\mathcal{O}_X$ linear extension of the Lie bracket to $\tilde{L} = \mathfrak{g} \otimes_{\mathbb{K}} \mathcal{O}_X$.

*Example* 10.2.5. Let $\theta_X$ be the tangent sheaf of $X$. Then the Lie bracket $b \in \mathrm{PDiff}(\theta_X, \theta_X; \theta_X)$ is a bidifferential operator which defines a Lie* structure on $\theta_{X,\mathcal{D}} = \theta_X \otimes_{\mathcal{O}_X} \mathcal{D}_X$.

*Remark* 10.2.6. By Remark 5.3.10, there is a natural functor $\mathrm{Alg}^{\mathrm{ch}}(X) \to \mathrm{Lie}^*(X)$. Concretely, it is the identity on the underlying $D$ module on $X$, and sends the chiral structure map

$$j_*^{(I)} j^{(I),*}(\boxtimes_{i \in I} A) \xrightarrow{\mu^I} \Delta_*^{(I)} A \qquad \text{to the composition} \qquad \boxtimes_{i \in I} A \to j_*^{(I)} j^{(I),*}(\boxtimes_{i \in I} A) \xrightarrow{\mu^I} \Delta_*^{(I)} L \ .$$

This is the desired analogue of the forgetful functor $\mathrm{Alg}_{\mathrm{Ass}} \to \mathrm{Alg}_{\mathrm{Lie}}$, as outlined in Equation 10.1.3. Given a chiral algebra $A \in \mathrm{Alg}^{\mathrm{ch}}(X)$, we denote its associated Lie* algebra by $A^{\mathrm{Lie}} \in \mathrm{Lie}^*(X)$.

10.3. **Comm! algebras.** We now recall the elementary definitions for the chiral analogue of the theory of commutative algebras.

*Definition* 10.3.1. A Comm! algebra on $X$ is a commutative algebra object in the symmetric monoidal category $D(X)^!$ of $D$ modules on $X$ with the $\otimes^!$ tensor structure of Definition A.2.6.

Let $\mathrm{Comm}^!(X) = \mathrm{Alg}_{\mathrm{Comm}}(D(X)^!)$ denote the category of Comm! algebras on $X$.

*Remark* 10.3.2. Concretely, a Comm! algebra on $X$ is given by a $D$ module $A \in D(X)$ together with commutative multiplication maps

$$m^i : \otimes_{i \in I}^! A_i \to A$$

in $D(X)$ for each $I \in \mathrm{fSet}$.

*Remark* 10.3.3. For each $A \in \mathrm{Alg}^{\mathrm{ch}}(X)$, the natural excision exact triangle

(10.3.1) $\qquad A \boxtimes A \xrightarrow{i} j_* j^*(A \boxtimes A) \xrightarrow{q'} \Delta_* \Delta^!(A \boxtimes A)[1] \qquad$ induces

Thus, we see that the chiral product map $\mu$ factors through a map $m : A \otimes^! A \to A$ as indicated if and only if the induced Lie* bracket map $b$ vanishes. This is precisely the condition in the definition 7.0.9 of commutative chiral algebra.

In particular, this defines an equivalence $\mathrm{Comm}^!(X) \xrightarrow{\cong} \mathrm{Alg}^{\mathrm{ch}}(X)^{\mathrm{Comm}}$ between the category of Comm! algebras and that of commutative factorization algebras on $X$. Together with the above discussion, this yields the desired proposition:

*Proposition* 10.3.4. The category $\mathrm{Comm}^!(X)$ is equivalent to the full subcategory of $\mathrm{Alg}^{\mathrm{ch}}(X)$ on objects whose image under the forgetful functor to $\mathrm{Lie}^*(X)$ have trivial Lie* structure maps.

10.4. **Coisson algebras and filtered quantizations.** Next, we recall the definition of Coisson algebra and the notion of a filtered quantization of a Coisson algebra to a chiral algebra.

*Definition* 10.4.1. A Coisson algebra $R \in \mathrm{Cois}(X)$ on $X$ is a $\mathrm{Comm}^!$ algebra together with a Lie* bracket $\{\cdot, \cdot\} : R \boxtimes R \to \Delta_* R$ which is a derivation of the $\mathrm{Comm}^!$ product $m : R \otimes^! R \to R$.

Given a one parameter family of chiral algebras $A_\hbar \in \mathrm{Alg}^{\mathrm{ch}}(X)_{/\mathbb{K}[\hbar]}$ with central fibre given by $A_0 = A_\hbar|_{\{\hbar=0\}} \in \mathrm{Comm}^!(X)$ a commutative chiral algebra, there is a canonical Coisson structure on $A_0$, defined by

$$(10.4.1) \qquad\qquad \{\cdot, \cdot\} := \frac{1}{\hbar}(b_\hbar)|_{\{\hbar=0\}} \ ,$$

where $b_\hbar : A_0 \boxtimes A_0 \to \Delta_* A_\hbar$ is defined using the induced family of Lie* brackets, as the composition

$$A_0 \boxtimes A_0 \hookrightarrow A_\hbar \boxtimes A_\hbar \xrightarrow{\iota_\hbar} j_* j^* (A_\hbar \boxtimes A_\hbar) \xrightarrow{\mu_\hbar} \Delta_* A_\hbar \ .$$

The composition vanishes to first order in $\hbar$, so that the expression in Equation 10.4.1 is indeed well defined; this follows immidiately from commutativity of $A_0$.

*Definition* 10.4.2. A family of chiral algebras $A_\hbar \in \mathrm{Alg}^{\mathrm{ch}}(X)_{/\mathbb{K}[\hbar]}$ with central fibre $A_0 \in \mathrm{Comm}^!(X)$ commutative is called a filtered quantization of the associated Coisson algebra.

10.5. **Lie, commutative, and Poisson vertex algebras.** Now, as in Section 8, we restrict to the case $X = \mathbb{A}^1$ and describe the objects of the preceding subsection in the weakly $\mathbb{G}_a$ equivariant.

*Warning* 10.5.1. Throughout this section, all of the objects will be of cohomological degree zero (in the heart of the relevant t-structure) and all the functors non-derived, in contrast with our general conventions.

Recall the diagram of Equation 9.2.2 and the surrounding discussion in Example 9.2.9, which establishes the equivalence between the chiral product map

$$\mu : j_* j^* (A \boxtimes A) \to \Delta_* A \qquad \text{and the map} \qquad \tilde{\mu} : A \boxtimes A \to \tilde{\Delta}_* A \ ,$$

which is the global generalization of the operator product expansion of a vertex algebra. In particular, the diagram of Equation 8.0.2 which specializes that of Equation 9.2.2 in the special case $X = \mathbb{A}^1$, together with the discussion in the proof of Theorem 8.0.5, explains the equivalence between the chiral product map and the usual vertex algebra operator product structure map.

We now extend the vertical arrows of each of the diagrams in equations 9.2.2 and 8.0.2 by both of the short exact sequences of equations 10.3.1 and 9.1.1, which yields the following:

$$(10.5.1)$$

These diagrams summarize the interaction of the relations between chiral, Lie* and Comm$^!$ algebras discussed in the preceding subsection, and the passage from chiral algebras to vertex algebras. In particular, we use them to deduce the analogues of Theorem 8.0.5 in the commutative and Lie case: applying the equivalence of Proposition 8.0.4 to weakly $\mathbb{G}_a$ equivariant Comm$^!$ and Lie* algebras, we find the categories of such are equivalent to commutative vertex algebras and vertex Lie algebras, respectfully:

*Proposition* 10.5.2. There is an equivalence of categories between the category of weakly translation equivariant commutative chiral algebras on $\mathbb{A}^1$ and the category of commutative vertex algebras

$$\mathrm{Alg}_{\mathrm{un}}^{\mathrm{ch}}(\mathbb{A}^1)^{\mathrm{Comm},(\mathbb{G}_a,w)} \xrightarrow{\;\cong\;} \mathrm{VOA}_{\mathbb{K}}^{\mathrm{Comm}} \qquad \text{defined by} \qquad A \mapsto A_0 \qquad \mathbb{V} \otimes_{\mathbb{K}} \mathbb{K}[x] \hookleftarrow \mathbb{V},$$

where $A \in D(\mathbb{A}^1)^{\mathbb{G}_a,w}$ is such that $\Delta_*^{\mathrm{main}} A = \mathcal{L} \in D(\mathrm{Ran}_{\mathbb{A}^1})_{\mathbb{A}^1}^{\mathbb{G}_a,w}$ and $A_0 = \iota_0^! A$, as in Proposition 8.0.4.

*Proof.* We apply the proof of Theorem 8.0.5 in the commutative case: From the commutative diagram 10.5.1, we see that the condition of commutativity of the chiral algebra $\mu \circ \iota = 0$ is equivalent to the condition $q \circ \tilde{\mu} = 0$. Since the operator product structure map $Y : \mathbb{V}^{\otimes 2} \to \mathbb{V} \otimes_{\mathbb{K}} \mathbb{K}((z))$ is just the data of such a map $\tilde{\mu}$ in the $\mathbb{G}_a$ equivariant case, $q \circ \tilde{\mu} = 0$ if and only if the operator product structure map is nonsingular. This latter condition is precisely the definition of a commutative vertex algebra. $\qquad\square$

*Remark* 10.5.3. Alternatively, the category Comm$^!(X)$ is equivalent to the category of affine $D$ schemes on $X$. In the case $X = \mathbb{A}^1$, the $\mathbb{G}_a$ equivariant $D$ schemes are by definition commutative algebra objects in $D(\mathbb{A}^1)^{!,(\mathbb{G}_a,w)}$. Applying Proposition 8.0.4, we obtain an equivalence with the category of commutative algebra objects in $\mathbb{K}[\partial]$-Mod, which is equivalent to the category of commutative vertex algebras by Proposition E.2.2.

*Proposition* 10.5.4. There is an equivalence of categories between the category of weakly translation equivariant Lie* algebras on $\mathbb{A}^1$ and the category of vertex Lie algebras

$$\mathrm{Lie}^*(\mathbb{A}^1)^{\mathbb{G}_a,w} \xrightarrow{\;\cong\;} \mathrm{VLA}_{\mathbb{K}} \qquad \text{defined by} \qquad L \mapsto L_0 \qquad \mathbb{V} \otimes_{\mathbb{K}} \mathbb{K}[x] \hookleftarrow \mathbb{V},$$

where $L \in D(\mathbb{A}^1)^{\mathbb{G}_a,w}$ is the $D$ module underlying the Lie* algebra and $L_0 = \iota_0^! L$, as in Proposition 8.0.4.

*Proof.* Again, the proof is essentially the same as that of Theorem 8.0.5. The bracket map

$$b : \mathbb{V}^{\otimes 2} \otimes_{\mathbb{K}} \mathbb{K}[x, y] \to \mathbb{V} \otimes_{\mathbb{K}} \mathbb{K}[x] \otimes_{\mathbb{K}} \delta_{x-y} \qquad \text{is determined on generators by} \qquad Y_- : \mathbb{V}^{\otimes 2} \to \mathbb{V} \otimes_{\mathbb{K}} \delta_{x-y},$$

by translation invariance. The latter is precisely the vertex Lie structure map, as desired. The remaining properties are checked as in the proof of Theorem 8.0.5, though the computation appears to be more involved to match the particular conventions from [FBZ04] which are listed in Definition E.3.1. $\qquad\square$

*Remark* 10.5.5. From the diagram of Equation 10.5.1 together with the above discussion, it is apparent that the forgetful functor $\mathrm{Alg}^{\mathrm{ch}}(X) \to \mathrm{Lie}^*(X)$ corresponds to remembering only the singular part of the OPE, which is the usual forgetful functor from vertex algebras to vertex Lie algebras.

*Example* 10.5.6. Suppose $L = \tilde{L}_{\mathcal{D}}$ an induced $D$ module on a translation invariant $\mathcal{O}_X$ module $\tilde{L}$ over $X = \mathbb{A}^1$, with fibre $\tilde{L}_0$ at the point $0 \in \mathbb{A}^1$. Then we have

$$L_0 = \tilde{L}_0 \otimes_{\mathbb{K}} \mathbb{K}[\partial]$$

is a free $\mathbb{K}[\partial]$ module, so that the bracket map is determined on $\tilde{L}_0$ by a map

$$\tilde{b} : \tilde{L}_0^{\otimes 2} \to \tilde{L}_0 \otimes \mathbb{K}[\partial] \otimes \delta_{x-y} \ .$$

This is equivalent to the data of a translation invariant bidifferential operator $\tilde{b}_0 \in \mathrm{PDiff}(\tilde{L}, \tilde{L}; \tilde{L})^{\mathbb{G}_a, w}$, as follows from Example 10.2.3.

*Example* 10.5.7. Let $\tilde{L}_0 = \mathfrak{g}$ a finite dimensional Lie algebra with Lie bracket $b_{\mathfrak{g}} : \mathfrak{g}^{\otimes 2} \to \mathfrak{g}$. Then the map

$$\tilde{b}_0 = b_{\mathfrak{g}} \otimes 1 : \tilde{L}_0^{\otimes 2} \to \tilde{L}_0 \otimes 1 \otimes \delta_{x-y}^{(0)} \hookrightarrow \tilde{L}_0 \otimes \mathbb{K}[\partial] \otimes \delta_{x-y}$$

defines a vertex Lie algebra, which corresponds to the Lie* algebra $\mathfrak{g} \otimes_{\mathbb{K}} \mathcal{D}_X \in \mathrm{Lie}^*(X)$ of Example 10.2.4 under Proposition 10.5.4. Rewritten in terms of the vertex Lie structure map $Y_- : L_0^{\otimes 2} \to L_0 \otimes_{\mathbb{K}} \delta_{x-y}$, this reads

$$(10.5.2) \qquad Y_-(J_{-1}^a, z)(J_{-1}^b) = \frac{J_{-1}^{[a,b]}}{z}$$

where $J_{-1}^a = a \otimes 1 \in \mathfrak{g} \otimes_{\mathbb{K}} \mathbb{K}[\partial] = L_0$ for each $a \in \mathfrak{g}$, recalling $\delta_{x-y} \cong z^{-1}\mathbb{K}[z^{-1}]$ so that $\delta_{x-y}^0 \mapsto z^{-1}$.

*Example* 10.5.8. Let $\tilde{L} = \theta_X$ be the tangent sheaf of $X = \mathbb{A}^1$. The Lie bracket of vector fields defines a bidifferential operator $\tilde{b} \in \mathrm{PDiff}(\theta_X, \theta_X; \theta_X)$ so that $L = \theta_{X, \mathcal{D}}$ defines a Lie* algebra on $X$.

For $X = \mathbb{A}^1$, we have $\theta_X \cong \mathbb{K}[x] \cdot \partial_x$ so that the fiber of $\tilde{L}$ over a fixed point $\tilde{L}_0 = \mathbb{K}$ is one dimensional. The Lie bracket is given by

$$\tilde{b} : \mathbb{K}[x]^{\otimes 2} \to \mathbb{K}[x] \otimes \mathbb{K}[\partial] \otimes \delta_{x-y} \qquad\qquad f \otimes g \mapsto ((\partial_x f)g - f(\partial_y g)) \otimes \delta_{x-y}$$

Now, applying the usual vertex algebra convention of fixing coordinates $(x, x - y)$ on $X^2$, as in Example 9.1.3, we have

$$\partial_x \mapsto \partial_x + \partial_{x-y} \qquad \partial_y \mapsto -\partial_{x-y} \qquad \text{so that} \qquad \partial_x - \partial_y \mapsto \partial_x + 2\partial_{x-y}$$

and thus the corresponding vertex Lie structure map is given by

$$(10.5.3) \qquad Y_-(l_{-2}, z)(l_{-2}) = \frac{l_{-3}}{z} + \frac{2l_{-2}}{z^2} \ ,$$

where $l_{-2} = 1 \in L_0 \cong \mathbb{K}[\partial]$ and $l_{-3} = \partial \in L_0$; this is chosen to match the usual notation for the Virasoro algebra generator $l_{-2}$ and its image under the translation operator $l_{-3}$.

We can also combine the above results to deduce the analogous relation between Coisson algebras and Poisson vertex algebras:

*Proposition* 10.5.9. There is an equivalence of categories between the category of weakly translation equivariant Coisson algebras on $\mathbb{A}^1$ and the category of Poisson vertex algebras

$$\mathrm{Cois}(\mathbb{A}^1)^{\mathbb{G}_a, w} \underset{\xleftarrow{\hspace{1cm}}}{\overset{\cong}{\xrightarrow{\hspace{1cm}}}} \mathrm{PVA}_{\mathbb{K}} \qquad \text{defined by} \qquad R \mapsto R_0 \qquad \mathbb{V} \otimes_{\mathbb{K}} \mathbb{K}[x] \leftarrow\!\shortmid \mathbb{V},$$

where $R \in D(\mathbb{A}^1)^{\mathbb{G}_a, w}$ is the $D$ module underlying the Coisson algebra and $R_0 = \iota_0^! R$, as in Proposition 8.0.4.

## 11. From Chiral algebras to topological associative algebras

In this section, together with the complementary Appendix D, we summarize the results of Chapter 3.6 of [BD04] and the closely related paper [Bei07].

*Warning* 11.0.1. Throughout this section, all of the objects will be of cohomological degree zero (in the heart of the relevant t-structure) and all the functors non-derived, in contrast with our general conventions.

### 11.1. **Modifications and topologies at a point on $D$ modules.**

In this subsection, we summarize the results of Section 2.1.13 of [BD04]. Throughout, let $X$ be an algebraic curve, $x \in X(\mathbb{K})$ a smooth, closed point inducing the complementary closed and open embeddings

$$\iota_x : \{x\} \hookrightarrow X \qquad j : U_x \hookrightarrow X \qquad \text{where} \qquad U_x = X \backslash \{x\} \ ,$$

and let $M \in D(X)^\heartsuit$ be a $D$ module.

*Definition* 11.1.1. Define the space $\Xi_x(M)$ of modifications of $M \in D(X)^\heartsuit$ at the point $x \in X(\mathbb{K})$ by

$$\Xi_x(M) = \{M_\xi \hookrightarrow M \text{ a } D_X \text{ submodule} \mid \text{supp}(M/M_\xi) = \{x\}\} \ .$$

Following [BD04], we let $h = \mathrm{dR}^0 : D(X)^\heartsuit \to \mathrm{Sh}_c(X)$ denote the zeroth cohomology of the the de Rham sheaf functor, as in Definition A.4.4, and $h_x = \iota_x^! : D(X)^\heartsuit \to \mathbb{K}\text{-Mod}$, so that $h_x(M) = \iota_x^! h(M) = h(M)_x$.

*Proposition* 11.1.2. For each $\xi \in \Xi_x(M)$, we have a short exact sequence

$$h(M_\xi)_x \hookrightarrow h(M)_x \twoheadrightarrow h(M/M_\xi)_x \ .$$

*Definition* 11.1.3. The $\Xi$ topology is the topology on $h(M)_x$ with basis of neighbourhoods of 0 given by the subspaces $h(M_\xi)_x \subset h(M)$ for each $\xi \in \Xi_x(M)$.

Let $\hat{h}_x : D(X; \mathrm{Top}_x) \to \mathbb{K}\text{-Mod}_{\mathrm{Top}}$ denote the functor taking $M$ to the completion of $h(M)_x$ in the $\Xi$ topology. See Appendix D for a review and conventions regarding functional analysis.

*Remark* 11.1.4. The completion of $h(M)_x$ in the $\Xi$ topology is given by

$$\hat{h}_x(M) = \lim_{\xi \in \Xi_x(M)} h(M)_x / h(M_\xi)_x \cong \lim_{\xi \in \Xi_x(M)} h(M/M_\xi)_x \ ,$$

by Proposition 11.1.2.

*Proposition* 11.1.5. Let $V \in \mathbb{K}\text{-Mod} \hookrightarrow \mathbb{K}\text{-Mod}_{\mathrm{Top}}$ be a discrete vector space. Then

$$\mathrm{Hom}_{\mathbb{K}\text{-Mod}_{\mathrm{Top}}}(\hat{h}_x(M), V) = \mathrm{Hom}_{D(X)}(M, \iota_{x*}V) \ .$$

*Definition* 11.1.6. A topology on $M$ at $x$ is defined to be a topology $\Xi_x^?(M)$ on $h(M)_x$ that is coarser than the $\Xi$ topology.

Let $D(X; \mathrm{Top}_x)$ denote the category of $D$ modules $M \in D(X)^\heartsuit$ together with a choice of topology $\Xi_x^?(M)$ on $M$ at $x$, with morphisms given by maps of $D$ modules inducing continuous maps on the corresponding completed topological vector spaces. Let $\hat{h}_x : D(X; \mathrm{Top}_x) \to \mathbb{K}\text{-Mod}_{\mathrm{Top}}$ denote the functor of taking the completion $\hat{h}_x^?(M) \in \mathbb{K}\text{-Mod}_{\mathrm{Top}}$ of $h(M)_x$ with respect to $\Xi_x^?(M)$.

*Remark* 11.1.7. As in remark 11.1.4, we have

$$\hat{h}_x^?(M) = \lim_{\xi \in \Xi_x^?(M)} h_x(M/M_\xi) \ .$$

11.2. **From chiral algebras to topological associative algebras.** We now recall the main result facilitating the construction of topological associative algebras from chiral algebras, and its variants for Lie*, Comm[!] and Coisson algebras.

*Proposition* 11.2.1. The functor $\hat{h}_x : D(X; \mathrm{Top}_x) \to \mathbb{K}\text{-Mod}_{\mathrm{Top}}$ that assigned to $M$ the completion of $h(M)_x$ in the $\Xi$ topology lifts to maps of operads

$$\hat{h}_x : D(X; \mathrm{Top}_x)^* \to \mathbb{K}\text{-Mod}_{\mathrm{Top}}^*$$

$$\hat{h}_x : D(X; \mathrm{Top}_x)^! \to \mathbb{K}\text{-Mod}_{\mathrm{Top}}^!$$

$$\hat{h}_x : D(X; \mathrm{Top}_x)^{\mathrm{ch}} \to \mathbb{K}\text{-Mod}_{\mathrm{Top}}^{\mathrm{ch},s}$$

where the three operad structures on $\mathbb{K}\text{-Mod}_{\mathrm{Top}}$ are as defined in remarks D.2.2, D.2.11, and Definition D.2.8, respectively.

The desired result follows from the preceding proposition:

*Corollary* 11.2.2. Let $A \in \mathrm{Alg}^{\mathrm{ch}}(X)$ be a chiral algebra on $X$, and let

$$\Xi_x^{\mathrm{as}}(A) = \{A_\xi \hookrightarrow A \text{ a chiral subalgebra} \mid \mathrm{supp}(A/A_\xi) = \{x\}\} .$$

Then $A_x^{\mathrm{as}} := \hat{h}_x(A, \Xi_x^{\mathrm{as}}) \in \mathrm{Alg}_{\mathrm{Ass}}(\mathbb{K}\text{-Mod}_{\mathrm{Top}}^{\mathrm{ch},s})$ defines a topological associative algebra with respect to the $\otimes^{\mathrm{ch}}$ tensor structure.

We now state the analogues for Lie*, Comm[!] and Coisson algebras:

*Corollary* 11.2.3. Let $L \in \mathrm{Lie}^*(X)$ be a Lie* algebra on $X$, and let

$$\Xi_x^{\mathrm{Lie}}(L) = \{L_\xi \hookrightarrow L \text{ a Lie* subalgebra} \mid \mathrm{supp}(L/L_\xi) = \{x\}\} .$$

Then $L_x^{\mathrm{Lie}} := \hat{h}_x(A, \Xi_x^{\mathrm{Lie}}) \in \mathrm{Alg}_{\mathrm{Lie}}(\mathbb{K}\text{-Mod}_{\mathrm{Top}}^*)$ defines a topological Lie algebra with respect to the $\otimes^*$ tensor structure.

*Corollary* 11.2.4. Let $R \in \mathrm{Comm}^!(X)$ be a Comm[!] algebra on $X$, and let

$$\Xi_x^{\mathrm{Comm}}(R) = \{R_\xi \hookrightarrow R \text{ a Comm}^! \text{ subalgebra} \mid \mathrm{supp}(R/R_\xi) = \{x\}\} .$$

Then $R_x^{\mathrm{Comm}} := \hat{h}_x(R, \Xi_x^{\mathrm{Comm}}) \in \mathrm{Alg}_{\mathrm{Ass}}(\mathbb{K}\text{-Mod}_{\mathrm{Top}}^!)$ defines a topological associative algebra with respect to the $\otimes^!$ tensor structure.

## 12. Chiral enveloping algebras of Lie* algebras

In Section 10, we recalled the notion of Lie* algebra on $X$, and explained that there was a forgetful functor

$$\mathrm{Alg}^{\mathrm{ch}}(X) \xrightarrow{(\cdot)^{\mathrm{Lie}}} \mathrm{Lie}^*(X) \qquad \text{analogous to the functor} \qquad \mathrm{Alg}_{\mathrm{Ass}}(\mathcal{O}) \xrightarrow{(\cdot)^{\mathrm{Lie}}} \mathrm{Alg}_{\mathrm{Lie}}(\mathcal{O})$$

from associative to Lie algebras given by remembering only commutators; these are defined naturally internal to an ambient operad $\mathcal{O}$, which we take to be $\mathrm{Vect}_{\mathbb{K}}$ with its usual symmetric monoidal structure for the remaining discussion.

The universal enveloping algebra of a Lie algebra can be understood as providing a left adjoint $\mathcal{U} : \mathrm{Alg}_{\mathrm{Lie}}(\mathrm{Vect}_{\mathbb{K}}) \to \mathrm{Alg}_{\mathrm{Ass}}(\mathrm{Vect}_{\mathbb{K}})$ to the above forgetful functor; for any associative algebra $A \in \mathrm{Alg}_{\mathrm{Ass}}(\mathrm{Vect}_{\mathbb{K}})$, and Lie algebra $\mathfrak{g} \in \mathrm{Alg}_{\mathrm{Lie}}(\mathrm{Vect}_{\mathbb{K}})$, there are natural isomorphisms

$$\mathrm{Hom}_{\mathrm{Alg}_{\mathrm{Lie}}(\mathrm{Vect}_K)}(\mathfrak{g}, A^{\mathrm{Lie}}) \xrightarrow{\cong} \mathrm{Hom}_{\mathrm{Alg}_{\mathrm{Ass}}(\mathrm{Vect}_K)}(\mathcal{U}(\mathfrak{g}), A) .$$

In this subsection, we recall the chiral enveloping algebra functor $\mathcal{U}^{\mathrm{ch}} : \mathrm{Lie}^*(X) \to \mathrm{Alg}^{\mathrm{ch}}(X)$, as well as its unital and twisted variants, which satisfy the analogous adjunction

$$\mathrm{Hom}_{\mathrm{Lie}^*(X)}(L, A^{\mathrm{Lie}}) \xrightarrow{\cong} \mathrm{Hom}_{\mathrm{Alg}^{\mathrm{ch}}(X)}(\mathcal{U}^{\mathrm{ch}}(L), A) \ ,$$

naturally for each $A \in \mathrm{Alg}^{\mathrm{ch}}(X)$ and $L \in \mathrm{Lie}^*(X)$.

*Definition* 12.0.1. Let $L \in \mathrm{Lie}^*(X)$ be a Lie$^*$ algebra on $X$. The non-unital chiral enveloping algebra $\mathcal{U}^{\mathrm{ch}}(L) \in \mathrm{Alg}^{\mathrm{ch}}(X)$ of $L$ is the non-unital chiral algebra corresponding to the non-unital factorization algebra

$$\tilde{C}^{\mathrm{CE},\otimes*}_\bullet(\Delta^{\mathrm{main}}_* L) \in \mathrm{Alg}^{\mathrm{fact}}(X) \qquad \text{where} \qquad \tilde{C}^{\mathrm{CE},\otimes*}_\bullet : \mathrm{Alg}_{\mathrm{Lie}}(D(\mathrm{Ran}_X)^*) \to \mathrm{CoAlg}_{\mathrm{Comm}^{\mathrm{nu}}}(D(\mathrm{Ran}_X)^*)$$

denotes the reduced Chevalley-Eilenberg chains object internal to $D(\mathrm{Ran}_X)^*$, and $\Delta^{\mathrm{main}}_* L \in \mathrm{Alg}_{\mathrm{Lie}}(D(\mathrm{Ran}_X)^*_X)$ is the Lie algebra object corresponding to $L \in \mathrm{Lie}^*(X)$.

The (unital) chiral enveloping algebra $\mathcal{U}^{\mathrm{ch}}_{\mathrm{un}}(L) \in \mathrm{Alg}^{\mathrm{ch}}_{\mathrm{un}}(X)$ of $L$ is the unitalization of the above non-unital variant, in the sense of Proposition 7.0.6.

*Remark* 12.0.2. The object $\tilde{C}^{\mathrm{CE},\otimes*}_\bullet(\Delta^{\mathrm{main}}_* L) \in D(\mathrm{Ran}_X)$ carries a natural non-unital, cocommutative coalgebra structure internal to $D(\mathrm{Ran}_X)^*$ as it is the reduced Chevalley-Eilenberg chains on a Lie algebra object, which induces such a structure internal to $D(\mathrm{Ran}_X)^{\mathrm{ch}}$ by applying the forgetful functor on coalgebras given by Proposition 5.3.5. The fact that the resulting cocommutative coalgebra in the chiral tensor structure is factorizeable is demonstrated in the proof of Theorem 6.4.2 in [FG11].

*Remark* 12.0.3. The $D$ module on $X$ underlying the chiral algebra $\mathcal{U}^{\mathrm{ch}}(L)$ is given by

$$\mathcal{U}^{\mathrm{ch}}(L) = \Delta^{\mathrm{main},!}\tilde{C}^{\mathrm{CE},\otimes*}_\bullet(\Delta^{\mathrm{main}}_* L)[-1] \cong \bigoplus_{k \in \mathbb{Z}_{>0}} \Delta^{\mathrm{main}!}\mathrm{Sym}^{k,\otimes*}(\Delta^{\mathrm{main}}_* L) \cong \bigoplus_{k \in \mathbb{Z}_{>0}} \mathrm{Sym}^k_!(L)$$

where the first isomorphism follows from Proposition 5.2.3, and $\mathrm{Sym}^k_!$ denotes the symmetric power of $L \in D(X)^!$. Similarly, the $D$ module underlying the unital chiral enveloping algebra is given by

$$\mathcal{U}^{\mathrm{ch}}_{\mathrm{un}}(L) \cong \mathrm{Sym}^\bullet_!(L) \ .$$

We now state the anticipated universal property, which is proved in [BD04] and [FG11]:

*Proposition* 12.0.4. There is a natural equivalence

$$\mathrm{Hom}_{\mathrm{Lie}^*(X)}(L, A^{\mathrm{Lie}}) \xrightarrow{\cong} \mathrm{Hom}_{\mathrm{Alg}^{\mathrm{ch}}(X)}(\mathcal{U}^{\mathrm{ch}}(L), A) \ ,$$

for each $A \in \mathrm{Alg}^{\mathrm{ch}}(X)$ and $L \in \mathrm{Lie}^*(X)$. The analogous statement holds for $\mathcal{U}^{\mathrm{ch}}_{\mathrm{un}}(L) \in \mathrm{Alg}^{\mathrm{ch}}_{\mathrm{un}}(X)$.

*Example* 12.0.5. Suppose $L = \tilde{L}_\mathcal{D} = \tilde{L} \otimes_{\mathcal{O}_X} \mathcal{D}_X \in \mathrm{Lie}^*(\mathbb{A}^1)^{\mathbb{G}_a,w}$ an induced $D$ module on a translation invariant $\mathcal{O}_X$ module $\tilde{L}$ over $X = \mathbb{A}^1$, as in Example 10.5.6, so that the fibre at the point $0 \in \mathbb{A}^1$ is given by $L_0 = \tilde{L}_0 \otimes_\mathbb{K} \mathbb{K}[\partial]$.

The vertex algebra corresponding to the translation invariant chiral enveloping algebra $A = \mathcal{U}^{\mathrm{ch}}_{\mathrm{un}}(L) \in \mathrm{Alg}^{\mathrm{ch}}_{\mathrm{un}}(\mathbb{A}^1)^{\mathbb{G}_a,w}$ under Theorem 8.0.5 has underlying vector space given by

$$\mathbb{V} = A_0 = \iota^!_0(\mathrm{Sym}^\bullet_!(L)) \cong \mathrm{Sym}^\bullet_\mathbb{K}(\tilde{L}_0 \otimes_\mathbb{K} \mathbb{K}[\partial]) \cong \mathrm{Sym}^\bullet_\mathbb{K}(\tilde{L}_0 \otimes_\mathbb{K} z^{-1}\mathbb{K}[z^{-1}]) \ .$$

*Example* 12.0.6. Let $\mathfrak{g}$ be a finite type Lie algebra over $\mathbb{K}$ and $L = \mathfrak{g}_\mathcal{D} \in \mathrm{Lie}^*(X)$ be the induced Lie$^*$ algebra, as in Example 10.5.7. The untwisted affine chiral algebra $A_0(\mathfrak{g}) = \mathcal{U}^{\mathrm{ch}}_{\mathrm{un}}(L) \in \mathrm{Alg}^{\mathrm{ch}}_{\mathrm{un}}(X)$

is defined as the chiral enveloping algebra of $L$. In particular, the corresponding vertex algebra is given by

$$\mathbb{V}_0(\mathfrak{g}) = \mathrm{Sym}_{\mathbb{K}}^{\bullet}(z^{-1}\mathfrak{g}[z^{-1}]) \qquad \text{with} \qquad Y(J_{-1}^a, z)(J_{-1}^b) = \frac{J_{-1}^{[a,b]}}{z}$$

as the singular part of the operator product map, where $J_{-1}^a = a \otimes z^{-1} \in \mathbb{V}_0(\mathfrak{g})$ is the usual generator corresponding to $a \in \mathfrak{g}$. The defining property of the chiral envelope is that it is freely generated subject to the relations imposed by the fixed singular terms, so the preceeding expression follows by definition from Equation 10.5.2.

*Example* 12.0.7. Let $L = \theta_{X,\mathcal{D}} \in \mathrm{Lie}^*(X)$ be the Lie$^*$ algebra given by the induced $D$ module on the tangent sheaf $\theta_X$ under Lie bracket of vector fields, as in Example 10.5.8. The untwisted Virasoro chiral algebra $A = \mathcal{U}_{\mathrm{un}}^{\mathrm{ch}}(L) \in \mathrm{Alg}_{\mathrm{un}}^{\mathrm{ch}}(X)$ is defined as the chiral enveloping algebra of $L$. In particular, the corresponding vertex algebra is given by

$$\mathrm{Vir}_0 = \mathrm{Sym}_{\mathbb{K}}^{\bullet}(z^{-1}\mathbb{K}[z^{-1}]) \qquad \text{with} \qquad Y(l_{-2}, z)(l_{-2}) = \frac{l_{-3}}{z} + \frac{2l_{-2}}{z^2}$$

as the singular part of the operator product map, where $l_{-2} = z^{-1}$ and $l_{-3} = -2z^{-2}$ are the usual Virasoro generators. As in the preceding example, this follows from the analogous Equation 10.5.3.

**12.1. $\omega_X$ extensions of Lie$^*$ algebras and twisted chiral enveloping algebras.** We now recall the variant of the results of the preceding subsection twisted by a central extension of the relevant Lie algebra object.

*Definition* 12.1.1. Let $L \in \mathrm{Lie}^*(X)$ be a Lie$^*$ algebra on $X$. An $\omega_X$ extension of $L$ is a Lie$^*$ algebra $L^{\flat} \in \mathrm{Lie}^*(X)$ fitting into a short exact sequence of Lie$^*$ algebras $\omega_X \hookrightarrow L^{\flat} \twoheadrightarrow L$.

*Remark* 12.1.2. An $\omega_X$-extension is necessarily central, as the canonical map

$$\mathrm{Hom}_{D(X)^*}(L, L; L^{\flat}) \xrightarrow{\cong} \mathrm{Hom}_{D(X)^*}(L^{\flat}, L^{\flat}; L^{\flat})$$

is an equivalence.

*Definition* 12.1.3. Let $L \in \mathrm{Lie}^*(X)$ be a Lie$^*$ algebra on $X$ with $\omega_X$ extension $L^{\flat} \in \mathrm{Lie}^*(X)$ and let $1^{\flat} \in H^0_{\mathrm{dR}}(L^{\flat})$ be the unit section of $\omega_X \hookrightarrow L^{\flat}$. The twisted chiral envelope $A = \mathcal{U}_{\mathrm{tw}}^{\mathrm{ch}}(L) \in \mathrm{Alg}^{\mathrm{ch}}(X)$ is the unital chiral algebra quotient of $\mathcal{U}_{\mathrm{un}}^{\mathrm{ch}}(L^{\flat})$ with kernel the ideal generated by $1 - 1^{\flat} \in H^0_{\mathrm{dR}}(\mathcal{U}_{\mathrm{un}}^{\mathrm{ch}}(L^{\flat}))$, where $1$ is the unit section of $\mathcal{U}_{\mathrm{un}}^{\mathrm{ch}}(L^{\flat})$.

*Remark* 12.1.4. The $D$ modules underlying $A = \mathcal{U}_{\mathrm{un}}^{\mathrm{ch}}(L)$ and $\mathcal{U}_{\mathrm{tw}}^{\mathrm{ch}}(L)$ are isomorphic, though not canonically unless the extension is trivial.

The following is the corresponding universal property of the twisted chiral enveloping algebra construction:

*Proposition* 12.1.5. There is a natural equivalence

$$\{\varphi \in \mathrm{Hom}_{\mathrm{Lie}^*(X)}(L^{\flat}, A^{\mathrm{Lie}}) \mid \varphi(1^{\flat}) = 1_A\} \xrightarrow{\cong} \mathrm{Hom}_{\mathrm{Alg}^{\mathrm{ch}}(X)}(\mathcal{U}_{\mathrm{tw}}^{\mathrm{ch}}(L), A) ,$$

for each $A \in \mathrm{Alg}^{\mathrm{ch}}(X)$ and $L \in \mathrm{Lie}^*(X)$, and $\omega_X$ extension

*Example* 12.1.6. Let $\mathfrak{g}$ be a finite type Lie algebra with non-degenerate, ad-invariant bilinear form $\kappa : \mathfrak{g}^{\otimes 2} \to \mathbb{K}$, and let $L = \mathfrak{g}_{\mathcal{D}} \in \mathrm{Lie}^*(X)$ be the induced Lie$^*$ algebra, as in Example 10.5.7. There is a canonical $\omega_X$ extension of $L$, called the Kac-Moody extension, and we define the affine Kac-Moody chiral algebra by $A_c(\mathfrak{g}) = \mathcal{U}_{\text{tw}}^{\text{ch}}(L) \in \text{Alg}_{\text{un}}^{\text{ch}}(X)$ as the twisted chiral enveloping algebra corresponding to the $\omega_X$ extension determined by $c \cdot \kappa$ for $c \in \mathbb{K}$; see e.g. 2.5.9 in [BD04].

On $X = \mathbb{A}^1$, the corresponding vertex algebra has the same underlying vector space as in the case $c = 0$, which reduces to Example 12.0.6; we have

$$\mathbb{V}_c(\mathfrak{g}) = \text{Sym}_{\mathbb{K}}^{\bullet}(\mathfrak{g} \otimes_{\mathbb{K}} z^{-1}\mathbb{K}[z^{-1}]) \qquad \text{with} \qquad Y(J_{-1}^a, z)(J_{-1}^b) = \frac{J_{-1}^{[a,b]}}{z} + \frac{c \, \kappa(a,b)}{z^2}$$

as the singular part of the operator product map.

*Example* 12.1.7. Let $\mathfrak{h}$ be a finite type, abelian Lie algebra with non-degenerate bilinear form $\kappa : \mathfrak{h}^{\otimes 2} \to \mathbb{K}$. In this case, the affine Kac-Moody chiral algebra is called the Heisenberg algebra $\text{Heis}_c(X, \mathfrak{h}) \in \text{Alg}_{\text{un}}^{\text{fact}}(X)$, and its corresponding vertex algebra is given by

$$\text{Heis}_c(\mathfrak{h}) = \text{Sym}_{\mathbb{K}}^{\bullet}(\mathfrak{h} \otimes_{\mathbb{K}} z^{-1}\mathbb{K}[z^{-1}]) \qquad \text{with} \qquad Y(J_{-1}^a, z)(J_{-1}^b) = \frac{c \, \kappa(a,b)}{z^2}$$

as the singular part of the operator product map.

*Example* 12.1.8. Let $L = \theta_{X,\mathcal{D}} \in \text{Lie}^*(X)$ be the Lie* algebra induced from the tangent sheaf $\theta_X$ under Lie bracket of vector fields, as in Example 10.5.8. There is a canonical $\omega_X$ extension of $L$, called the Virasoro extension, and we define the Virasoro chiral algebra by $\text{Vir}_c(X) = \mathcal{U}_{\text{tw}}^{\text{ch}}(L) \in \text{Alg}_{\text{un}}^{\text{ch}}(X)$ as the twisted chiral enveloping algebra corresponding to the $\omega_X$ extension determined by $c \cdot \kappa$; see e.g. 2.5.10 in [BD04].

On $X = \mathbb{A}^1$, the corresponding vertex algebra has the same underlying vector space as in the case $c = 0$, which reduces to Example 12.0.7; we have

$$\text{Vir}_c = \text{Sym}_{\mathbb{K}}^{\bullet}(z^{-1}\mathbb{K}[z^{-1}]) \qquad \text{with} \qquad Y(l_{-2}, z)(l_{-2}) = \frac{l_{-3}}{z} + \frac{2l_{-2}}{z^2} + \frac{\frac{1}{2}c}{z^4}$$

as the singular part of the operator product map.

*Example* 12.1.9. Let $X$ be a smooth curve, $\omega_X^{\frac{1}{2}}$ be a spin structure on $X$, and $V$ be a symplectic vector space. Consider the $D$ module $V \otimes_{\mathbb{K}} \omega_X^{\frac{1}{2}} \otimes_{\mathcal{O}_X} \mathcal{D}_X \in D(X)$ as defining an abelian Lie* algebra $L_V \in \text{Lie}^*(X)$. Then the bidifferential operator

$$\omega_V \otimes_{\mathbb{K}} (-) \wedge (-) \in \text{PDiff}(\omega_X^{\frac{1}{2}} \otimes_{\mathbb{K}} V, \omega_X^{\frac{1}{2}} \otimes_{\mathbb{K}} V; \omega_X) \ ,$$

given by applying the symplectic form $\omega_V : V^{\otimes 2} \to \mathbb{K}$ together with the multiplication of algebraic densities, defines an $\omega_X$ extension of $L_V$.

The chiral Weyl algebra $\mathcal{W}^{\text{ch}}(X, V) = \mathcal{U}_{\text{tw}}^{\text{ch}}(L_V)$ is defined as the twisted chiral enveloping algebra corresponding to this $\omega_X$ extension; see e.g. 3.8.1 in [BD04].

On $X = \mathbb{A}^1$, the corresponding vertex algebra is given by

$$\mathcal{W}^{\text{ch}}(V) = \text{Sym}_{\mathbb{K}}^{\bullet}(V \otimes_{\mathbb{K}} z^{-1}\mathbb{K}[z^{-1}]) \qquad \text{with} \qquad Y(\varphi_{-1}^v; z)\varphi_{-1}^w = \frac{\omega_V(v,w)}{z}$$

as the singular part of the operator product map, where $\varphi_{-1}^v = v \otimes z^{-1} \in \mathcal{W}(V)$.

*Example* 12.1.10. For $V = T^*N = N \oplus N^{\vee}$ a cotangent vector space with its canonical symplectic form, a variant of the chiral Weyl algebra of the previous example can be defined independent of a choice of spin structure. Consider instead the trivial Lie* algebra $L_N = \tilde{L}_{N,\mathcal{D}}$ induced from $\tilde{L}_N = (N \otimes_{\mathbb{K}} \mathcal{O}_X) \oplus (N^{\vee} \otimes_{\mathbb{K}} \omega_X)$ together with the evident analogue of the above $\omega_X$ extension.

The resulting chiral algebra $\mathcal{W}^{\mathrm{ch}}(X, N, N^\vee) \in \mathrm{Alg}^{\mathrm{ch}}_{\mathrm{un}}(X)$ is called the '$\beta$-$\gamma$ system on $N$' or 'linear chiral differential operators on $N$'. On $X = \mathbb{A}^1$, the corresponding vertex algebra is given by

(12.1.1)
$$\mathcal{W}^{\mathrm{ch}}(N) = \mathrm{Sym}^\bullet_{\mathbb{K}}\big(\big[N \otimes_{\mathbb{K}} z^{-1}\mathbb{K}[z^{-1}]\big] \oplus \big[N^\vee \otimes_{\mathbb{K}} z^{-1}\mathbb{K}[z^{-1}]\big] \cdot dz\big) \qquad \text{with} \qquad Y(\varphi^n_{-1}; z)\varphi^{*,\xi}_0 = \frac{\xi(n)}{z},$$

and similarly with the roles of $n$ and $\xi$ exchanged, as the singular part of the operator product map, where $\varphi^n_{-1} = n \otimes z^{-1}, \varphi^{*,\xi}_0 = \xi \otimes z^{-1} \in \mathcal{W}(V)$ for $n \in N$ and $\xi \in N^\vee$.

More generally, a chiral algebra $\mathcal{W}^{\mathrm{ch}}(X, M) \in \mathrm{Alg}^{\mathrm{ch}}_{\mathrm{un}}(X)$ can be defined for any coherent $D$ module $M$, by taking $L_N = M \oplus M^\circ$ together with the canonical duality pairing $\langle \cdot, \cdot \rangle_M \in \mathrm{Hom}_{D(X)^*}(M, M^\circ; \omega_X)$, where $M^\circ = \mathbb{D}M \in D(X)$ is the dual $D$ module. The $D$ module underlying this chiral algebra is given by

$$\mathcal{W}^{\mathrm{ch}}(X, M) \cong \mathrm{Sym}^\bullet_!(M \oplus M^\circ) \cong \mathrm{Sym}^\bullet_!(M) \otimes^! \mathrm{Sym}^\bullet_!(M^\circ) \ .$$

The previous construction corresponds to the special case $M = N \otimes_{\mathbb{K}} \mathcal{D}_X$ so that $M^\circ = N^\vee \otimes_K \omega_X \otimes_{\mathcal{O}_X} \mathcal{D}_X$.

## 13. BRST REDUCTION OF CHIRAL ALGEBRAS AND VERTEX ALGEBRAS

In this section, we again restrict to the case where $X$ is a smooth algebraic curve.

*Warning* 13.0.1. Throughout this section, all of the objects will be of cohomological degree zero (in the heart of the relevant t-structure) and all the functors non-derived, in contrast with our general conventions.

Let $L \in \mathrm{Lie}^*(X)$ be a Lie* algebra on $X$ such that the underlying $D$ module of $L$ is coherent, and let $L^\circ = \mathbb{D}L \in D(X)$ denote the dual of the underlying $D$ module. For simplicity, we also assume $L$ is torsion-free as an $\mathcal{O}_X$ module, as this holds in all our examples of interest.

*Definition* 13.0.2. The chiral Clifford algebra $\mathcal{C}l^{\mathrm{ch}}(X, L) = \mathcal{W}^{\mathrm{ch}}(X, L[1]) \in \mathrm{Alg}^{\mathrm{ch}}(X)$ is defined as the (graded) chiral Weyl algebra on the $D$ module $L[1]$.

*Remark* 13.0.3. The $D$ module underlying this chiral algebra is given by

(13.0.1)
$$\mathcal{W}^{\mathrm{ch}}(X, L[1]) \cong \mathrm{Sym}^\bullet_!(L[1] \oplus L^\circ[-1]) \cong \mathrm{Sym}^\bullet_!(L[1]) \otimes^! \mathrm{Sym}^\bullet_!(L^\circ[-1]) \ .$$

*Example* 13.0.4. For $L = \tilde{L}_D$ on $X = \mathbb{A}^1$ an induced $D$ module on a translation invariant $\mathcal{O}_X$ module $\tilde{L}$ with fibre $\tilde{L}_0$ at $0 \in \mathbb{A}^1$, as in Example 12.0.5, the corresponding vertex algebra is given by

$$\mathcal{C}l^{\mathrm{ch}}(L)_0 = \mathrm{Sym}^\bullet_{\mathbb{K}}\big(\big[\tilde{L}_0 \otimes_{\mathbb{K}} z^{-1}\mathbb{K}[z^{-1}]\big][1] \oplus \big[\tilde{L}^\vee_0 \otimes_{\mathbb{K}} z^{-1}\mathbb{K}[z^{-1}]\big][-1] \cdot dz\big) \qquad \text{with} \qquad Y(\psi^a_{-1}; z)\psi^{*,\xi}_0 = \frac{\xi(a)}{z},$$

and similarly with the roles of $a$ and $\xi$ exchanged, as the singular part of the operator product map, where $\psi^a_{-1} = a \otimes z^{-1}, \psi^{*,\xi}_0 = \xi \otimes z^{-1} \in \mathcal{C}l^{\mathrm{ch}}(L)_0$ for $a \in \tilde{L}_0$ and $\xi \in \tilde{L}^\vee_0$.

Following [BD04], we denote by $\mathcal{C}l^{\mathrm{ch}}(L)^j_i$ the cohomological degree $j$ summand of the $i^{th}$ PBW filtration step of $\mathcal{C}l^{\mathrm{ch}}(L)$. Note this conflicts with the notation $\mathcal{C}l^{\mathrm{ch}}(L)_0$ used above for the corresponding vertex algebra, but we will not refer to the latter object again until Example 13.0.18 at the end of this section.

*Definition* 13.0.5. Let $M \in D_c(X)$ be a coherent $D$ module on $X$ which is torsion free as an $\mathcal{O}_X$ module. The Tate extension $\mathfrak{gl}(M)^\flat \in \mathrm{Lie}^*(X)$ is the $\omega_X$ extension

$$\omega_X \hookrightarrow \mathfrak{gl}(M)^\flat \twoheadrightarrow \mathfrak{gl}(M) \qquad \text{of} \qquad \mathfrak{gl}(M) = \mathcal{H}om_{D(X)^*}(M, M) = M \otimes^! M^\circ \in \mathrm{Lie}^*(M)$$

the endomorphism Lie* algebra $\mathfrak{gl}(M) \in \mathrm{Lie}^*(M)$, defined by

$$\mathfrak{gl}(M)^\flat := \Delta^! \mathrm{Cone}\left[ M \boxtimes M^\circ[1] \xrightarrow{\iota \oplus \langle \cdot, \cdot \rangle_M} j_* j^*(M \boxtimes M^\circ) \oplus \Delta_* \omega_X \right]$$

where $\iota : M \boxtimes M^\circ \to j_* j^*(M \boxtimes M^\circ)$ is the unit of the $(j_*, j^*)$ adjunction, and $\langle \cdot, \cdot \rangle_M : M \boxtimes M^\circ \to \Delta_* \omega_X$ is the $\otimes^*$ duality pairing.

*Remark* 13.0.6. The usual excision short exact sequence induces the sequence

$$\Delta_* \omega_X \hookrightarrow \mathrm{Cone}\left[ M \boxtimes M^\circ[1] \xrightarrow{\iota \oplus \langle \cdot, \cdot \rangle_M} j_* j^*(M \boxtimes M^\circ) \oplus \Delta_* \omega_X \right] \twoheadrightarrow \Delta_*(M \otimes^! M^\circ)$$

so that $\mathfrak{gl}(M)^\flat$ is an $\omega_X$ extension of $\mathfrak{gl}(M)$, by Kashiwara's lemma

The Tate extension $L^\flat \in \mathrm{Lie}^*(X)$ of $L \in \mathrm{Lie}^*(X)$ is the $\omega_X$ extension of $L$ pulled back from the extension of $\mathfrak{gl}(L)^\flat$ defined above under the adjoint action map of Lie* algebras $L \to \mathfrak{gl}(L)$.

*Remark* 13.0.7. By construction, the adjoint action extends to a map of Lie* algebras $L^\flat \to \mathfrak{gl}(L)^\flat$.

*Remark* 13.0.8. Agreement with usual Tate extension on topological Lie algebras

*Proposition* 13.0.9. The restriction of the Lie* bracket on $\mathcal{C}l^{\mathrm{ch}}(L)^{\mathrm{Lie}}$ to $\mathcal{C}l^{\mathrm{ch}}(L)_2^0 \in \mathrm{Lie}^*(X)$ defines an $\omega_X$ extension given by

$$\mathcal{C}l^{\mathrm{ch}}(L)_0^0 = \omega_X \hookrightarrow \mathcal{C}l^{\mathrm{ch}}(L)_2^0 \twoheadrightarrow \mathcal{C}l^{\mathrm{ch}}(L)_2^0 / \mathcal{C}l^{\mathrm{ch}}(L)_0^0 \cong L \otimes L^\circ .$$

Moreover, $\mathcal{C}l^{\mathrm{ch}}(L)_2^0 = \mathfrak{gl}(L)^\flat \in \mathrm{Lie}^*(X)$ is canonically equivalent to the Tate extension.

*Corollary* 13.0.10. There is a natural morphism of Lie* algebras $\beta : L^\flat \to \mathcal{C}l^{\mathrm{ch}}(L)^{\mathrm{Lie}}$ given by the composition $L^\flat \hookrightarrow \mathfrak{gl}(L)^\flat \cong \mathcal{C}l^{\mathrm{ch}}(L)_2^0 \hookrightarrow \mathcal{C}l^{\mathrm{ch}}(L)$.

Throughout the remainder of this section, let $A \in \mathrm{Alg}_{\mathrm{un}}^{\mathrm{ch}}(X)$ be a chiral algebra on $X$.

*Definition* 13.0.11. A BRST datum for $A$ with respect to $L$ is a map $\alpha : L^\flat \to A^{\mathrm{Lie}}$ of Lie* algebras on $X$ such that $\alpha(1^\flat) = -1_A$.

*Remark* 13.0.12. Let $\alpha : L^\flat \to A^{\mathrm{Lie}}$ be a BRST datum for $A$ with respect to $L$. Then there is a morphism of Lie* algebras

$$l^0 := \alpha + \beta : L \to A \otimes \mathcal{C}l^{\mathrm{ch}}(L)$$

where the sum descends to $L$ since $\alpha(1^\flat) + \beta(1^\flat) = 0$; we also abuse notation throughout by omitting the superscript Lie where we have applied the forgetful functor $(\cdot)^{\mathrm{Lie}} : \mathrm{Alg}^{\mathrm{ch}}(X) \to \mathrm{Lie}^*(X)$. The image $\mathrm{im}(l^0) \subset \mathcal{C}l^{\mathrm{ch}}(L)^0$ is concentrated in cohomological degree 0.

The map of $D$ modules

$$l^{-1} : L[1] \to A \otimes \mathcal{C}l^{\mathrm{ch}}(L) \qquad \text{defined by} \qquad L[1] \hookrightarrow 1_A \otimes \mathcal{C}l^{\mathrm{ch}}(L)_1^{-1} \hookrightarrow A \otimes \mathcal{C}l^{\mathrm{ch}}(L)$$

extends $l^0$ above to a map of Lie* algebras

$$l = l^{-1} \oplus l^0 : L[1] \rtimes L \to A \otimes \mathcal{C}l^{\mathrm{ch}}(L)$$

where the former is the semidirect product of $L$ with the abelian Lie* algebra $L[1]$.

*Remark* 13.0.13. Note that the tensor factor $\mathrm{Sym}_!^\bullet(L^\circ[-1])$ of the $D$ module underlying $\mathfrak{Cl}^{\mathrm{ch}}(L)$, as in Equation 13.0.1, can be identified with the underlying $D$ module of the Chevalley-Eilenberg cochains $C_{\mathrm{CE}}^\bullet(L) \in \mathrm{Comm}^!(X)$ on $L \in \mathrm{Lie}^*(X)$, and in particular is equipped with a canonical differential

$$\delta^{\mathrm{CE}} : \mathrm{Sym}_!^\bullet(L^\circ[-1]) \to \mathrm{Sym}_!^\bullet(L^\circ[-1])[1] \qquad \text{defined by} \qquad \delta^{\mathrm{CE}}|_{L^\circ[-1]} = b^\circ \in \mathrm{Hom}_{D(X)^!}(L^\circ[-1]; L^\circ[-1], L^\circ[-1])^1$$

the cohomological degree 1, arity one-to-two operation in the cooperad $D(X)^{!,\mathrm{op}}$, corresponding to the two-to-one operation $b \in \mathrm{Hom}_{D(X)^*}(L, L; L)$ underlying the Lie* bracket.

Now, we construct the so-called BRST charge as follows: Let

$$\tilde{\chi} = \mu \circ (l^0 \otimes \mathbb{1}_{L^\circ[-1]}) - \mu \circ (l^{-1} \otimes \delta^{\mathrm{CE}}|_{L^\circ[-1]}) \ \in \mathrm{Hom}_{D(X)^{\mathrm{ch}}}(L, L^\circ; A \otimes \mathfrak{Cl}^{\mathrm{ch}}(L)^1[1])$$

denote the arity two chiral operation defined by the given composition. We have the following key lemma from 3.8.9 in [BD04]:

*Lemma* 13.0.14. The following arity two operations in the operad $D(X)^*$ agree:

$$b \circ (l^0 \otimes \mathbb{1}_{L^\circ[-1]}) = b \circ (l^{-1} \otimes \delta^{\mathrm{CE}}|_{L^\circ[-1]}) \ \in \mathrm{Hom}_{D(X)^*}(L, L^\circ[-1]; A \otimes \mathfrak{Cl}^{\mathrm{ch}}(L)^1) \ .$$

The preceding lemma implies that $\tilde{\chi}$ induces an arity two operation in the operad $D(X)^!$

$$\chi \in \mathrm{Hom}_{D(X)^!}(L, L^\circ; A \otimes \mathfrak{Cl}^{\mathrm{ch}}(L)^1[1]) = \mathrm{Hom}_{D(X)}(L \otimes L^\circ; A \otimes \mathfrak{Cl}^{\mathrm{ch}}(L)^1[1]) \ .$$

This allows us to make the following key definition:

*Definition* 13.0.15. The BRST charge corresponding to the BRST datum $\alpha$ is defined by

$$\mathfrak{d}_\alpha := \chi(\mathbb{1}_L) \in \Gamma_{\mathrm{dR}}(X, A \otimes \mathfrak{Cl}^{\mathrm{ch}}(L)^1[1]) \ .$$

Further, the BRST differential corresponding to $\alpha$ is defined by

$$d_\alpha := b(\chi(\mathbb{1}_L) \otimes (\cdot)) : A \otimes \mathfrak{Cl}^{\mathrm{ch}}(L) \to A \otimes \mathfrak{Cl}^{\mathrm{ch}}(L)[1] \ .$$

*Theorem* 13.0.16. [BD04] The BRST charge satisfies $b(\mathfrak{d}_\alpha, \mathfrak{d}_\alpha) = 0$ and thus $d_\alpha^2 = 0$.

*Definition* 13.0.17. The BRST reduction of $A$ by $L$ via the BRST datum $\alpha$ is the DG chiral algebra

$$\mathrm{C}_{\mathrm{BRST}}(L; A) := (A \otimes \mathfrak{Cl}^{\mathrm{ch}}(L), d_\alpha) \ .$$

*Example* 13.0.18. Concretely, suppose that $X = \mathbb{A}^1$ and all the objects in the construction are weakly $\mathbb{G}_a$ equivariant, as in sections 8 and 10.5. Then $A$ is equivalent to a vertex algebra $\mathbb{V}$, $L$ is equivalent to a vertex Lie algebra $L_0$, and the vertex algebra corresponding to $\mathfrak{Cl}^{\mathrm{ch}}(L)$ is given by

$$\mathfrak{Cl}^{\mathrm{ch}}(L)_0 = \mathrm{Sym}_{\mathbb{K}}^\bullet\big(\big[\tilde{L}_0 \otimes_{\mathbb{K}} z^{-1}\mathbb{K}[z^{-1}]\big][1] \oplus \big[\tilde{L}_0^\vee \otimes_{\mathbb{K}} z^{-1}\mathbb{K}[z^{-1}]\big][-1] \cdot dz\big) \ ,$$

as in Example 13.0.4. Further, suppose for simplicity that $L_0$ is as in Example 10.5.7, defined by $\tilde{L}_0 = \mathfrak{g}$ a finite type Lie algebra.

Then the the BRST charge $\mathfrak{d}_\alpha$ is given concretely by

$$Q_\alpha = \sum_i e_{-1}^i \otimes \psi_0^{*,i} - \frac{1}{2} \sum_{i,j,k} 1 \otimes c_k^{ij} \psi_0^{*,i} \psi_0^{*,j} \psi_{k,-1} \ ,$$

where $i, j, k$ are indices labelling a basis for $\mathfrak{g}$ and $\psi^a_{-1} = a \otimes z^{-1}, \psi^{*,\xi}_0 = \xi \otimes z^{-1} \in \mathcal{C}l^{\mathrm{ch}}(L)_0$ for $a \in \tilde{L}_0$ and $\xi \in \tilde{L}^\vee_0$. Heuristically, this is computed by making the following identifications:

$$\mu \circ (l^0 \otimes \mathbb{1}_{L^\circ[-1]})(\mathbb{1}_L) = \sum_i e^i_{-1} \otimes \psi^{*,i}_0 \qquad + \text{ singular terms} \qquad \text{, and}$$

$$\mu \circ (l^{-1} \otimes \delta^{\mathrm{CE}}|_{L^\circ[-1]})(\mathbb{1}_L) = \frac{1}{2} \sum_{i,j,k} 1 \otimes c^{ij}_k \psi^{*,i}_0 \psi^{*,j}_0 \psi_{k,-1} \qquad + \text{ singular terms} \qquad ,$$

where the singular terms in each of the above expressions are equal so that the difference constitutes a non-singular section, in keeping with Lemma 13.0.14.

## 14. Francis-Gaitsgory chiral Koszul duality: $\mathrm{Alg}^{\mathrm{fact}}(X) \cong \mathrm{Alg}^{\mathrm{ch}}(X)$

**14.1. Overview.** In this section, we explain the correspondence between factorization algebras $\mathcal{A} \in \mathrm{Alg}^{\mathrm{fact}}(X)$ and chiral algebras $A \in \mathrm{Alg}^{\mathrm{ch}}(X)$, following [BD04] and [FG11] throughout. The main idea is that the $D$ module $A$ underlying the chiral algebra is defined by

$$A[1] = A_1 = \Delta^{\mathrm{main},!}\mathcal{A} \in D(X)$$

the restriction of $\mathcal{A} \in D(\mathrm{Ran}_X)$ to the first stratum of the Ran space of $X$, and the chiral product $\mu : j_* j^* A^{\boxtimes 2} \to \Delta_* A$ is equivalent to the data required to extend $A_1$ to a factorizable sheaf $\mathcal{A}$ on $\mathrm{Ran}_X$. We begin with an outline of the equivalence in geometric terms, before describing the more structured algebraic perspective that facilitates the proof.

The data of the $D$ module on $\mathrm{Ran}_X$ underlying a factorization algebra $\mathcal{A} \in \mathrm{Alg}^{\mathrm{fact}}(X)$ is almost completely specified by $A_1 = \Delta^{\mathrm{main},!}\mathcal{A}$, since the factorization data gives an identification of the restriction of $A_2 = \Delta^{\{1,2\},!}\mathcal{A}$ to the complement of the diagonal with that of $A^{\boxtimes 2}_1$:

$$(14.1.1) \qquad\qquad\qquad j^*(A_2) \cong j^*(A^{\boxtimes 2}_1) ,$$

and similarly on higher cardinality products. From the excision sequence for $A_2$,

$$\Delta_* \Delta^! A_2 \to A_2 \to j_* j^* A_2 ,$$

we see that the additional information required to reconstruct $A_2$ from $A_1$ is equivalent to the boundary map

$$j_* j^*(A^{\boxtimes 2})[2] \cong j_* j^* A_2 \to \Delta_* \Delta^! A_2[1] \cong \Delta_* A[2] ,$$

which is precisely the desired chiral product map

$$\mu : j_* j^*(A^{\boxtimes 2}) \to \Delta_* A$$

after identifying the restrictions of $A_2$ to the diagonal and its complement with their descriptions in terms of $A$ as shown. Indeed, we can recover $A_2$ as

$$(14.1.2) \qquad\qquad\qquad A_2 = \ker\left[j_* j^*(A^{\boxtimes 2})[2] \to \Delta_* A[1][1]\right] .$$

Moreover, as we explain below, this sheaf extends coherently to $A_3 \in D(X^3)$ satisfying the required gluing and factorizability conditions if and only if $\mu$ satisfies the Jacobi identity.

14.2. **Cocommutative-Lie Koszul duality for** $D(\mathbf{Ran}_X)^{\mathrm{ch}}$. We now explain the more structured algebraic perspective on this equivalence. Chiral algebras and factorization algebras are formally defined as certain Lie and cocommutative coalgebra objects internal to the category $D(\mathrm{Ran}_X)^{\mathrm{ch}}$ of $D$ modules on the Ran space of $X$ with respect to the chiral tensor structure $\otimes^{\mathrm{ch}}$, respectively. For a general well-behaved symmetric monoidal category $\mathcal{C}^{\otimes}$, there is a canonical equivalence between Lie algebra and cocommutative coalgebra objects in it, given by the Chevalley-Eilenberg chains functor

$$\mathrm{C}^{\bullet}(\cdot) : \mathrm{Alg}_{\mathrm{Lie}}(\mathcal{C}^{\otimes}) \to \mathrm{Alg}_{\mathrm{CoComm}}(\mathcal{C}^{\otimes}) \qquad L \mapsto \mathrm{C}^{\bullet}(L) = (\mathrm{Sym}^{\bullet}(L[1]), d_{\mathrm{CE}})$$

where the differential $d_{\mathrm{CE}}$ is generated by the map $L[1] \otimes L[1] \to L[1][1]$ given by the Lie bracket.

In [FG11] it is shown that the category $D(\mathrm{Ran}_X)^{\mathrm{ch}}$ satisfies the hypotheses required to construct such an equivalence, and moreover that the resulting functor induces an equivalence between the full subcategories of chiral algebras and factorization algebras:

*Theorem* 14.2.1. [BD04, FG11] The Chevalley-Eilenberg chains functor on $D(\mathrm{Ran}_X)^{\mathrm{ch}}$ induces equivalences

$$
\begin{array}{ccc}
\mathrm{Alg}^{\mathrm{ch}}(X) & \xrightarrow{\;\cong\;} & \mathrm{Alg}^{\mathrm{fact}}(X) \\
\downarrow & & \downarrow \\
\mathrm{Alg}_{\mathrm{Lie}}(D(\mathrm{Ran}_X)^{\mathrm{ch}}) & \xrightarrow{\;\cong\;} & \mathrm{Alg}_{\mathrm{CoComm}}(D(\mathrm{Ran}_X)^{\mathrm{ch}})
\end{array}
$$

such that the preceding diagram commutes.

We now outline the proof of the Theorem, emphasizing that the resulting equivalence reproduces the geometric arguements given in the overview above. Let $\mathcal{L} \in \mathrm{Alg}_{\mathrm{Lie}}(D(\mathrm{Ran}_X)^{\mathrm{ch}})$ be a Lie algebra object with respect to the chiral tensor product. The free graded cocommutative coalgebra object in $D(\mathrm{Ran}_X)$ generated by $\mathcal{L}$ is given by the $S_n$ coinvariants of

$$\tilde{\mathrm{C}}_{\bullet}(\mathcal{L}) = \bigoplus_{n \in \mathbb{N}} \mathcal{L}[1]^{\otimes n} \ .$$

For notational simplicity, we describe the construction omitting the $S_n$ coinvariants throughout. In the case at hand, recall from Proposition 5.3.3 that

$$(\otimes^{\mathrm{ch}}_{j \in J} \mathcal{L})_I = \bigoplus_{\pi : I \twoheadrightarrow J} j(\pi)_* j(\pi)^* (\boxtimes_{j \in J} L_{I_j}) \qquad\qquad \text{so that}$$

$$\mathrm{C}_{\bullet}(\mathcal{L})_I = \bigoplus_{n \in \mathbb{N}} (\mathcal{L}[1]^{\otimes n})_I = \bigoplus_{n=0}^{|I|} (\mathcal{L}^{\otimes n})_I [n] = \bigoplus_{[\pi : I \twoheadrightarrow S] \in \mathrm{fSet}_{I/}} j(\pi)_* j(\pi)^* (\boxtimes_{s \in S} L_{I_s})[|S|] \qquad .$$

The Chevalley-Eilenberg differential on $\mathrm{C}_{\bullet}(\mathcal{L})$ is defined over each $I$ as follows: Fix $\pi \in \mathrm{Hom}_{\mathrm{fSet}_{I/}}(T, S)$, that is, $\pi_T : I \twoheadrightarrow T, \pi_S : I \twoheadrightarrow S$ and $\pi : T \twoheadrightarrow S$ and such that $\pi_S = \pi \circ \pi_T$. We are interested in the case $|T| = |S| + 1$, so for concreteness say $\pi(t_0) = \pi(t_1) = s_1$ and $\pi(t_i) = s_i$ for $i = 2, ..., |S|$. Then we have a map

$$
\begin{aligned}
j(\pi_T)_* j(\pi_T)^* (\boxtimes_{t \in T} L_{I_t}) &= j(\pi_T)_* j(\pi_T)^* \big( (L_{I_{t_0}} \boxtimes L_{I_{t_1}}) \boxtimes (\boxtimes_{s \neq s_1} L_{I_s}) \big) \\
&\cong j(\pi_S)_* j(\pi_S)^* \big( j(\tilde{\pi}) \big)_* j(\tilde{\pi})^* (L_{I_{t_0}} \boxtimes L_{I_{t_1}}) \boxtimes (\boxtimes_{s \neq s_1} L_{I_s}) \big) \\
&\to j(\pi_S)_* j(\pi_S)^* (\boxtimes_{s \in S} L_{I_s})
\end{aligned}
$$

where $\tilde{\pi} : I_{t_0} \sqcup I_{t_1} \to \{1, 2\}$ and the map is given by $j(\pi_S)_* j(\pi_S)^* (b(\tilde{\pi}) \boxtimes \mathbb{1})$, where

$$b(\tilde{\pi}) : j(\tilde{\pi})_* j(\tilde{\pi})^* (L_{I_{t_0}} \boxtimes L_{I_{t_1}}) \to L_{I_{s_1}}$$

is the chiral Lie algebra structure map. For fixed $\pi : I \to T$, summing over all such $\pi : T \twoheadrightarrow S$ defines the component of the Chevalley-Eilenberg differential on the summand $j(\pi_T)_* j(\pi_T)^*(\boxtimes_{t \in T} L_{I_t})$ of $\mathrm{C}_\bullet(\mathcal{L})_I$.

Now, suppose $\mathcal{L} \in \mathrm{Alg}^{\mathrm{ch}}(X)$, so that there exists $A \in D(X)$ such that $\mathcal{L} = \Delta_*^{\mathrm{main}} A \in D(\mathrm{Ran}_X)$, or equivalently

$$L_I = \Delta_*^{(I)} A \ ,$$

for each $I \in \mathrm{fSet}_{\mathrm{surj}}$. Then we have

$$j(\pi)_* j(\pi)^*(\boxtimes_{s \in S} L_{I_s}) = j(\pi)_* j(\pi)^*(\boxtimes_{s \in S} \Delta_*^{(I_s)} A) \cong j(\pi)_* j(\pi)^* \Delta(\pi)_* A^{\boxtimes |S|} \cong \Delta(\pi)_* j_*^{(S)} j^{(S)*} A^{\boxtimes |S|}$$

so that

$$\tag{14.2.1} \mathrm{C}_\bullet(\mathcal{L})_I = \bigoplus_{[\pi : I \,\twoheadrightarrow\, S] \in \mathrm{fSet}_{I/}} \Delta(\pi)_* j_*^{(S)} j^{(S)*} A^{\boxtimes |S|}[|S|] \ .$$

On the first stratum $X \hookrightarrow \mathrm{Ran}_X$ of the Ran space, we have the desired equality

$$A_1 := \mathrm{C}(\mathcal{L})_1 = A[1] \ ,$$

since the tensor powers of arity greater than one with respect to the chiral tensor structure vanish when restricted to the main diagonal, by Remark 5.3.4. Similarly, over $X^2$ there are two non-vanishing terms in the expression of Equation 14.2.1 for the homological Chevalley-Eilenberg complex, given by

$$A_2 := \mathrm{C}(\mathcal{L})_2 = \left[ j_* j^*(A^{\boxtimes 2})[2] \to \Delta_* A[1] \right] \ ,$$

in agreement with Equation 14.1.2 from our geometric explanation. Finally, over $X^3$ the Chevalley-Eilenberg complex is given by

$$A_3 := \mathrm{C}(\mathcal{L})_3 = \left[ j_*^{\{1,2,3\}} j^{\{1,2,3\},*}(A^{\boxtimes 3})[3] \to \bigoplus_{\substack{i \neq j \\ i,j=1,2,3}} \Delta_*^{x_i = x_j} j_* j^*(A^{\boxtimes 2})[2] \to \Delta_* A[1] \right] \ .$$

This construction manifestly defines $A_3 \in D(X^3)$ compatibly extending $A_2 \in D(X^2)$ as defined above, and the requirement that the differential squares to zero so that it actually gives a well defined complex of $D$ modules is equivalent to the Jacobi identity for the chiral product map $\mu$.

# Chapter 2

# Equivariant factorization algebras and the localization theorem

In this chapter, we develop an analogous theory of equivariant factorization algebras, as outlined in Section 1.3.2 of the introduction.

## 15. A REVIEW OF EQUIVARIANT $D$ MODULES

We begin with a brief overview of the theory of equivariant $D$ modules, paralleling the theory of equivariant constructible sheaves recalled in Appendix B.1. Let $G$ be an algebraic group, $X$ a finite type scheme over $\mathbb{K}$, and fix an action of $G$ on $X$. Let $m : G \times G \to G$ denote the multiplication map, $a : G \times X \to X$ the action map, $\mathrm{p}_X : G \times X \to X$ the projection to $X$ and $\mathrm{p}_G : X \times G \to G$ the projection to $G$.

*Definition* 15.0.1. A $G$ equivariant structure on $\mathcal{F} \in \mathrm{QCoh}(X)$ is an isomorphism $\alpha : a^\bullet \mathcal{F} \xrightarrow{\cong} p_X^\bullet \mathcal{F}$ in $\mathrm{QCoh}(G \times X)$, together with commutativity of the diagram in $\mathrm{QCoh}(G \times G \times X)$ defined by

$$(15.0.1) \qquad \begin{array}{ccc} (a \circ (\mathbb{1}_G \times a))^\bullet M & \xrightarrow{\cong} & (a \circ (m \times \mathbb{1}_X))^\bullet M \\ {\scriptstyle \cong} \downarrow {\scriptstyle (\mathbb{1}_G \times a)^\bullet \alpha} & & {\scriptstyle \cong} \downarrow {\scriptstyle (m \times \mathbb{1}_X)^\bullet \alpha} \\ (\mathbb{1}_G \times a)^\bullet \mathrm{p}_X^\bullet M & \xrightarrow[\cong]{\mathbb{1}_{\mathrm{QCoh}(G)} \boxtimes \alpha} & \tilde{\mathrm{p}}_X^\bullet M \end{array} \quad ,$$

where $\tilde{\mathrm{p}}_X : G \times G \times X \to X$ is the projection to $X$.

*Definition* 15.0.2. A weak $G$ equivariant structure on $M \in D(X)$ is an isomorphism $\alpha : a^! M \xrightarrow{\cong} \mathrm{p}_X^! M$ of complexes of $O_G \boxtimes D_X$ modules, together with commutativity data for the analogue of the diagram 15.0.1 in the category of complexes of $\mathcal{O}_{G \times G} \boxtimes D_X$ modules.

Let $D(X)^{G,w}$ denote the category of weakly $G$ equivariant $D$ modules on $X$.

*Remark* 15.0.3. Heuristically, a weak equivariant structure on a $D$ module is an equivariant structure on the underlying $\mathcal{O}_X$ module such that for each $g \in G$, the induced equivalence $\alpha_g : a_g^! M \xrightarrow{\cong} M$ lifts to an isomorphism in $D(X)$.

*Definition* 15.0.4. A (strong) $G$ equivariant structure on $M \in D(X)$ is an isomorphism $\alpha : a^! M \xrightarrow{\cong} \mathrm{p}_X^! M$ in $D(G \times X)$, together with commutativity data for the analogue of the diagram 15.0.1 in the category $D(G \times G \times X)$.

Let $D(X)^G$ denote the category of (strongly) $G$ equivariant $D$ modules on $X$.

*Remark* 15.0.5. Heuristically, a strong equivariant structure on a $D$ module is a weak equivariant structure such that the induced equivalences $\alpha_g : a_g^! M \xrightarrow{\cong} M$ in $D(X)$ are locally constant along $G$. For $G$ connected, the equivariant structure appears to be uniquely determined by the local constancy condition, since it is fixed by the requirement $\alpha_e = \mathbb{1}_M$, so that strong equivariance is a property of the underlying $D$ module, rather than a structure. This statement is true for a strict $D$ module, but not for a complex, as we explain more carefully below.

*Example* 15.0.6. There is evidently a forgetful functor $D(X)^G \to D(X)^{G,w}$ coming from the forgetful functor $D(G) \to \mathrm{QCoh}(G)$. The additional data required to lift a weak equivariant structure to a strong equivariant structure is as follows:

Let $M \in D(X)^{G,w}$ and define the Lie derivative of $M$ with respect to the $G$ equivariant structure by

$$\mathcal{L}_{(\cdot)} : \mathfrak{g} \to \mathrm{End}_{\mathrm{Sh}(X;\mathbb{K})}(M) \qquad \text{by} \qquad \mathcal{L}_X(m) = \partial_t \alpha_{\exp(tX)}(m)|_{t=0}$$

for each $m \in M$.

Let $da : \mathfrak{g} \to \Gamma(X, T_X)$ denote the infinitesimal action map and $\nabla^M : \Gamma(X, T_X) \to \mathrm{End}_{\mathrm{Sh}(X)}(M)$ the connection underlying the $D$ module structure which compose to define a map $\nabla_{da(\cdot)} : \mathfrak{g} \to \mathrm{End}_{\mathrm{Sh}(X)}(M)$.

For $M \in D(X)^{\heartsuit,(G,w)}$ a weakly $G$ equivariant $D$ module concentrated in a single cohomological degree, the weak equivariant structure defines a strong equivariant structure if $\mathcal{L}_{(\cdot)} = \nabla_{da(\cdot)}$. The map $\nabla_{da(\cdot)}$ is defined independent of the $G$ equivariant structure, and the $G$ equivariant structure necessarily integrates the map $\mathcal{L}_{(\cdot)}$, so that a strong equivariant structure on $M \in D(X)^{\heartsuit}$ is unique for $G$ connected. The condition of its existence is the integrability of the representation $\mathcal{L}_{(\cdot)} : \mathfrak{g} \to \mathrm{End}_{\mathrm{Sh}(X)}(M)$.

Now, for a general object $M \in D(X)^{G,w}$, let $(M^\bullet, d_M)$ denote the underlying complex of $D$ modules. Then a lift to a strong $G$ equivariant structure on $M$ is equivalent to a homotopy trivializing the difference of these endomorphisms, that is, a map $h : \mathfrak{g} \to \mathrm{End}_{D(X)}^{-1}(M^\bullet)$ such that $\mathcal{L}_{(\cdot)} - \nabla_{da(\cdot)} = [h, d_M]$. Note that the difference $\mathcal{L}_{(\cdot)} - \nabla_{da(\cdot)} \in \mathrm{End}_{D(X)}(M)$.

*Example* 15.0.7. The dualizing sheaf $\omega_X \in D^r(X)^G$ and the constant sheaf $\omega_X[-2d_X] \in D^r(X)^G$ admit canonical strong $G$ equivariant structures for any action of $G$ on $X$, given by the identifications $a^! \omega_X = \omega_{G \times X} = \mathrm{p}_X^! \omega_X$.

Under the quasiisomorphism $\omega_X[-2d_X] \cong \Omega_{X,\mathcal{D}}^\bullet$ of Proposition A.4.2, the induced strong equivariant structure on $\Omega_{X,\mathcal{D}}^\bullet \in D^r(X)$ is given by

$$h = \iota_{da(\cdot)} : \mathfrak{g} \to \mathrm{End}_{D^r(X)}^{-1}(\Omega_{X,\mathcal{D}}^\bullet) \cong \mathrm{Diff}(\Omega_X^\bullet, \Omega_X^\bullet[-1]) \qquad \text{where} \qquad \begin{cases} da : \mathfrak{g} \to \Gamma(X, T_X) \\ \iota_{(\cdot)} : \Gamma(X, T_X) \to \mathrm{End}^{-1}(\Omega_{X,\mathcal{D}}^\bullet) \end{cases}$$

are the infinitesimal action map and the interior product operation. The compatibility follows from the Cartan formula, as the endomorphism $\mathcal{L}_X - \Delta_{da(\cdot)}$ is given by the usual Lie derivative.

*Example* 15.0.8. For $X = \mathbb{A}^1$ we have $\omega_X = \mathbb{K}[x]$ and

$$\Omega_{\mathbb{A}^1,\mathcal{D}}^\bullet = \mathbb{K}[x, \partial_x] \xrightarrow{m_{\partial_x}} \mathbb{K}[x, \partial_x][-1] \ ,$$

where $m_{\partial_x}$ denotes the left multiplication map. For the action of $G = \mathbb{G}_a$, the homotopy is given by $h = \mathbb{1} : \mathbb{K}[x, \partial_x][-1] \to \mathbb{K}[x, \partial_x]$ the identity map. For the action of $G = \mathbb{G}_m$, the homotopy is given by $h = m_x : \mathbb{K}[x, \partial_x][-1] \to \mathbb{K}[x, \partial_x]$.

*Example* 15.0.9. The category $D(\mathrm{pt})^G$ is the category of complexes of $G$ representations $(V, d) \in \mathrm{Rep}(G)_{\mathbb{K}}$ together with a homotopy $h : \mathfrak{g} \to End^{-1}(V)$ trivializing the infinitesimal action $d\rho : \mathfrak{g} \to \mathrm{End}_{\mathrm{Vect}_{\mathbb{K}}}(V)$, that is, such that $[d, h] = d\rho$. There is a natural functor

$$D(\mathrm{pt})^G \to \mathrm{D}(H_G^\bullet(\mathrm{pt})) \qquad \text{defined by} \qquad (V, d, h) \mapsto (V \otimes_{\mathbb{K}} \mathrm{Sym}^\bullet(\mathfrak{g}^\vee[-2]))^G \ , \ d_u = d \otimes \mathbb{1}_{\mathfrak{g}^\vee} + h_{\xi^i} \otimes m_{u_i}$$

where $u_i \in g^\vee[-2]$ are some choice of linear generators of cohomological degree 2, and $\xi^i \in \mathfrak{g}$ are the corresponding dual basis.

*Proposition* 15.0.10. The above functor $D(\mathrm{pt})^G \to D(H_G^\bullet(\mathrm{pt}))$ induces an equivalence $D_{\mathrm{c}}(\mathrm{pt})^G \cong D_{\mathrm{fg}}^b(H_G^\bullet(\mathrm{pt}))$, in keeping with Theorems B.1.10 and B.2.4.

*Theorem* 15.0.11. Let $X$ be a smooth, finite type variety over $\mathbb{C}$. For $H$ a subgroup of $G$, there are restriction and induction adjunctions:

$$\mathrm{Res}_H^G : D_{\mathrm{rh}}(X)^G \rightleftarrows D_{\mathrm{rh}}(X)^H : \mathrm{Ind}_{H,*}^G \qquad \mathrm{Ind}_{H,!}^G : D_{\mathrm{rh}}(X)^H \rightleftarrows D_{\mathrm{rh}}(X)^G : \mathrm{Res}_H^G$$

Moreover, there are natural functors as in A.3.1 and A.5.3 between the corresponding $G$ equivariant categories $D_{G,\mathrm{rh}}$, satisfying the same adjunctions adjunctions and relations, defined for $G$ equivariant maps $f : X \to Y$. These functors all commute with $\mathrm{Res}_H^G$, while $f_*$ and $f^!$ commute with $\mathrm{Ind}_{h,*}^G$, and $f_!$ and $f^*$ commute with $\mathrm{Ind}_{H,!}^G$.

*Definition* 15.0.12. The $G$ equivariant cochains functor is defined by $C_G^\bullet = \pi_* : D(X)^G \to D(\mathrm{pt})^G$. The $G$ equivariant chains functor is defined by $C_\bullet^G = \pi_! : D(X)^G \to D(\mathrm{pt})^G$. The $G$ equivariant (Borel-Moore) de Rham (co)chains and cohomology are defined as in B.1.4 and B.1.7.

*Example* 15.0.13. Computing $C_G^\bullet(X)$ in terms of the de Rham model as in Example A.4.3 and Example 15.0.7, we find its image under the equivalence of Proposition 15.0.10 is

$$C_G^\bullet(X) = (\Omega_X^\bullet \otimes_{\mathbb{K}} \mathrm{Sym}^\bullet(\mathfrak{g}^\vee[-2])^G , \ d = d_{\mathrm{dR}} \otimes \mathbb{1} + \iota_{\xi_i} \otimes m_{u_i} .$$

This is the usual Cartan model for equivariant cohomology of $X$. More generally, $C_G^\bullet(X; A)$ is computed by the Cartan model with coefficients in the equivariant complex $A \in D(X)^G$.

In this case, the homotopy $h$ corresponds to the $C_\bullet(G; \mathbb{K})$ module structure on $C_{\mathrm{dR}}^\bullet(X; \mathbb{K})$, and the complex $C_G^\bullet(X; \mathbb{K})$ above is equivalent to the image of $C_{\mathrm{dR}}^\bullet(X; \mathbb{K})$ under the functor of Remark B.2.3.

*Example* 15.0.14. Suppose $G$ acts on $X$ trivially. Then $D(X)^G \cong D(X) \otimes D(\mathrm{pt})^G$ and thus by Example 15.0.9 above, there is a natural functor $D(X)^G \to D(X) \otimes D(H_G^\bullet(\mathrm{pt}))$ from equivariant $D$ modules to families of $D$ modules on $X$ over $\mathrm{Spec}\ H_G^\bullet(\mathrm{pt})$.

## 16. The category $D(\mathrm{Ran}_X)^G$

In this section, we define the category $D(\mathrm{Ran}_X)^G$ of $G$ equivariant $D$ modules on $\mathrm{Ran}_X$, and breifly outline the analogues of various structures on $D(\mathrm{Ran}_X)$ in the equivariant setting.

For each $I \in \mathrm{fSet}$ there is a diagonal action of $G$ on $X^I$, and for each $\pi : I \twoheadrightarrow J$ the corresponding diagonal embedding $\Delta(\pi) : X^J \hookrightarrow X^I$ is a morphism of $G$ varieties.

*Remark* 16.0.1. The equivariance of the diagonal embeddings under $G$ implies that the diagram defining $\mathrm{Ran}_X$ can be understood in the category of $G$ schemes, and thus the colimit $\mathrm{Ran}_X$ has a natural action of $G$. Heuristically, this is simply the action of $G$ on the space of finite subsets of $X$ by $g \cdot \{x_i\}_{i \in I} = \{g \cdot x_i\}$, which is evidently modelled by the diagonal action as above.

In analogy with Definition 4.2.1 and the discussion of that section, we make the following definition:

*Definition* 16.0.2. An object $\mathcal{A} \in D(\mathrm{Ran}_X)^G$ is an assignment

$$I \mapsto A_I \in D(X^I)^G \qquad\qquad [\pi : I \twoheadrightarrow J] \mapsto [\Delta(\pi)^! A_I \xrightarrow{\cong} A_J]$$

defined for each finite set $I \in \mathrm{fSet}$ and surjection $\pi : I \twoheadrightarrow J$, where the isomorphism $\Delta(\pi)^! A_I \xrightarrow{\cong} A_J$ is required to be in $D(X^J)^G$.

A morphism $f : \mathcal{A} \to \mathcal{B}$ between $\mathcal{A}, \mathcal{B} \in D(\mathrm{Ran}_X)^G$ is given by an assignment

$$I \mapsto [f_I : A_I \to B_I] \qquad [\pi : I \twoheadrightarrow J] \mapsto \qquad
\begin{array}{ccc}
\Delta(\pi)^! A_I & \longrightarrow & A_J \\
{\scriptstyle \Delta(\pi)^!(f_I)} \downarrow \ \ {\scriptstyle \cong} \!\!\!\!\!\!\! \diagup & & \downarrow {\scriptstyle f_J} \\
\Delta(\pi)^! B_I & \longrightarrow & B_J
\end{array}$$

defined for each finite set $I \in \mathrm{fSet}$ and surjection $\pi : I \twoheadrightarrow J$, where all required morphisms are in $D(X^I)^G$.

An object $\mathcal{A} \in D(\mathrm{Ran}_{X,\mathrm{un}})^G$ and morphism of such is defined similarly, analogously following Definition 4.3.1.

An object $\mathcal{A} \in D(\mathrm{Ran}_X)^G$ is called coherent, holonomic, ... if $A_I \in D(X^I)^G$ is so for each $I \in \mathrm{fSet}$.

*Remark* 16.0.3. Following Remark 4.2.4, the category $D(\mathrm{Ran}_X)^G$ can be equivalently defined as the category $G$ equivariant $D$ modules on the pseudo indscheme $\mathrm{Ran}_X$.

*Remark* 16.0.4. The definition 16.0.2 is stated exactly as in Definition 4.2.1 by replacing all the objects and morphisms of $D$ with their $G$ equivariant analogues. This is possible because all of the underlying geometric maps involved are $G$ equivariant so that there are natural lifts of the resulting functors to the $G$ equivariant category. In what follows, we list the various structures induced on $D(\mathrm{Ran}_X)^G$ following this pattern:

*Remark* 16.0.5. The category $D(\mathrm{Ran}_X)^{G,w}$ is defined as in definitions 4.2.1 and 16.0.2, with $A_I \in D(X^I)^{G,w}$ weakly $G$ equivariant for each $I \in \mathrm{fSet}$ and all required morphisms in the relevant weakly equivariant categories.

*Remark* 16.0.6. As in Remark 4.2.5, there are canonical functors $\Delta_*^I : D(X^I)^G \xrightleftharpoons{\hspace{1cm}} D(\mathrm{Ran}_X)^G : \Delta^{I,!}$ for each $I \in \mathrm{fSet}$. For $I = \{\mathrm{pt}\}$ these functors induce an equivalence $D(X)^G \cong D(\mathrm{Ran}_X)_X^G$.

*Example* 16.0.7. The object $\omega_{\mathrm{Ran}_X} \in D(\mathrm{Ran}_X)$ naturally lifts to $\omega_{\mathrm{Ran}_X} \in D(\mathrm{Ran}_X)^G$, as there is a canonical equivariant structure $\omega_{X^I} \in D(X^I)^G$ for each $I \in \mathrm{fSet}$, as in Example 15.0.7.

*Definition* 16.0.8. The monoidal structures $\otimes^!, \otimes^*, \otimes^{\mathrm{ch}} : \times_{j \in J} D(\mathrm{Ran}_X)^G \to D(\mathrm{Ran}_X)^G$ are presented by

$$\times_{j \in J} D(X^{I_j})^G \to D(X^I)^G \qquad\qquad (M_{I_j}) \mapsto \otimes^!(\Delta_*^{I_j} M_{I_j})_I$$
$$\times_{j \in J} D(X^{I_j})^G \to D(X^I)^G \qquad\qquad (M_{I_j}) \mapsto \boxtimes_{j \in J} M_{I_j}$$
$$\times_{j \in J} D(X^{I_j})^G \to D(X^I)^G \qquad\qquad (M_{I_j}) \mapsto j(\pi)_* j(\pi)^!(\boxtimes_{j \in J} M_{I_j}) \ ,$$

defined for each $\pi : I \twoheadrightarrow J$.

*Remark* 16.0.9. The the monoidal structures of Definition 16.0.8 are the natural lifts of the definitions 5.1.1 ,5.2.1, and 5.3.1 to the G equivariant category, in keeping with Remark 16.0.4 above.

*Definition* 16.0.10. The $\otimes^!$, $\otimes^*$, and $\otimes^{\mathrm{ch}}$ operad structures on $D(X)$ are defined by $D(X)^G \hookrightarrow D(\mathrm{Ran}_X)^G$, where the latter is equipped with the corresponding monoidal structure, following Example C.1.11.

*Corollary* 16.0.11. The functors $\Delta_*^{\mathrm{main}} : D(X)^G \to D(\mathrm{Ran}_X)^G$ and $\Delta^{\mathrm{main},!} : D(\mathrm{Ran}_X)^G \to D(X)^G$ of 4.2.5 define an equivalence of operads between $D(\mathrm{Ran}_X)_X^G$ and $D(X)^G$ under $\otimes^!$, $\otimes^*$ or $\otimes^{\mathrm{ch}}$.

## 17. Equivariant factorization algebras and equivariant chiral algebras

In this section, we define equivariant factorization algebras and equivariant chiral algebras, closely following the usual definition of factorization algebras in [BD04, FG11].

*Definition* 17.0.1. A non-unital factorization algebra on $X$ is a non-unital cocommutative coalgebra object $\mathcal{A} \in D(\mathrm{Ran}_X)^{G,\mathrm{ch}}$ such that the induced maps

$$j(\pi)^* A_I \xrightarrow{\cong} j(\pi)^*(\boxtimes_{j \in J} A_{I_j})$$

are equivalences in $D(U(\pi))^G$ for each $I, J$ and $\pi : I \twoheadrightarrow J$.

A unital factorization algebra on $X$ is an object $\mathcal{A} \in D(\mathrm{Ran}_{X,\mathrm{un}})^G$ with a non-unital factorization algebra structure on its image in $D(\mathrm{Ran}_X)^G$, and compatibility data with the unital structure on $\mathcal{A}$, as in Definition 4.2.1.

Let $\mathrm{Alg}^{\mathrm{fact}}(X)^G$ denote the category of non unital $G$ equivariant factorization algebras, defined as the full subcategory of $\mathrm{CoComm}^{\mathrm{nu}}(D(\mathrm{Ran}_X)^{G,\mathrm{ch}})$. Similarly, let $\mathrm{Alg}^{\mathrm{fact}}_{\mathrm{un}}(X)^G$ denote the category of unital $G$ equivariant factorization algebras.

*Example* 17.0.2. The dualizing sheaf $\omega_{\mathrm{Ran}_X} \in D(\mathrm{Ran}_X)^G$ of Example 16.0.7 defines a $G$ equivariant factorization algebra, with structure maps given by the natural $G$ equivariant lifts of those of Example 6.0.4.

*Definition* 17.0.3. A (non-unital) chiral algebra on $X$ is a Lie algebra object in $\mathcal{L} \in \mathrm{Lie}(D(\mathrm{Ran}_X)^{G,\mathrm{ch}})$ such that underlying object $\mathcal{L} \in D(\mathrm{Ran}_X)^G_X$ is supported on the image of the main diagonal $\Delta^{\mathrm{main}} : X \to \mathrm{Ran}_X$.

A unital chiral algebra on $X$ is an object $\mathcal{L} \in D(\mathrm{Ran}_{X,\mathrm{un}})^G$ with a non unital chiral algebra structure on its image in $D(\mathrm{Ran}_X)$, and compatibility data with the unital structure on $\mathcal{L}$, as in Definition 4.3.1.

Let $\mathrm{Alg}^{\mathrm{ch}}(X)^G$ denote the category of non unital $G$ equivariant chiral algebras, defined as the full subcategory of $\mathrm{Lie}(D(\mathrm{Ran}_X)^{\mathrm{ch},G})$. Similarly, let $\mathrm{Alg}^{\mathrm{ch}}_{\mathrm{un}}(X)^G$ denote the category of unital $G$ equivariant chiral algebras.

*Corollary* 17.0.4. A non unital, $G$-equivariant chiral algebra is equivalent to a Lie algebra object $L \in \mathrm{Lie}(D(X)^{G,\mathrm{ch}})$ internal to the operad $D(X)^{G,\mathrm{ch}}$, by Corollary 16.0.11.

*Example* 17.0.5. In particular, the equivariant chiral algebra structure maps are given by

$$b_I \in \mathrm{Hom}_{D(X)^{G,\mathrm{ch}}}(\{L\}_{i \in I}, L) = \mathrm{Hom}_{D(X^I)^G}(j^{(I)}_* j^{(I),*}(\boxtimes_{i \in I} L), \Delta^{(I)}_* L) \ .$$

There is also a weakly equivariant analogue of these definitions, which is used in Section 8 to relate chiral algebras to vertex algebras.

*Definition* 17.0.6. The category of weakly $G$ equivariant chiral algebras $\mathrm{Alg}^{\mathrm{ch}}(X)^{G,w}$ is the full subcategory of $\mathrm{Lie}(D(\mathrm{Ran}_X)^{\mathrm{ch},(G,w)})$ on algebras with underlying object $A \in D(\mathrm{Ran}_X)^{G,w}_X$ supported on the main diagonal $X \hookrightarrow \mathrm{Ran}_X$; see also Remark 16.0.5.

The main result of [FG11], recalled in Section 14, generalizes to the $G$ equivariant setting, by the same arguement used in the proof of *loc. cit.* lifted to the $G$ equivariant category:

*Corollary* 17.0.7. There is a canonical equivalence of categories $\mathrm{Alg}^{\mathrm{fact}}(X)^G \cong \mathrm{Alg}^{\mathrm{ch}}(X)^G$ between $G$ equivariant factorization algebras and chiral algebras.

*Example* 17.0.8. Suppose $G$ acts trivially on $X$. Then applying Example 15.0.14 to $X^I$ for each $I \in$ fSet, we obtain a natural functor $\mathrm{Alg}^{\mathrm{fact}}(X)^G \to \mathrm{Alg}^{\mathrm{fact}}(X)_{/H^\bullet_G(\mathrm{pt})}$ from $G$ equivariant factorization algebras on $X$ to families of factorization algebras over $H^\bullet_G(\mathrm{pt})$.

## 18. Equivariant factorization homology and the localization theorem

In this section, we define equivariant factorization homology, as well as the pullback of factorization algebras, and prove an analogue of the classical equivariant localization theorem in the setting of factorization homology.

### 18.1. Factorization homology.

Factorization homology is one of the primary invariants of factorization algebras, analogous to sheaf cohomology of sheaves, which generalizes (the dual space to) the spaces of conformal blocks of vertex algebras; it was originally introduced in [BD04]. We now give the definition of factorization homology of equivariant factorization algebras, following the standard definition given in [FG11]. In summary, the functor of factorization homology is given by the composition

$$\mathrm{Alg}^{\mathrm{fact}}(X)^G \xrightarrow{\mathrm{oblv}} D(\mathrm{Ran}_X)^G \xrightarrow{\mathrm{p}_{\mathrm{Ran}_X*}} D(\mathrm{pt})^G$$

where oblv denotes the forgetful functor to $G$ equivariant $D$ modules on $\mathrm{Ran}_X$, $\mathrm{p}_{\mathrm{Ran}_X} : \mathrm{Ran}_X \to \mathrm{pt}$ is the unique such map and $\mathrm{p}_{\mathrm{Ran}_X*}$ denotes the induced pushforward functor on equivariant $D$ modules defined below. In particular, taking $G = \{e\}$ to be the trivial group, this gives the usual definition of factorization homology.

*Remark* 18.1.1. Recall from Remark 4.2.3 that the category of ($G$ equivariant) $D$ modules on $\mathrm{Ran}_X$ can be defined formally as

$$D(\mathrm{Ran}_X)^G = \lim_{I \in \mathrm{fSet}^{\mathrm{surj}}} D^!(X^I)^G = \lim \left[ \ldots \to D(X^I)^G \xrightarrow{\Delta(\pi)^!} D(X^J)^G \to \ldots \right],$$

so that an object $\mathcal{A} \in D(\mathrm{Ran}_X)^G$ is given, as in Definition 16.0.2, by an assignment

$$I \mapsto A_I \in D(X^I)^G \qquad [\pi : I \twoheadrightarrow J] \mapsto [\Delta(\pi)^! A_I \xrightarrow{\cong} A_J]$$

defined for each finite set $I \in$ fSet and surjection $\pi : I \twoheadrightarrow J$. Alternatively, passing to left adjoints yields the description

$$D(\mathrm{Ran}_X)^G = \mathrm{colim}_{I \in \mathrm{fSet}^{\mathrm{surj}}} D^*(X^I)^G = \mathrm{colim} \left[ \ldots \leftarrow D(X^I)^G \xleftarrow{\Delta(\pi)_*} D(X^J)^G \leftarrow \ldots \right];$$

concretely, applying the $(\Delta(\pi)_*, \Delta(\pi)^!)$ adjunctions to the equivalences

$$(18.1.1) \qquad \Delta(\pi)^! A_I \xrightarrow{\cong} A_J \qquad \text{gives maps} \qquad \Delta(\pi)_* A_J \cong \Delta(\pi)_* \Delta(\pi)^! A_I \to A_I$$

for each $\pi : I \twoheadrightarrow J$.

From the latter description in the preceding remark, the functor $\mathrm{p}_{\mathrm{Ran}_X*} : D(\mathrm{Ran}_X)^G \to D(\mathrm{pt})^G$ is induced by the system of equivariant de Rham cohomology functors

$$\mathrm{p}_{I*} : D(X^I)^G \to D(\mathrm{pt})^G \qquad \text{noting} \qquad \mathrm{p}_{I*} \Delta(\pi)_* = \mathrm{p}_{J*} : D(X^J)^G \to D(\mathrm{pt})^G .$$

Concretely, for $\mathcal{A} = (A_I)_{I \in \mathrm{fSet}^{\mathrm{surj}}} \in D(\mathrm{Ran}_X)^G$ presented in terms of the limit description, we have

$$\mathrm{p}_{\mathrm{Ran}_X*} \mathcal{A} = \mathrm{colim}_{I \in \mathrm{fSet}^{\mathrm{surj}}} \mathrm{p}_{I*} A_I \qquad \text{with diagram structure maps} \qquad \mathrm{p}_{J*} A_J = \mathrm{p}_{I*} \Delta(\pi)_* A_J \to \mathrm{p}_{I*} A_I ,$$

given by the image under $\mathrm{p}_{I*}$ of the maps of Equation 18.1.1 above.

*Definition* 18.1.2. The functor of equivariant factorization homology over $X$ is defined by

$$\int_X^G (\cdot) := \mathrm{p}_{\mathrm{Ran}_X *} \circ \mathrm{oblv} : \mathrm{Alg}^{\mathrm{fact}}(X)^G \to D(\mathrm{pt})^G \qquad \mathcal{A} \mapsto \int_X^G \mathcal{A} = \mathrm{p}_{\mathrm{Ran}_X *}\mathcal{A} = \operatorname*{colim}_{I \in \mathrm{fSet}^{\mathrm{surj}}} \mathrm{p}_{I*}A_I \;.$$

18.2. **Pullback of factorization algebras.** Towards the statement of the equivariant localization theorem in factorization homology, we need to formulate the notion of pullback of factorization algebras. Let $f : X \to Y$ be an equivariant map of smooth algebraic varieties with $G$ action, and let $f^I : X^I \to Y^I$ and $\mathrm{Ran}(f) : \mathrm{Ran}_X \to \mathrm{Ran}_Y$ be the induced maps on products and on Ran spaces.

*Definition* 18.2.1. The pullback of equivariant $D$ modules on the Ran space

$$\mathrm{Ran}(f)^! : D(\mathrm{Ran}_Y)^G \to D(\mathrm{Ran}_X)^G \qquad \text{is defined by} \qquad \mathcal{A} = (A_I)_{I \in \mathrm{fSet}^{\mathrm{surj}}} \mapsto \mathrm{Ran}(f)^!\mathcal{A} = ((f^I)^! A_I)_{I \in \mathrm{fSet}^{\mathrm{surj}}} \;,$$

with gluing data for $\mathrm{Ran}(f)^!\mathcal{A} \in D(\mathrm{Ran}_X)^G$ given by

$$\Delta_X(\pi)^!(\mathrm{Ran}(f)^!\mathcal{A})_I = \Delta_X(\pi)^!(f^I)^! A_I \cong (f^J)^!\Delta_Y(\pi)^! A_I \xrightarrow{\cong} (f^J)^! A_J = (\mathrm{Ran}(f)^!\mathcal{A})_J \;,$$

where the arrow is given by the image of the gluing data $\Delta_Y(\pi)^! A_I \xrightarrow{\cong} A_J$ for $\mathcal{A}$ under $(f^J)^!$.

Now, suppose $\mathcal{A} \in \mathrm{Alg}^{\mathrm{fact}}(Y)^G$ is an equivariant factorization algebra on $Y$, with factorization structure maps

$$A_I \to j(\pi)_* j(\pi)^*(\boxtimes_{j \in J} A_{I_j})$$

for each $I, J$ and $\pi : I \twoheadrightarrow J$. Further, suppose $f : X \to Y$ is a closed embedding, and note that the commutative diagram

$$\begin{array}{ccc} U_X(\pi) & \xrightarrow{j_X(\pi)} & X^I \\ {\scriptstyle f^I}\big\downarrow & & \big\downarrow{\scriptstyle f^I} \\ U_Y(\pi) & \xrightarrow{j_Y(\pi)} & Y^I \end{array}$$

is cartesian for $f$ injective. Then we have:

*Proposition* 18.2.2. The pullback $\mathrm{Ran}(f)^!\mathcal{A} \in \mathrm{Alg}^{\mathrm{fact}}(X)^G$ is canonically an equivariant factorization algebra on $X$.

*Proof.* The structure maps are given by

$$\begin{aligned} (f^I)^! A_I &\to (f^I)^! j_Y(\pi)_* j_Y(\pi)^*(\boxtimes_{j \in J} A_{I_j}) \\ &\cong j_X(\pi)_*(f^I)^! j_Y(\pi)^*(\boxtimes_{j \in J} A_{I_j}) \\ &\cong j_X(\pi)_* j_X(\pi)^*(f^J)^!(\boxtimes_{j \in J} A_{I_j}) \\ &= j_X(\pi)_* j_X(\pi)^*(\boxtimes_{j \in J}(f^{I_j})^! A_{I_j}) \end{aligned}$$

where the arrow is given by the image of the structure maps for $A$ under $(f^I)^!$. These commutative coalgebra structure maps satisfy the factorization property since the map which is required to be a homotopy equivalence is given by

$$j_X(\pi)^*(f^I)^! A_I \cong (f^I)^! j_Y(\pi)^* A_I \to (f^I)^! j_Y(\pi)^*(\boxtimes_{j \in J} A_{I_j}) \cong j_X(\pi)^*(\boxtimes_{j \in J}(f^{I_j})^! A_{I_j})$$

where the arrow is given by the image of the equivalence $j(\pi)^* A_I \to j(\pi)^*(\boxtimes_{j \in J} A_{I_j})$ under $(f^I)^!$.   $\square$

*Remark* 18.2.3. Throughout the remainder of the text we will denote the pullback factorization algebra $\mathrm{Ran}(f)^!\mathcal{A} \in \mathrm{Alg}^{\mathrm{fact}}(X)^G$ by simply $f^!\mathcal{A}$.

Finally, in preparation for the statement of the localization theorem, we note the following property of the pullback of factorization algebras:

*Proposition* 18.2.4. Let $f : X \to Y$ be an equivariant, closed embedding of smooth $G$ varieties, and $\mathcal{A} \in \mathrm{Alg}^{\mathrm{fact}}(X)^G$ an equivariant factorization algebra on $Y$. There is a canonical map

$$\int_X^G f^! \mathcal{A} \to \int_Y^G \mathcal{A} \qquad \text{in } D(\mathrm{pt})^G.$$

*Proof.* For each $I \in \mathrm{fSet}$, the $(f_*^I, (f^I)^!)$ adjunction gives a canonical map
(18.2.1)
$$f_*^I (f^! \mathcal{A})_I = f_*^I (f^I)^! A_I \to A_I \qquad \text{and thus} \qquad \mathrm{p}_{X^I *} (f^! \mathcal{A})_I \cong \mathrm{p}_{Y^I *} f_*^I (f^! \mathcal{A})_I = \mathrm{p}_{Y^I *} f_*^I (f^I)^! A_I \to \mathrm{p}_{Y^I *} A_I .$$

These maps are evidently compatible with the structure maps for the colimit over $I \in \mathrm{fSet}^{\mathrm{surj}}$ and thus induce the desired map

$$\int_X^G f^! \mathcal{A} = \underset{I \in \mathrm{fSet}^{\mathrm{surj}}}{\mathrm{colim}} \, \mathrm{p}_{X^I *} (f^! \mathcal{A})_I \to \underset{I \in \mathrm{fSet}^{\mathrm{surj}}}{\mathrm{colim}} \, \mathrm{p}_{Y^I *} A_I = \int_Y^G \mathcal{A} \ .$$

$\square$

18.3. **The equivariant localization theorem for factorization homology.** We now formulate and prove the analogue of the equivariant localization theorem for factorization homology. Let $G$ be a connected, reductive algebraic group, $X$ a smooth $G$ variety, and $\iota : X^G \hookrightarrow X$ the inclusion of the variety of $G$-fixed points $X^G$. Further, recall from Example 15.0.9 that there is a canonical functor $D(\mathrm{pt})^G \to \mathrm{D}(H_G^\bullet(\mathrm{pt}))$; throughout this section we abuse notation and identify $\int_X^G \mathcal{A} \in D(\mathrm{pt})^G$ with its image under this functor. Finally, following Appendix B.3, for simplicity we restrict to the case that $G = (\mathbb{C}^\times)^n$ is given by an algebraic torus, and choose $\{f_i \in H_G^\bullet(\mathrm{pt})\}$ generators of an ideal whose corresponding subvariety contains the union of the stabilizer subalgebras $\mathfrak{g}_x[2] \hookrightarrow \mathfrak{g}[2] = \mathrm{Spec} \, H_G^\bullet(\mathrm{pt})$ over all non-fixed points $x \in X \backslash X^G$.

*Theorem* 18.3.1. Let $\mathcal{A} \in \mathrm{Alg}^{\mathrm{fact}}(X)^G$ be an equivariant factorization algebra on $X$. The canonical map of Proposition 18.2.4 induces an isomorphism

$$\int_{X^G}^G \iota^! \mathcal{A} \xrightarrow{\cong} \int_X^G \mathcal{A} \qquad \text{over } H_G^\bullet(\mathrm{pt})[f_i^{-1}].$$

*Proof.* First, we note that for each $I \in \mathrm{fSet}$ the fixed points $(X^I)^G = (X^G)^I$ in $X^I$ are given by the $I$-fold product of the fixed points variety $X^G$, and no additional fixed points can occur in the partial colimits $\mathrm{Ran}_{X, \leqslant n}$, since $G$ is connected. Further, the set of possible stabilizer subtori $G_{\underline{x}} \hookrightarrow G$ of points $\underline{x} \in X^I$ is exhausted by those occuring in $X$.

Thus, for each $I \in \mathrm{fSet}$, we can apply Theorem B.3.1 to the map of Equation 18.2.1 to conclude

$$\mathrm{p}_{(X^G)^I *} (\iota^! \mathcal{A})_I = \mathrm{p}_{(X^I)^G *} \iota_I^! A_I \xrightarrow{\cong} \mathrm{p}_{X^I *} A_I \qquad \text{is an isomorphism over } H_G^\bullet(\mathrm{pt})[f_i^{-1}] \ ,$$

where $\iota_I : (X^I)^G \hookrightarrow X^I$ is the inclusion of the $G$ fixed points in $X^I$. It follows that the map induced on colimits as in Proposition 18.2.4 is itself an isomorphism, as claimed. $\square$

## 19. Equivariant and topological vertex algebras

In this section, we define equivariant topological vertex algebras in terms of equivariant chiral algebras, and recover the notion of topological vertex algebra in [Hua94] in the case of $\mathbb{G}_m$ equivariant vertex algebras in dimension 1.

Let $X = \mathbb{A}^n$ be $n$ dimensional affine space, let $G$ act linearly on $\mathbb{A}^n$ by $a : G \times \mathbb{A}^n \to \mathbb{A}^n$, and let $\mathfrak{a} = \mathrm{Lie}(\mathbb{G}_a^n) \cong \mathbb{A}^n$ so that we can interpret the infinitesimal action as $da : \mathfrak{g} \to \mathrm{End}_{\mathbb{K}}(\mathfrak{a})$.

*Definition* 19.0.1. An $n$ dimensional vertex algebra $\mathbb{V}$ is a weakly $\mathbb{G}_a^n$ equivariant chiral algebra $A \in \mathrm{Alg}_{\mathrm{un}}^{\mathrm{ch}}(\mathbb{A}^n)^{\mathbb{G}_a^n, w}$ on $\mathbb{A}^n$.

A topological vertex algebra is an $n$ dimensional vertex algebra $\mathbb{V}$ together with a lift of the weak $\mathbb{G}_a^n$ equivariant structure on the corresponding chiral algebra $A \in \mathrm{Alg}_{\mathrm{un}}^{\mathrm{ch}}(\mathbb{A}^n)^{\mathbb{G}_a^n, w}$ to a strong $\mathbb{G}_a^n$ equivariant structure.

A $G$ equivariant topological vertex algebra is an $n$ dimensional vertex algebra $\mathbb{V}$ together with a lift of the weak $\mathbb{G}_a^n \ltimes \mathbb{G}_a^n$ equivariant structure to a strong $G \ltimes \mathbb{G}_a^n$ equivariant structure $A \in \mathrm{Alg}_{\mathrm{un}}^{\mathrm{ch}}(\mathbb{A}^n)^{G \times \mathbb{G}_a^n}$.

A framed topological vertex algebra is a $\mathbb{G}_a^n \rtimes \mathfrak{so}(2n; \mathbb{K})$ equivariant vertex algebra.

*Remark* 19.0.2. In terms of the vertex algebra data underlying the $n$ dimensional vertex algebra $\mathbb{V} \in \mathrm{VOA}_n$, a $G$ equivariant structure gives the following data:

- The weak $G \ltimes \mathbb{G}_a^n$ equivariant structure yields a $G$ representation $\rho_{\mathbb{V}} : G \to \mathrm{Aut}_{\mathbb{K}}(\mathbb{V})$, such that

  (19.0.1)          $d\rho_{\mathbb{V}} \circ T = T \circ d\rho_{\mathbb{V}} + T \circ da$      as maps      $\mathfrak{g} \times \mathfrak{a} \to \mathrm{End}_{\mathbb{K}}(\mathbb{V})$

  where $d\rho_{\mathbb{V}} : G \to \mathrm{End}_{\mathbb{K}}(\mathbb{V})$ is the corresponding Lie algebra representation, $T : \mathfrak{a} \to \mathrm{End}_{\mathbb{K}}(\mathbb{V})$ is the translation operator and $da : \mathfrak{g} \to \mathrm{End}_{\mathbb{K}}(\mathfrak{a})$ is the infinitesimal action map.

- The strong $\mathbb{G}_a^n$ equivariant structure yields a Lie algebra map

  $g_{-1} : \mathfrak{a} \to \mathrm{Der}_{\mathrm{VOA}_n}^{-1}(\mathbb{V})$        such that        $[d_{\mathbb{V}}, g_{-1}] = T : \mathfrak{a} \to \mathrm{Der}_{\mathrm{VOA}_n}(\mathbb{V})$ .

  This is interpreted as a homotopy trivializing the translation operator.

- The compatible strong $G$ equivariant structure yields a Lie algebra map

  $h_{\mathbb{V}} : \mathfrak{g} \to \mathrm{End}_{\mathbb{K}}^{-1}(\mathbb{V})$        such that        $[d_{\mathbb{V}}, h_{\mathbb{V}}] = d\rho_{\mathbb{V}} : \mathfrak{g} \to \mathrm{End}_{\mathbb{K}}(\mathbb{V})$ .

  This is interpreted as a homotopy trivializing the infinitesimal action $d\rho_{\mathbb{V}} : \mathfrak{g} \to \mathrm{End}_{\mathbb{K}}(\mathbb{V})$. The endomorphisms $d\rho_{\mathbb{V}}$ and $h_{\mathbb{V}}$ do not act by vertex algebra derivations, as is apparent from equation 19.0.1, but act by derivations twisted by $T \circ da$.

*Example* 19.0.3. Concretely, a topological vertex algebra in dimensional 1 is just a DG vertex algebra $(\mathbb{V}, d)$ together with a degree $-1$ derivation $g_{-1} : \mathrm{Der}^{-1}(\mathbb{V})$ trivializing the translation operator, that is, such that $[d, g_{-1}] = T$. This is equivalent to a particular subset of the structure of a strong topological vertex algebra in [Hua94], as we explain in Example 19.0.4 below. See also Section 20 below for the relation to $\mathbb{E}_{2n}$ algebras.

*Example* 19.0.4. A framed topological vertex algebra in dimension 1 is a graded DG vertex algebra $(\mathbb{V}, d)$ together with degree $-1$ endomorphisms $g_0, g_{-1} \in \mathrm{End}^{-1}(\mathbb{V})$ of graded degrees 0 and 1, such that

- $[d, g_{-1}] = T$,
- $[d, g_0] = L_0$,
- $[T, g_{-1}] = 0$ and moreover $g_{-1}$ acts by derivations of $\mathbb{V}$, and
- $[T, g_0] = -g_{-1}$ and moreover $g_0$ acts by twisted derivations of $\mathbb{V}$.

This agrees with the notion of topological vertex algebra in [Hua94].

## 20. $\mathbb{G}_a^n$ EQUIVARIANT FACTORIZATION ALGEBRAS ON $\mathbb{A}^n$ AND $E_{2n}$ ALGEBRAS

In this section, we sketch a proof of the folklore result that translation invariant factorization algebras on $\mathbb{A}^n$ over $\mathbb{K} = \mathbb{C}$ are equivalent to $\mathbb{E}_{2n}$ algebras. Let $X = \mathbb{A}_{\mathbb{C}}^n$ be $n$ dimensional complex affine space and let $G = \mathbb{G}_a^n$ act on $\mathbb{A}^n$ by translation. The action of $\mathbb{G}_a^n$ on $\mathbb{A}^n$ is free and transitive, so we have an equivalence of categories

$$(20.0.1) \qquad D(\mathbb{A}^n)^{\mathbb{G}_a^n} \underset{\longleftarrow}{\overset{\longrightarrow}{\rule{1cm}{0pt}}} \mathrm{Vect}_{\mathbb{K}} \qquad \text{defined by} \qquad M \mapsto C_{\mathrm{dR}}^\bullet(\mathbb{A}^n, M) \qquad V \otimes \omega_{\mathbb{A}^n} \hookleftarrow V \ .$$

*Proposition* 20.0.1. There is an equivalence of categories

$$\mathrm{Alg}_{\mathrm{un}}^{\mathrm{ch}}(\mathbb{A}_{\mathbb{C}}^n)^{\mathbb{G}_a^n} \overset{\cong}{\longrightarrow} \mathrm{Alg}_{\mathbb{E}_{2n}}(\mathrm{Vect}_{\mathbb{K}}) \qquad \text{defined by} \qquad A \mapsto C_{\mathrm{dR}}^\bullet(\mathbb{A}^n, A) \ .$$

*Proof.* Let $A \in \mathrm{Alg}_{\mathrm{un}}^{\mathrm{ch}}(\mathbb{A}^n)^{\mathbb{G}_a^n}$ and recall from Warning 7.0.4 that this notation refers to the underlying object $A \in D(\mathbb{A}^n)$. Then $A = \omega_{\mathbb{A}_n} \otimes V$ where $V = C_{\mathrm{dR}}^\bullet(\mathbb{A}^n, A)$, and we exhibit an equivalence between the chiral algebra structure maps on $A$ and $\mathbb{E}_{2n}$ algebra structure maps on $V$, natural in $A \in D(\mathbb{A}^n)^{\mathbb{G}_a^n}$ and correspondingly $V \in \mathrm{Vect}_{\mathbb{K}}$.

The data of a chiral algebra structure on $A$ is given by compatible structure maps

$$(20.0.2) \qquad b_I \in \mathrm{Hom}_{D(X)^{\mathrm{ch}}}(\{A\}_{i \in I}, A)^{\mathbb{G}_a^n} \qquad \text{defined by maps} \qquad b_I : j_*^{(I)} j^{(I),*}(\boxtimes_{i \in I} A) \to \Delta_*^{(I)} A$$

in $D(X^I)^{\mathbb{G}_a^n}$ as in Example 17.0.5. Applying the equivalence of Equation 20.0.1 above, we find

$$\mathrm{Hom}_{D(X)^{\mathrm{ch}}}(\{A\}_{i \in I}, A)^{\mathbb{G}_a^n} \cong \mathrm{Hom}_{D(X^I)^{\mathbb{G}_a^n}}(j_*^{(I)} j^{(I),*} \omega_{X^I}, \Delta_*^{(I)} \omega_X) \otimes_{\mathbb{K}} \mathrm{Hom}_{\mathrm{Vect}_{\mathbb{K}}}(V^{\otimes I}, V)$$

Moreover, we have

$$\mathrm{Hom}_{D(X^I)^{\mathbb{G}_a^n}}(j_*^{(I)} j^{(I),*} \omega_{X^I}, \Delta_*^{(I)} \omega_X) \cong C_c^\bullet(\mathrm{Conf}^I(\mathbb{A}^n)) \ .$$

Thus, the required structure maps of Equation 20.0.2 are equivalent to structure maps

$$V^{\otimes I} \to V \otimes_{\mathbb{K}} C_c^\bullet(\mathrm{Conf}^I(\mathbb{A}^n)) \qquad \text{or equivalently} \qquad C_\bullet(\mathrm{Conf}^I(\mathbb{A}^n); \mathbb{K}) \to \mathrm{Hom}_{\mathrm{Vect}_{\mathbb{K}}}(V^{\otimes I}, V) \ .$$

defining $V \in \mathrm{Alg}_{\mathbb{E}_{2n}}(\mathrm{Vect}_{\mathbb{K}})$. Similarly, one checks that morphisms in $\mathrm{Alg}_{\mathrm{un}}^{\mathrm{ch}}(\mathbb{A}^n)^{\mathbb{G}_a^n}$ of such chiral algebras are equivalent to maps of the corresponding $\mathbb{E}_{2n}$ algebras. $\qquad\square$

*Example* 20.0.2. Consider the arity 2 chiral structure map

$$\mu_2 : j_* j^*(A^{\boxtimes 2}) \to \Delta_* A \ ,$$

where $\Delta : \mathbb{A}^n \to \mathbb{A}^{2n}$ is the diagonal embedding and $j : \mathbb{A}^{2n} \backslash \Delta \to \mathbb{A}^{2n}$ is the complementary open embedding, and the map is in the category $\mathrm{D}_G(\mathbb{A}^{2n})$. Then applying

$$\mathrm{Hom}(\cdot, \Delta_* A) : D(\mathbb{A}^{2n}) \to \mathrm{Vect}_{\mathbb{K}} \qquad \text{to the exact triangle} \qquad A^{\boxtimes 2} \hookrightarrow j_* j^*(A^{\boxtimes 2}) \twoheadrightarrow \Delta_* A[1]$$

induces cochain maps to $\mathrm{Hom}_{\mathrm{Vect}_{\mathbb{K}}}(V^{\otimes 2}, V)$ from the exact sequence

$$C^\bullet(\mathrm{pt}) \hookrightarrow C_c^\bullet(\mathbb{A}^n \backslash \{0\})[1] \to C_c^\bullet(\mathbb{A}^n)[1] \qquad \text{or concretely} \qquad \mathbb{K}_m \hookrightarrow \mathbb{K}_m \oplus \mathbb{K}_\pi[2n-1] \twoheadrightarrow \mathbb{K}_\pi[2n-1] \ .$$

Thus, we see that the data of the Lie$^*(X)$ algebra underlying a chiral algebra, which also determines the Poisson vertex structure on the associated graded, corresponds to the shifted Poisson bracket of the corresponding homology $\mathbb{P}_n$ algebra, and that the induced $\mathrm{Comm}^!(X)$ structure coming from the necessarily non-singular chiral bracket corresponds to the commutative multiplication underlying the $\mathbb{P}_n$ algebra. Recall that these structures on a chiral algebra were discussed in Section 10 and summarized in the diagram of Equation 10.5.1, and the relevant descriptions of the $\mathbb{E}_n$ and $\mathbb{P}_n$ operads are summarized in appendices C.4 and C.5.

*Example* 20.0.3. Let $(\mathbb{V}, d_{\mathbb{V}}, T, g_{-1})$ be a topological vertex algebra as in Definition 19.0.1. The corresponding translation invariant chiral algebra $A \in \mathrm{Alg}_{\mathrm{un}}^{\mathrm{ch}}(\mathbb{A}^1)^{\mathbb{G}_a}$ is given by

$$A \cong \mathbb{V} \otimes_{\mathbb{K}} (\mathbb{K}[x, \partial_x] \oplus \mathbb{K}[x, \partial_x][-1]) \ , \quad d = d_{\mathbb{V}} \otimes \mathbb{1} + \mathbb{1} \otimes d_{\mathrm{dR}}$$

with $h = g_{-1} \otimes \mathbb{1} + \mathbb{1} \otimes \iota_{\partial_x}$; see also Example 15.0.8. The above gives an equivalence $\mathbb{V} \mapsto C_{\mathrm{dR}}^\bullet(\mathbb{A}^n, A)$ between topological vertex algebras of dimension $n$ and $\mathbb{E}_{2n}$ algebras.

## 21. Equivariant and semidirect product operads

In this section, we recall the formalism of equivariant and semidirect product operads following [SW03]. Similar results are discussed in [Wes07] in the homotopy setting. Let $\mathcal{C}$ be a cartesian symmetric monoidal category, $G \in \mathrm{Grp}(\mathcal{C})$ be a group object of $\mathcal{C}$, and let $G\text{-Mod}(\mathcal{C})$ denote the category of objects $C \in \mathcal{C}$ with an action of $G$ on $C$ and morphisms those in $\mathcal{C}$ equipped with $G$ equivariant structure.

*Remark* 21.0.1. The category $G\text{-Mod}(\mathcal{C})$ is naturally symmetric monoidal with respect to the underlying monoidal structure on $\mathcal{C}$.

*Example* 21.0.2. Let $\mathcal{C} = \mathrm{Top}$ be the category of spaces. A group object $G \in \mathrm{Grp}(\mathrm{Top})$ is a topological group and $G\text{-Mod}(\mathcal{C}) = \mathrm{Top}_G$ is the category of $G$ spaces. Similarly, for $\mathcal{C} = \mathrm{Sch}_{\mathbb{K}}$ the categoy of schemes, a group object is an algebraic group over $\mathbb{K}$ and $G\text{-Mod}(\mathrm{Sch}_{\mathbb{K}})$ is the category of $G$ schemes.

*Example* 21.0.3. Let $\mathcal{C} = \mathrm{CoComm}(\mathrm{Vect}_{\mathbb{K}})$ be the category of cocommutative coalgebras in $\mathrm{Vect}_{\mathbb{K}}$. The structure maps of a group object $G \in \mathrm{Grp}(\mathcal{C})$ define a compatible product, unit, and antipode on $G$, so that $G$ itself is naturally a cocommutative Hopf algebra $\Lambda \in \mathrm{Hopf}^{\mathrm{co}}(\mathrm{Vect}_{\mathbb{K}})$. There is a natural equivalence $G\text{-Mod}(\mathcal{C}) \cong \Lambda\text{-Mod}(\mathrm{CoComm}(\mathrm{Vect}_{\mathbb{K}}))$ and the induced symmetric monoidal structure on $G\text{-Mod}(\mathcal{C})$ is that corresponding to the coproduct on $\Lambda$.

*Proposition* 21.0.4. Let $F : \mathcal{C} \to \mathcal{C}'$ be a symmetric monoidal functor of cartesian categories. Then $F$ induces a functor $F : \mathrm{Grp}(\mathcal{C}) \to \mathrm{Grp}(\mathcal{C}')$ and symmetric monoidal functors $F_G : G\text{-Mod}(\mathcal{C}) \to F(G)\text{-Mod}(\mathcal{C}')$ for each $G \in \mathrm{Grp}(\mathcal{C})$.

*Example* 21.0.5. The functor $C_\bullet(\cdot; \mathbb{K}) : \mathrm{Top} \to \mathrm{CoComm}(\mathrm{Vect}_{\mathbb{K}})$ of Remark C.2.7 is symmetric monoidal. Thus, each group object $G \in \mathrm{Top}$ defines a cocommutative Hopf algebra $\Lambda = C_\bullet(G; \mathbb{K})$, and the induced symmetric monoidal functor $C_\bullet(\cdot; \mathbb{K}) : \mathrm{Top}_G \to \Lambda\text{-Mod}(\mathrm{CoComm}(\mathrm{Vect}_{\mathbb{K}}))$ restricts to that of Remark B.2.2.

*Definition* 21.0.6. A $G$ equivariant operad in $\mathcal{C}$ is an operad in the category $G\text{-Mod}(\mathcal{C})$ with its induced symmetric monoidal structure.

Let $\mathrm{Op}_G(\mathcal{C}) = \mathrm{Op}(G\text{-Mod}(\mathcal{C}))$ denote the category of $G$ equivariant operads.

*Definition* 21.0.7. Let $\mathcal{O} \in \mathrm{Op}_G(\mathcal{C})$ be a G operad. The semidirect product operad $\mathcal{O} \rtimes G \in \mathrm{Op}(\mathcal{C})$ is defined by

$$\mathrm{col}\,(\mathcal{O} \rtimes G) = \mathrm{col}\,\mathcal{O} \qquad (\mathcal{O} \rtimes G)(\{c_i\}, d) = \mathcal{O}(\{c_i\}, d) \otimes G^{\otimes |I|}$$

together with composition maps for each $\pi : I \to J$ given by

$$(21.0.1) \quad \bigotimes_{j \in J} \mathcal{O}(\{c_i\}_{i \in I_j}, d_j) \otimes G^{\otimes |I_j|} \otimes \mathcal{O}(\{d_j\}_{j \in J}, e) \otimes G^{\otimes |J|} \qquad \to \qquad \mathcal{O}(\{c_I\}_{i \in I}, e) \otimes G^{\otimes |I|}$$

$$(21.0.2) \quad \otimes_j (\beta_j, (g_i)_{i \in I_j}) \otimes (\alpha, (g_j)_{j \in J}) \qquad \mapsto \qquad \otimes_j (g_j \cdot \beta_j) \circ \alpha, (g_i \cdot g_{\pi(i)})_{i \in I}$$

and units defined by $\mathbb{1}_c^G = \mathbb{1}_c \otimes e \in \mathcal{O}(c, c) \otimes G$, for $e : \mathrm{u}_{\mathcal{C}} \to G$ the identity structure map.

*Proposition* 21.0.8. Let $\mathcal{O} \in \mathrm{Op}_G(\mathcal{C})$ be a $G$ operad. There is a natural symmetric monoidal equivalence

$$\mathrm{Alg}_{\mathcal{O} \rtimes G}(\mathcal{C}) \xrightarrow{\cong} \mathrm{Alg}_{\mathcal{O}}(G\text{-}\mathrm{Mod}(\mathcal{C})) \ .$$

*Proof.* The composition $\mathcal{O} \to \mathcal{O} \rtimes G \to \mathcal{C}^{\otimes}$ defines an object of $\mathrm{Alg}_{\mathcal{O}}(\mathcal{C})$. The composition maps for $\mathcal{O} \rtimes G$ define a lift to a map $\mathcal{O} \to G\text{-}\mathrm{Mod}(\mathcal{C})^{\otimes}$ in $\mathrm{Op}(\mathcal{C})$. The associativity data for the composition law in $\mathcal{O} \rtimes G$ defines equivariance data lifting the map $\mathcal{O} \to G\text{-}\mathrm{Mod}(\mathcal{C})^{\otimes}$ to $\mathrm{Op}(G\text{-}\mathrm{Mod}(\mathcal{C}))$.

Note that although we have used the language of $\mathcal{C}^{\otimes} \in \mathrm{Op}(\mathcal{C})$ which can only be interpretted literally for $\mathcal{C}$ closed, the arguement extends naturally to the general setting via hom tensor adjunction. $\square$

## 22. The $K$ equivariant little $d$-cubes operad $\mathbb{E}_d^K$

In this section, we recall the construction of the $K$ equivariant little $d$-cubes operad, following Section 5.4.2 of [Lur12], and references therein. Throughout, let $K$ be a connected topological group and $\mathrm{Top}(d) = \mathrm{Aut}_{\mathrm{Top}}(\mathbb{R}^d)$ denote the topological automorphism group of $\mathbb{R}^d$, which naturally defines $\mathrm{Top}(d) \in \mathrm{Grp}(\mathrm{Top})$ a topological group.

The action of $\mathrm{Top}(d)$ on $\mathbb{R}^d$ induces an action on $\mathrm{Conf}^I(\mathbb{R}^d)$ for each finite set $I$, so that the little $d$-cubes operad $\mathbb{E}_d \in \mathrm{Op}(\mathrm{Top}(d)\text{-}\mathrm{Mod}(\mathrm{Top}))$ is naturally a $\mathrm{Top}(d)$ equivariant operad in $\mathrm{Top}$. More generally, for any map of topological groups $\rho : K \to \mathrm{Top}(d)$, we obtain the structure of a $K$ equivariant operad in $\mathrm{Top}$ on $\mathbb{E}_d \in \mathrm{Op}(K\text{-}\mathrm{Mod}(\mathrm{Top}))$.

*Definition* 22.0.1. The $K$ equivariant little $d$-cubes operad $\mathbb{E}_d^K = \mathbb{E}_d \rtimes K \in \mathrm{Op}(\mathrm{Top})$ is the semidirect product of $\mathbb{E}_d$ with $K$ under the action of $\rho$.

*Remark* 22.0.2. More concretely, the $K$ equivariant little $d$-cubes operad $\mathbb{E}_d^K$ is presented by

$$\mathbb{E}_d^K(I) = \mathrm{Conf}^I(\mathbb{R}^d) \times K^I \qquad \text{with} \qquad (\times_{j \in J}\mathrm{Conf}^{I_j}(\mathbb{R}^d) \times K^{I_j}) \times \mathrm{Conf}^J(\mathbb{R}^d) \times K^J \to \mathrm{Conf}^I(\mathbb{R}^d) \times K^I$$

specified up to homotopy equivalence by group multiplication along $\pi : I \twoheadrightarrow J$ in the $K$ factors, together with the (homotopy equivalence class of) composition map on the $\mathrm{Conf}(\mathbb{R}^d)$ factors determined by the operad structure on $\mathbb{E}_d$, twisted by the action of $K$ on the configuration spaces according to the formula 21.0.2. A strict model for this operad is given by the skew little cubes operad of [DHK18], for example.

*Remark* 22.0.3. The little $d$-cubes operad together with the $K$ equivariant structure above defines

$$C_{\bullet}^{\rho}(\mathbb{E}_d) \in \mathrm{Op}(H_{\bullet}(K)\text{-}\mathrm{Mod}(\mathrm{CoComm}(\mathrm{Vect}_{\mathbb{K}}))) \ ,$$

by Example 21.0.5. We abuse notation and denote by $C_{\bullet}^{\rho}(\mathbb{E}_d) \in \mathrm{Op}(H_{\bullet}(K)\text{-}\mathrm{Mod}(\mathrm{Vect}_{\mathbb{K}}))$ its image under the forgetful functor to $H_{\bullet}(K)\text{-}\mathrm{Mod}(\mathrm{Vect}_{\mathbb{K}})$.

*Example* 22.0.4. Let $K = \mathrm{SO}(d)$ and $\rho : \mathrm{SO}(d) \to \mathrm{Top}(d)$ the canonical inclusion. The operad $\mathbb{E}_d^{\mathrm{fr}} = \mathbb{E}_d^{\mathrm{SO}(d)}$ is the *framed* or *oriented* little $d$-cubes operad.

*Example* 22.0.5. The framed little 2-cubes operad $\mathbb{E}_2^{S^1} \in \mathrm{Op}(\mathrm{Top})$ was introduced in [Get94a]. The operad $\mathbb{E}_2^{S^1}$ is formal [GS10], and its homology operad $H_{\bullet}(\mathbb{E}_2^{S^1})$ is the Batalin-Vilkovisky operad $\mathrm{BV} \in \mathrm{Op}(\mathrm{Vect}_{\mathbb{Z}})$ [Get94a], which is generated in arity 1 and 2 by

$$\mathrm{BV}(1) = \mathbb{K}_{\Delta}[-1]\langle 1 \rangle \qquad \mathrm{BV}(2) = \mathbb{P}_2(2)$$

subject to the usual relations of the $\mathbb{P}_2$ operad, as in Definition C.5.1, together with the additional relation

$$(22.0.1) \qquad \Delta \circ m - m \circ (\Delta \otimes \mathbb{1}) - m \circ (\mathbb{1} \otimes \Delta) = \pi \; ,$$

where $m : A^{\otimes 2} \to A$, and $\pi : A^{\otimes 2} \to A[-1]$ are the commutative multiplication and Poisson bracket structure maps.

The derived category $\mathrm{D}^b_{\mathrm{fg}}(\Lambda)$ of graded modules over $\Lambda = H_\bullet(S^1) = \mathbb{C}[\varepsilon]/\varepsilon^2$ is equivalent to the category of mixed complexes. The conclusion of Proposition 21.0.8 in this example is that there is an equivalence between BV algebras $A \in \mathrm{Alg}_{\mathrm{BV}}(\mathrm{Vect}_{\mathbb{K}})$ and $\mathbb{P}_2$ algebras in the category of mixed complexes $A \in \mathrm{Alg}_{\mathbb{P}_2}(\mathrm{D}(\Lambda)_{\mathbb{Z}})$ such that the mixed differential $\Delta = \rho([S^1]) : A \to A[-1]$ satisfies the relation 22.0.1 above. This was observed in [Get94b], for example.

The Dunn additivity Theorem, recalled in Theorem C.4.9, admits the following equivariant enhancement, which was established in Remark 5.4.2.14 following Theorem 5.1.2.2 in [Lur12]:

*Theorem* 22.0.6. There is a natural equivalence of operads

$$\mathbb{E}_d^K \star \mathbb{E}_{d'}^{K'} \cong \mathbb{E}_{d+d'}^{K \times K'} \; .$$

## 23. Goresky-Kottwitz-MacPherson Koszul duality for equivariant operads

In this section, we explain an application of the Goresky-Kottwitz-MacPherson result describing equivariant cohomology in terms of Koszul duality, in the context of equivariant operads following Section 21. Let $\mathcal{C} = \mathrm{CoComm}(\mathrm{Vect}_{\mathbb{K}})$ as in Example 21.0.3, so that $G \in \mathrm{Grp}(\mathcal{C})$ naturally defines $\Lambda \in \mathrm{Hopf}^{\mathrm{co}}(\mathrm{Vect}_{\mathbb{K}})$ and we identify $G\text{-Mod}(\mathcal{C}) \cong \Lambda\text{-Mod}(\mathrm{CoComm}(\mathrm{Vect}_{\mathbb{K}}))$. In this case, the Proposition 21.0.8 gives for each $\mathcal{O} \in \mathrm{Op}_G(\mathcal{C})$ a natural a symmetric monoidal equivalence

$$\mathrm{Alg}_{\mathcal{O} \rtimes G}(\mathrm{CoComm}(\mathrm{Vect}_{\mathbb{K}})) \cong \mathrm{Alg}_{\mathcal{O}}(\Lambda\text{-Mod}(\mathrm{CoComm}(\mathrm{Vect}_{\mathbb{K}}))) \; .$$

In particular, this equivalence identifies algebras in the essential images of the free functor to $\mathrm{CoComm}(\mathrm{Vect}_{\mathbb{K}})$, inducing an equivalence

$$(23.0.1) \qquad \mathrm{Alg}_{\mathcal{O} \rtimes G}(\mathrm{Vect}_{\mathbb{K}}) \cong \mathrm{Alg}_{\mathcal{O}}(\Lambda\text{-Mod}(\mathrm{Vect}_{\mathbb{K}})_{\mathbb{Z}}) \; .$$

Now, for simplicity we restrict to bounded, finitely generated, derived categories as in the statements of the summary theorem B.2.4. Then applying the results of *loc. cit.* together with the above discussion, we obtain:

*Proposition* 23.0.1. Let $G$ be a connected Lie group, and consider the graded algebras $\Lambda = H_\bullet(G; \mathbb{K})$ and $\mathbb{S} = H^\bullet_G(\mathrm{pt}; \mathbb{K})$. Further, let $\mathcal{O} \in \mathrm{Op}(\Lambda\text{-Mod}(\mathrm{CoComm}(\mathrm{Vect}_{\mathbb{K}})))$ be a $G$ equivariant operad in $\mathrm{CoComm}(\mathrm{Vect}_{\mathbb{K}})$. Then there are natural symmetric monoidal equivalences

$$(23.0.2) \qquad \mathrm{Alg}_{\mathcal{O} \rtimes \Lambda}(\mathrm{Perf}_{\mathbb{K}}) \xrightarrow{\cong} \mathrm{Alg}_{\mathcal{O}}(\mathrm{D}^b_{\mathrm{fg}}(\Lambda)) \xrightarrow{\cong} \mathrm{Alg}_{t(\mathcal{O})}(\mathrm{D}^b_{\mathrm{fg}}(\mathbb{S}))$$

where $t : \mathrm{D}^b_{\mathrm{fg}}(\Lambda) \to \mathrm{D}^b_{\mathrm{fg}}(\mathbb{S})$ is the Koszul duality functor B.2.1 extended as in Proposition C.1.14.

In particular, if $\mathcal{O} \in \mathrm{Op}(G\text{-Mod}(\mathrm{Top}_f))$ is a $G$ equivariant operad in $\mathrm{Top}_f$, there are natural symmetric monoidal equivalences

$$(23.0.3) \qquad \mathrm{Alg}_{C_\bullet(\mathcal{O}) \rtimes \Lambda}(\mathrm{Perf}_{\mathbb{K}}) \xrightarrow{\cong} \mathrm{Alg}_{C_\bullet(\mathcal{O})}(\mathrm{D}^b_{\mathrm{fg}}(\Lambda)) \xrightarrow{\cong} \mathrm{Alg}_{C^G_\bullet(\mathcal{O})}(\mathrm{D}^b_{\mathrm{fg}}(\mathbb{S}))$$

where $C^G_\bullet(\mathcal{O}) = t(C_\bullet(\mathcal{O})) \in \mathrm{Op}(\mathrm{D}^b_{\mathrm{fg}}(\mathbb{S}\text{-Mod}))$ denotes the $G$ equivariant chains on $\mathcal{O}$, as in Definition 23.0.2 below.

*Proof.* The first equivalence of 23.0.2 is just the restriction of 23.0.1 to bounded, finitely generated derived categories. The second equivalence of 23.0.2 follows from applying Proposition C.1.14 to the symmetric monoidal equivalence of Theorem B.2.4. □

*Definition* 23.0.2. Let $\mathcal{O} \in \mathrm{Op}(G\text{-Mod}(\mathrm{Top}))$ be an equivariant operad. Then equivariant chains operad on $\mathcal{O}$ is the operad

$$C_\bullet^G(\mathcal{O}) \in \mathrm{Op}(\mathrm{D}_\mathrm{fg}^b(\mathbb{S})) \qquad \text{defined by} \qquad C_\bullet^G(\mathcal{O})(I) = C_\bullet^G(\mathcal{O}(I)) \ .$$

*Example* 23.0.3. Let $\rho : K \to \mathrm{Top}(d)$ and $C_\bullet^\rho(\mathbb{E}_d) \in \mathrm{Op}(\mathrm{D}_\mathrm{fg}^b(\Lambda))$ as in Definition 22.0.1 and Remark 22.0.3. Then the above proposition gives an equivalence between $K$ framed $\mathbb{E}_d$ algebras $A \in \mathrm{Alg}_{C_\bullet(\mathbb{E}_d^K)}(\mathrm{Vect}_\mathbb{K})$ in $\mathrm{Vect}_\mathbb{K}$ and algebras $A \in \mathrm{Alg}_{C_\bullet^K(\mathbb{E}_d)}(\mathrm{D}_\mathrm{fg}^b(\mathbb{S}))$ over the operad

$$C_\bullet^K(\mathbb{E}_d) \in \mathrm{Op}(\mathrm{D}_\mathrm{fg}^b(\mathbb{S})) \qquad \text{given by} \qquad C_\bullet^K(\mathbb{E}_d)(I) = C_\bullet^K(\mathrm{Conf}^I(\mathbb{R}^d)) \ .$$

*Example* 23.0.4. In particular, consider the framed little 2-cubes operad $\mathbb{E}_2^{S^1} \in \mathrm{Op}(\mathrm{Top})$, and let $\mathbb{S} = H_{S^1}^\bullet(\mathrm{pt}) \cong \mathbb{K}[u]$ where $u$ is the cohomological degree 2 generator. The corresponding equivariant chains operad

$$\mathbb{BD}_0^u := C_\bullet^{S^1}(\mathbb{E}_2) \in \mathrm{Op}(\mathrm{D}_\mathrm{fg}^b(\mathbb{K}[u]))$$

defines a 2 periodic analogue of the operad $\mathbb{BD}_0^h \in \mathrm{Op}(\mathrm{D}_\mathrm{fg}^b(\mathbb{K}[\hbar]))$ of Definition C.6.2. In particular, $\mathbb{BD}_0^u$ interpolates between the $\mathbb{P}_2$ operad and the $\mathbb{E}_0$ operad: it is generated in degree 2 by

$$\mathbb{BD}_0^u(2) := C_\bullet^{S^1}(\mathbb{E}_2)(2) \cong \left[ \mathbb{K}[u]_m \xrightarrow{m_u} \mathbb{K}[u]_\pi[1]\langle 1 \rangle \right] \quad \in \quad \mathrm{D}_\mathrm{fg}^b(\mathbb{K}[u][S_2])$$

where $\mathbb{K}[u]_m$ is the trivial representation and $\mathbb{K}[u]_\pi$ is the sign, subject to the relations of the $\mathbb{P}_2$ operad of Definition C.5.1 extended linearly to $\mathbb{K}[u]$. This can be understood explicitly via formality by applying the Koszul duality functor of Example B.2.3 to the explicit presentation from Definition C.5.1 of the $\mathbb{P}_2$ operad.

Thus, applied to this example, Proposition 23.0.1 gives a symmetric monoidal equivalence

$$(23.0.4) \qquad \mathrm{Alg}_{\mathbb{E}_2^{S^1}}(\mathrm{Perf}_\mathbb{K}) \cong \mathrm{Alg}_{\mathbb{BD}_0^u}(\mathrm{D}_\mathrm{fg}^b(\mathbb{K}[u])) \ .$$

Motivated by the strong Poisson additivity theorem of Rozenblyum, we also make the following definition:

*Definition* 23.0.5. The operad $\mathbb{BD}_n^u \in \mathrm{Op}(\mathrm{D}^b(\mathbb{K}[u]))$ is defined as

$$\mathbb{BD}_n^u := \mathbb{E}_n \circ \mathbb{BD}_0^u \in \mathrm{Op}(\mathrm{D}^b(\mathbb{K}[u]))$$

the Boardman-Vogt tensor product of the operad $\mathbb{BD}_0^u \in \mathrm{Op}(\mathrm{D}^b(\mathbb{K}[u]))$ defined in Example 23.0.4 above, with the operad $\mathbb{E}_n \in \mathrm{Op}(\mathrm{Perf}_\mathbb{K})$ of Definition C.4.2.

## 24. $\mathbb{G}_a^n \rtimes G$ EQUIVARIANT FACTORIZATION ALGEBRAS ON $\mathbb{A}^n$ AND $\mathbb{E}_{2n}^K$ ALGEBRAS

In this section, we extend the identification of Section 20 to identify $\mathbb{G}_a^n \rtimes G$ equivariant factorization algebras on complex affine space with equivariant $\mathbb{E}_{2n}$ algebras for $K$ the maximal compact of $G$. Let $X = \mathbb{A}_\mathbb{C}^n$ be $n$ dimensional complex affine space and let $G = \mathbb{G}_a^n$ act on $\mathbb{A}^n$ by translation. Let $G$ be a complex reductive group with maximal compact subgroup $K$, and $\rho : G \to \mathrm{Aut}(\mathbb{A}_\mathbb{C}^n)$ a linear action of $G$ on $\mathbb{A}_\mathbb{C}^n$. Then we have:

*Proposition* 24.0.1. There is an equivalence of categories

$$\mathrm{Alg}_\mathrm{un}^\mathrm{ch}(\mathbb{A}_\mathbb{C}^n)^{\mathbb{G}_a^n \rtimes G} \xrightarrow{\cong} \mathrm{Alg}_{\mathbb{E}_{2n}^K}(\mathrm{Vect}_\mathbb{K}) \qquad \text{defined by} \qquad A \mapsto C_\mathrm{dR}^\bullet(\mathbb{A}^n, A) \ .$$

*Proof.* Following the proof of Proposition 20.0.1, the data of a $\mathbb{G}_a^n \rtimes G$ equivariant chiral algebra structure on $A$ is given by compatible structure maps

$$b_I \in \operatorname{Hom}_{D(X)^{\mathrm{ch}}}(\{A\}_{i \in I}, A)^{\mathbb{G}_a^n \rtimes G} . \tag{24.0.1}$$

Applying the equivalence of Equation 20.0.1 above, we find

$$\operatorname{Hom}_{D(X)^{\mathrm{ch}}}(\{A\}_{i \in I}, A)^{\mathbb{G}_a^n \rtimes G} \cong \operatorname{Hom}_{D(X^I)^{\mathbb{G}_a^n \rtimes G}}(j_*^{(I)} j^{(I),*} \omega_{X^I}, \Delta_*^{(I)} \omega_X) \otimes_{\mathbb{K}} \operatorname{Hom}_{\mathrm{Vect}_{\mathbb{K}}}(V^{\otimes I}, V)$$

Moreover, we have

$$\operatorname{Hom}_{D(X^I)^{\mathbb{G}_a^n \rtimes G}}(j_*^{(I)} j^{(I),*} \omega_{X^I}, \Delta_*^{(I)} \omega_X) \cong C_{G,c}^\bullet(\operatorname{Conf}^I(\mathbb{A}^n)) .$$

Thus, the required structure maps of Equation 20.0.2 are equivalent to structure maps

$$V^{\otimes I} \to V \otimes_{\mathbb{K}} C_{G,c}^\bullet(\operatorname{Conf}^I(\mathbb{A}^n)) \qquad \text{or equivalently} \qquad C_\bullet^\bullet(\operatorname{Conf}^I(\mathbb{A}^n); \mathbb{K}) \to \operatorname{Hom}_{\mathrm{Vect}_{\mathbb{K}}}(V^{\otimes I}, V) .$$

defining $V \in \operatorname{Alg}_{\mathbb{E}_{2n}^K}(\mathrm{Vect}_{\mathbb{K}})$, as desired. $\qquad \square$

## 25. Deformation quantization in the Omega background

In this section, we explain the relationship between $\mathbb{G}_m$ equivariant factorization algebras and quantization.

### 25.1. **Quantization of $\mathbb{E}_n$ algebras.**
To begin, we explain the interpretation of the equivalence in Equation 23.0.4 of Example 23.0.4 as relating $S^1$ equivariance data on $\mathbb{E}_{n+2}$ algebras to 'two-periodic graded quantizations' of their homology $\mathbb{P}_{n+2}$ algebras to $\mathbb{E}_n$ algebras. Throughout, we again let $\mathbb{K}[u] = H_{S^1}^\bullet(\mathrm{pt})$ be the $S^1$ equivariant cohomology if a point.

The main result of this subsection is the following:

*Proposition* 25.1.1. There is an equivalence of categories

$$\operatorname{Alg}_{\mathbb{E}_{n+2}^{S^1}}(\operatorname{Perf}_{\mathbb{K}}) \xrightarrow{\cong} \operatorname{Alg}_{\mathbb{BD}_n^u}(\mathrm{D}_{\mathrm{fg}}^b(\mathbb{K}[u])) ,$$

intertwining the functor of taking homology $H_\bullet : \operatorname{Alg}_{\mathbb{E}_{n+2}^{S^1}}(\operatorname{Perf}_{\mathbb{K}}) \to \operatorname{Alg}_{\mathbb{P}_{n+2}}(\operatorname{Perf}_{\mathbb{K}})$ and the specialization to the central fibre $(\cdot)|_{\{0\}} : \operatorname{Alg}_{\mathbb{BD}_n^u}(\mathrm{D}_{\mathrm{fg}}^b(\mathbb{K}[u])) \to \operatorname{Alg}_{\mathbb{P}_{n+2}}(\operatorname{Perf}_{\mathbb{K}})$, the two periodic analogue of the functor of Proposition C.6.10.

*Remark* 25.1.2. A similar result was obtained by explicit calculation in [BBZB+20] in the case $n = 1$, and the analogous statement for general $n$ was announced there as to appear in [BZN].

*Proof.* Applying the equivariant Dunn-Lurie theorem [Lur12], recalled in Theorem 22.0.6, there is an equivalence

$$\operatorname{Alg}_{\mathbb{E}_{n+2}^{S^1}}(\operatorname{Perf}_{\mathbb{K}}) \xrightarrow{\cong} \operatorname{Alg}_{\mathbb{E}_n}(\operatorname{Alg}_{\mathbb{E}_2^{S^1}}(\operatorname{Perf}_{\mathbb{K}})) .$$

Further, the equivalence of Equation 23.0.4 induces an equivalence

$$\operatorname{Alg}_{\mathbb{E}_n}(\operatorname{Alg}_{\mathbb{E}_2^{S^1}}(\operatorname{Perf}_{\mathbb{K}})) \xrightarrow{\cong} \operatorname{Alg}_{\mathbb{E}_n}(\operatorname{Alg}_{\mathbb{BD}_0^u}(\mathrm{D}_{\mathrm{fg}}^b(\mathbb{K}[u]))) .$$

The strong Poisson additivity theorem of Rozenblyum gives the final desired equivalence

$$\operatorname{Alg}_{\mathbb{E}_n}(\operatorname{Alg}_{\mathbb{BD}_0^u}(\mathrm{D}_{\mathrm{fg}}^b(\mathbb{K}[u]))) \xrightarrow{\cong} \operatorname{Alg}_{\mathbb{BD}_n^u}(\mathrm{D}_{\mathrm{fg}}^b(\mathbb{K}[u])) .$$

$\square$

*Remark* 25.1.3. The $\mathbb{B}\mathbb{D}_n^u$ operad controls two-periodic graded quantizations of $\mathbb{P}_{n+2}$ algebras to $\mathbb{E}_n$ algebras, in the sense of Proposition C.6.10. Thus, the above gives a correspondence between $S^1$ equivariant structures on a $\mathbb{E}_{n+2}$ algebra and two-periodic graded quantizations of its homology $\mathbb{P}_{n+2}$ algebra.

*Remark* 25.1.4. Concretely, for $A \in \mathrm{Alg}_{\mathbb{E}_{d+2}^{S^1}}(\mathrm{Perf}_\mathbb{K})$ an algebra over the little $n$ disks operad equipped with an $S^1$ equivariant structure, we obtain $t(A) \in \mathrm{Alg}_{\mathbb{B}\mathbb{D}_n^u}(\mathrm{D}_{\mathrm{fg}}^b(\mathbb{K}[u]))$ such that the central fibre $t(A)|_{\{0\}} \cong H_\bullet(A) \in \mathrm{Alg}_{\mathbb{P}_{n+2}}(\mathrm{Vect})$ is equivalent to the homology $\mathbb{P}_{n+2}$ algebra of $A$. Thus, we can interpret the $S^1$ equivariance data on $A$ as defining a deformation $t(A)$ of $H_\bullet(A)$ to an $\mathbb{E}_n$ algebra $t(A)|_{\{1\}} \in \mathrm{Alg}_{\mathbb{E}_n}(\mathrm{Perf}_\mathbb{K})$.

*Example* 25.1.5. The special case $n = 1$ of the above gives an equivalence

$$\mathrm{Alg}_{\mathbb{E}_3^{S^1}}(\mathrm{Perf}_\mathbb{K}) \xrightarrow{\cong} \mathrm{Alg}_{\mathbb{B}\mathbb{D}_1^u}(\mathrm{D}_{\mathrm{fg}}^b(\mathbb{K}[u])) \ .$$

Heuristically, this result identifies $S^1$ equivariant structures on an $\mathbb{E}_3$ algebra with deformation quantizations of its homology $\mathbb{P}_3$ algebra to an $\mathbb{E}_1$ algebra. Subsections II-11.3 and II-13.5 explain examples of this phenomenon.

## 25.2. The equivariant cigar reduction principle for $\mathbb{E}_n$ algebras.

*Example* 25.2.1. Let $A \in \mathrm{Alg}_{\mathbb{E}_{n+2}}(\mathrm{Perf}_\mathbb{K})$ be an $\mathbb{E}_{n+2}$ algebra in $\mathrm{Perf}_\mathbb{K}$, and consider its image

$$A_0 = \mathrm{oblv}_{\mathbb{E}_{n+2}}^{\mathbb{E}_n}(A) \in \mathrm{Alg}_{\mathbb{E}_n}(\mathrm{Perf}_\mathbb{K}) \qquad \text{under} \qquad \mathrm{oblv}_{\mathbb{E}_{n+2}}^{\mathbb{E}_n} : \mathrm{Alg}_{\mathbb{E}_{n+2}}(\mathrm{Perf}_\mathbb{K}) \to \mathrm{Alg}_{\mathbb{E}_n}(\mathrm{Perf}_\mathbb{K})$$

the forgetful functor of Example C.4.12. Then $A_0$ is canonically a module over $A$ in the $\mathbb{E}_2$ sense, that is, the pair $(A, A_0)$ canonically define a $\mathrm{Disk}_{n \subset n+2}^{\mathrm{fr}}$-algebra in the sense of [AFT17]. Equivalently, by Proposition 4.8 of [AFT16], $A_0$ is canonically a module over the Hochschild chains algebra $\mathrm{CH}_\bullet(A)$ in the $\mathbb{E}_1$ sense, that is, there is a canonical map

(25.2.1)          $\mathrm{CH}_\bullet(A) \to \mathrm{CH}^\bullet(A_0)$      in the category      $\mathrm{Alg}_{\mathbb{E}_{n+1}}(\mathrm{Perf}_\mathbb{K}) \ .$

In terms of the Dunn additivity equivalence

$$A \in \mathrm{Alg}_{\mathbb{E}_{n+2}}(\mathrm{Perf}_\mathbb{K}) \cong \mathrm{Alg}_{\mathbb{E}_2}(\mathrm{Alg}_{\mathbb{E}_n}(\mathrm{Perf}_\mathbb{K})) \ ,$$

the map in Equation 25.2.1 encodes the fact that $A$ is canonically a module over itself (in the $\mathbb{E}_2$ sense) internal to $\mathrm{Alg}_{\mathbb{E}_n}(\mathrm{Perf}_\mathbb{K})$, and $A_0$ is the underlying object of this module. Equivalently, analogously identifying

$$\mathrm{CH}_\bullet(A) \in \mathrm{Alg}_{\mathbb{E}_{n+1}}(\mathrm{Perf}_\mathbb{K}) \cong \mathrm{Alg}_{\mathbb{E}_1}(\mathrm{Alg}_{\mathbb{E}_n}(\mathrm{Perf}_\mathbb{K})) \ ,$$

the object $A_0$ admits a canonical module structure

$$A_0 \in \mathrm{CH}_\bullet(A)\text{-}\mathrm{Mod}(\mathrm{Alg}_{\mathbb{E}_n}(\mathrm{Perf}_\mathbb{K})) \ .$$

Further, we have:

*Example* 25.2.2. The negative cyclic chains defines a canonical deformation

$$\mathrm{CC}_\bullet^-(A) \in \mathrm{Alg}_{\mathbb{E}_{n+1}}(\mathrm{D}_{\mathrm{fg}}^b(\mathbb{K}[u]))$$

with central fibre

$$\mathrm{CC}_\bullet^-(A)|_{\{0\}} = \mathrm{CH}_\bullet(A) \ \in \mathrm{Alg}_{\mathbb{E}_{n+1}}(\mathrm{Perf}_\mathbb{K})$$

the Hochschild chains algebra.

Now, we let $A \in \mathrm{Alg}_{\mathbb{E}_{n+2}}(\mathrm{Perf}_{\mathbb{K}})$ and $A_0 = \mathrm{oblv}_{\mathbb{E}_{n+2}}^{\mathbb{E}_n}(A) \in \mathrm{Alg}_{\mathbb{E}_n}(\mathrm{Perf}_{\mathbb{K}})$ be as in Example 25.2.1 above, and state the main result of this subsection:

*Proposition* 25.2.3. An $S^1$ equivariant structure on $A$ in the $\mathbb{E}_2$ direction, that is, a lift to

$$A \in \mathrm{Alg}_{\mathbb{E}_2^{S^1}}(\mathrm{Alg}_{\mathbb{E}_n}(\mathrm{Perf}_{\mathbb{K}})) \ ,$$

is equivalent to a deformation

$$A_u \in \mathrm{CC}_\bullet^-(A)\text{-}\mathrm{Mod}(\mathrm{Alg}_{\mathbb{E}_n}(\mathrm{D}_{\mathrm{fg}}^b(\mathbb{K}[u])))$$

such that the central fibre

$$A_u|_{\{0\}} = A_0 \ \in \mathrm{CH}_\bullet(A)\text{-}\mathrm{Mod}(\mathrm{Alg}_{\mathbb{E}_n}(\mathrm{Perf}_{\mathbb{K}}))$$

is equivalent to $A_0$ equipped with the $\mathrm{CH}_\bullet(A)$ module structure of Example 25.2.1 above.

*Remark* 25.2.4. The results of Subsection II-11.4 provide an example of the above phenomenon.

25.3. **Quantization of factorization $\mathbb{E}_n$ algebras.** We now give the analogue of the above discussion for factorization $\mathbb{E}_n$ algebras on $X$, in the sense defined in Section II-7, which describes quantization in the $\Omega$-background for mixed holomorphic-topological field theories:

*Proposition* 25.3.1. There is a natural equivalence of categories

$$\mathrm{Alg}_{\mathbb{E}_{n+2}^{S^1}, \mathrm{un}}^{\mathrm{fact}}(X) \xrightarrow{\cong} \mathrm{Alg}_{\mathbb{BD}_n^u, \mathrm{un}}^{\mathrm{fact}}(X) \ ,$$

intertwining the functor of taking homology $H_\bullet : \mathrm{Alg}_{\mathbb{E}_{n+2}^{S^1}, \mathrm{un}}^{\mathrm{fact}}(X) \to \mathrm{Alg}_{\mathbb{P}_{n+2}, \mathrm{un}}^{\mathrm{fact}}(X)$ and the specialization to the central fibre $(\cdot)|_{\{0\}} : \mathrm{Alg}_{\mathbb{BD}_n^u, \mathrm{un}}^{\mathrm{fact}}(X) \to \mathrm{Alg}_{\mathbb{P}_{n+2}, \mathrm{un}}^{\mathrm{fact}}(X)$.

*Proof.* □

*Remark* 25.3.2. Concretely, for $\mathcal{A} \in \mathrm{Alg}_{\mathbb{E}_{n+2}^{S^1}, \mathrm{un}}^{\mathrm{fact}}(X)$ a factorization $\mathbb{E}_{n+2}$ algebra equipped with an $S^1$ equivariant structure, we obtain $t(\mathcal{A}) \in \mathrm{Alg}_{\mathbb{BD}_n^u, \mathrm{un}}^{\mathrm{fact}}(X)$, so that the central fibre $t(\mathcal{A})|_{\{0\}} \in \mathrm{Alg}_{\mathbb{P}_{n+2}, \mathrm{un}}^{\mathrm{fact}}(X)$ is a factorization $\mathbb{P}_{n+2}$ algebra, which is identified with an $(n+1)$-shifted Coisson algebra by the chiral Poisson additivity theorem of Rozenblyum. Thus, we can interpret the $S^1$ equivariance data on $\mathcal{A}$ as defining a deformation $t(\mathcal{A})$ of this shifted Coisson algebra to a factorization $\mathbb{E}_n$ algebra $t(\mathcal{A})|_{\{1\}} \in \mathrm{Alg}_{\mathbb{E}_n, \mathrm{un}}^{\mathrm{fact}}(X)$.

*Example* 25.3.3. The special case $n = 0$ of the above gives an equivalence

$$\mathrm{Alg}_{\mathbb{E}_2^{S^1}, \mathrm{un}}^{\mathrm{fact}}(X) \xrightarrow{\cong} \mathrm{Alg}_{\mathbb{BD}_0^u, \mathrm{un}}^{\mathrm{fact}}(X) \ .$$

The strong chiral Poisson additivity theorem of Rozenblyum identifies the latter category with that of (two-periodic) filtered quantizations of factorization algebras, in the sense of Definition 10.4.2. Thus, this result interprets $S^1$ equivariant structures on a factorization $\mathbb{E}_2$ algebra as two-periodic graded quantizations of the corresponding shifted Coisson algebra. The results of subsection II-19.3 are an example of this phenomenon.

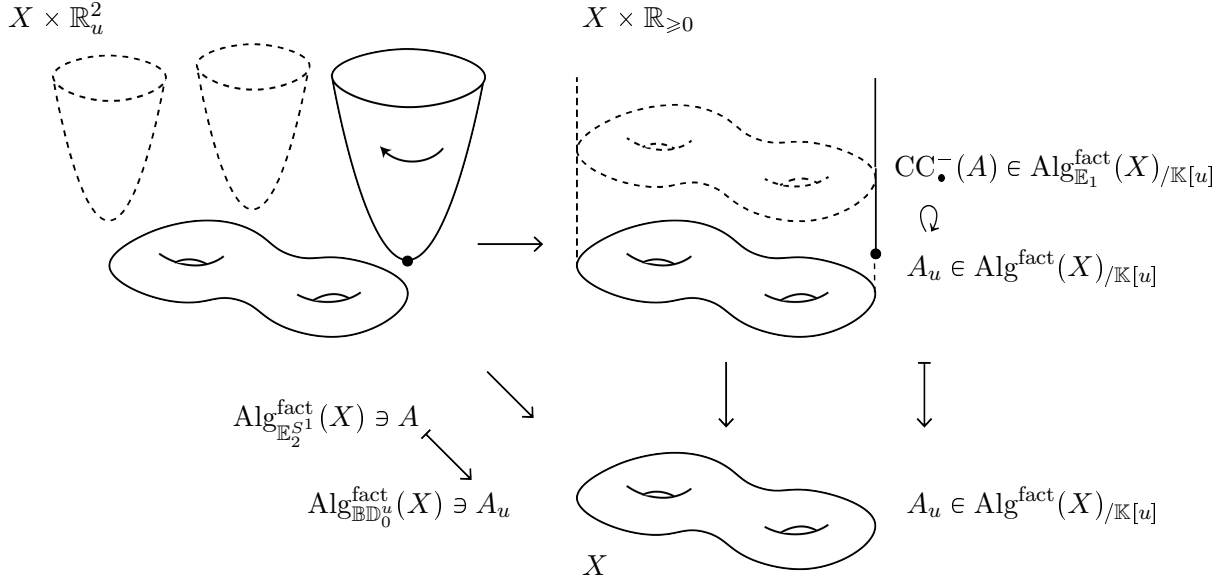

FIGURE 4. The equivariant cigar reduction principle pictured in the case $n = 0$

25.4. **The equivariant cigar reduction principle for factorization $\mathbb{E}_n$ algebras.** We now state the analogues of the results of Subsection 25.2 above for factorization $\mathbb{E}_n$ algebras. The main result, Proposition 25.4.3 below, is illustrated in Figure 4.

As in Example 25.2.1, we have:

*Example* 25.4.1. Let $\mathcal{A} \in \mathrm{Alg}^{\mathrm{fact}}_{\mathbb{E}_{n+2},\mathrm{un}}(X)$ be a factorization $\mathbb{E}_{n+2}$ algebra on $X$, and

$$\mathcal{A}_0 = \mathrm{oblv}^{\mathbb{E}_n}_{\mathbb{E}_{n+2}} \mathcal{A} \in \mathrm{Alg}^{\mathrm{fact}}_{\mathbb{E}_n,\mathrm{un}}(X) \ .$$

Then by the equivalence induced by Dunn additivity together with Corollary II-7.1.8, we have

$$\mathcal{A} \in \mathrm{Alg}^{\mathrm{fact}}_{\mathbb{E}_{n+2},\mathrm{un}}(X) \cong \mathrm{Alg}_{\mathbb{E}_2}(\mathrm{Alg}^{\mathrm{fact}}_{\mathbb{E}_n,\mathrm{un}}(X)) \ ,$$

and analogously for the Hochschild chains algebra

$$\mathrm{CH}_\bullet(\mathcal{A}) \in \mathrm{Alg}^{\mathrm{fact}}_{\mathbb{E}_{n+1},\mathrm{un}}(X) \cong \mathrm{Alg}_{\mathbb{E}_1}(\mathrm{Alg}^{\mathrm{fact}}_{\mathbb{E}_n,\mathrm{un}}(X)) \ .$$

Moreover, $\mathcal{A}_0$ admits a canonical module structure

$$\mathcal{A}_0 \in \mathrm{CH}_\bullet(\mathcal{A})\text{-Mod}(\mathrm{Alg}^{\mathrm{fact}}_{\mathbb{E}_n,\mathrm{un}}(X)) \ .$$

Further, as in Example 25.2.2, we have:

*Example* 25.4.2. The negative cyclic chains define a canonical deformation

$$\mathrm{CC}^-_\bullet(\mathcal{A}) \in \mathrm{Alg}^{\mathrm{fact}}_{\mathbb{E}_{n+1},\mathrm{un}}(X)_{/(\mathbb{A}^1/\mathbb{G}_m)} \cong \mathrm{Alg}_{\mathbb{E}_1}(\mathrm{Alg}^{\mathrm{fact}}_{\mathbb{E}_n,\mathrm{un}}(X)_{/(\mathbb{A}^1/\mathbb{G}_m)})$$

with central fibre given by the Hochschild chains algebra

$$\mathrm{CC}^-_\bullet(\mathcal{A})|_{\{0\}} = \mathrm{CH}_\bullet(\mathcal{A}) \ \in \mathrm{Alg}^{\mathrm{fact}}_{\mathbb{E}_{n+1},\mathrm{un}}(X) \ .$$

Now, we let $\mathcal{A} \in \mathrm{Alg}^{\mathrm{fact}}_{\mathbb{E}_{n+2},\mathrm{un}}(X)$ and $\mathcal{A}_0 = \mathrm{oblv}^{\mathbb{E}_n}_{\mathbb{E}_{n+2}} \mathcal{A} \in \mathrm{Alg}^{\mathrm{fact}}_{\mathbb{E}_n,\mathrm{un}}(X)$ be as in Example 25.4.1 above, and state the main result of this subsection:

*Proposition* 25.4.3. An $S^1$ equivariant structure on $\mathcal{A}$ in the $\mathbb{E}_2$ direction, that is, a lift to

$$\mathcal{A} \in \mathrm{Alg}^{\mathrm{fact}}_{\mathbb{E}^{S^1}_{n+2}, \mathrm{un}}(X) \ ,$$

is equivalent to a deformation

$$\mathcal{A}_u \in \mathrm{CC}^-_\bullet(\mathcal{A})\text{-}\mathrm{Mod}(\mathrm{Alg}^{\mathrm{fact}}_{\mathbb{E}_n, \mathrm{un}}(X)_{/(\mathbb{A}^1/\mathbb{G}_m)})$$

such that the central fibre

$$\mathcal{A}_u|_{\{0\}} = \mathcal{A}_0 \ \in \mathrm{CH}_\bullet(\mathcal{A})\text{-}\mathrm{Mod}(\mathrm{Alg}^{\mathrm{fact}}_{\mathbb{E}_n, \mathrm{un}}(X))$$

is equivalent to $\mathcal{A}_0$ equipped with the $\mathrm{CH}_\bullet(\mathcal{A})$ module structure of Example 25.4.1 above.

*Remark* 25.4.4. The results of Subsection II-19.3 provide an example of the above phenomenon.

# Chapter 3
# Appendices

## Appendix A. Sheaf theory

Let $X$ be a smooth variety of dimension $d_X$ over $\mathbb{K} = \mathbb{C}$ or a field of characteristic 0. We write $\mathcal{O}_X$ for the sheaf of regular functions, $\mathcal{D}_X$ for the sheaf of differential operators, $\Theta_X$ for the tangent sheaf, $\Omega_X^1$ for the sheaf of Kähler differentials, $\Omega_X^{d_X}$ for the sheaf of sections of the canonical bundle, and $\omega_X = \Omega_X^{d_X}[d_X]$ for the dualizing sheaf on $X$. Let $\mathrm{Sh}_z(X)$ denote the category of complexes of sheaves of $\mathbb{K}$-modules on $X$ in the Zariski topology.

*Warning* A.0.1. Note that we assume $X$ is a smooth variety throughout, and only define the category of $D$ modules on more general spaces in Supappendix II-A.5.

A.1. $\mathcal{O}$-**module conventions.** Let $\mathrm{D}(\mathcal{O}_X)$ be the DG category of complexes of $\mathcal{O}_X$-modules, $\mathrm{QCoh}(X)$ and $\mathrm{Coh}(X)$ be the full subcategories of complexes with quasi-coherent and coherent cohomology sheaves, and $\mathrm{Perf}(X)$ the subcategory of bounded complexes with finitely generated cohomology sheaves. The category $\mathrm{D}(\mathcal{O}_X)$ is symmetric monoidal with respect to the tensor product $\otimes_{\mathcal{O}_X}$, with unit object $\mathcal{O}_X$, and $\mathrm{QCoh}(X)$, $\mathrm{Coh}(X)$, and $\mathrm{Perf}(X)$ are monoidal subcategories.

*Definition* A.1.1. Let $f : X \to Y$ a map of schemes. The inverse and direct image functors are

$$f^{\bullet} : \mathrm{D}(\mathcal{O}_Y) \to \mathrm{D}(\mathcal{O}_X) \qquad f^{\bullet}\mathcal{F} = f^{-1}F \otimes_{f^{-1}\mathcal{O}_Y} \mathcal{O}_X \qquad \text{and} \qquad f_{\bullet} : \mathrm{D}(\mathcal{O}_X) \to \mathrm{D}(\mathcal{O}_Y) \qquad f_{\bullet}\mathcal{F} = f_{\bullet}\mathcal{F},$$

where $f_{\bullet} : \mathrm{Sh}_z(X) \to \mathrm{Sh}_z(Y)$ and $f^{-1} : \mathrm{Sh}_z(Y) \to \mathrm{Sh}_z(X)$ are the usual direct and inverse image functors on sheaves of $\mathbb{K}$-modules.

*Remark* A.1.2. Note that $f^{\bullet}$ preserves quasicoherence, as does $f_{\bullet}$ for quasicompact, quasiseperated maps. We define the global sections functor by $\Gamma = \pi_{\bullet} : \mathrm{D}(\mathcal{O}_X) \to \mathrm{Vect}$ where $\pi : X \to \mathrm{pt}$.

*Definition* A.1.3. Let $\mathcal{F}, \mathcal{G} \in \mathrm{D}(\mathcal{O}_X)$. The internal hom object in $\mathrm{D}(\mathcal{O}_X)$ is

$$\underline{\mathrm{Hom}}_{\mathcal{O}_X}(\mathcal{F}, \mathcal{G}) \in \mathrm{D}(\mathcal{O}_X) \qquad \text{by} \qquad \underline{\mathrm{Hom}}_{\mathcal{O}_X}(\mathcal{F}, \mathcal{G})(U) := \mathrm{Hom}_{\mathcal{O}_X|_U}(\mathcal{F}|_U, \mathcal{G}|_U) .$$

*Remark* A.1.4. For $\mathcal{H} \in \mathrm{D}(\mathcal{O}_X)$, we have

$$\mathrm{Hom}(\mathcal{H}, \underline{\mathrm{Hom}}_{\mathcal{O}_X}(\mathcal{F}, \mathcal{G})) \cong \mathrm{Hom}(\mathcal{H} \otimes_{\mathcal{O}_X} \mathcal{F}, \mathcal{G}) .$$

In particular, the space of homomorphisms is given by the space of sections of the internal hom object

$$\mathrm{Hom}(\mathcal{F}, \mathcal{G}) = \mathrm{Hom}(\mathcal{O}_X, \underline{\mathrm{Hom}}_{\mathcal{O}_X}(\mathcal{F}, \mathcal{G})) = \Gamma(X, \underline{\mathrm{Hom}}_{\mathcal{O}_X}(\mathcal{F}, \mathcal{G})).$$

*Remark* A.1.5. For $\mathcal{F} \in \mathrm{Coh}(X)$ coherent and $\mathcal{G} \in \mathrm{QCoh}(X)$ quasi-coherent, the object $\underline{\mathrm{Hom}}_{\mathcal{O}_X}(\mathcal{F}, \mathcal{G}) \in \mathrm{QCoh}(X)$ is quasi-coherent. If $\mathcal{F}, \mathcal{G} \in \mathrm{Coh}(X)$ are both coherent, then $\underline{\mathrm{Hom}}_{\mathcal{O}_X}(\mathcal{F}, \mathcal{G}) \in \mathrm{Coh}(X)$ is also coherent.

*Definition* A.1.6. The duality functor on coherent $\mathcal{O}_X$-modules is defined by

$$(-)^{\vee} : \mathrm{Coh}(X) \to \mathrm{Coh}(X) \qquad \text{by} \qquad \mathcal{F} \mapsto \mathcal{F}^{\vee} := \underline{\mathrm{Hom}}_{\mathcal{O}_X}(\mathcal{F}, \mathcal{O}_X).$$

*Remark* A.1.7. There are canonical isomorphisms $\underline{\mathrm{Hom}}_{\mathcal{O}_X}(\mathcal{F}, \mathcal{G}) \cong \mathcal{G} \otimes_{\mathcal{O}_X} \mathcal{F}^{\vee}$ and $(\mathcal{F}^{\vee})^{\vee} \cong \mathcal{F}$.

### A.2. $D$-module conventions.

Let $D^l(X)$ and $D^r(X)$ be the concrete DG categories of complexes of left and right $\mathcal{D}_X$-modules which are quasicoherent as $\mathcal{O}_X$-modules, and let $D^l_{\mathrm{c}}(X)$ and $D^r_{\mathrm{c}}(X)$ denote the full sub DG categories of complexes with cohomology that is coherent as a module over $\mathcal{D}_X$.

*Example* A.2.1. The sheaf of regular functions $\mathcal{O}_X \in D^l(X)^\heartsuit$ has the structure of a left $D$ module, given by the defining action of the sheaf of differential operators $\mathcal{D}_X$ on $\mathcal{O}_X$.

More generally, a left $D$ module (or a complex of such) $M \in D^l(X)$ on $X$ is given by a quasicoherent sheaf (or a complex of such) $M \in \mathrm{QCoh}(X)$, together with a flat connection, that is, $\nabla \in \mathrm{Hom}_{\mathrm{Sh}_{\mathbb{K}}(X)}(M, \Omega^1_X \otimes_{\mathcal{O}_X} M)$ such that

- $\nabla_\theta(fs) = \theta(f)s + f\nabla_\theta(s)$ , and
- $\nabla_{[\theta_1,\theta_2]}s = [\nabla_{\theta_1}, \nabla_{\theta_2}]s$ ,

where $\theta, \theta_1, \theta_2 \in \Theta_X$, $f \in \mathcal{O}_X$, and $s \in M$. The first condition is that $\nabla$ defines a connection, and the second that $\nabla$ is flat.

*Example* A.2.2. The sheaf of sections of the canonical bundle $\Omega^{d_X}_X \in D^r(X)^\heartsuit$ is the protypical example of a right $D_X$ module, with action of vector fields given by $\theta(\eta) = -\mathrm{Lie}_\theta(\eta)$ for $\theta \in \Theta_X$ and $\eta \in \Omega^{d_X}_X$.

*Remark* A.2.3. There is a canonical equivalence of the categories $D^l(X)$ and $D^r(X)$

$$D^l(X) \underset{(-)^l}{\overset{(-)^r}{\rightleftarrows}} D^r(X) \qquad \text{defined by} \qquad \begin{cases} M \mapsto M^l := M \otimes_{\mathcal{O}_X} \omega^\vee_X & \text{for } M \in D^r(X) \quad \text{and} \\ L \mapsto L^r := \omega_X \otimes_{\mathcal{O}_X} L & \text{for } L \in D^l(X). \end{cases}$$

We write $D(X)$ for the abstract DG category given by the common value of $D^r(X)$ and $D^l(X)$ under this identification, and $D_{\mathrm{c}}(X)$ for the full sub DG category corresponding to $D^r_{\mathrm{c}}(X)$ and $D^l_{\mathrm{c}}(X)$, which are also identified under this equivalence. $D^r(X)$ and $D^l(X)$ both have natural forgetful functors to $\mathrm{QCoh}(X)$, which are intertwined by tensoring with $\omega_X$. This perspective is summarized in the following diagram:

$$\begin{array}{ccc} D^l(X) & \xrightarrow{\ \omega_X\ }_{\simeq} & D^r(X) \\ \downarrow{\scriptstyle o^l} & & \downarrow{\scriptstyle o^r} \\ \mathrm{QCoh}^l(X) & \xrightarrow{\ \omega_X\ }_{\simeq} & \mathrm{QCoh}^r(X) \end{array} \qquad \text{so that} \qquad \begin{array}{c} D(X) \\ {\scriptstyle o^l}\downarrow \quad \searrow{\scriptstyle o^r} \\ \mathrm{QCoh}(X) \xrightarrow[\simeq]{\ \omega_X\ } \mathrm{QCoh}(X) \end{array} \quad ,$$

where $\mathrm{QCoh}^l(X)$ and $\mathrm{QCoh}^r(X)$ are just the category $\mathrm{QCoh}(X)$

*Remark* A.2.4. Throughout, when defining a functor involving (potentially several copies of) the category $D(X)$, we will prescribe the values of the functor in terms of a particular choice of realization $D^r(X)$ or $D^l(X)$ for each copy of $D(X)$, with the extension to all other choices of concrete realizations of $D(X)$ implicitly specified via the above equivalence.

*Remark* A.2.5. Note that the above equivalence is exact up to a cohomological degree shift of $d_X = \dim_{\mathbb{K}} X$, so that the category $D(X)$ inherits two different t-structures, which differ only by this shift. We choose to preference the right t structure, and all statements about exactness of functors involving $D(X)$ will be given in these terms. This t-structure will be the one which corresponds to the perverse t-structure on constructible sheaves under the Riemann-Hilbert correspondence. In particular, under this identification $\omega_X \in D(X)$ is the dualizing sheaf, $\omega_X[-d_X] \in D(X)^\heartsuit$ is the IC sheaf, and $\underline{\mathbb{K}}_X := \omega_X[-2d_X] \in D(X)$ is the constant sheaf.

*Definition* A.2.6. The $\otimes^!$ monoidal structure on $D(X)$ is $\otimes^! : D(X)^{\otimes 2} \to D(X)$ defined by

$$\otimes^! : D^l(X) \times D^l(X) \to D^l(X) \qquad M \otimes^! N = M \otimes_{\mathcal{O}_X} N \quad \text{with} \quad P(m \otimes n) = Pm \otimes n + m \otimes Pn \ ,$$

for $P \in D_X$.

*Remark* A.2.7. This formula agrees with the usual definition of the tensor product of connections, and tensor products of flat connections are flat. The corresponding functor $\otimes^! : D^r(X) \otimes D^r(X) \to D^r(X)$ is given by $M \otimes^! N = M \otimes_{\mathcal{O}_X} N \otimes_{\mathcal{O}_X} \omega_X^\vee$.

Let $\mathbb{1} \in D(X)$ denote the tensor unit, and note $o^l(\mathbb{1}) = \mathcal{O}_X$ and $o^r(\mathbb{1}) = \omega_X$. We will often use just $\otimes$ to denote this symmetric monoidal structure on $D(X)$.

*Definition* A.2.8. Let $f : X \to Y$ be a map of smooth varieties. The inverse image functor $f^! : D(Y) \to D(X)$ is defined by

$$f^! : D^l(Y) \to D^l(X) \qquad f^!(M) = f^\bullet(M) \quad \text{equipped with the pullback flat connection.}$$

*Remark* A.2.9. This functor is symmetric monoidal with respect to $\otimes^!$, and in particular maps the tensor unit $\mathbb{1}_Y$ to $\mathbb{1}_X$.

*Remark* A.2.10. The corresponding functor $f^! : D^r(Y) \to D^r(X)$ is given by

$$f^!(M) = f^\bullet(M \otimes_{\mathcal{O}_Y} \omega_Y^\vee) \otimes_{\mathcal{O}_X} \omega_X \cong f^\bullet M \otimes_{\mathcal{O}_X} \omega_{X/Y} \ .$$

*Definition* A.2.11. The exterior product is defined by

$$\boxtimes : D(X) \times D(Y) \to D(X \times Y) \qquad \text{by} \qquad M \boxtimes N = \pi_X^! M \otimes \pi_Y^! N \ ,$$

for $\pi_X : X \times Y \to X, \pi_Y : X \times Y \to Y$.

*Remark* A.2.12. Note that

$$M \otimes N = \Delta^!(M \boxtimes N)$$

for $M, N \in D(X)$ and $\Delta : X \to X \times X$ the diagonal embedding.

*Definition* A.2.13. Let $f : X \to Y$ again be a map of smooth varieties. The direct image functor is

$$f_* : D^r(X) \to D^r(Y) \qquad f_*(M) = f_\bullet(M \otimes_{\mathcal{D}_X} \mathcal{D}_{X \to Y}) \qquad \text{for} \qquad \mathcal{D}_{X \to Y} := f^! \mathcal{D}_Y \in (\mathcal{D}_X, f^{-1}\mathcal{D}_Y)\text{-Mod}$$

where $\mathcal{D}_{X \to Y} = f^! \mathcal{D}_Y \in D^l(X)$ is defined in terms of $\mathcal{D}_Y \in D^l(Y)$ as a left module, so that the additional $(\mathcal{D}_Y, \mathcal{D}_Y)$-bimodule structure on $\mathcal{D}_Y$ equips $\mathcal{D}_{X \to Y}$ with the structure of a $(\mathcal{D}_X, f^{-1}\mathcal{D}_Y)$-bimodule.

*Definition* A.2.14. The de Rham cochains functor is $C_{\mathrm{dR}}^\bullet := \pi_* : D(X) \to \mathrm{Vect}_{\mathbb{K}}$, where $\pi : X \to \mathrm{pt}$. The de Rham chains functor is $C_\bullet^{\mathrm{dR}} := \pi_! : D(X) \to \mathrm{Vect}_{\mathbb{K}}$.

*Remark* A.2.15. Note that the de Rham cochain and chains functors are calculated as

$$C_{\mathrm{dR}}^\bullet : D^r(X) \to \mathrm{Vect} \qquad\qquad C_{\mathrm{dR}}^\bullet(X; M) = \pi_\bullet(M \otimes_{D_X} \mathcal{O}_X)$$

$$C_{\mathrm{dR}}^\bullet : D^l(X) \to \mathrm{Vect} \qquad\qquad C_{\mathrm{dR}}^\bullet(X; M) = \pi_\bullet(\omega_X \otimes_{D_X} M) \ .$$

*Definition* A.2.16. The sheaf internal hom functor

$$\underline{\mathrm{Hom}}_{D(X)}(\cdot, \cdot) : D(X)^{\mathrm{op}} \times D(X) \to \mathrm{Sh}_z(X) \qquad \text{by} \qquad \underline{\mathrm{Hom}}_{D(X)}(M, N)(U) = \mathrm{Hom}_{D(U)}(j^! M, j^! N) \ ,$$

for $M, N \in D(X)$, where $j : U \to X$ is the open embedding.

*Remark* A.2.17. Note that

$$\Gamma \circ \underline{\mathrm{Hom}} = \mathrm{Hom} : D(X)^{\mathrm{op}} \times D(X) \to \mathrm{Vect} \qquad\qquad \Gamma(X, \underline{\mathrm{Hom}}_{D(X)}(M, N)) = \mathrm{Hom}_{D(X)}(M, N) \ .$$

*Definition* A.2.18. The duality functor $\mathbb{D} : D_c(X)^{\mathrm{op}} \to D_c(X)$ is defined by

$$\mathbb{D} : D_c^r(X)^{\mathrm{op}} \to D_c^l(X) \qquad\qquad \mathbb{D}(M) = \underline{\mathrm{Hom}}_{D^r(X)}(M, \mathcal{D}_X) \ ,$$

where $\mathcal{D}_X \in D^r(X)$ is considered as a $(\mathcal{D}_X, \mathcal{D}_X)$-bimodule so that $\mathbb{D}(M)$, which is a priori an object in $\mathrm{Sh}_z(X)$, defines an object of $D^l(X)$ as desired.

*Remark* A.2.19. Note that $\mathbb{D}$ preserves coherence, but if $M$ is not coherent, then the resulting object of $\mathcal{D}_X$-Mod is not in general quasicoherent as an object of $\mathrm{D}(\mathcal{O}_X)$.

*Definition* A.2.20. The genuine internal hom functor $\mathcal{H}om(\cdot, \cdot) : D_c(X)^{\mathrm{op}} \otimes D(X) \to D(X)$ is defined by

$$\mathcal{H}om(\cdot, \cdot) : D_c^r(X)^{\mathrm{op}} \times D^l(X) \to D^l(X) \qquad\qquad \mathcal{H}om(M, N) = \underline{\mathrm{Hom}}_{D^r(X)}(M, N \otimes_{\mathcal{O}_X} \mathcal{D}_X) \ ,$$

where $N \otimes_{\mathcal{O}_X} \mathcal{D}_X \in D^r(X)$ is considered as a $(\mathcal{D}_X, \mathcal{D}_X)$-bimodule so that $\mathcal{H}om(M, N) \in D^l(X)$ as above.

*Remark* A.2.21. Note that

$$C_{\mathrm{dR}}^\bullet \circ \mathcal{H}om = \mathrm{Hom} : D_c(X)^{\mathrm{op}} \times D(X) \to \mathrm{Vect} \qquad\qquad C_{\mathrm{dR}}^\bullet(X, \mathcal{H}om_{D(X)}(M, N)) = \mathrm{Hom}_{D(X)}(M, N) \ .$$

Further, we have

$$\mathcal{H}om(\cdot, \cdot) = \mathbb{D}(\cdot) \otimes^! (\cdot) : D_c(X)^{\mathrm{op}} \times D(X) \to D(X) \ ,$$

and in particular $\mathcal{H}om(\cdot, \mathbb{1}) = \mathbb{D} : D(X)^{\mathrm{op}} \to D(X)$; we could equivalently take this as the definition of $\mathcal{H}om$.

*Remark* A.2.22. The pushforward and pullback functors $f_*$ and $f^!$ above were defined on the entire category $D(X)$, but their putative adjoints can not always be defined. In general, the best we can do is the following: Let $f : X \to Y$ again be a map of smooth varieties, and let $D_c^{f^!}(Y)$ be the full subcategory of objects $M \in D_c(Y)$ such that $f^! \mathbb{D} M \in D_c(X)$ is coherent, and similarly $D_c^{f_*}(X)$ be the full subcategory of objects $M \in D_c(X)$ such that $f_* \mathbb{D} M \in D_c(Y)$ is coherent. Then we define

$$f^* := \mathbb{D} f^! \mathbb{D} : D_c^{f^!}(Y) \to D_c(X) \qquad\qquad f_! := \mathbb{D} f_* \mathbb{D} : D_c^{f_*}(X) \to D_c(Y) \ .$$

In various situations, these definitions simplify to more useful ones, as in the following propositions.

*Proposition* A.2.23. Let $f : X \to Y$ be a smooth map of relative dimension $d = d_X - d_Y$ of smooth varieties. Then $f^! : D_c(Y) \to D_c(X)$ preserves coherence, so that $f^* : D_c(Y) \to D_c(X)$ is defined. Moreover, in this case $f^* = f^![-2d]$, and we have a natural isomorphism

$$\mathrm{Hom}_{D(X)}(f^* M, N) \cong \mathrm{Hom}_{D(Y)}(M, f_* N)$$

of functors $D_c(X) \times D(Y) \to \mathrm{Vect}$.

*Proposition* A.2.24. Let $f : X \to Y$ be a proper map of smooth varieties. Then $f_* : D_c(X) \to D_c(Y)$ preserves coherence, so that $f_! : D_c(X) \to D_c(Y)$ is defined. Moreover, in this case $f_! = f_*$ and we have a natural isomorphism

$$\mathrm{Hom}_{D(Y)}(f_! M, N) \cong \mathrm{Hom}_{D(X)}(M, f^! N)$$

of functors $D_c(X) \times D(Y) \to \mathrm{Vect}$.

A.3. **The six functors formalism.** Let $D_{\mathrm{rh}}(X)$ the full subcategory of $D(X)$ on bounded complexes with regular holonomic cohomology modules.

*Theorem* A.3.1. There are functors

$$\otimes^! : D_{\mathrm{rh}}(X) \times D_{\mathrm{rh}}(X) \to D_{\mathrm{rh}}(X) \qquad \mathbb{D} : D_{\mathrm{rh}}(X)^{\mathrm{op}} \xrightarrow{\cong} D_{\mathrm{rh}}(X) \qquad \mathcal{H}\mathrm{om}_X : D_{\mathrm{rh}}(X)^{\mathrm{op}} \times D_{\mathrm{rh}}(X) \to D_{\mathrm{rh}}(X) \,,$$

and for $f : X \to Y$ natural adjunctions,

$$f^* : D_{\mathrm{rh}}(Y) \rightleftharpoons D_{\mathrm{rh}}(X) : f_* \qquad f_! : D_{\mathrm{rh}}(X) \rightleftharpoons D_{\mathrm{rh}}(Y) : f^! \ .$$

Moreover, these satisfy:

- for $f : X \to Y$ smooth of relative dimension $d$, $f^* = f^![-2d]$ as in A.2.23 above;
- for $f : X \to Y$ proper, $f_* = f_!$ as in A.2.24 above;
- for $f : X \to Y$, there are natural equivalences $\mathbb{D}_Y f_* \cong f_! \mathbb{D}_Y$, $\mathbb{D}_X f^* \cong f^! \mathbb{D}_Y$;
- $\otimes^!$ defines a symmetric monoidal structure on $D_{\mathrm{rh}}(X)$; and
- for $f : X \to Y$, there are natural equivalences

$$f_!(M \otimes f^* N) \cong f_!(M) \otimes N \qquad \mathcal{H}\mathrm{om}_Y(f_! M, N) \cong f_* \mathcal{H}\mathrm{om}_X(M, f^! N) \qquad f^! \mathcal{H}\mathrm{om}_Y(M, N) \cong \mathcal{H}\mathrm{om}_X(f^* M, f^! N) \,.$$

- For a Cartesian square

(A.3.1)
$$\begin{array}{ccc} Z & \xrightarrow{\tilde{g}} & X \\ {\scriptstyle \tilde{f}}\downarrow & & \downarrow{\scriptstyle f} \\ Y & \xrightarrow{g} & W \end{array} \qquad \text{there is a natural isomorphism} \qquad \tilde{f}_* \circ \tilde{g}^! \cong g^! \circ f_* \ .$$

*Definition* A.3.2. The $\otimes^*$ tensor structure $\otimes^* : D_{\mathrm{rh}}(X)^{\times 2} \to D_{\mathrm{rh}}(X)$ is defined by $(M, N) \mapsto \Delta^*(M \boxtimes N)$.

*Definition* A.3.3. The constant $D$ module on $X$ is $\underline{\mathbb{K}}_X = \pi^* \underline{\mathbb{K}}_{\mathrm{pt}} \in D_{\mathrm{rh}}(X)$, where $\pi : X \to \mathrm{pt}$ and $\underline{\mathbb{K}}_{\mathrm{pt}} \in D(\mathrm{pt})$ is the object corresponding to $\mathbb{K} \in \mathrm{Vect}_{\mathbb{K}} \cong D(\mathrm{pt})$.

*Example* A.3.4. If $X$ is smooth of dimension $d_X$, then $\underline{\mathbb{K}}_X \cong \omega_X[-2d_X]$ is a shift of the dualizing sheaf, as in Remark A.2.5.

*Definition* A.3.5. The de-Rham (Borel-Moore) (co)chains on $X$ are

$$C_\bullet(X) = \pi_! \pi^! \underline{\mathbb{K}}_{\mathrm{pt}} \qquad C_\bullet^{\mathrm{BM}}(X) = \pi_* \pi^! \underline{\mathbb{K}}_{\mathrm{pt}} \qquad C^\bullet(X) = \pi_* \pi^* \underline{\mathbb{K}}_{\mathrm{pt}} \qquad C_{\mathrm{c}}^\bullet(X) = \pi_! \pi^* \underline{\mathbb{K}}_{\mathrm{pt}} \qquad \in \quad D(\mathrm{pt})$$

where $\pi : X \to \mathrm{pt}$ and $\underline{\mathbb{K}}_{\mathrm{pt}} \in D(\mathrm{pt})$ are as in the preceding definition.

*Remark* A.3.6. Note that we use cochain complexes throughout, and do not use the convention of reversing the grading on homology. Thus, classes in homology which are correspond to higher dimensional cycles geometrically contribute to the homology groups of lower (cohomological) degree.

*Example* A.3.7. For $X$ a smooth variety of dimension $d_X$, we have an isomorphism $C_\bullet^{\mathrm{BM}}(X) \cong C_{\mathrm{dR}}^\bullet(X)[2d_X]$.

*Remark* A.3.8. The unit and counit of the above adjunctions for $f : X \to Y$ give canonical maps

$$A \to f_* f^* A \qquad \text{and} \qquad f_! f^! A \to A \ .$$

Applying these to $A = \underline{\mathbb{K}}_Y$ and $A = \omega_Y$, respectively, and composing with $\pi_*$ for $\pi : Y \to \mathrm{pt}$, we obtain maps

$$f^* : C^\bullet(Y) \to C^\bullet(X) \qquad \text{and} \qquad f_* : C_\bullet(X) \to C_\bullet(Y)$$

of objects in $D(\mathrm{pt}) = \mathrm{Vect}$, as expected for usual chains and cochains.

If $f$ is proper, then we similarly have maps

$$f^* : C_c^\bullet(Y) \to C_c^\bullet(X) \qquad \text{and} \qquad f_* : C_\bullet^{\mathrm{BM}}(X) \to C_\bullet^{\mathrm{BM}}(Y) \ ,$$

while if $f$ is smooth of relative dimension $d = d_X - d_Y$, then we have maps

$$f^* : C_\bullet^{\mathrm{BM}}(Y) \to C_\bullet^{\mathrm{BM}}(X)[-2d] \qquad \text{and} \qquad f_* : C_c^\bullet(X) \to C_c^\bullet(Y)[-2d] \ .$$

Finally, if $f$ is proper and smooth of relative dimension $d$, then we have maps

$$f^* : C_\bullet(Y) \to C_\bullet(X)[-2d] \qquad \text{and} \qquad f_* : C^\bullet(X) \to C^\bullet(Y)[-2d] \ .$$

A.4. **The de Rham functor and Riemann-Hilbert correspondence.** Throughout this section, let $X$ be a smooth, finite dimensional variety over $\mathbb{K} = \mathbb{C}$, and let $\Omega_X^\bullet \in \mathrm{Sh}_{\mathrm{z}}(X)$ denote the algebraic de Rham complex, viewed as a complex of sheaves with the usual de Rham differential.

*Remark* A.4.1. Each $\Omega_X^i \in \mathcal{O}_X$-Mod is a coherent $\mathcal{O}_X$ module, but the differential on $\Omega_X^\bullet$ is not $\mathcal{O}_X$ linear. Rather, the de Rham differential $d_{\mathrm{dR}} \in \mathrm{Diff}(\Omega^i, \Omega^{i+1})$ is a differential operator, so the de Rham complex can equivalently be described in terms of the induced complex of $D$ modules $\Omega_{X,\mathcal{D}}^\bullet = \Omega_X^\bullet \otimes_{\mathcal{O}_X} \mathcal{D}_X \in D^r(X)$, recalling $\mathrm{Diff}(\mathcal{F}, \mathcal{G}) = \mathrm{Hom}_{D(X)}(\mathcal{F}_D, \mathcal{G}_D)$.

*Proposition* A.4.2. There is a natural quasiisomorphism $\Omega_{X,\mathcal{D}}^\bullet \xrightarrow{\cong} \omega_X[-2d_X] = \underline{\mathbb{K}}_X \in D^r(X)$.

*Example* A.4.3. For $M \in D(X)$, applying this resolution to the calculation of de Rham cochains of $M^l \in D^l(X)$ following A.2.15 yields

$$C_{\mathrm{dR}}^\bullet(X; A) = \pi_\bullet(\omega_X \otimes_{\mathcal{D}_X} M) \cong \Gamma(X, \Omega_X^\bullet \otimes_{\mathcal{O}_X} M^l)[2d_X]$$

where $\Omega_X^\bullet \otimes_{\mathcal{O}_X} M \in \mathrm{Sh}_{\mathrm{z}}(X)$ denotes the usual de Rham complex with coefficients in a complex of $\mathcal{O}_X$ modules with flat connection. For $X$ smooth and projective, this is calculated by $\Gamma(X^{\mathrm{an}}, \Omega_{X^{\mathrm{an}}}^\bullet \otimes_{\mathcal{O}_{X^{\mathrm{an}}}} M^{\mathrm{an}})$ where $\Omega_{X^{\mathrm{an}}}^\bullet \otimes_{\mathcal{O}_{X^{\mathrm{an}}}} M^{\mathrm{an}} \in \mathrm{Sh}(X^{\mathrm{an}}; \mathbb{C})$ denotes the analytic variant of the above de Rham complex.

*Definition* A.4.4. The analytic de Rham functor is

$$\mathrm{dR} : D(X) \to \mathrm{D}^b(X) \qquad \text{defined by} \qquad \mathrm{dR}(M) = \Omega_{X^{\mathrm{an}}}^\bullet \otimes_{\mathcal{O}_{X^{\mathrm{an}}}} M^{l,\mathrm{an}}[2d_x] \ ,$$

for each $M \in D(X)$, where $\mathrm{D}^b(X) = \mathrm{D}^b(X(\mathbb{C}))$ denotes the derived category of sheaves on $X$ in the analytic topology, as in subappendix A.5.

Let $\mathrm{D}_c^b(X) = \mathrm{D}_c^b(X(\mathbb{C}); \mathbb{C})$ denote the bounded derived category constructible sheaves on $X$ in the analytic topology, as defined in *loc. cit.*.

*Theorem* A.4.5. The de Rham functor restricts to a derived equivalence $\mathrm{dR} : D_{\mathrm{rh}}(X) \to \mathrm{D}_c^b(X)$.

Moreover, it naturally intertwines the six functors operations stated in Theorems A.3.1 and A.5.3, with the caveat that it intertwines the $\otimes^*$ tensor structure on $D_{\mathrm{rh}}(X)$ of Definition A.3.2 with that of Theorem A.5.3.

*Remark* A.4.6. The object $\mathcal{O}_X \in D^l(X)$ or equivalently $\omega_X \in D^r(X)$ corresponds to

$$\mathrm{dR}(\mathcal{O}_X) = \Omega_X^\bullet[2n] \cong \underline{\mathbb{K}}_X[2n]$$

which is the dualizing sheaf in $D_c^b(X)$. Equivalently, the objects $\underline{\mathbb{K}}_X = \mathcal{O}_X[-2n]$ and $\underline{\mathrm{IC}}_X = \mathcal{O}_X[-n]$ correspond to the constant sheaf $\mathbb{K}_X$ and the intersection cohomology sheaf $\mathrm{IC}_X = \mathbb{K}_X[n]$, in keeping with Remark A.2.5.

A.5. **Constructible Sheaves.** Let $X$ be a quasiprojective algebraic variety over $\mathbb{C}$; we use the same notation to denote $X(\mathbb{C})$ in the analytic topology. Let $\mathrm{Sh}(X; \mathbb{K})$ denote the category of sheaves of $\mathbb{K}$ vector spaces on $X$ and $\mathrm{D}^b(X)$ and $\mathrm{D}^+(X)$ the bounded and bounded below derived categories; we fix the base coefficient field $\mathbb{K}$ once and for all, and supress it from the notation throughout.

*Definition* A.5.1. A *stratification* of $X$ is a finite collection $(X_s)_{s \in \mathcal{S}}$ of disjoint, smooth, connected, locally closed subvarieties such that $X = \cup_{s \in \mathcal{S}} X_s$, and for any $s, t \in \mathcal{S}$ the intersection $\overline{X}_s \cap X_t$ is either empty or $X_t$.

The set $\mathcal{S}$ is equipped with the *closure partial order*, defined by $t \leqslant s$ if $X_t \subset \overline{X}_s$.

*Definition* A.5.2. A sheaf $\mathcal{F} \in \mathrm{Sh}(X, \mathbb{K})$ is *constructible with respect to a stratification* $(X_s)_{s \in \mathcal{S}}$ if its restriction $\mathcal{F}|_{X_s}$ to each stratum $X_s$ is a finite rank local system; $\mathcal{F}$ is *constructible* if it is constructible with respect to some stratification of $X$.

A complex $\mathcal{F} \in \mathrm{D}^b(X, \mathbb{K})$ is *constructible with respect to a stratification* if its cohomology sheaves $H^k(\mathcal{F})$ are constructible with respect to that stratification for each $k$; $\mathcal{F}$ is *constructible* if it is constructible with respect to some stratification of $X$.

Let $\mathrm{Sh}_c(X)$ denote the abelian category of constructible sheaves, $\mathrm{D}^b_c(X)$ its bounded derived category, and similarly $\mathrm{Sh}_{\mathcal{S}}(X)$ and $\mathrm{D}^b_{\mathcal{S}}(X)$ the categories of sheaves constructible with respect to a fixed stratification $(X_s)_{s \in \mathcal{S}}$. The canonical functor $\mathrm{D}^b_c(X) \to \mathrm{D}^b(X)$ is fully faithful, so that $\mathrm{D}^b_c(X)$ is equivalent to the full subcategory of $\mathrm{D}^b(X)$ of constructible objects, and similarly for $\mathrm{D}^b_{\mathcal{S}}(X)$. Let $\mathrm{D}^+_c(X)$ and $\mathrm{D}^+_{\mathcal{S}}(X)$ denote the analogous bounded below derived categories. We identify the $\mathrm{D}^b_c(\mathrm{pt}) = \mathrm{D}^b_{\mathrm{fg}}(\mathbb{K}\text{-Mod}) = \mathrm{Perf}_{\mathbb{K}}$ with the bounded derived category of complexes with finite dimensional cohomology throughout.

*Theorem* A.5.3. There are functors

$$\otimes : \mathrm{D}^b_c(X) \times \mathrm{D}^b_c(X) \to \mathrm{D}^b_c(X) \qquad \mathbb{D}_X : \mathrm{D}^b_c(X)^{\mathrm{op}} \to \mathrm{D}^b_c(X) \qquad \underline{\mathrm{Hom}}_X : \mathrm{D}^b_c(X)^{\mathrm{op}} \times \mathrm{D}^b_c(X) \to \mathrm{D}^b_c(X) \,,$$

and for $f : X \to Y$ natural adjunctions,

$$f^* : \mathrm{D}^b_c(Y) \rightleftarrows \mathrm{D}^b_c(X) : f_* \qquad f_! : \mathrm{D}^b_c(X) \rightleftarrows \mathrm{D}^b_c(Y) : f^! \quad .$$

These satisfy the standard six functor formalism compatibilities, as in A.3.1 above.

*Definition* A.5.4. The constructible (Borel-Moore) (co)chains on $X$ are

$$C_\bullet(X) = \pi_! \pi^! \underline{\mathbb{K}}_{\mathrm{pt}} \qquad C_\bullet^{\mathrm{BM}}(X) = \pi_* \pi^! \underline{\mathbb{K}}_{\mathrm{pt}} \qquad C^\bullet(X) = \pi_* \pi^* \underline{\mathbb{K}}_{\mathrm{pt}} \qquad C_c^\bullet(X) = \pi_! \pi^* \underline{\mathbb{K}}_{\mathrm{pt}} \quad \in \quad D^b_c(\mathrm{pt})$$

where $\pi : X \to \mathrm{pt}$ is the unique map. More generally, for $A \in \mathrm{D}^b_c(X)$ we define $C^\bullet(X; A) = \pi_* A$ and $C_\bullet(X; A) = \pi_! A$.

*Remark* A.5.5. The above objects are equivalent to those from Definition A.3.5, by the Riemann-Hilbert correspondence A.4.5, so there is no ambiguity in the notation.

*Remark* A.5.6. The above objects satisfy the same functoriality as in Remark A.3.8, again by the Riemann-Hilbert correspondence A.4.5.

*Remark* A.5.7. The diagonal map $\Delta : X \to X \times X$ endows $C^\bullet(X) \in \mathrm{Comm}(\mathrm{Perf}_{\mathbb{K}})$ with the structure of a commutative algebra and $C_\bullet(X) \in \mathrm{CoComm}(\mathrm{Perf}_{\mathbb{K}})$ a cocommutative coalgebra.

*Definition* A.5.8. The (Borel-Moore) (co)*homology* groups of $X$ are defined as the images of the objects in Definition A.5.4 under $H^\bullet : \mathrm{D}(\mathbb{K}\text{-Mod}) \to \mathbb{K}\text{-Mod}_{\mathbb{Z}}$.

*Remark* A.5.9. The functor $H^\bullet$ is lax monoidal, so that $H^\bullet(X) \in \mathrm{Comm}(\mathbb{K}\text{-}\mathrm{Mod}_\mathbb{Z})$ is again a commutative algebra; similarly, but via universal coefficients, $H_\bullet(X) \in \mathrm{CoComm}(\mathbb{K}\text{-}\mathrm{Mod}_\mathbb{Z})$ is a cocommutative algebra.

## Appendix B. Equivariant Cohomology

B.1. **Equivariant sheaves.** In this subsection, we explain the formalism of equivariant sheaf theory, and how it gives rise to equivariant cohomology. We review the theory of equivariant $D$ modules in Section 15 of the main body of the present work. We give the exposition here in terms of constructible sheaves.

Let $X$ be a quasiprojective algebraic variety over $\mathbb{C}$, $G$ a connected reductive algebraic group, and $K$ the compact real form of $G$; we use the same notation to denote $X(\mathbb{C}), G(\mathbb{C})$ and $K(\mathbb{R})$ in the analytic topology. In particular, the topological space $X$ together with the action of $G$ or $K$ satisfies the hypotheses of the references [BL94] and [GKM97], which we follow closely throughout.

Let $\mathrm{D}^b_G(X)$ and $\mathrm{D}^+_G(X)$ denote the bounded, and bounded below, derived categories of $G$ equivariant sheaves on $X$, as defined in 2.2 and 2.8 of [BL94].

*Remark* B.1.1. The category $\mathrm{D}^b_G(X)$ can be presented as the category $\mathrm{D}^b(X/G)$ of sheaves on the quotient stack $X/G$, viewed as a simplicial space, and in particular we have canonical functors

$$\mathrm{for}_G : \mathrm{D}^b_G(X) \to \mathrm{D}^b(X) \qquad q^* : \mathrm{D}^b(\bar{X}) \to \mathrm{D}^b_G(X)$$

defined by forgetting the equivariant structure, and by pullback along the canonical map $q$ to the quotient topological space $\bar{X}$, respectively. The latter is an equivalence if the action of $G$ is free.

*Definition* B.1.2. An equivariant sheaf $\mathcal{F} \in \mathrm{D}^b_G(X)$ is constructible with respect to a stratification if the underlying sheaf in $\mathrm{D}^b(X)$ is constructible.

We denote the bounded derived category of $G$ equivariant constructible sheaves by $\mathrm{D}^b_{G,c}(X)$, and similarly for $\mathrm{D}^+_{G,c}(X)$.

*Theorem* B.1.3. For $H$ a subgroup of $G$, there are restriction and induction adjunctions:

$$\mathrm{Res}^G_H : \mathrm{D}^b_{G,c}(X) \rightleftarrows \mathrm{D}^b_{H,c}(X) : \mathrm{Ind}^G_{H,*} \qquad \mathrm{Ind}^b_{H,!} : \mathrm{D}^b_{H,c}(X) \rightleftarrows \mathrm{D}^b_{G,c}(X) : \mathrm{Res}^G_H$$

Moreover, there are functors as in A.5.3, satisfying the same adjunctions and relations, defined for $G$ equivariant maps $f : X \to Y$. These functors all commute with $\mathrm{Res}^G_H$, while $f_*$ and $f^!$ commute with $\mathrm{Ind}^G_{h,*}$, and $f_!$ and $f^*$ commute with $\mathrm{Ind}^G_{H,!}$.

*Definition* B.1.4. The equivariant (Borel-Moore) (co)chains on $X$ are

$$C^G_\bullet(X) = \pi_! \pi^! \underline{\mathbb{K}}_{\mathrm{pt}} \qquad C^{G,\mathrm{BM}}_\bullet(X) = \pi_* \pi^! \underline{\mathbb{K}}_{\mathrm{pt}} \qquad C^\bullet_G(X) = \pi_* \pi^* \underline{\mathbb{K}}_{\mathrm{pt}} \qquad C^\bullet_{G,\mathrm{BM}}(X) = \pi_! \pi^* \underline{\mathbb{K}}_{\mathrm{pt}} \quad \in \quad \mathrm{D}^b_{G,c}(\mathrm{pt}),$$

where $\pi : X \to \mathrm{pt}$ is the unique map. More generally, for $A \in \mathrm{D}^b_{G,c}(X)$, we define $C^\bullet(X; A) = \pi_* A$ and $C_\bullet(X; A) = \pi_! A$.

*Remark* B.1.5. The functoriality properties outlined in Remark A.3.8 hold for equivariant maps. In particular, by Remark A.5.7, $C^\bullet_G(X) \in \mathrm{Comm}(\mathrm{D}^b_{G,c}(\mathrm{pt}))$ is a commutative algebra and $C^G_\bullet(X) \in \mathrm{CoComm}(\mathrm{D}^b_{G,c}(\mathrm{pt}))$ is a cocommutative coalgebra.

*Remark* B.1.6. There is a functor $H^\bullet : \mathrm{D}^b_{G,c}(\mathrm{pt}) \to \mathbb{K}\text{-Mod}_\mathbb{Z}$ given by forgetting the equivariant structure and applying the usual cohomology object functor. This evidently does not depend on the equivariant structure, and for example $H^\bullet(C^\bullet_G(X)) = H^\bullet(X)$ is just the usual cohomology.

There is another functor on $\mathrm{D}^b_{G,c}(\mathrm{pt})$ which is the equivariant analogue of the cohomology object functor: the data of the equivariant structure defines a functor $\mathrm{D}^b_{G,c}(\mathrm{pt}) \to \mathrm{D}^b(BG)$ and composing with the global sections functor $\pi_* : \mathrm{D}^b(BG) \to \mathrm{D}^b(\mathrm{pt})$ gives the desired functor $\mathrm{D}^b_{G,c}(\mathrm{pt}) \to$ $\mathbb{K}$-Mod$_{\mathbb{Z}}$. This functor is lax monoidal, so that the image of $\underline{\mathbb{K}}_{pt} = C_G^\bullet(pt)$ under it defines a commutative algebra object, which we denote by $H_G^\bullet(pt) \in \mathrm{Comm}(\mathbb{K}\text{-Mod}_{\mathbb{Z}})$.

The object $C_G^\bullet(pt)$ is the monoidal unit of $\mathrm{D}_{G,c}^b(pt)$, so that every object is canonically a module object for it, and thus the above functor lifts to define a functor $H_G^\bullet : \mathrm{D}_{G,c}^b(pt) \to H_G^\bullet(pt)\text{-Mod}_{\mathbb{Z}}$, which we call the equivariant cohomology object functor.

More generally, the above functor $\pi_* : \mathrm{D}_{G,c}^b(pt) \to \mathrm{D}^+(pt)$ lifts to a functor $G : \mathrm{D}_{G,c}^b(pt) \to \mathrm{D}_{fg}^+(H_G^\bullet(pt))$.

*Definition* B.1.7. The equivariant (Borel-Moore) (co)homology groups of a $G$ space $X$ are the images of the objects in Definition B.1.4 under $H_G^\bullet : \mathrm{D}_{G,c}^b(pt) \to H_G^\bullet(pt)\text{-Mod}_{\mathbb{Z}}$.

*Remark* B.1.8. The image $H_G^\bullet(C_G^\bullet(pt)) = H_G^\bullet(pt)$ is given by the object defined previously, so the notation is consistent.

For the remainder of this section, let $\mathbb{S} = H_G^\bullet(pt) \in \mathrm{Comm}(\mathbb{K}\text{-Mod}_{\mathbb{Z}})$. Let $C^+(\mathbb{S})$ denote the DG category of bounded below DG modules over $\mathbb{S}$, $K^+(\mathbb{S})$ the homotopy category obtained by quotienting by homotopy equivalences, $\mathrm{D}^+(\mathbb{S})$ the derived category obtained by localizing at quasi-isomorphisms, and similarly $C_{fg}^+(\mathbb{S})$, $K_{fg}^+(\mathbb{S})$ and $\mathrm{D}_{fg}^+(\mathbb{S})$ those with finitely generated cohomology modules. The canonical functor $\mathrm{D}_{fg}^+(\mathbb{S}) \to \mathrm{D}^+(\mathbb{S})$ is fully faithful, so that $\mathrm{D}_{fg}^+(\mathbb{S})$ is equivalent to the full subcategory of objects with finitely generated cohomology.

*Remark* B.1.9. There are standard functors on categories of DG modules

$$\otimes_{\mathbb{S}} : \mathrm{D}^+(\mathbb{S})^{\times 2} \to \mathrm{D}^+(\mathbb{S}) \qquad H^\bullet : \mathrm{D}^+(\mathbb{S}) \to \mathbb{S}\text{-Mod}_{\mathbb{Z}} \qquad \mathbb{D}_{\mathbb{S}} : \mathrm{D}_{fg}^+(\mathbb{S}) \to \mathrm{D}_{fg}^+(\mathbb{S}) \qquad \underline{\mathrm{Hom}}_{\mathbb{S}} : \mathrm{D}_{fg}^+(\mathbb{S}) \times \mathrm{D}^+(\mathbb{S}) \to \mathrm{D}^+(\mathbb{S})$$

preserving the subcategories $\mathrm{D}_{fg}^+(\mathbb{S})$.

*Theorem* B.1.10. There is a canonical triangulated equivalence $\mathcal{L}_G : \mathrm{D}^+(\mathbb{S}) \xrightarrow{\cong} \mathrm{D}_G^+(pt)$, inverse to a functor $\mathrm{D}_G^+(pt) \to \mathrm{D}^+(\mathbb{S})$ generalizing $G$ defined in Remark B.1.6, and intertwining the equivariant cohomology and tensor product functors. Further, this induces an equivalence $\mathrm{D}_{fg}^+(\mathbb{S}) \cong \mathrm{D}_{G,c}^+(pt)$ of full triangulated subcategories, inverse to $G$ of *loc cit.*, and intertwining the duality and internal Hom functors.

The intertwining conditions are stated precisely in Theorem B.2.4 below, alongside those for another model of $\mathrm{D}_{G,c}^b(pt)$, which we discuss in the following section.

## B.2. Goresky-Kottwitz-MacPherson Koszul Duality.

Let $\Lambda = H_\bullet(G) \in \mathrm{CoAss}(\mathbb{K}\text{-Mod}_{\mathbb{Z}})$ denote the homology of $G$, considered as a graded cocommutative coalgebra over $\mathbb{K}$. The group structure maps define a compatible unit, antipode, and associative product on $\Lambda$, making it into a cocommutative Hopf algebra over $\mathbb{K}$; see also Section 23. Let $\mathrm{D}^b(\Lambda), \mathrm{D}_{fg}^b(\Lambda), \mathrm{D}^+(\Lambda), \mathrm{D}_{fg}^+(\Lambda)$ the bounded, bounded below and/or with finite dimensional cohomology derived categories, as above. Our convention is such that $H_\bullet(G)$ is non-positively graded so that it acts on modules by non-positive degree endomorphisms of a complex.

*Remark* B.2.1. The cocommutative coalgebra structure on $\Lambda$ gives a lift of the tensor product over $\mathbb{K}$ to $\otimes_{\mathbb{K}} : \mathrm{D}^+(\Lambda)^{\times 2} \to \mathrm{D}^+(\Lambda)$ which defines a symmetric monoidal structure on $\mathrm{D}^+(\Lambda)$. Similarly, the antipode on $\Lambda$ gives a lift of the dual over $\mathbb{K}$ to $\mathbb{D}_\Lambda = \mathrm{Hom}_{\mathbb{K}}(\cdot, \mathbb{K}) : \mathrm{D}_{fg}^+(\Lambda) \to \mathrm{D}_{fg}^+(\Lambda)$ and together these define an internal Hom functor $\underline{\mathrm{Hom}}_\Lambda : \mathrm{D}_{fg}^+(\Lambda) \times \mathrm{D}^+(\Lambda) \to \mathrm{D}^+(\Lambda)$.

*Remark* B.2.2. The action map $G \times X \to X$ on a $G$ space $X$ defines a natural $\Lambda$ module structure on $H_\bullet(X)$, and dually $H^\bullet(X)$, so that they naturally lift to objects $H_\bullet(X), H^\bullet(X) \in \Lambda\text{-Mod}_{\mathbb{Z}}$.

More generally, the forgetful functor $D_{G,c}^b(\mathrm{pt}) \to D_c^b(\mathrm{pt}) = \mathrm{Perf}_{\mathbb{K}}$ lifts to $E: D_{G,c}^b(\mathrm{pt}) \to D_{\mathrm{fg}}^+(\Lambda)$, and similarly for $D^+$.

*Remark* B.2.3. The graded associative algebras $\mathbb{S} = H_G^\bullet(\mathrm{pt})$ and $\Lambda = H_\bullet(G)$ are Koszul dual in the sense of [BGG71, BGS96]. In particular, there is a canonical functor

(B.2.1)
$$t: C^+(\Lambda) \to C^+(\mathbb{S}) \qquad (N, d_N) \mapsto t(N, d_N) = \left( \mathbb{S} \otimes_{\mathbb{K}} N \ , \ d_{t(N)}(s,n) = \sum_i \xi_i s \otimes x_i n + s \otimes d_N(n) \right)$$

where $(x_i)$ denotes a basis for the generators of $\Lambda$ over $\mathbb{K}$ and $(\xi_i)$ the dual basis for the generators of $\mathbb{S}$. There is a functor $h: C^+(\mathbb{S}) \to C^+(\Lambda)$ defined similarly, and these induce inverse equivalences on $D^+$ and $D_{\mathrm{fg}}^b$.

The object $t(N) \in C^+(\mathbb{S})$ above can also be understood as the total complex of the double complex:

(B.2.2)

In the case $N = C^\bullet(X; \mathbb{K})$ for a $G$ space $X$ as in Remark B.2.2, this is precisely the double complex presentation that induces the Serre spectral sequence for the cohomology of the fibration $X \hookrightarrow X/G \twoheadrightarrow BG$.

The results of this section on models for $D_{G,c}^b(\mathrm{pt})$ and their compatibilities are summarized in the following theorem from [GKM97], following [BL94] and the Koszul duality results from [BGG71, BGS96].

*Theorem* B.2.4. Let $G$ be a connected Lie group, $\Lambda = H_\bullet(G; \mathbb{K})$, and $\mathbb{S} = H_G^\bullet(\mathrm{pt}; \mathbb{K})$. There exist commuting triangulated equivalences

such that we have compatible commutativity of the following diagrams:

In particular, for a $G$ space $X$ and $A \in \mathrm{D}^b_{G,c}(X)$ we have

$$H^\bullet(X; A) = H^\bullet \circ E \circ \pi_*(A) \qquad\qquad H^\bullet_G(X; A) = H^\bullet \circ G \circ \pi_*(A) \quad \text{, and}$$

$$H_\bullet(X; A) = H^\bullet \circ E \circ \pi_!(A) \qquad\qquad H^G_\bullet(X; A) = H^\bullet \circ G \circ \pi_!(A) \quad .$$

*Remark* B.2.5. As we explain in Example 15.0.13, in the $D$ module setting the preceding theorem gives rise to the usual Cartan model for equivariant de Rham cohomology, by applying the functor from Equation B.2.1 to the de Rham complex together with its canonical equivariant structure.

B.3. **The equivariant localization theorem.** In this section, we recall the equivariant localization theorem, originally proved in [AB95], in the setting of sheaf cohomology, following Section 6.2 of [GKM97]. In fact, we recall a variant that uses the homological excision sequence rather than the more commonly stated version, which relies on excision in cohomology; both variants follow readily from the results of *loc. cit.*. For simplicity, we restrict to the case that $G = (\mathbb{C}^\times)^n$ is an algebraic torus, so that we identify the equivariant cohomology of a point

$$\mathbb{S} = H^\bullet_G(\mathrm{pt}) = \mathbb{K}[\mathfrak{g}[2]] \ ,$$

with the coordinate ring of the vector space underlying the Lie algebra $\mathfrak{g}$, shifted in cohomological degree by $-2$.

For each point $x \in X$, let $G_x$ be (the connected component of the identity in) the stabilizer in $G$ of $x$, and let $\mathfrak{g}_x = \mathrm{Lie}(G_x)$ be its Lie algebra. Further, let $\iota : Z \hookrightarrow X$ be the inclusion of a closed, $G$-invariant subvariety, $j : U = X \backslash Z \hookrightarrow X$ the inclusion of the complementary open, and recall that for each $A \in \mathrm{D}^b_G(X)$ the homological excision exact triangle

$$\iota_* \iota^! A \to A \to j_* j^! A \qquad \text{induces} \qquad C^G_\bullet(Z; A) \to C^G_\bullet(X; A) \to C^G_\bullet(X, Z; A) \ ,$$

in the category $\mathrm{D}^b_{\mathrm{fg}}(\mathbb{S})$, by applying $E \circ \pi_!$. We now state the main result of this section

*Theorem* B.3.1. [GKM97] Let $Z \hookrightarrow X$ be a closed, $G$-invariant subvariety of $X$ containing the $G$-fixed points $X^G \subset Z$. Then $H^G_\bullet(X, Z; A)$ is a torsion module over $\mathbb{S}$, with support

$$\mathrm{supp}(H^G_\bullet(X, Z; A)) \subset \bigcup_{x \in X \backslash Z} \mathfrak{g}_x$$

contained in the union of the (finitely many distinct) stabilizer subalgebras $\mathfrak{g}_x$ of points $x \in X \backslash Z$. In particular, if $\{f_i \in \mathbb{S}\}$ generate an ideal whose corresponding subvariety of $\mathrm{Spec}\,\mathbb{S}$ contains $\cup_{x \in X \backslash Z} \mathfrak{g}_x$, then the natural map

$$H^G_\bullet(Z; A) \xrightarrow{\cong} H^G_\bullet(X; A) \qquad\qquad \text{is an isomorphism over } H^G_\bullet(\mathrm{pt})[f_i^{-1}].$$

## Appendix C. Operads

C.1. **Operads and algebras.** Let $\mathcal{C}$ be a symmetric monoidal category with monoidal structure $\otimes : \mathcal{C} \times \mathcal{C} \to \mathcal{C}$, tensor unit $u_{\mathcal{C}} \in \mathcal{C}$.

*Definition* C.1.1. A (symmetric, coloured) *operad* $\mathcal{O}$ in $\mathcal{C}$ is:

- A collection col $\mathcal{O}$, elements of which are called *colours* or *objects* of $\mathcal{O}$
- For each finite set $I$ and indexed collection of objects $\{c_i\}_{i \in I}$ and $d$ of $\mathcal{O}$, an object $\mathcal{O}(\{c_i\}_{i \in I}, d) \in \mathcal{C}$ called the *multilinear operations* in $\mathcal{O}$.
- For each map of finite sets $\pi : I \to J$, and indexed collections of objects $\{c_i\}_{i \in I}$, $\{d_j\}_{j \in J}$ and $e$, a morphism

$$\bigotimes_{j \in J} \mathcal{O}(\{c_i\}_{i \in I_j}, d_j) \otimes \mathcal{O}(\{d_j\}_{j \in J}, e) \to \mathcal{O}(\{c_I\}_{i \in I}, e)$$

  of objects of $\mathcal{C}$ called the *composition law* in $\mathcal{O}$.
- For each object $c \in \mathcal{C}$, a morphism $\mathbb{1}_c \in \mathcal{O}(c, c)$ called the *identity map* on $c$, which is both a left and right unit for the composition law.
- For each sequence of maps $I \xrightarrow{\pi} J \xrightarrow{\varphi} K$, and indexed collections of objects $\{c_i\}_{i \in I}$, $\{d_j\}_{j \in J}$, $\{e_k\}_{k \in K}$ and $f$, the commutativity of the diagram

$$\bigotimes_{j \in J} \mathcal{O}(\{c_i\}_{i \in I_j}, d_j) \otimes \bigotimes_{k \in K} \mathcal{O}(\{d_j\}_{j \in J_k}, e_k) \otimes \mathcal{O}(\{e_k\}_{k \in K}, f) \longrightarrow \bigotimes_{k \in K} \mathcal{O}(\{c_i\}_{i \in I_k}, e_k) \otimes \mathcal{O}(\{e_k\}_{k \in K}, f)$$

$$\downarrow \qquad\qquad\qquad\qquad\qquad\qquad\qquad\qquad\qquad\qquad\qquad\qquad\qquad\qquad \downarrow$$

$$\bigotimes_{j \in J} \mathcal{O}(\{c_i\}_{i \in I_j}, d_j) \otimes \mathcal{O}(\{d_j\}, f) \longrightarrow \mathcal{O}(\{c_i\}_{i \in I}, f)$$

A map $\varphi : \mathcal{O} \to \mathcal{O}'$ of operads in $\mathcal{C}$ is:

- A map col $\mathcal{O} \to$ col $\mathcal{O}'$
- For each map of finite sets $\pi : I \to J$, and indexed collections of objects $\{c_i\}_{i \in I}$ and $d$ of $\mathcal{O}$, a morphism

$$\mathcal{O}(\{c_i\}_{i \in I_j}, d) \to \mathcal{O}'(\{\varphi(c_i)\}_{i \in I_j}, \varphi(d))$$

  such that $\mathbb{1}_c$ maps to $\mathbb{1}_{\varphi(c)}$ for each $c \in$ col $\mathcal{O}$.
- For each map of finite sets $\pi : I \to J$, and indexed collections of objects $\{c_i\}_{i \in I}$, $\{d_j\}_{j \in J}$ and $e$, the commutativity of the diagram

$$\bigotimes_{j \in J} \mathcal{O}(\{c_i\}_{i \in I_j}, d_j) \otimes \mathcal{O}(\{d_j\}_{j \in J}, e) \longrightarrow \mathcal{O}(\{c_i\}_{i \in I}, e)$$

$$\downarrow \qquad\qquad\qquad\qquad\qquad\qquad\qquad\qquad\qquad \downarrow$$

$$\bigotimes_{j \in J} \mathcal{O}'(\{\varphi(c_i)\}_{i \in I_j}, \varphi(d_j)) \otimes \mathcal{O}'(\{\varphi(d_j)\}_{j \in J}, \varphi(e)) \longrightarrow \mathcal{O}'(\{\varphi(c_i)\}_{i \in I}, \varphi(e))$$

The collection of operads in $\mathcal{C}$ thus defines a category, denoted $\mathrm{Op}(\mathcal{C})$. We write simply $\mathrm{Op}$ in the case that $\mathcal{C} = \mathrm{Set}$, and $\mathrm{Op}_{\mathbb{K}}$ in the case that $\mathcal{C} = \mathbb{K}\text{-Mod}$.

Let $\mathrm{Alg}_{\mathcal{O}}(\mathcal{O}')$ denote the space $\mathrm{Hom}_{\mathrm{Op}(\mathcal{C})}(\mathcal{O}, \mathcal{O}')$ of map of operads $\mathcal{O} \to \mathcal{O}'$.

*Remark* C.1.2. For an operad $\mathcal{O}$ with a single object col $\mathcal{O} = \{c\}$, we use the notation $\mathcal{O}(I) = \mathcal{O}(\{c\}_{i \in I}, c)$ and $\mathcal{O}(n) = \mathcal{O}(I)$ for $I = \{1, ..., n\}$. Note that $\mathcal{O}(I)$ has the structure of an $\mathrm{Aut}_{\mathrm{fSet}}(I)$ module, and similarly $\mathcal{O}(n)$ an $S_n$ module in $\mathcal{C}$.

*Example* C.1.3. Let $D$ be a symmetric monoidal category enriched over $\mathcal{C}$. Then $D$ defines an operad $\mathcal{O}_D \in \mathrm{Op}(\mathcal{C})$ in $\mathcal{C}$, with objects given by objects of $D$ and multilinear operations defined by

$$\mathcal{O}_D(\{c_i\}_{i \in I}, d) = \mathrm{Hom}_D(\otimes_{i \in I} c_i, d) \ .$$

In this case, we abreviate $\mathrm{Alg}_{\mathcal{O}}(\mathcal{O}_D) = \mathrm{Alg}_{\mathcal{O}}(D)$.

In particular, if $\mathcal{C}$ is a closed monoidal category, then for any operad $\mathcal{O}$ in $\mathcal{C}$, we have a canonical category $\mathrm{Alg}_{\mathcal{O}}(\mathcal{C}) := \mathrm{Alg}_{\mathcal{O}}(\mathcal{O}_{\mathcal{C}})$ of algebras over $\mathcal{O}$ internal to $\mathcal{C}$. This definition generalizes to arbitrary symmetric monoidal $\mathcal{C}$ by hom-tensor adjunction, though $\mathcal{O}_{\mathcal{C}}$ no longer defines an operad in $\mathcal{C}$.

*Example* C.1.4. Suppose $\mathcal{C}$ has initial object $\varnothing_{\mathcal{C}}$. The trivial operad $\mathrm{triv}_{\mathcal{C}}$ in $\mathcal{C}$ is defined as having a single object, with multilinear operations given by $\mathrm{triv}_{\mathcal{C}}(\varnothing) = \mathrm{u}_{\mathcal{C}}$, $\mathrm{triv}_{\mathcal{C}}(\{\mathrm{pt}\}) = \mathrm{u}_{\mathcal{C}}$ and $\mathrm{triv}(I) = \varnothing_{\mathcal{C}}$ for $|I| \neq 0, 1$. The category of triv algebras $\mathrm{Alg}_{\mathrm{triv}_{\mathcal{C}}}(\mathcal{C}) = \mathcal{C}$ is equivalent to the underlying category $\mathcal{C}$.

*Example* C.1.5. The commutative operad $\mathrm{Comm}_{\mathbb{K}}$ in $\mathrm{Vect}_{\mathbb{K}}$ is defined as having a single object, with multilinear operations given by $\mathrm{Comm}_{\mathbb{K}}(I) = \mathbb{K}$ for all $I \in \mathrm{fSet}$. For $\mathcal{C}$ a $\mathbb{K}$ linear symmetric monoidal category, the category of $\mathrm{Comm}_{\mathbb{K}}$ algebras $\mathrm{Comm}_{\mathbb{K}}(\mathcal{C}) := \mathrm{Alg}_{\mathrm{Comm}_{\mathbb{K}}}(\mathcal{C})$ in $\mathcal{C}$ is the usual category of (unital) commutative algebra objects in $\mathcal{C}$.

*Example* C.1.6. The associative operad $\mathrm{Ass}_{\mathbb{K}}$ in $\mathrm{Vect}_{\mathbb{K}}$ is defined as having a single object, with multilinear operations given by $\mathrm{Ass}_{\mathbb{K}}(I) = \mathbb{K}[S_I]$ the regular representation of the symmetric group $S_I = \mathrm{Aut}_{\mathrm{fSet}}(I)$ on $I$ for each $I \in \mathrm{fSet}$. For $\mathcal{C}$ a $\mathbb{K}$ linear symmetric monoidal category, the category of $\mathrm{Ass}_{\mathbb{K}}$ algebras $\mathrm{Ass}_{\mathbb{K}}(\mathcal{C}) := \mathrm{Alg}_{\mathrm{Ass}_{\mathbb{K}}}(\mathcal{C})$ in $\mathcal{C}$ is the usual category of (unital) associative algebra objects in $\mathcal{C}$.

*Example* C.1.7. Let $M = (M_n)_{n \in \mathbb{N}}$ with $M_n \in \mathbb{K}[S_n]$-Mod be a sequence of symmetric group modules in $\mathrm{Vect}_{\mathbb{K}}$. The free operad $\mathcal{F}(M) \in \mathrm{Op}(\mathrm{Vect}_{\mathbb{K}})$ on $M$ is characterized by the property that

$$\mathrm{Alg}_{\mathcal{F}(M)}(\mathcal{O}) = \oplus_{n \in \mathbb{N}} \mathrm{Hom}_{\mathbb{K}[S_n]}(M_n, \mathcal{O}(n))$$

for each $\mathcal{O} \in \mathrm{Op}(\mathrm{Vect}_{\mathbb{K}})$.

More generally, we say an operad $\mathcal{O}$ is generated over $M$ if there exist $R = (R_n)_{n \in \mathbb{N}}$ with $R_n \in \mathbb{K}[S_n]$-Mod together with maps $R_n \hookrightarrow \mathcal{F}(M)(n)$ defining an operadic ideal of $\mathcal{F}(M)$, such that $\mathcal{O} = \mathcal{F}(M)/R$.

*Example* C.1.8. The associative operad $\mathrm{Ass}_{\mathbb{K}}$ is the free operad on $\mathrm{Ass}_{\mathbb{K}}(2) = \mathbb{K}_m \oplus \mathbb{K}_{m^{\mathrm{op}}} \in \mathbb{K}[S_2]$-Mod given by the regular representation, subject to the single relation

$$\mathbb{K}_{\mathrm{Ass}(m)} \hookrightarrow \mathcal{F}(3) \qquad \text{defined by} \qquad 1 \mapsto \mathrm{Ass}(m) = m \circ (m \otimes \mathbb{1}) - m \circ (\mathbb{1} \otimes m) \in \mathcal{F}(3) \ .$$

The commutative operad $\mathrm{Comm}$ is generated by the trivial representation $\mathrm{Comm}(2) = \mathbb{K}_m \in \mathbb{K}[S_2]$-Mod and subject to the same single relation.

*Example* C.1.9. The Lie operad $\mathrm{Lie}_{\mathbb{K}} \in \mathrm{Op}(\mathrm{Vect}_{\mathbb{K}})$ is the operad generated by $\mathrm{Lie}(2) = \mathbb{K}_b \in \mathbb{K}[S_2]$-Mod given by the sign representation, subject to the relation

$$\mathbb{K}_{\mathrm{Jac}(\pi)} \hookrightarrow \mathcal{F}(3) \qquad \text{defined by} \qquad 1 \mapsto \mathrm{Jac}(\pi) = \pi \circ (\mathbb{1} \otimes \pi) - \pi \circ (\pi \otimes \mathbb{1}) - \pi \circ (\mathbb{1} \otimes \pi) \circ \sigma_{12} \ .$$

*Example* C.1.10. Let $D$ be a category enriched over $\mathcal{C}$. Then $D$ defines an operad in $\mathcal{C}$, with objects given by those of $D$ and multilinear operations defined by

$$\mathcal{O}_D(\{c_i\}_{i \in I}, d) = \begin{cases} \mathrm{Hom}_D(c, d) & \text{if } |I| = 1 \\ \varnothing & \text{otherwise} \end{cases} \ .$$

In fact, there is an equivalence of categories between the category $\mathrm{Cat}(\mathcal{C})$ of categories enriched in $\mathcal{C}$, and the category $\mathrm{Op}(\mathcal{C})_{/\mathrm{triv}_{\mathcal{C}}}$ of operads in $\mathcal{C}$ over $\mathrm{triv}_{\mathcal{C}}$, as long as the initial object of $\mathcal{C}$ is strict.

*Example* C.1.11. Let $D$ be a symmetric monoidal category enriched over $\mathcal{C}$ and $D' \hookrightarrow D$ be a subcategory. Then $D'$ is not necessarily closed under the symmetric monoidal structure on $D$, and thus does not necessarily define a symmetric monoidal subcategory. However, $D'$ still defines a suboperad $\mathcal{O}_{D'} \hookrightarrow \mathcal{O}_D$ given by

$$\mathcal{O}_{D'}(\{c_i\}_{i \in I}, d) = \mathrm{Hom}_D(\otimes_{i \in I} c_i, d) \ .$$

In this sense, operads are a natural generalization of symmetric monoidal categories, and are sometimes called multicategories or pseudo-tensor categories, as in C.6.5. Again, we abreviate $\mathrm{Alg}_{\mathcal{O}}(\mathcal{O}_{D'}) = \mathrm{Alg}_{\mathcal{O}}(D')$.

*Example* C.1.12. Let $\mathcal{C}$ be a symmetric monoidal category. Then $\mathcal{C}^{\mathrm{op}}$ is canonically symmetric monoidal, and we define $\mathrm{CoAlg}_{\mathcal{O}}(\mathcal{C}) = \mathrm{Alg}_{\mathcal{O}}(\mathcal{C}^{\mathrm{op}})^{\mathrm{op}}$. In particular, we define coassociative coalgebras and cocommutative coalgebras in $\mathcal{C}$ by $\mathrm{CoAss}(\mathcal{C}) = \mathrm{CoAlg}_{\mathrm{Ass}}(\mathcal{C})$ and $\mathrm{CoComm}(\mathcal{C}) = \mathrm{CoAlg}_{\mathrm{Comm}}(\mathcal{C})$.

*Remark* C.1.13. Let $\varphi : \mathcal{O} \to \mathcal{O}'$ a map of operads in $\mathcal{C}$. There is a functor $\varphi^* : \mathrm{Alg}_{\mathcal{O}'}(\mathcal{O}'') \to \mathrm{Alg}_{\mathcal{O}}(\mathcal{O}'')$.

*Proposition* C.1.14. Let $F : \mathcal{C} \to \mathcal{C}'$ be a lax symmetric monoidal functor. Then $F$ naturally defines a functor $F : \mathrm{Op}(\mathcal{C}) \to \mathrm{Op}(\mathcal{C}')$ and in particular defines $F_{\mathcal{O}} : \mathrm{Alg}_{\mathcal{O}}(\mathcal{O}') \to \mathrm{Alg}_{F(\mathcal{O})}(F(\mathcal{O}'))$ for each $\mathcal{O}, \mathcal{O}' \in \mathrm{Op}(\mathcal{C})$. Further in particular, we obtain functors $\mathrm{Alg}_{\mathcal{O}}(\mathcal{C}) \to \mathrm{Alg}_{F(\mathcal{O})}(F(\mathcal{C})) \to \mathrm{Alg}_{F(\mathcal{O})}(\mathcal{C}')$.

*Example* C.1.15. The functor $C_{\bullet}(\cdot; \mathbb{K}) : \mathrm{Top} \to \mathrm{Vect}_{\mathbb{K}}$ is symmetric monoidal, and thus defines a functor $C_{\bullet}(\cdot; \mathbb{K}) : \mathrm{Op}(\mathrm{Top}) \to \mathrm{Op}(\mathrm{Vect}_{\mathbb{K}})$, as well as $\mathrm{Alg}_{\mathcal{O}}(\mathrm{Top}) \to \mathrm{Alg}_{F(\mathcal{O})}(\mathrm{Vect}_{\mathbb{K}})$.

*Example* C.1.16. The functor $H^{\bullet} : \mathrm{Vect}_{\mathbb{K}} \to \mathbb{K}\text{-}\mathrm{Mod}_{\mathbb{Z}}$ of taking the cohomology object is lax symmetric monoidal, and thus defines a functor $H^{\bullet} : \mathrm{Op}(\mathrm{Vect}_{\mathbb{K}}) \to \mathrm{Op}(\mathbb{K}\text{-}\mathrm{Mod}_{\mathbb{Z}})$. The precomposition of this functor with $C_{\bullet}(\cdot; \mathbb{K}) : \mathrm{Op}(\mathrm{Top}) \to \mathrm{Op}(\mathrm{Vect}_{\mathbb{K}})$ from C.1.15 defines $H_{\bullet}(\cdot; \mathbb{K}) : \mathrm{Op}(\mathrm{Top}) \to \mathrm{Op}(\mathbb{K}\text{-}\mathrm{Mod}_{\mathbb{Z}})$ the homology operad functor.

## C.2. The Hadamard tensor product and Hopf operads.

Let $\mathcal{C}$ be a symmetric monoidal category. Define the *Hadamard tensor product*

$$\otimes^H : \mathrm{Op}(\mathcal{C}) \times \mathrm{Op}(\mathcal{C}) \to \mathrm{Op}(\mathcal{C}) \qquad \text{by} \qquad (\mathcal{O} \otimes^H \mathcal{O}')(\{c_i\}, d) = \mathcal{O}(\{c_i\}, d) \otimes \mathcal{O}'(\{c_i\}, d)$$

for each finite set $I$ and indexed collections $\{c_i\}_{i \in I}$ and $d$ of objects of $\mathcal{O}$. The composition morphisms are defined as the tensor products of those for $\mathcal{O}$ and $\mathcal{O}'$.

*Proposition* C.2.1. The category $\mathrm{Op}(\mathcal{C})$ is symmetric monoidal with respect to $\otimes^H$, and tensor unit given by $\mathrm{Comm}_{\mathcal{C}}$.

*Proposition* C.2.2. Let $\mathcal{O}, \mathcal{O}' \in \mathrm{Op}(\mathcal{C})$. The symmetric monoidal structure on $\mathcal{C}$ lifts to a bifunctor

$$\mathrm{Alg}_{\mathcal{O}}(\mathcal{C}) \times \mathrm{Alg}_{\mathcal{O}'}(\mathcal{C}) \to \mathrm{Alg}_{\mathcal{O} \otimes^H \mathcal{O}'}(\mathcal{C}) \ .$$

*Definition* C.2.3. A *Hopf operad* in $\mathcal{C}$ is a coassociative coalgebra object in the category $\mathrm{Op}(\mathcal{C})$. Concretely, a Hopf operad is an operad $\mathcal{O}$ together with morphisms $\Delta : \mathcal{O} \to \mathcal{O}^{\otimes^H 2}$ and $\varepsilon : \mathcal{O} \to \mathrm{Comm}_{\mathcal{C}}$ satisfying the usual relations of a coalgebra.

A Hopf operad is called *cocommutative* if it is cocommutative as a coalgebra object. We denote the category of Hopf operads in $\mathcal{C}$ by $\mathrm{HOp}(\mathcal{C}) = \mathrm{CoAss}(\mathrm{Op}(\mathcal{C}))$ and the full subcategory of cocommutative Hopf operads by $\mathrm{HOp}^{\mathrm{co}}(\mathcal{C}) = \mathrm{CoComm}(\mathrm{Op}(\mathcal{C}))$.

*Proposition* C.2.4. Let $\mathcal{O}$ be a Hopf operad with structure maps $\Delta : \mathcal{O} \to \mathcal{O}^{\otimes^{H}2}$ and $\varepsilon : \mathcal{O} \to \mathrm{Comm}_{\mathcal{C}}$ as above. Then the bifunctor

$$\otimes : \mathrm{Alg}_{\mathcal{O}}(\mathcal{C}) \times \mathrm{Alg}_{\mathcal{O}}(\mathcal{C}) \to \mathrm{Alg}_{\mathcal{O}}(\mathcal{C}) \qquad \text{defined by} \qquad (A, B) \mapsto \Delta^*(A \otimes B)$$

defines a monoidal structure on the category $\mathrm{Alg}_{\mathcal{O}}(\mathcal{C})$ of $\mathcal{O}$ algebras in $\mathcal{C}$. If $\mathcal{O}$ is cocommutative, then this defines a symmetric monoidal structure.

*Proposition* C.2.5. Let $\mathcal{O}$ be a cocommutative Hopf operad in $\mathcal{C}$. Then there is a natural symmetric monoidal equivalence

$$\mathrm{Alg}_{\mathcal{O}}(\mathrm{Op}(\mathcal{C})) \cong \mathrm{Op}(\mathrm{Alg}_{\mathcal{O}}(\mathcal{C})) \ ,$$

preserving the forgetful functor to $\mathrm{Op}(\mathcal{C})$, where $\mathrm{Alg}_{\mathcal{O}}(\mathrm{Op}(\mathcal{C}))$ is defined using the Hadamard monoidal structure on $\mathrm{Op}(\mathcal{C})$ and equipped with the monoidal structure coming from the Hopf structure on $\mathcal{O}$, while $\mathrm{Op}(\mathrm{Alg}_{\mathcal{O}}(\mathcal{C}))$ is defined using the latter, and is equipped with the former.

*Corollary* C.2.6. There is an equivalence $\mathrm{HOp}(\mathcal{C}) \cong \mathrm{Op}(\mathrm{CoAss}(\mathcal{C}))$ inducing $\mathrm{HOp}^{\mathrm{co}}(\mathcal{C}) \cong \mathrm{Op}(\mathrm{CoComm}(\mathcal{C}))$.

*Example* C.2.7. There is a natural symmetric monoidal lift $C_\bullet(\cdot; \mathbb{K}) : \mathrm{Top} \to \mathrm{CoComm}(\mathrm{Vect}_{\mathbb{K}})$ of the functor of Example C.1.15, by Remark A.5.7. Thus, there is a natural lift $C_\bullet(\cdot; \mathbb{K}) : \mathrm{Op}(\mathrm{Top}) \to \mathrm{HOp}^{\mathrm{co}}(\mathrm{Vect}_{\mathbb{K}})$.

*Example* C.2.8. Let $\mathcal{C}$ be a cartesian monoidal category. Then every object of $\mathcal{C}$ is canonically a cocommutative colagebra, determining a canonical equivalence $\mathcal{C} \cong \mathrm{CoComm}(\mathcal{C})$, with structure map given by the diagonal. This induces a canonical equivalence $\mathrm{Op}(\mathcal{C}) \cong \mathrm{HOp}^{\mathrm{co}}(\mathcal{C})$.

C.3. **The Boardman-Vogt tensor product.** There is an alternate symmetric monoidal structure called the Boardman-Vogt tensor product [BV73], which is defined on cocommutative Hopf operads $\mathrm{HOp}^{\mathrm{co}}(\mathcal{C})$ by the following:

*Proposition* C.3.1. There is a unique symmetric monoidal structure $\star : \mathrm{HOp}^{\mathrm{co}}(\mathcal{C})^{\times 2} \to \mathrm{HOp}^{\mathrm{co}}(\mathcal{C})$ equipped with natural isomorphisms

$$\mathrm{Alg}_{\mathcal{O}\star\mathcal{P}}(\mathcal{C}) \cong \mathrm{Alg}_{\mathcal{O}}(\mathrm{Alg}_{\mathcal{P}}(\mathcal{C})) \ .$$

The preceding proposition can be interpreted as the statement that the Boardman-Vogt tensor product makes cocommutative Hopf operads into a closed cartesian symmetric monoidal category, with internal Hom objects $\mathcal{H}om(\mathcal{O}, \mathcal{P}) = \mathrm{Alg}_{\mathcal{O}}(\mathcal{P})$.

C.4. **The little $d$-cubes operad $\mathbb{E}_d$.** Let $\square^d = (-1, 1)^d$ denote the open cube of dimension $d$, coordinatized as a submanifold of $\mathbb{R}^d$.

*Definition* C.4.1. A map $f : \square^d \to \square^d$ is called a rectilinear embedding if it is defined by

$$f(x_1, ..., x_d) = (a_1 x_1 + b_1, ..., a_d x_d + b_d)$$

for some $a_1, ..., a_d, b_1, ..., b_d \in \mathbb{R}$ with $a_i > 0$. More generally, an embedding $f : I \times \square^d \to \square^d$ for some finite set $I$ is called rectilinear if it its restriction to $\{i\} \times \square^d$ is rectilinar for each $i \in I$.

Let $\mathrm{Rect}_d^I \in \mathrm{Top}$ denote the space of rectilinear embeddings, topologized as an open subset of $(\mathbb{R}^{2d})^I$, or equivalently as a subspace of $\mathrm{Emb}(I \times \square^d, \square^d)$ with the compact-open topology.

*Definition* C.4.2. The little $d$-cubes operad $\mathbb{E}_d \in \mathrm{Op}(\mathrm{Top})$ is the single coloured operad in $\mathrm{Top}$ defined by

$$\mathbb{E}_d(I) = \mathrm{Rect}_d^I \qquad \text{with} \qquad \times_{j \in J} \mathrm{Rect}_d^{I_j} \times \mathrm{Rect}_d^J \to \mathrm{Rect}_d^I$$

given by the obvious composition of rectilinear embeddings.

*Remark* C.4.3. The little $d$-cubes operad defines an operad $C_\bullet(\mathbb{E}_d; \mathbb{K}) \in \mathrm{Op}(\mathrm{Vect}_\mathbb{K})$ in $\mathrm{Vect}_\mathbb{K}$, as in Examples C.1.15 and C.2.7.

*Definition* C.4.4. Let $I$ be a finite set and $X \in \mathrm{Top}$ a topological space. The configuration space $\mathrm{Conf}^I(X)$ of configurations of $I$ points in $X$ is the space

$$\mathrm{Conf}^I(X) = \mathrm{Emb}(I, X) = (X)^I \backslash \{(x_i) | x_i = x_j \text{ for } i \neq j\} \quad \in \quad \mathrm{Top}$$

topologized as an open subset of $(\mathbb{R}^d)^I$, or equivalently with the compact-open topology.

Note there is a canonical map $\mathrm{ev} : \mathrm{Rect}_d^I \to \mathrm{Conf}^I(\square^d)$, given by evaluation at the origin $\{0\} \in \square^d$ for each $i \in I$.

*Proposition* C.4.5. The evaluation map $\mathrm{ev} : \mathrm{Rect}_d^I \xrightarrow{\sim} \mathrm{Conf}^I(\square^d)$ defines a homotopy equivalence.

*Remark* C.4.6. Fix a homeomorphism $\square^d \xrightarrow{\cong} \mathbb{R}^d$. This induces homemorphisms $\mathrm{Conf}^I(\square^d) \xrightarrow{\cong} \mathrm{Conf}^I(\mathbb{R}^d)$ for each $I$, and thus the configuration spaces of points $\mathrm{Conf}^I(\mathbb{R}^d)$ in $\mathbb{R}^d$ define a homotopy equivalent model of $\mathbb{E}_d$, by an enhancement of C.4.5 above.

The space of choices of homeomorphism $\square^d \xrightarrow{\cong} \mathbb{R}^d$ is a torsor for the topological group $\mathrm{Top}(d) = \mathrm{Aut}_{\mathrm{Top}}(\mathbb{R}^d)$, and the structure of an algebra in Top over $\mathbb{E}_d$ does not determine equivariance data for the structure maps with respect to this action of $\mathrm{Top}(d)$. Such additional data is equivalent to a lift to an algebra over the *unoriented $d$-cubes* operad $\mathbb{E}_d^{\mathrm{Top}(d)}$, defined in 22.0.1.

*Example* C.4.7. The $\mathbb{E}_0$ operad. The tensor unit for unital cocommutative Hopf operads under Boardman Vogt tensor product.

*Proposition* C.4.8. There is a homotopy equivalence of operads $C_\bullet(\mathbb{E}_1; \mathbb{K}) \cong \mathrm{Ass}_\mathbb{K} \in \mathrm{Op}(\mathrm{Vect}_\mathbb{K})$.

The following was initally proved in [Dun88], and in the context of quasioperads in Theorem 5.1.2.2 in [Lur12]:

*Theorem* C.4.9. The Boardman-Vogt tensor product of the $\mathbb{E}_n$ and $\mathbb{E}_m$ operads is canonically equivalent to the $\mathbb{E}_{n+m}$ operad:

$$\mathbb{E}_n \star \mathbb{E}_m \cong \mathbb{E}_{n+m} \ .$$

*Remark* C.4.10. Concretely, for $n = 2$ for example, the preceding Theorem identifies $\mathbb{E}_2$ algebras with associative algebra objects in the category of associative algebras, or equivalently vector spaces equipped with two compatible associative algebra structures.

*Remark* C.4.11. In the non-derived setting, the Eckmann-Hilton arguement implies that $\mathbb{E}_n$ algebras for $n \geqslant 2$ are necessarily commutative, and moreover that the various compatible associative algebra structures are all canonically equivalent. However, this arguement fails to extend homotopy coherently, and is for example obstructed by the induced $\mathbb{P}_n$ algebra structure on homology.

*Remark* C.4.12. There are canonical maps of operads $\mathbb{E}_n \to \mathbb{E}_{n+1}$ and induced forgetful functors $\mathrm{Alg}_{\mathbb{E}_{n+1}}(\mathrm{Vect}) \to \mathrm{Alg}_{\mathbb{E}_n}(\mathrm{Vect})$, which correspond to forgetting one of the various compatible associative algebra structures.

*Definition* C.4.13. The $\mathbb{E}_\infty$ operad is the colimit $\mathbb{E}_\infty = \mathrm{colim}_n \mathbb{E}_n$ of the $\mathbb{E}_n$ operads.

*Proposition* C.4.14. There is a canonical homotopy equivalence $\mathbb{E}_\infty \xrightarrow{\cong} \mathrm{Comm}$.

*Remark* C.4.15. Concretely, an algebra over the $\mathbb{E}_\infty$ operad is equivalent to a coherent system of $\mathbb{E}_n$ algebras for each $n \in \mathbb{N}$. A commutative algebra evidently induces such a system, which defines the above map, and the statement is that up to homotopy all $\mathbb{E}_\infty$ algebras are of this form.

## C.5. The $(d-1)$-shifted Poisson operad $\mathbb{P}_d$.

*Definition* C.5.1. The $(d-1)$ shifted Poisson operad is the operad $\mathbb{P}_d \in \mathrm{Op}(\mathrm{Vect}_{\mathbb{K}})$ generated by

$$\mathbb{P}_d(2) = \mathbb{K}_m \oplus \mathbb{K}_\pi[d-1] \quad \in \quad \mathbb{K}[S_2]\text{-Mod}$$

where $\mathbb{K}_m$ is the trivial representation and $\mathbb{K}_\pi$ is the sign, subject to the relations

$$\mathbb{K}_{\mathrm{Jac}(\pi)}[2d-2] \hookrightarrow \mathcal{F}(3) \qquad \text{defined by} \qquad 1 \mapsto \qquad \mathrm{Jac}(\pi) = \pi \circ (\mathbb{1} \otimes \pi) - \pi \circ (\pi \otimes \mathbb{1}) - \pi \circ (\mathbb{1} \otimes \pi) \circ \sigma_{12} \ ,$$

$$\mathbb{K}_{\mathrm{Ass}(m)}[\phantom{2d-2}] \hookrightarrow \mathcal{F}(3) \qquad \text{defined by} \qquad 1 \mapsto \qquad \mathrm{Ass}(m) = m \circ (m \otimes \mathbb{1}) - m \circ (\mathbb{1} \otimes m)$$

$$\mathbb{K}_{\mathrm{Dist}(m,\pi)}[d-1] \hookrightarrow \mathcal{F}(3) \qquad \text{defined by} \qquad 1 \mapsto \quad \mathrm{Dist}(m,\pi) = \pi \circ (\mathbb{1} \otimes m) - m \circ (\pi \otimes \mathbb{1}) - m \circ (\mathbb{1} \otimes \pi) \ .$$

where $\mathcal{F} = \mathcal{F}(\mathbb{K}_m \oplus \mathbb{K}_\pi[d-1]) \in \mathrm{Op}(\mathrm{Vect}_{\mathbb{K}})$ is the free operad on these generators.

*Example* C.5.2. The category $\mathrm{Alg}_{\mathbb{P}_d}(\mathrm{Vect}_{\mathbb{K}})$ is evidently given by the usual category of $(d-1)$-shifted Poisson algebras, that is, commutative algebras $A \in \mathrm{Comm}(\mathrm{Vect}_{\mathbb{K}})$ together with a Lie bracket $\pi : A^{\otimes 2} \to A[1-d]$ such that $\pi$ is a derivation of the product $m$.

*Example* C.5.3. Consider the operad $H_\bullet(\mathbb{E}_n) \in \mathrm{Op}(\mathrm{Vect}_{\mathbb{K}})$, as defined in Example C.1.16. We have

$$H_\bullet(\mathbb{E}_n)(I) = H_\bullet(\mathrm{Conf}^I(\mathbb{R}^d); \mathbb{K})$$

for each $I \in \mathrm{fSet}$. For each $i, j \in I$, define the map

$$F_{ij} : \mathrm{Conf}^I(\mathbb{R}^d) \to S^{d-1} \qquad (x_i)_{i \in I} \mapsto (x_i - x_j)/|x_i - x_j|$$

and let $\omega_{ij} = F_{ij}^* \Omega \in H^{d-1}(\mathrm{Conf}^I(\mathbb{R}^d); \mathbb{K})$ where $\Omega \in H^{d-1}(S^{d-1}; \mathbb{K})$ is a fixed choice of generator. Note that the natural $\mathrm{Aut}_{\mathrm{fSet}}(I)$ action on $H^{d-1}(\mathrm{Conf}^I(\mathbb{R}^d); \mathbb{K})$ satisfies $\pi \cdot \omega_{ij} = \omega_{\pi(i)\pi(j)}$.

The following theorem was proved by Arnold in the case $d = 2$ and by F. Cohen for $d \geqslant 2$:

*Theorem* C.5.4. The cohomology ring $H^\bullet(\mathrm{Conf}^I(\mathbb{R}^d); \mathbb{K})$ is generated by the classes $\omega_{ij} \in H^{d-1}(\mathrm{Conf}^I(\mathbb{R}^d); \mathbb{K})$ for each $i, j \in I$, subject to the relations:

- $\omega_{ij} = (-1)^d \omega_{ji}$ ,
- $\omega_{ij}\omega_{jk} + \omega_{jk}\omega_{ki} + \omega_{ki}\omega_{ij} = 0$ , and
- $\omega_{ij}^2 = 0$, for $n$ odd.

*Remark* C.5.5. In fact, it was proved in *loc. cit.* that $H_\bullet(\mathrm{Conf}^I(\mathbb{R}^d); \mathbb{Z})$ is torsion free, and that the above result remains true over $\mathbb{Z}$.

Note that the map $F_{12} : \mathrm{Conf}^2(\mathbb{R}^d) \xrightarrow{\sim} S^{d-1}$ is a homotopy equivalence, inducing an identification

$$H_\bullet(\mathbb{E}_d; \mathbb{K})(2) \cong H_\bullet(S^{d-1}; \mathbb{K}) = \mathbb{K} \oplus \mathbb{K}_{\omega_{12}}^\vee[d-1]$$

Thus, the above result implies that for $d \geqslant 2$, the operad $H_\bullet(\mathbb{E}_d; \mathbb{K})$ is generated by its arity two operations, and these have the same symmetric group action and relations as the generators $\mathbb{P}_2(2)$ above, and we obtain:

*Corollary* C.5.6. The $(d-1)$-shifted Poisson operad $\mathbb{P}_d \cong H_\bullet(\mathbb{E}_d; \mathbb{K}) \in \mathrm{Op}(\mathrm{Vect}_{\mathbb{K}})$ is equivalent to the homology of the litte $d$-disks operad $\mathbb{E}_d$.

C.6. **The Beilinson-Drinfeld operad $\mathbb{BD}_d$ and quantization of $\mathbb{P}_d$ algebras.** Throughout, let $\mathbb{K}[\hbar] = \mathrm{Sym}^\bullet(\mathbb{K}_\hbar\langle 1 \rangle)$ denote a graded polynomial ring over $\mathbb{K}$ with $\hbar$ of weight $+1$. The grading is interpreted as defining a $\mathbb{G}_m$ action on $\mathbb{A}^1_\hbar = \mathrm{Spec}\,\mathbb{K}[\hbar]$. Let $\mathrm{D}^b_{\mathrm{fg}}(\mathbb{K}[\hbar])$ denote the bounded derived category of graded modules over $\mathbb{K}[\hbar]$ with finitely generated cohomology.

*Remark* C.6.1. The grading ensures an equivalence of the localization $\mathrm{D}^b_{\mathrm{fg}}(\mathbb{K}[\hbar^{\pm 1}]) \cong \mathrm{Perf}_\mathbb{K}$ with the bounded derived category of complexes over $\mathbb{K}$ with finite dimensional cohomology. Thus, specialization over $\mathbb{A}^1/\mathbb{G}_m$ to the central and generic fibres defines symmetric monoidal functors

(C.6.1) $\qquad (\cdot)|_{\{0\}} : \mathrm{D}^b_{\mathrm{fg}}(\mathbb{K}[\hbar]) \to \mathrm{Perf}_\mathbb{K} \qquad$ and $\qquad (\cdot)|_{\{1\}} : \mathrm{D}^b_{\mathrm{fg}}(\mathbb{K}[\hbar]) \to \mathrm{D}^b_{\mathrm{fg}}(\mathbb{K}[\hbar^{\pm 1}]) \cong \mathrm{Perf}_\mathbb{K}$ .

*Definition* C.6.2. The dimension 0 Beilinson-Drinfeld operad is the operad $\mathbb{BD}^\hbar_0 \in \mathrm{Op}(\mathrm{D}^b_{\mathrm{fg}}(\mathbb{K}[\hbar]))$ generated by

$$\mathbb{BD}_0(2) = \left[ \mathbb{K}[\hbar]_m \xrightarrow{\hbar} \mathbb{K}[\hbar]_\pi[-1]\langle 1 \rangle \right] \quad \in \quad \mathrm{D}^b_{\mathrm{fg}}(\mathbb{K}[\hbar][S_2])$$

where $\mathbb{K}[\hbar]_m$ is the trivial representation and $\mathbb{K}[\hbar]_\pi$ is the sign, subject to the relations of the $\mathbb{P}_0$ operad C.5.1 extended linearly to $\mathbb{K}[\hbar]$.

*Example* C.6.3. Concretely, an object $A \in \mathrm{Alg}_{\mathbb{BD}_0}(\mathrm{D}^b_{\mathrm{fg}}(\mathbb{K}[\hbar]))$ is given by

- a complex $(A, d) \in \mathrm{D}^b_{\mathrm{fg}}(\mathbb{K}[\hbar])$ of graded $\mathbb{K}[\hbar]$ modules,
- a commutative multiplication $\cdot : A^{\otimes 2} \to A$, and
- a Lie bracket $\{,\} : A^{\otimes 2} \to A[-1]$ of degree $-1$ ,

such that $\{,\}$ is a biderivation, as for a usual Poisson algebra, and moreover for each $a, b \in A$,

$$d(a \cdot b) = d(a) \cdot b + (-1)^{|a|} a \cdot d(b) + \hbar\{a, b\} \ .$$

Note that the specialization at $\hbar = 0$ of such an algebra is just a usual $\mathbb{P}_0$ algebra in $\mathrm{Perf}_\mathbb{K}$, while for $\hbar \neq 0$ the complex of generators is acyclic so that operations on $(A, d)$ are compatibly trivializeable up to homotopy so that the resulting object defines an $\mathbb{E}_0$ algebra.

Thus, the $\mathbb{BD}_0$ operad controls quantizations $\mathbb{P}_0$ algebras to $\mathbb{E}_0$ algebras, in the following sense:

*Proposition* C.6.4. There are canonical equivalences of operads

$\qquad \mathbb{BD}_0|_{\{0\}} \cong \mathbb{P}_0 \quad \in \quad \mathrm{Op}(\mathrm{Perf}_\mathbb{K}) \qquad$ and $\qquad \mathbb{BD}_0|_{\{1\}} \cong \mathbb{E}_0 \quad \in \quad \mathrm{Op}(\mathrm{D}^b_{\mathrm{fg}}(\mathbb{K}[\hbar^{\pm 1}])) \cong \mathrm{Op}(\mathrm{Perf}_\mathbb{K})$ .

In particular, specialization over $\mathbb{A}^1/\mathbb{G}_m$ as in Equation C.6.1 defines symmetric monoidal functors

$\qquad (\cdot)|_{\{0\}} : \mathrm{Alg}_{\mathbb{BD}_0}(\mathrm{D}^b_{\mathrm{fg}}(\mathbb{K}[\hbar])) \to \mathrm{Alg}_{\mathbb{P}_0}(\mathrm{Perf}_\mathbb{K}) \qquad$ and $\qquad (\cdot)|_{\{1\}} : \mathrm{Alg}_{\mathbb{BD}_0}(\mathrm{D}^b_{\mathrm{fg}}(\mathbb{K}[\hbar])) \to \mathrm{Alg}_{\mathbb{E}_0}(\mathrm{Perf}_\mathbb{K})$ .

*Definition* C.6.5. The dimension 1 Beilinson-Drinfeld operad is the operad $\mathbb{BD}_1 \in \mathrm{Op}(\mathrm{D}^b_{\mathrm{fg}}(\mathbb{K}[\hbar]))$ generated by

$$\mathbb{BD}_1(2) = (\mathbb{K}[\hbar]_m \oplus \mathbb{K}[\hbar]_{m^{\mathrm{op}}}) \oplus \mathbb{K}[\hbar]_\pi\langle 1 \rangle \quad \in \quad \mathrm{D}^b_{\mathrm{fg}}(\mathbb{K}[\hbar][S_2])$$

where $\mathbb{K}[\hbar]_m \oplus \mathbb{K}[\hbar]_{m^{\mathrm{op}}}$ is the regular representation and $\mathbb{K}[\hbar]_\pi$ is the sign, subject to the relations of the $\mathbb{P}_1$ operad C.5.1 extended linearly to $\mathbb{K}[\hbar]$, together with the relation

$$\mathbb{K}[\hbar]_{\mathrm{BD}(m,\pi)} \hookrightarrow \mathbb{BD}_0(2) \subset \mathcal{F}(\mathbb{BD}_1(2)) \qquad \text{defined by} \qquad 1 \mapsto \mathrm{BD}(m, \pi) = m - m^{\mathrm{op}} - \hbar\pi \ .$$

*Example* C.6.6. Concretely, an object $A \in \mathrm{Alg}_{\mathbb{BD}_0}(\mathrm{D}^b_{\mathrm{fg}}(\mathbb{K}[\hbar]))$ is given by an associative algebra $A$, together with a Lie bracket $\{,\} : A^{\otimes 2} \to A$ which is a biderivation of the associative product, and satisfies

$$ab - (-1)^{|a||b|} ba = \hbar\{a, b\} \ .$$

Note that at $\hbar = 0$ the product is commutative and thus $A$ defines a usual Poisson algebra, while for $\hbar \neq 0$ the operation $\{,\}$ is determined by the associative product, so that $A$ is just a usual associative algebra.

Thus, the $\mathbb{BD}_1$ operad classifies quantizations of $\mathbb{P}_1$ algebras to $\mathbb{E}_1$ algebras, in the following sense:

*Proposition* C.6.7. There are canonical equivalences of operads

$$\mathbb{BD}_1|_{\{0\}} \cong \mathbb{P}_1 \quad \in \quad \mathrm{Op}(\mathrm{Vect}_{\mathbb{K}}) \qquad \text{and} \qquad \mathbb{BD}_1|_{\{1\}} \cong \mathbb{E}_1 \quad \in \quad \mathrm{Op}(\mathrm{D}^b_{\mathrm{fg}}(\mathbb{K}[\hbar^{\pm 1}])) \cong \mathrm{Op}(\mathrm{Perf}_{\mathbb{K}}) \ .$$

In particular, specialization over $\mathbb{A}^1/\mathbb{G}_m$ as in Equation C.6.1 defines symmetric monoidal functors

$$(\cdot)|_{\{0\}} : \mathrm{Alg}_{\mathbb{BD}_1}(\mathrm{D}^b_{\mathrm{fg}}(\mathbb{K}[\hbar])) \to \mathrm{Alg}_{\mathbb{P}_1}(\mathrm{Perf}_{\mathbb{K}}) \qquad \text{and} \qquad (\cdot)|_{\{1\}} : \mathrm{Alg}_{\mathbb{BD}_1}(\mathrm{D}^b_{\mathrm{fg}}(\mathbb{K}[\hbar])) \to \mathrm{Alg}_{\mathbb{E}_1}(\mathrm{Perf}_{\mathbb{K}}) \ .$$

More generally, for $n \geqslant 2$, we make the following definition:

*Definition* C.6.8. The dimension $n$ Beilinson-Drinfeld operad is the operad $\mathbb{BD}_n \in \mathrm{Op}(\mathrm{D}^b_{\mathrm{fg}}(\mathbb{K}[\hbar]))$ defined as the image under the Rees construction of the operad $\mathbb{E}_n \in \mathrm{Op}(\mathrm{Perf}_{\mathbb{K}})$ together with the Postnikov filtration.

*Remark* C.6.9. This definition does *not* agree with the definitions given above when applied to the cases $n = 0, 1$; the definitions stated above are the correct ones.

Generalizing Propositions C.6.4 and C.6.7 above, we have:

*Proposition* C.6.10. There are canonical equivalences of operads

$$\mathbb{BD}_n|_{\{0\}} \cong \mathbb{P}_n \quad \in \quad \mathrm{Op}(\mathrm{Vect}_{\mathbb{K}}) \qquad \text{and} \qquad \mathbb{BD}_n|_{\{1\}} \cong \mathbb{E}_n \quad \in \quad \mathrm{Op}(\mathrm{D}^b_{\mathrm{fg}}(\mathbb{K}[\hbar^{\pm 1}])) \cong \mathrm{Op}(\mathrm{Perf}_{\mathbb{K}}) \ .$$

In particular, specialization over $\mathbb{A}^1/\mathbb{G}_m$ as in Equation C.6.1 defines symmetric monoidal functors

$$(\cdot)|_{\{0\}} : \mathrm{Alg}_{\mathbb{BD}_n}(\mathrm{D}^b_{\mathrm{fg}}(\mathbb{K}[\hbar])) \to \mathrm{Alg}_{\mathbb{P}_n}(\mathrm{Perf}_{\mathbb{K}}) \qquad \text{and} \qquad (\cdot)|_{\{1\}} : \mathrm{Alg}_{\mathbb{BD}_n}(\mathrm{D}^b_{\mathrm{fg}}(\mathbb{K}[\hbar])) \to \mathrm{Alg}_{\mathbb{E}_n}(\mathrm{Perf}_{\mathbb{K}}) \ .$$

## Appendix D. Functional Analysis

D.1. **Topological Vector Spaces.**

*Definition* D.1.1. A topological vector space over $\mathbb{K}$ is a vector space $V \in \mathbb{K}$-Mod with a complete, separated, linear topology.

Let $\mathbb{K}$-Mod$_{\text{Top}}$ denote the category of topological vector spaces with continuous linear maps.

*Example* D.1.2. Any vector space $V$ equipped with the discrete topology defines a topological vector space $V$, with the propery that any linear map $V \to W$ is continuous.

All finite dimensional vector spaces will be considered with the discrete topology by default.

*Proposition* D.1.3. A discrete vector space $V$ is canonically equivalent to the (filtered) colimit $V = \text{colim}_k V_k$ of its finite dimensional subspaces $V_k \in \mathbb{K}$-Mod$_{\text{fg}}$ in the category $\mathbb{K}$-Mod$_{\text{Top}}$.

*Remark* D.1.4. The preceding proposition defines a fully faithful embedding

$$\mathbb{K}\text{-Mod} = \text{Ind}(\mathbb{K}\text{-Mod}_{\text{fg}}) \hookrightarrow \mathbb{K}\text{-Mod}_{\text{Top}} \ .$$

*Example* D.1.5. The discrete vector space $t^{-1}\mathbb{K}[t^{-1}]$, or more generally $\mathbb{K}((t))/t^n\mathbb{K}[[t]]$, is presented as the colimit

$$\mathbb{K}((t))/t^n\mathbb{K}[[t]] = \text{colim}_k t^{-k}\mathbb{K}[[t]]/t^n\mathbb{K}[[t]] \ .$$

*Example* D.1.6. Let $V = \lim_i V_i$ with $V_i \in \mathbb{K}$-Mod$_{\text{fg}}$ be a pro-finite dimensional vector space. Then $V$ has a canonical profinite topology defined by the basis of neighbourhoods of $0 \in V$ given by the subspaces $\ker(\pi_i)$, where $\pi_i : V \to V_i$ is the canonical projection.

All pro-finite dimensional vector spaces will be considered with the pro-finite topology by default.

*Proposition* D.1.7. A pro-finite dimensional vector space $V$ is canonically equivalent to the limit $V = \lim_i V_i$ of its finite dimensional quotients $V_i \in \mathbb{K}$-Mod$_{\text{fg}}$ in the category $\mathbb{K}$-Mod$_{\text{Top}}$.

*Remark* D.1.8. The preceding proposition defines a fully faithful embedding $\text{Pro}(\mathbb{K}\text{-Mod}_{\text{fg}}) \hookrightarrow \mathbb{K}\text{-Mod}_{\text{Top}}$.

*Example* D.1.9. The vector space $\mathbb{K}[[t]]$ of formal power series is pro-finite, as it is given by the limit

$$\mathbb{K}[[t]] = \lim_n \mathbb{K}[[t]]/t^n\mathbb{K}[[t]] \ .$$

*Definition* D.1.10. A Tate vector space is a topological vector space $V$ that admits a direct sum decomposition $V = U \oplus W$ for $U$ a discrete vector space and $V$ a pro-finite dimensional vector space, as topological vector spaces.

*Example* D.1.11. The prototypical example of a Tate vector space is the field of Laurent series $\mathbb{K}((t))$, presented for example as $\mathbb{K}((t)) = t^{-1}\mathbb{K}[t^{-1}] \oplus \mathbb{K}[[t]]$.

D.2. **Tensor structures on topological vector spaces.** In this appendix, we summarize the main results of the paper [Bei07] of Beilinson, building on Chapter 3.6 of [BD04]. All of the objects will be of cohomological degree zero and all the functors non-derived, in contrast with our general conventions. However, we remark that [Ras20b] establishes some analogous results in the derived setting, which we will also need.

Let $\{V_i\}_{i \in I}$ denote a finite collection of topological vector spaces $V_i \in \mathbb{K}$-Mod$_{\text{Top}}$, and consider the algebraic tensor product $\otimes_i V_i \in \mathbb{K}$-Mod.

*Definition* D.2.1. The $\otimes^*$ symmetric monoidal structure $\otimes^* : \mathbb{K}\text{-Mod}_{\text{Top}}^{\times I} \to \mathbb{K}\text{-Mod}_{\text{Top}}$ is defined by setting $\otimes_i^* V_i \in \mathbb{K}\text{-Mod}_{\text{Top}}$ to be the completion of $\otimes_i V_i$ with respect to the topology defined as follows: a subspace $Q \subset \otimes_i V_i$ is open if for every $J \subset I$ and $v \in \otimes_{i \in I \setminus J} V_i$ there exist open subspaces $P_j \subset V_j$ for each $j \in J$ such that

$$\otimes_j P_j \otimes v \quad \subset \quad Q .$$

*Remark* D.2.2. The induced operad $\mathbb{K}\text{-Mod}_{\text{Top}}^*$ is defined, following Example C.1.3, by

$$\text{Hom}_{\mathbb{K}\text{-Mod}_{\text{Top}}^*}(\{V_i\}, W) := \text{Hom}_{\mathbb{K}\text{-Mod}_{\text{Top}}}(\otimes_i^* V_i, W) = \{F : \times_i V_i \to W \mid F \text{ is continuous and multilinear}\} .$$

In particular, an associative algebra object in the category $\mathbb{K}\text{-Mod}_{\text{Top}}^*$ is given by a topological vector space $A \in \mathbb{K}\text{-Mod}_{\text{Top}}$ together with an associative, bilinear product $\mu : A \otimes A \to A$ such that the corresponding map $A \times A \to A$ is continuous.

*Definition* D.2.3. The $\otimes^{\text{ch}}$ monoidal structure $\otimes^{\text{ch}, \tau} : \mathbb{K}\text{-Mod}_{\text{Top}}^{\times 2} \to \mathbb{K}\text{-Mod}_{\text{Top}}$ is defined for each linear order $\tau : \{1, ..., n\} \to I$ by setting

$$\otimes_i^{\text{ch}, \tau} V_i = V_{\tau(1)} \otimes^{\text{ch}} ... \otimes^{\text{ch}} V_{\tau(n)} \quad \in \quad \mathbb{K}\text{-Mod}_{\text{Top}}$$

to be the completion of $\otimes_i V_i$ with respect to the topology defined as follows: a subspace $Q \subset \otimes_i V_i$ is open if for every $a \in \{1, ..., n\}$ and $v \in V_{\tau(a+1)} \otimes ... \otimes V_{\tau(n)}$ there exists an open subspace $P_a \subset V_a$ such that

$$V_{\tau(1)} \otimes ... \otimes V_{\tau(a-1)} \otimes P_a \otimes v \quad \subset \quad Q .$$

*Remark* D.2.4. Equivalently, the $\otimes^{\text{ch}}$ monoidal structure is defined iteratively for $V = \lim_n V_n$ with each $V_n = \text{colim}_k V_{n,k}$ by

$$U \otimes^{\text{ch}} V = \lim_n \text{colim}_k U \otimes V_{n,k} ,$$

for any $U \in \mathbb{K}\text{-Mod}_{\text{Top}}$.

*Remark* D.2.5. Note that our notation differs slightly from that of [Bei07], [BD04] and [Ras20b], as we use $\otimes^{\text{ch}}$ in place of $\otimes^{\to}$.

*Remark* D.2.6. An associative algebra object in the category $\mathbb{K}\text{-Mod}_{\text{Top}}^{\text{ch}}$ is given by a topological vector space $A \in \mathbb{K}\text{-Mod}_{\text{Top}}$ together with an associative, bilinear product $\mu : A \otimes A \to A$ such that the corresponding map $A \times A \to A$ is continuous, and the open left ideals of $A$ form a basis for the topology of $A$.

*Remark* D.2.7. Although $\otimes^{\text{ch}}$ is evidently not symmetric, it defines a natural (symmetric) operad stucture as follows:

*Definition* D.2.8. The induced operad $\mathbb{K}\text{-Mod}_{\text{Top}}^{\text{ch}, s}$ is defined by

$$\text{Hom}_{\mathbb{K}\text{-Mod}_{\text{Top}}^{\text{ch}, s}}(\{V_i\}, W) = \bigoplus_{\tau \in S_I} \text{Hom}_{\mathbb{K}\text{-Mod}_{\text{Top}}}(\otimes_i^{\text{ch}, \tau} V_i, W) ,$$

where $S_I$ is the $S_n$ torsor of linear orders $\tau : \{1, ..., n\} \xrightarrow{\cong} I$ and $n = |I|$.

*Definition* D.2.9. The $\otimes^!$ symmetric monoidal structure $\otimes^! : \mathbb{K}\text{-Mod}_{\text{Top}}^{\times I} \to \mathbb{K}\text{-Mod}_{\text{Top}}$ is defined by setting $\otimes_i^! V_i \in \mathbb{K}\text{-Mod}_{\text{Top}}$ to be the completion of $\otimes_i V_i$ with respect to the topology with basis of neighbourhoods at 0 given by subspaces of the form

$$\sum_{i \in I} P_i \otimes \left( \otimes_{i' \in I \setminus \{i\}} V_{i'} \right)$$

for $P_i \subset V_i$ an open subspace.

*Remark* D.2.10. Equivalently, the $\otimes^!$ tensor product is defined on $V_i = \lim_{n_i} V/P_{n_i}$ by

$$\otimes_i^! V_i = \lim_{(n_i)_{i \in I}} \otimes_i (V_i/P_{n_i}) \ .$$

Thus, the $\otimes^!$ tensor product is simply the usual, completed tensor product.

*Remark* D.2.11. The induced operad $\mathbb{K}\text{-Mod}_{\text{Top}}^!$ is defined, following Example C.1.3, by

$$\text{Hom}_{\mathbb{K}\text{-Mod}_{\text{Top}}^!}(\{V_i\}, W) := \text{Hom}_{\mathbb{K}\text{-Mod}_{\text{Top}}}(\otimes_i^! V_i, W) \ .$$

An associative algebra object in the category $\mathbb{K}\text{-Mod}_{\text{Top}}^!$ is given by a topological vector space $A \in \mathbb{K}\text{-Mod}_{\text{Top}}$ together with an associative, bilinear product $\mu : A \otimes A \to A$ such that the corresponding map $A \times A \to A$ is continuous, and the open two-sided ideals of $A$ form a basis for the topology of $A$.

*Remark* D.2.12. In general, the topology on $\otimes_i V_i$ underlying the $\otimes^!$ monoidal structure is strictly coarser than that underlying the $\otimes^{\text{ch}}$ monoidal structure (for each fixed $\tau$), which is strictly coarser than that underlying the $\otimes^*$ monoidal structure. Thus, we have natural maps

$$\otimes_i^* V_i \to \otimes_i^{\text{ch},\tau} V_i \to \otimes_i^! V_i \qquad \text{and} \qquad \text{Hom}_{\mathbb{K}\text{-Mod}_{\text{Top}}}(\otimes_i^! V_i, W) \to \text{Hom}_{\mathbb{K}\text{-Mod}_{\text{Top}}}(\otimes_i^{\text{ch},\tau} V_i, W) \to \text{Hom}_{\mathbb{K}\text{-Mod}_{\text{Top}}}(\otimes_i^* V_i, W) \ ,$$

for any $W \in \mathbb{K}\text{-Mod}_{\text{Top}}$. These induce natural maps of operads

$$\mathbb{K}\text{-Mod}_{\text{Top}}^! \otimes^H \text{Ass} \to \mathbb{K}\text{-Mod}_{\text{Top}}^{\text{ch},s} \to \mathbb{K}\text{-Mod}_{\text{Top}}^* \otimes^H \text{Ass} \ ,$$

where Ass denotes the associative operad and $\otimes^H$ the Hadamard tensor product; see Appendix C for a review and conventions regarding operads.

Composing with the projection Ass $\twoheadrightarrow$ Comm and precomposing with the inclusion Lie $\hookrightarrow$ Ass, we also have maps of operads

$$(\text{D.2.1}) \qquad\qquad \mathbb{K}\text{-Mod}_{\text{Top}}^! \otimes^H \text{Lie} \to \mathbb{K}\text{-Mod}_{\text{Top}}^{\text{ch},s} \to \mathbb{K}\text{-Mod}_{\text{Top}}^* \ .$$

*Proposition* D.2.13. For any $U, V \in \mathbb{K}\text{-Mod}_{\text{Top}}$ there is a short exact sequence of topological vector spaces

$$U \otimes^* V \hookrightarrow U \otimes^{\text{ch}} V \oplus V \otimes^{\text{ch}} U \twoheadrightarrow U \otimes^! V \ ,$$

where the left map is given by the diagonal inclusion, and the right map is given by the difference of projections.

*Corollary* D.2.14. In the special case of arity 2 operations, the sequence of maps in Equation D.2.1 gives a left exact sequence

$$\text{Hom}_{\mathbb{K}\text{-Mod}_{\text{Top}}^!}(U, V; W) \otimes_{\mathbb{K}} \text{Lie}(2) \to \text{Hom}_{\mathbb{K}\text{-Mod}_{\text{Top}}^{\text{ch},s}}(U, V; W) \to \text{Hom}_{\mathbb{K}\text{-Mod}_{\text{Top}}^*}(U, V; W) \ .$$

## Appendix E. Vertex Algebras

### E.1. **Vertex Algebras.**

*Definition* E.1.1. An element $A \in \mathrm{End}(V)[[z^{\pm 1}]]$ is a field if $A(v) \in V((z))$ for each $v \in V$.

*Remark* E.1.2. More explicitly, $A = \sum_{n \in \mathbb{Z}} a_n z^{-n-1} \in \mathrm{End}(V)[[z^{\pm 1}]]$ is a field if for each $v \in V$ there exists $N \in \mathbb{Z}$ such that $a_n(v) = 0$ for all $n > N$.

*Definition* E.1.3. A vertex algebra is a tuple $(V, \varnothing, T, Y)$ of:

- a vector space $V$, the state space
- an element $\varnothing \in V$, the vacuum
- a linear map $T \in \mathrm{End}(V)$, the translation operator
- a linear map $Y(\cdot, z) : V \to \mathrm{End}(V)[[z^{\pm 1}]]$, the vertex operator

such that:

- $Y(\varnothing, z) = \mathbb{1}_V$
- $Y(a, z) = \sum_{n \in \mathbb{Z}} a_n z^{-n-1} \in \mathrm{End}(V)[[z^{\pm 1}]]$ is a field for each $a \in V$ and $v \in V$
- $Y(a, z)(\varnothing) \in V[[z]] \subset V((z))$ for each $a \in V$, and the resulting evaluation satisfies $Y(a, z)(\varnothing)|_{z=0} = a$.
- $[T, Y(a, z)] = \partial_z Y(a, z)$ for each $a \in V$
- $T\varnothing = 0$
- For each $a, b \in V$, the fields $Y(a, z), Y(b, z) \in \mathrm{End}(V)[[z^{\pm 1}]]$ are local with respect to one another.

A $\mathbb{Z}$ grading on a vertex algebra is a $\mathbb{Z}$ grading on its underlying vector space $V$ such that

- $\varnothing \in V_0$
- $T : V \to V[1]$
- $a_n : V \to V[m - n - 1]$ for each $a \in V_m$

A dg vertex algebra is a $\mathbb{Z}$ graded vertex algebra $(V, \varnothing, T, Y)$ together with a linear map $d : V \to V[1]$ such that

- $d^2 = 0$
- $[d, T] = 0$
- $dY(A, z)(b) = Y(dA, z)(B) + Y(A, z)(dB)$ for each $a, b \in V$.

*Proposition* E.1.4. Let $(V, \varnothing, T, Y, d)$ a dg vertex algebra. Then $H^\bullet(V, d)$ is a $\mathbb{Z}$ graded vertex algebra.

*Example* E.1.5. Let $(V, \varnothing, T, Y)$ be a $\mathbb{Z}$ graded vertex algebra and fix $a \in V$ of degree and let

$$d_a = \int Y(A, z)dz := a_0 : V \to V[1]$$

Then $[d_a, T] = 0$ and moreover by corollary 3.3.8 in FBZ we have

$$[a_0, Y(b, z)] = Y(a_0(b), z)$$

Thus, we must only require

$$d_a^2 = (a_0)^2 = \int : Y(a, z)Y(a, z) : dz = 0$$

to ensure $d_a$ defines a differential making $V$ into a dg vertex algebra.

## E.2. **Commutative Vertex Algebras.**

*Definition* E.2.1. A vertex algebra $(V, \varnothing, T, Y)$ is called commutative if the field $Y(A, z) \in \mathrm{End}(V)[z]$ is non-singular, for each $a \in V$.

*Proposition* E.2.2. The following are equivalent:

- A commutative vertex algebra
- A (unital) commutative algebra with a derivation
- A (unital) commutative algebra object in $\mathbb{K}[T]$-Mod

## E.3. **Vertex Lie Algebras.**

*Definition* E.3.1. A vertex Lie algebra is a tuple $(L_0, T, Y_-)$ of:

- a vector space $L_0$
- a linear operator $T \in \mathrm{End}(L_0)$
- a linear map $Y_- : L_0 \to \mathrm{End}(L_0) \otimes z^{-1}\mathbb{C}[[z^{-1}]]$

such that

- $Y_-(a, z) \in \mathrm{End}(L_0) \otimes z^{-1}\mathbb{C}[[z^{-1}]] \subset \mathrm{End}(L_0)[[z^{\pm 1}]]$ is a field for each $a \in z$.
- $Y_-(Ta, z) = \partial_z Y_-(a, z)$ for each $a \in L_0$
- $Y_-(a, z)b = (e^{zT}Y_-(b, -z)a)_-$ for each $a, b \in L_0$
- for any $a, b \in L_0$ with $Y_-(a, z) = \sum_{n \geqslant 0} a_n z^{-n-1}$, we have

$$[a_m, Y_-(b, w)] = \sum_{n \geqslant 0} \binom{m}{n} (w^{m-n} Y_-(a_n(b), w))_-$$

*Email address*: dbutson@perimeterinstitute.ca